\documentclass[12pt,a4paper]{article}
\usepackage{mdframed,mathrsfs,amsthm,graphicx,xcolor,verbatim,bbm,amsmath,amsfonts,amssymb,newclude,nicefrac,hyperref,bm,geometry,mathtools,cite,mleftright,enumerate,xparse}
\usepackage{xurl}
\usepackage[shortlabels]{enumitem}
\usepackage[english]{babel}
\usepackage[capitalise,sort,nameinlink]{cleveref} 

\crefname{enumi}{item}{items}

\usepackage{etoolbox}
\makeatletter
\apptocmd{\cref@getref}{\xdef\@lastusedlabel{#1}}{}{error}
\crefname{lemma}{Lemma\adddep{\loc}{\@lastusedlabel}}{Lemmas}
\crefname{prop}{Proposition\adddep{\loc}{\@lastusedlabel}}{Propositions}
\crefname{cor}{Corollary\adddep{\loc}{\@lastusedlabel}}{Corollaries}
\crefname{theorem}{Theorem\adddep{\loc}{\@lastusedlabel}}{Theorems}
\crefname{setting}{Setting\adddep{\loc}{\@lastusedlabel}}{Settings}
\makeatother
\crefname{equation}{}{}
\crefname{figure}{Figure}{Figures}

\newcommand{\SetStatementReference}[3]{%
	\crefname{#1}{#2}{#3}%
	\Crefname{#1}{#2}{#3}%
	\AddToHook{env/#1/begin}{%
		\crefalias{theorem}{#1}%
	}%
}

\SetStatementReference{theorem}{Theorem}{Theorems}
\SetStatementReference{lemma}{Lemma}{Lemmas}
\SetStatementReference{prop}{Proposition}{Propositions}
\SetStatementReference{cor}{Corollary}{Corollaries}
\SetStatementReference{definition}{Definition}{Definitions}
\SetStatementReference{setting}{Setting}{Settings}
\SetStatementReference{conjecture}{Conjecture}{Conjectures}
\SetStatementReference{example}{Example}{Examples}
\SetStatementReference{remark}{Remark}{Remarks}

\usepackage{tikz}
\usetikzlibrary{matrix,chains,positioning,decorations.pathreplacing,arrows}
\usetikzlibrary{shapes,arrows}
\tikzset{
	font={\fontsize{9pt}{12}\selectfont}}
\usepackage{adjustbox}

\tikzset{
	font={\fontsize{9pt}{12}\selectfont}}

\usepackage{environ}
\makeatletter
\newsavebox{\measure@tikzpicture}
\NewEnviron{scaletikzpicturetowidth}[1]{%
	\def\tikz@width{#1}%
	\begin{lrbox}{\measure@tikzpicture}%
		\BODY
	\end{lrbox}%
	\pgfmathparse{#1/\wd\measure@tikzpicture}%
	\BODY
}
\makeatother

\hypersetup{
	colorlinks,
	linkcolor={blue!80!black},
	citecolor={green},
	urlcolor={blue!80!black}
}

\theoremstyle{plain}
\newtheorem{theorem}{Theorem}[section]

\newtheorem{prop}[theorem]{Proposition}

\newtheorem{setting}[theorem]{Setting}

\theoremstyle{remark}

\theoremstyle{definition}
\newtheorem{definition}[theorem]{Definition}
\numberwithin{equation}{section}

\newcommand{\wdt}{B}
\newcommand{\dom}{R}

\newmdenv[
skipabove=\topsep,
skipbelow=\topsep,
backgroundcolor=gray!10,
linecolor=gray!80,
linewidth=1pt,
innertopmargin=4pt,
innerbottommargin=0.5\baselineskip,
innerleftmargin=0.5\baselineskip,
innerrightmargin=0.5\baselineskip
]{boxedthm}

\newmdenv[
skipabove=\topsep,
skipbelow=\topsep,
backgroundcolor=blue!10,
linecolor=gray!80,
linewidth=1pt,
innertopmargin=4pt,
innerbottommargin=0.5\baselineskip,
innerleftmargin=0.5\baselineskip,
innerrightmargin=0.5\baselineskip
]{boxeddef}

\DeclareMathAlphabet{\mathpzc}{OT1}{pzc}{m}{it}

\DeclareFontEncoding{LS1}{}{}
\DeclareFontSubstitution{LS1}{stix}{m}{n}
\DeclareMathAlphabet{\mathscr}{LS1}{stixscr}{m}{n}

\renewcommand{\P}{\mathbb{P}}

\newcommand{\R}{\mathbb{R}}
\newcommand{\N}{\mathbb{N}}

\newcommand{\Z}{\mathbb{Z}}

\renewcommand{\emptyset}{\varnothing}

\newcommand{\fa}[1]{\forall\, #1 \in}

\DeclarePairedDelimiter{\norm}{\lVert}{\rVert_{2}\cfadd{def:Euclidean_norm}}
\DeclarePairedDelimiter{\intronorm}{\lVert}{\rVert}
\DeclarePairedDelimiter{\supn}{\lVert}{\rVert_{\infty}\cfadd{def:matrixnorm}}
\DeclarePairedDelimiter{\vass}{\lvert}{\rvert}
\DeclarePairedDelimiter{\abs}{\lvert}{\rvert}
\DeclarePairedDelimiter{\pr}{(}{)}
\DeclarePairedDelimiter{\PR}{[}{]}
\DeclarePairedDelimiter{\pR}{\{}{\}}
\DeclarePairedDelimiter{\ceil}{\lceil}{\rceil\cfadd{def:ceiling}}
\DeclarePairedDelimiter{\floor}{\lfloor}{\rfloor\cfadd{def:ceiling}}

\makeatletter
\newcommand{\vast}{\bBigg@{4}}
\newcommand{\Vast}{\bBigg@{5.5}}
\makeatother

\NewDocumentCommand{\fabs}{sO{}m}{{\IfBooleanTF{#1}
		{\fabsaux{\left|}{\right|}{#3}}
		{\fabsaux{#2|}{#2|}{#3}}}
}
\makeatletter
\newcommand{\fabsaux}[3]{\mathpalette\fabsaux@i{{#1}{#2}{#3}}}
\newcommand{\fabsaux@i}[2]{\fabsaux@ii#1#2}
\newcommand{\fabsaux@ii}[4]{%
	\sbox\z@{$\m@th#1#2#4#3$}%
	\sbox\tw@{$\m@th\|$}%
	\mathopen{\hbox to\wd\tw@{\hss\vrule height \ht\z@ depth \dp\z@ width .3\wd\tw@\hss}}%
	\mkern-2mu #4 \mkern-2mu
	\mathclose{\hbox to\wd\tw@{\hss\vrule height \ht\z@ depth \dp\z@ width .3\wd\tw@\hss}}%
}
\makeatother

\newcommand{\vertiii}[1]{{\left\vert\kern-0.25ex\left\vert\kern-0.25ex\left\vert #1 
		\right\vert\kern-0.25ex\right\vert\kern-0.25ex\right\vert}}

\newcommand{\prb}[1]{\pr[\big]{ #1 }}
\newcommand{\PRb}[1]{\PR[\big]{ #1 }}
\newcommand{\pRb}[1]{\pR[\big]{ #1 }}
\newcommand{\prbbb}[1]{\pr[\bigg]{ #1 }}
\newcommand{\PRbbb}[1]{\PR[\bigg]{ #1 }}

\newcommand{\prbb}[1]{\pr[\Big]{ #1 }}
\newcommand{\PRbb}[1]{\PR[\Big]{ #1 }}

\newcommand{\Indfct}[1]{\mathbbm{1}_{#1}}

\newcommand{\forallDist}{\forall\,}

\newcommand{\radius}{R}

\newcommand{\Cost}{\cfadd{def:cost}{\operatorname{Cost}}}
\newcommand{\COST}{\cfadd{def:cost} {\operatorname{Cost}}}
\newcommand{\cost}{\cfadd{def:cost} {\operatorname{cost}}}
\newcommand{\ent}{\cfadd{def:ent}\mathrm{ent}}
\newcommand{\ENT}{\cfadd{def:ent}\mathrm{Ent}}

\newcommand{\edgy}{\mathfrak{s}}

\newcommand{\normmm}[1]{{\left\vert\kern-0.25ex\left\vert\kern-0.25ex\left\vert #1 
		\right\vert\kern-0.25ex\right\vert\kern-0.25ex\right\vert}} 

\newcommand{\qandq}{\qquad\text{and}\qquad}

\newcommand{\id}{\operatorname{id}}

\newcommand{\param}{\mathcal{P}}

\newcommand{\realisation}{\cfadd{def:ANNrealization}\mathcal{R}}
\newcommand{\size}{\cfadd{def:size}\mathcal{S}}

\newcommand{\indicator}[1]{\mathbbm{1}_{#1}}

\newcommand{\ANNSET}{\cfadd{def:ANN}\mathbf{N}}

\newcommand{\paramANN}{\cfadd{def:ANN}\mathcal{P}}

\newcommand{\paramnotzero}{\mathfrak{P}\cfadd{def:paramnonzero}}

\newcommand{\lengthANN}{\cfadd{def:ANN}\mathcal{L}}
\newcommand{\weight}[2]{\cfadd{def:ANN}\mathcal{W}_{#1,#2}}
\newcommand{\bias}[2]{\cfadd{def:ANN}\mathcal{B}_{#1,#2}}

\newcommand{\hidlengthANN}{\cfadd{def:ANN}\mathcal{H}}
\newcommand{\inDimANN}{\cfadd{def:ANN}\mathcal{I}}
\newcommand{\compANN}[2]{{#1 \bullet #2\cfadd{def:ANNcomposition}}}

\newcommand{\outDimANN}{\cfadd{def:ANN}\mathcal{O}}

\newcommand{\ReLUidANN}[1]{\cfadd{def:ReLU_identity}\mathbb{I}_{#1}}

\newcommand{\dims}{\cfadd{def:ANN}\mathcal{D}}
\newcommand{\singledims}{\cfadd{def:ANN}\mathbb{D}}
\newcommand{\insize}{\cfadd{def:size}\mathbb{S}_0}
\newcommand{\outsize}{\cfadd{def:size}\mathbb{S}_1}
\newcommand{\rsize}[1]{\cfadd{def:size}\mathbb{S}_{#1}}

\newcommand{\parallelizationSpecial}{\cfadd{def:simpleParallelization}\mathbf{P}}

\newcommand{\RELU}{\mathfrak{R}\cfadd{def:RELU}}
\newcommand{\Hyp}{\operatorname{Hyp}\cfadd{def:hyper2}}
\newcommand{\Lin}{\mathfrak{L}\cfadd{def:Lin}}

\newcommand{\convex}[1]{\mathfrak{C}\cfadd{def:convex}}

\renewenvironment{pmatrix}{\mleft(\begin{matrix}}{\end{matrix}\mright)}

\DeclareMathOperator*{\ssum}{\textstyle\sum}

\DeclareMathOperator{\ssssum}{\textstyle\sum}

\DeclareMathOperator*{\sprod}{\textstyle\prod}

\newcommand{\is}{\curvearrowleft}

\newcommand{\eps}{\varepsilon}

\ExplSyntaxOn

\seq_new:N \g_cflist_loaded
\seq_new:N \g_cflist_pending

\NewDocumentCommand{\cfadd}{ m }
{
	\seq_if_in:NnF \g_cflist_loaded { #1 } {
		\seq_if_in:NnF \g_cflist_pending { #1 } {
			\seq_gput_right:Nn \g_cflist_pending { #1 }
		}
	}
}

\NewDocumentCommand{\cfconsiderloaded}{ m }{
	\seq_gput_right:Nn \g_cflist_loaded {#1}
}

\NewDocumentCommand{\cfremove}{ m }
{
	\seq_gremove_all:Nn \g_cflist_pending { #1 }
}

\NewDocumentCommand{\cfload}{ o }
{
	\seq_if_empty:NTF \g_cflist_pending {\unskip} {
		(cf.\ \cref{\seq_use:Nn \g_cflist_pending {,}})\IfValueTF{#1}{#1~}{\unskip}
		\seq_gconcat:NNN \g_cflist_loaded \g_cflist_loaded \g_cflist_pending
		\seq_gclear:N \g_cflist_pending
	}
}

\NewDocumentCommand{\cfclear} {} {
	\seq_gclear:N \g_cflist_loaded
	\seq_gclear:N \g_cflist_pending
}

\NewDocumentCommand{\cfout}{ o }
{
	\seq_if_empty:NTF \g_cflist_pending {\unskip} {
		(cf.\ \cref{\seq_use:Nn \g_cflist_pending {,}})\IfValueTF{#1}{#1~}{\unskip}
		\seq_gclear:N \g_cflist_pending
	}
}

\NewDocumentCommand{\ifnocf} { m } {
	\seq_if_empty:NT \g_cflist_pending { #1 }
}

\ExplSyntaxOff

\ExplSyntaxOn

\bool_new:N \g_noteobserve

\NewDocumentCommand{\setnote}{}{
	\bool_gset_true:N \g_noteobserve
}

\NewDocumentCommand{\setobserve}{}{
	\bool_gset_false:N \g_noteobserve
}

\NewDocumentCommand{\nobs}{ o }{
	\IfValueT{#1}{
		\str_if_eq:noTF {note} {#1} {
			\bool_gset_true:N \g_noteobserve
		} {
			\str_if_eq:noTF {Note} {#1} {
				\bool_gset_true:N \g_noteobserve
			} {
				\bool_gset_false:N \g_noteobserve
			}
		}
	}
	\bool_if:nTF { \g_noteobserve } {
		\bool_gset_false:N \g_noteobserve 
		note
	} {
		\bool_gset_true:N \g_noteobserve 
		observe
	}
	\IfValueF{#1}{~}
}

\NewDocumentCommand{\Nobs}{ o }{
	\IfValueT{#1}{
		\str_if_eq:noTF {note} {#1} {
			\bool_gset_true:N \g_noteobserve
		} {
			\str_if_eq:noTF {Note} {#1} {
				\bool_gset_true:N \g_noteobserve
			} {
				\bool_gset_false:N \g_noteobserve
			}
		}
	}
	\bool_if:nTF { \g_noteobserve } {
		\bool_gset_false:N \g_noteobserve 
		Note
	} {
		\bool_gset_true:N \g_noteobserve 
		Observe
	}
	\IfValueF{#1}{~}
}

\bool_new:N \g_hencetherefore

\NewDocumentCommand{\hence}{ o }{
	\IfValueT{#1}{
		\str_if_eq:noTF {hence} {#1} {
			\bool_gset_true:N \g_hencetherefore
		} {
			\str_if_eq:noTF {Hence} {#1} {
				\bool_gset_true:N \g_hencetherefore
			} {
				\bool_gset_false:N \g_hencetherefore
			}
		}
	}
	\bool_if:nTF { \g_hencetherefore } {
		\bool_gset_false:N \g_hencetherefore 
		hence
	} {
		\bool_gset_true:N \g_hencetherefore 
		therefore
	}
	\IfValueF{#1}{~}
}

\NewDocumentCommand{\Hence}{ o }{
	\IfValueT{#1}{
		\str_if_eq:noTF {hence} {#1} {
			\bool_gset_true:N \g_hencetherefore
		} {
			\str_if_eq:noTF {Hence} {#1} {
				\bool_gset_true:N \g_hencetherefore
			} {
				\bool_gset_false:N \g_hencetherefore
			}
		}
	}
	\bool_if:nTF { \g_hencetherefore } {
		\bool_gset_false:N \g_hencetherefore 
		Hence
	} {
		\bool_gset_true:N \g_hencetherefore 
		Therefore
	}
	\IfValueF{#1}{~}
}

\int_new:N \g_furthermore

\NewDocumentCommand{\Moreover}{ o o }{
	\IfValueT{#1}{
		\str_case:nn {#1} {
			{Furthermore} {\int_set:Nn {\g_furthermore} {0}}
			{Moreover} {\int_set:Nn {\g_furthermore} {1}}
			{In~addition} {\int_set:Nn {\g_furthermore} {2}}
			{note} {\bool_gset_true:N \g_noteobserve}
			{observe} {\bool_gset_false:N \g_noteobserve}
		}
		\IfValueT{#2}{
			\str_case:nn {#2} {
				{Furthermore} {\int_set:Nn {\g_furthermore} {0}}
				{Moreover} {\int_set:Nn {\g_furthermore} {1}}
				{In~addition} {\int_set:Nn {\g_furthermore} {2}}
				{note} {\bool_gset_true:N \g_noteobserve}
				{observe} {\bool_gset_false:N \g_noteobserve}
			}
		}
	}
	\int_case:nn { \int_mod:nn {\g_furthermore} {3} } {
		{ 0 } { Furthermore,~\nobs that}
		{ 1 } { Moreover,~\nobs that}
		{ 2 } { In~addition,~\nobs that}
	}
	\int_incr:N \g_furthermore
	\IfValueF{#1}{~}
}

\ExplSyntaxOff

\ExplSyntaxOn

\bool_new:N \g_forexample

\NewDocumentCommand{\eg}{ o }{
	\IfValueT{#1}{
		\str_if_eq:noTF {fe} {#1} {
			\bool_gset_true:N \g_forexample
		} {\bool_gset_false:N \g_forexample}
	}
	\bool_if:nTF { \g_forexample } {
		\bool_gset_false:N \g_forexample
		for~example
	}{
		\bool_gset_true:N \g_forexample
		for~instance
	}
}

\ExplSyntaxOn
\global\def\loc{dummy}

\NewDocumentEnvironment {athm} {m m o} {%
		\IfValueT{#3}{\begin{#1}[#3]}
			\IfValueF{#3}{\begin{#1}}
				\label{#2}\global\def\loc{#2}%
			}{%
			\end{#1}%
		}
		
		\NewDocumentEnvironment{aproof} {} {%
			\begin{proof}[Proof~of~\cref{\loc}]%
			}{%
				\global\def\loc{dummy}\end{proof}%
		}
		
		\ExplSyntaxOff

		\newcommand{\lref}[1]{\cref{\loc.#1}}
		
		\newcommand{\llabel}[1]{\label{\loc.#1}}
		
		\newcommand{\finishproofthus}{The proof of \cref{\loc} is thus complete.}

		\ExplSyntaxOn
		
		\seq_new:N \g_deplist
		\seq_new:N \l_done
		\tl_new:N \l_first_tl
		\tl_new:N \l_second_tl
		\int_new:N \l_count_int
		\iow_new:N \l_logfile
		
		\NewDocumentCommand{\adddep}{ m m }
		{
			\seq_gput_right:Nx \g_deplist { #1 }
			\seq_gput_right:Nx \g_deplist { #2 }
		}
		
		\NewDocumentCommand{\listdeps}{}
		{
			\int_set:Nn \l_count_int {\seq_count:N {\g_deplist}}
			\seq_clear:N \l_done
			\int_while_do:nNnn {\l_count_int} > {0} {
				\iow_open:Nn \l_logfile {dependencies.dot}
				\seq_gpop_left:NN \g_deplist \l_first_tl
				\seq_gpop_left:NN \g_deplist \l_second_tl
				\cs_if_eq:NNTF { \l_first_tl } { \l_second_tl } {} {
					\tl_if_eq:NnTF { \l_first_tl } { dummy } {} {
						\seq_if_in:NxTF { \l_done } {\l_first_tl \l_second_tl} {} {
							\iow_now:Nn \l_logfile {\cref*{\tl_use:N \l_second_tl}}
							\iow_now:Nn \l_logfile {->}
							\iow_now:Nn \l_logfile {\cref*{\tl_use:N \l_first_tl}}
							{\cref{\tl_use:N \l_second_tl}}~
							~$\to$~
							{\cref{\tl_use:N \l_first_tl}}\\
							\seq_put_left:Nx {\l_done} {\l_first_tl \l_second_tl}
						}
					}
				}
				\int_decr:N {\l_count_int}
				\int_decr:N {\l_count_int}
			}
			\iow_close:N \l_logfile
		}
		
		\ExplSyntaxOff

\ExplSyntaxOn
\prop_new:N \l__verbs
\prop_put:Nnn \l__verbs {show} {shows}
\prop_put:Nnn \l__verbs {imply} {implies}
\prop_put:Nnn \l__verbs {demonstrate} {demonstrates}
\prop_put:Nnn \l__verbs {prove} {proves}
\prop_put:Nnn \l__verbs {establish} {establishes}
\prop_put:Nnn \l__verbs {ensure} {ensures}
\prop_put:Nnn \l__verbs {assure} {assures}

\seq_const_from_clist:Nn \g_prove_mru {
	establish,
	demonstrate,
	prove,
	show,
	imply,
	ensure
}

\tl_new:N \g_wordtmp
\seq_new:N \l_mytmps

\cs_generate_variant:Nn \str_if_in:nnTF { nVTF }

\NewDocumentCommand{\prove}{ o }{
	\IfValueTF{#1}{
		\seq_clear:N \l_mytmps
		\seq_map_inline:Nn \g_prove_mru {
			\str_if_eq:nnTF {##1} {ensure} {
				\str_set:Nn \l_temps {n}
			} {
				\str_set:Nx \l_temps {\str_head_ignore_spaces:n {##1}}
			}
			\str_if_in:nVTF {#1} \l_temps {
				\seq_put_right:Nn \l_mytmps {##1}
			} { }
		}
		\seq_get_right:NN \l_mytmps \g_wordtmp
	} {
		\seq_get_right:NN \g_prove_mru \g_wordtmp
	}
	\tl_use:N \g_wordtmp
	\seq_gput_left:NV \g_prove_mru \g_wordtmp
	\seq_gremove_duplicates:N \g_prove_mru
}

\NewDocumentCommand{\proves}{ o }{
	\IfValueTF{#1}{
		\seq_clear:N \l_mytmps
		\seq_map_inline:Nn \g_prove_mru {
			\str_if_eq:nnTF {##1} {ensure} {
				\str_set:Nn \l_temps {n}
			} {
				\str_set:Nx \l_temps {\str_head_ignore_spaces:n {##1}}
			}
			\str_if_in:nVTF {#1} \l_temps {
				\seq_put_right:Nn \l_mytmps {##1}
			} { }
		}
		\seq_get_right:NN \l_mytmps \g_wordtmp
	} {
		\seq_get_right:NN \g_prove_mru \g_wordtmp
	}
	\str_set:NV \l_tmpa_str \g_wordtmp
	\prop_get:NVN \l__verbs \l_tmpa_str \l_tmpa_tl
	\tl_use:N \l_tmpa_tl
	\seq_gput_left:NV \g_prove_mru \g_wordtmp
	\seq_gremove_duplicates:N \g_prove_mru
}

\ExplSyntaxOff

\makeatletter
\ExplSyntaxOn

\seq_new:N \g_abbrs
\prop_new:N \g_abbr_counts
\tl_new:N \l_abbr_count_tl

\NewDocumentCommand{\abbr}{m m O{#1} m m O{#4} m}{
	\expandafter\newcommand\csname#3\endcsname[1][]{
		\seq_if_in:NnTF \g_abbrs {#1} {
			\prop_get:NnN \g_abbr_counts {#1} \l_abbr_count_tl
			\prop_gput:Nnx \g_abbr_counts {#1} {\int_eval:n {\l_abbr_count_tl + 1}}
			\hyperref[#1]{#7}
		} {
			\seq_gput_left:Nn \g_abbrs {#1}
			\prop_gput:Nnn \g_abbr_counts {#1} {1}
			\expandafter\gdef\csname#1@def\endcsname{#2}
			\phantomsection\label{#1}
			\str_if_eq:nnTF{##1}{}{\emph{#2}}{##1}~(\hyperref[#1]{#7})
		}
	}
	\expandafter\newcommand\csname#6\endcsname[1][]{
		\seq_if_in:NnTF \g_abbrs {#1} {
			\prop_get:NnN \g_abbr_counts {#1} \l_abbr_count_tl
			\prop_gput:Nnx \g_abbr_counts {#1} {\int_eval:n {\l_abbr_count_tl + 1}}
			\hyperref[#1]{#4}
		} {
			\expandafter\gdef\csname#1@def\endcsname{#5}
			\seq_gput_left:Nn \g_abbrs {#1}
			\prop_gput:Nnn \g_abbr_counts {#1} {1}
			\phantomsection\label{#1}
			\str_if_eq:nnTF{##1}{}{\emph{#5}}{##1}~(\hyperref[#1]{#4})
		}
	}
}

\ExplSyntaxOff
\makeatother

\abbr{ANN}{Artificial neural network}{ANNs}{artificial neural networks}{ANN}
\abbr{DNN}{deep artificial neural network}{DNNs}{deep artificial neural networks}{DNN}
\abbr{SGD}{Stochastic gradient descent}{SGDs}{Stochastic gradient descent}{SGD}
\abbr{GD}{gradient descent}{GDs}{gradient descent}{GD}
\abbr{RMSprop}{root mean square propagation}{RMSprops}{root mean square propagation}{RMSprop}
\abbr{Adam}{adaptive moment estimation}{Adams}{Adaptive moment estimation}{Adam}
\abbr{iid}{independent and identically distributed}{i.i.d}{independent and identically distributed}{i.i.d.}
\abbr{DL}{Deep learning}{DLs}{Deep learning}{DL}
\abbr{AI}{artificial intelligence}{AIs}{artificial intelligence}{AI}
\abbr{LLM}{large language model}{LLMs}{large language models}{LLM}
\abbr{PDE}{partial differential equation}{PDEs}{partial differential equations}{PDE}
\abbr{as}{almost surely}{ass}{almost surely}{a.s.}
\abbr{pas}{$\P$-almost surely}{pass}{almost surely}{$\P$-a.s.}
\abbr{OP}{optimization problem}{OPs}{optimization problems}{OP}
\abbr{Adagrad}{adaptive gradient}{Adagrads}{Adaptive gradients}{Adagrad}
\abbr{MUON}{momentum orthogonalized by Newton-Schulz}{MUONs}{momentum orthogonalized by Newton-Schulzs}{MUON}
\abbr{ReLU}{rectified linear unit}{ReLUs}{rectified linear units}{ReLU}
\abbr{BERT}{Bidirectional Encoder Representations from
	Transformers}{BERTs}{Bidirectional Encoder Representations from
	Transformers}{BERT}
\abbr{IBC}{information based complexity}{IBCs}{information based complexity}{IBC}
\abbr{GLUE}{General Language Understanding Evaluation}{GLUEs}{General Language Understanding Evaluation}{GLUE}
\abbr{MultiNLI}{Multi-Genre Natural Language Inference}{MultiNLIs}{Multi-Genre Natural Language Inferences}{MultiNLI}
\abbr{ILSVRC}{ImageNet Large Scale Visual Recognition Challenge}{ILSVRCs}{ImageNet Large Scale Visual Recognition Challenges}{ILSVRC}
\abbr{MSCOCO}{Microsoft Common Objects in Context}
{MS~COCO datasets}{Microsoft Common Objects in Context datasets}
[MSCOCOdatasets]{MS~COCO}

\makeatletter

\renewcommand{\setminus}{\backslash}
\makeatother

\title{
	The necessity of depth for artificial neural networks
	\\ to approximate certain classes of smooth 
	and bounded\\  functions without the curse of dimensionality}

\author{Lukas Gonon$^{1}$, Robin Graeber$^{2,3}$, and Arnulf Jentzen$^{4,5}$
	\bigskip
	\\
	\small{$^1$ Faculty of Natural Sciences,
		Department of Mathematics,}
	\vspace{-0.1cm}\\
	\small{Imperial College London, United Kingdom}
	\vspace{-0.1cm}\\
	\small{e-mail: \texttt{l.gonon@imperial.ac.uk}}
	\smallskip
	\\
    \small{$^2$ School of Data Science,	The Chinese University of Hong Kong, Shenzhen}
    \vspace{-0.1cm}\\
    \small{(CUHK-Shenzhen), China, e-mail:
    	\texttt{223040041@link.cuhk.edu.cn}}
    \smallskip
    \\
    \small{$^3$ Applied Mathematics: Institute for Analysis and Numerics,}
    \vspace{-0.1cm}\\
    \small{Faculty of Mathematics and Computer Science,}
    \vspace{-0.1cm}\\
    \small{University of M\"unster, Germany}
    \vspace{-0.1cm}\\
    \small{e-mail: \texttt{r\_grae02@uni-muenster.de}}
    \smallskip
    \\
	\small{$^4$ School of Data Science and Shenzhen Research Institute of Big Data,}
	\vspace{-0.1cm}\\
	\small{The Chinese University of Hong Kong, Shenzhen (CUHK-Shenzhen), China}
	\vspace{-0.1cm}\\
	\small{e-mail: \texttt{ajentzen@cuhk.edu.cn}}
	\smallskip
	\\
	\small{$^5$ Applied Mathematics: Institute for Analysis and Numerics,}
	\vspace{-0.1cm}\\
	\small{Faculty of Mathematics and Computer Science,}
	\vspace{-0.1cm}\\
	\small{University of M{\"u}nster, Germany}
	\vspace{-0.1cm}\\
	\small{e-mail: \texttt{ajentzen@uni-muenster.de}}
	\smallskip
	\\
}

\begin{document}

\maketitle
\newpage
\begin{abstract}
In this article we study high-dimensional approximation capacities of shallow and deep artificial neural networks (ANNs) with the rectified linear unit (ReLU) activation.
In particular, it is a key contribution of this work to reveal that for all $a\in\R$, $b\in(a,\infty)$ with $\max\{\vass{a},\vass{b}\}>1$ we have that the functions 
$[a,b]^d\ni x=(x_1,\ldots,x_d)\mapsto\prod_{i=1}^d x_i\in\R$
for $d\in\N=\{1,2,3,\ldots\}$ 
as well as the functions 
$[a,b]^d\ni x=(x_1,\ldots,x_d)\mapsto\sin(\prod_{i=1}^d x_i)\in\R$
for $d\in\N$ 
can neither be approximated without the \textit{curse of dimensionality} by means of shallow ANNs nor insufficiently deep ANNs with ReLU activation but can be approximated without the \emph{curse of dimensionality} by sufficiently deep ANNs with ReLU activation.
More specifically, we prove that in the case of shallow ReLU ANNs or deep ReLU ANNs with
 a fixed number of hidden layers we have that the number of ANN parameters \emph{must grow at least exponentially} in the dimension $d\in\N$
while in the case of deep ReLU ANNs with the number of hidden layers growing in the dimension $d\in\N$ we have that the number of scalar real parameters of the approximating ANNs \emph{grows at most polynomially} in both the inverse of the prescribed approximation accuracy $\eps>0$ and the dimension $d\in\N$,
even if the absolute values of the ANN parameters are assumed to be uniformly bounded by $1$.
We thus show that the product functions and the sine of the product functions are \emph{polynomially tractable} approximation problems among the approximating class of deep ReLU ANNs with the number of hidden layers being allowed to grow in the dimension $d\in\N$.
We establish the above outlined statement not only for the product functions and the sine of the product functions but also for other classes of target functions, in particular, for classes of uniformly globally bounded $C^{\infty}$-functions with compact support.
Roughly speaking, in this work we lay open that simple approximation problems such as approximating the sine or cosine of products cannot be solved in standard implementation frameworks by shallow or insufficiently deep ANNs with ReLU activation in polynomial time, 
but can be approximated by sufficiently deep ReLU ANNs with the number of bits on a computer storing the approximating deep ANNs growing at most polynomially in both the inverse of the prescribed approximation accuracy and the dimension of the problem.
\end{abstract}

\newpage
\tableofcontents

\newpage
\section{Introduction}
\label{sect:intro}
\newcommand{\f}{\mathscr{f}}
\newcommand{\g}{\mathscr{g}}
\newcommand{\h}{\mathscr{h}}
\newcommand{\W}{\mathfrak{w}}
\newcommand{\tgt}{g}
\newcommand{\shft}{\lambda}
\newcommand{\strch}{A}
\newcommand{\substrch}{B}

\newcommand{\dimsbene}{Beneventano et al.~\cite[Lemma~2.17]{beneventano21} (see, \eg, also Grohs et al.~\cite[Proposition~2.6]{GrohsHornungJentzen2019})}
\newcommand{\benecorlip}{Beneventano et al.~\cite[Lemma~3.22]{beneventano21}}
\newcommand{\benecorcomp}{Beneventano et al.~\cite[Lemma~6.5]{beneventano21}}
\newcommand{\benetwoprod}{Beneventano et al.~\cite[Lemma~6.3]{beneventano21}}
\newcommand{\benepairprod}{Beneventano et al.~\cite[Lemma~6.4]{beneventano21}}
\newcommand{\benepotprod}{Beneventano et al.~\cite[Lemma~6.7]{beneventano21}}
\newcommand{\beneprod}{Beneventano et al.~\cite[Proposition~6.8]{beneventano21}}
\ANN\ approximations are ubiquitous in our digital world and appear in diverse areas, whether 
language processing (cf., \eg, Devlin et al.~\cite{BERT}), 
image classification (cf., \eg, Chen et al.~\cite{ChenEtAl19}), 
predictive models for cancer diagnosis (cf., \eg, Sidey-Gibbons  \& Sidey-Gibbons \cite{sidey2019machine}), 
risk assessment (cf., \eg, Paltrinieri et al.~\cite{PALTRINIERI2019475}), 
or 
biomedical imaging and signal processing (cf., \eg, Min et al.~\cite{min2017deep}).
In such learning problems, \ANNs\ are employed to approximate the target function which, roughly speaking, describes the best relationship of the input data to the output data in the considered learning problem.
There are a large number of numerical simulation results which indicate that \ANNs\ are comparatively well suited to approximate the target functions in such learning problems.
The success of \ANN\ approximations becomes even more remarkable if one takes into account that the target functions in the above-mentioned learning problems are usually extremely high-dimensional functions.

For example, in an object recognition problem, the 10000-dimensional unit cube $[0,1]^{10000}$ can be employed to represent the set of all grey-scale images with $100\times 100$ pixels and the target function $f\colon[0,1]^{10000}\to[0,1]$ of the considered learning problem is then a function from the 10000-dimensional unit cube $[0,1]^{10000}$ to the interval $[0,1]$ modelling for every image the probability that it contains a certain object, say, a car.
%
Loosely speaking, it is impossible to approximate such target functions by classical deterministic approximation methods (such as finite differences or finite elements in the context of \PDEs; cf., \eg, Jovanovi\'{c} \& Süli \cite{jovanovic2014analysis} and Tadmor \cite{TadmorEitan12}), as such classical approximations suffer from the \textit{curse of dimensionality} in the sense that the number of parameters to describe such approximations grows at least exponentially in the input dimension (cf., \eg, Bellman \cite{Bellman1957}, 
cf., \eg, Hinrichs et al.~\cite{CoDHinrichs},
Novak \& Wo\'zniakowski \cite[Chapter~1]{NovakWozniakowski2008}, and
Novak \& Wo\'zniakowski \cite[Section~9.7]{NovakWozniakowski2010}).

In many cases numerical simulations for \ANNs\ suggest that \ANN\ approximations are capable of approximating such extremely high-dimensional input-output data relationships and, in particular, numerical simulations suggest that \ANN\ approximations for such problems seem to overcome the curse of dimensionality in the sense that the amount of real numbers used to describe those approximations seems to grow at most polynomially in the reciprocal $\eps^{-1}$ of the approximation precision $\eps>0$ and the dimension $d\in\N$ of the domain of the target function of the considered learning problem.
In the \IBC\ literature this polynomial growth estimate in both the reciprocal of the approximation precision and the problem dimension is also often referred to as \textit{polynomial tractability} (cf., \eg,  Novak \& Wo\'zniakowski \cite[Chapter~1]{NovakWozniakowski2008}
 and
Novak \& Wo\'zniakowski \cite[Section~9.7]{NovakWozniakowski2010}).
In the most simple form, an \ANN\ describes a function (the so-called realization function of the \ANN) which is given by iterated compositions of affine linear functions 
(with the entries of the multiplicative matrix in the affine linear function and the entries of the additive vector in the affine linear function described through a parameter vector of the \ANN)
and certain fixed nonlinear functions (referred to as activation functions).
\color{black}
Roughly speaking, the result of such iterated composition after each nonlinear function represents a hidden layer of the \ANN\ and
\ANNs\ with one (or none) hidden layers are referred to as shallow \ANNs\
while \ANNs\ with two or more hidden layers are called deep \ANNs\ with the number of hidden layers representing the depth of the \ANN\ (see also \cref{figure1} below for a graphical illustration of the architecture of an \ANN).

Successful implementations in the above named learning problems usually employ deep \ANNs\ with a large number of hidden layers.
In particular, the modern language processing framework \BERT\ introduced in Devlin et al.~\cite{BERT} set new benchmarks in several natural language processing tasks like \GLUE\ (see Wang et al.~\cite{wang2018glue}) and \MultiNLI\ (see Williams et al.~\cite{Williams18}) using \ANNs\ with 12 and 24 Transformer encoder layers. 
The methods described in He et al.~\cite{He_2016_CVPR} won several image recognition competitions in 2015 like \ILSVRC\ (see Russakovsky et al.~\cite{ILSVRC15}) and \MSCOCO\
(see Lin et al.~\cite{Tsung-Yiet2014}) 
by successfully implementing and training 152-layer residual \ANNs\ and
in 2017 the DenseNets in Huang et al.~\cite{Huang_2017_CVPR} consisting of 38 to 248 hidden layers outperformed state of the art techniques in visual object recognition.

 The large number of numerical simulations in the above named learning problems also indicate that shallow or insufficiently deep \ANNs\ might not be able to approximate the considered high-dimensional target functions without the curse of dimensionality.
 Taking this into account, a natural topic of research is to develop a mathematical theory which rigorously explains why (and for which classes of target functions) deep \ANNs\ seem to be capable of overcoming the curse of dimensionality while shallow or insufficiently deep \ANNs\ seem to fail to do so in many relevant learning problems.
 In the scientific literature there are also a few mathematical research articles which contribute or have strong connections to this area of research.

In particular, 
we refer to Daniely \cite{daniely2017depth} for a class of functions which can be approximated in $L^2$ with respect to the uniform probability measure on a product of spheres without the curse of dimensionality by \textit{\ReLU\ \ANNs\ with two hidden layers} but not by \textit{shallow \ReLU\ \ANNs} under suitable exponential bounds on the absolute values of the parameters of the shallow \ANNs,
we refer to Chui et al.~\cite{chui2019deep} for approximation
of radial functions by \ANNs\ with suitable smooth bounded
sigmoidal activation functions, where the required number of parameters is bounded polynomially
in the reciprocal precision with a dimension-independent exponent
for four hidden layers, while shallow \ANNs\ require exponents
growing with the input dimension,
we refer to Eldan \& Shamir \cite{eldan2016power} for a sequence of two hidden layer \ANNs\ (with the number of parameters growing at most polynomially in the input dimension) which cannot be approximated by shallow \ANNs\ without the curse of dimensionality (teacher-student setup; cf., \eg, Saad \& Solla \cite{SaadSolla95} and Riegler \& Biehl \cite{RieglerBiehl95}),
and we refer to Venturi et al.~\cite{Venturi21} for a family of oscillating complex-valued functions which can be approximated in the $L^2$-sense with respect to a certain absolutely continuous probability measure without the curse of dimensionality by \ANNs\ with two hidden layers but not by shallow \ANNs.

We refer to Telgarsky \cite{Telgarsky15,Telgarsky16} and Yu et al.~\cite{YuAnnan21} for suitable families of \textit{deep} \ANNs\ indexed over an external parameter with at most polynomially many \ANN\ parameters (with respect to the external parameter) which can only be approximated by \textit{insufficiently deep} \ANNs\ if the number of \ANN\ parameters in the insufficiently deep \ANNs\ grows at least exponentially in the external parameter (teacher-student setup; cf., \eg, Saad \& Solla \cite{SaadSolla95} and Riegler \& Biehl \cite{RieglerBiehl95})
and
we refer to Chen \& Wu \cite{ChenWu2019} for families of functions whose Fourier transformations can be approximated on cubes by deep \ANNs\ with the number of parameters growing at most logarithmically in the length of the edges of the cubes but which can only be approximated on cubes by insufficiently deep \ANNs\ with the number of parameters growing at least linearly in the length of the edges of the cubes.
We refer to Safran \& Shamir \cite{safran17a} for
$L^2([0,1]^d)$-approximation of nonaffine twice continuously
differentiable functions by \ReLU\ \ANNs, where fixed depth
requires the number of parameters to grow at least polynomially
in the reciprocal precision, while uniform approximation
with polylogarithmic depth and parameter counts is possible
for target functions admitting uniform approximations using
polylogarithmically many weighted additions and multiplications
with polynomially bounded intermediate values.
Grohs et al.~\cite{GrohsIbrgimovJentzen2021} establish
a depth-separation result in the input dimension for
approximation in $L^2$ with respect to the standard
normal distribution under polynomial bounds on the
absolute values of the \ANN\ parameters.
The respective target functions can be approximated by deep \ANNs\ 
with the number of parameters growing at most cubically
in the input dimension, while approximation by
insufficiently deep \ANNs\ suffers from the curse
of dimensionality.

In this article we investigate the necessity of increasing
depth for uniform approximation on cubes.
More specifically, we construct families of uniformly globally
bounded and infinitely often differentiable target functions
for which every fixed bound on the number of hidden layers
leads to exponential lower bounds in the input dimension,
even at a fixed approximation precision and even if the
parameters of the approximating \ReLU\ \ANNs\
are arbitrarily large.
The same target functions can be approximated with polynomially
many parameters in both the input dimension and the reciprocal
of the prescribed approximation precision by sufficiently deep
\ReLU\ \ANNs\ whose parameter sizes are uniformly bounded by $1$.
In \cref{Thm:gen}, the target functions are additionally
chosen to have compact support.
The sine of the product functions provides an explicit
family of target functions with the depth-separation
properties established in \cref{Thm:sin}.
%
This is particularly relevant as the number and the size of the real valued parameters in the approximating \ANN\ are strongly linked to the amount of memory needed on a computer to store a solution of the approximation problem. 
In particular, we show 
\begin{enumerate}[(I)]
	\item 
that
 the number of bits needed to approximate the considered globally bounded and infinitely often differentiable functions via \ANNs\ with the depth being bounded by any arbitrarily large fixed number \emph{must grow at least exponentially in the dimension of the function}
and 
\item
that
it is possible to approximate the considered globally bounded and infinitely often differentiable functions via sufficiently deep \ANNs\ with the number of bits needed to store the approximating \ANNs\ \emph{growing at most polynomially in both the dimension of the function and the reciprocal of the prescribed approximation accuracy};
\end{enumerate}
see \cref{cor:bit:both:bounds} in \cref{subsection:entropy_results} below for details.
%
In this sense, we show that 
the amount of memory needed to store such shallow or insufficiently deep \ANN\ approximations \emph{must grow at least exponentially in the problem dimension} and
that the amount of memory needed to store sufficiently deep \ANN\ approximations \emph{grows at most polynomially in the problem dimension and the reciprocal of the prescribed approximation accuracy} (see \cref{rmk:1} in \cref{subsection:entropy_results} below).

\subsection{Mathematical framework and approximation complexity}
\label{intro:framework}
To illustrate the findings of this work in more details, we now depict in this introductory section three representative key \ANN\ approximation results of this article in a precise and self-contained way (see \cref{Thm:gen}, \cref{Thm:sin}, and \cref{Thm:prod} below).
Each of these three \ANN\ approximation results employs the mathematical description of standard fully-connected feedforward \ANNs\ with the \ReLU\ activation which is the subject of the following mathematical framework; see \eqref{setting:eq1}, \eqref{setting:eq2}, and \eqref{setting:eq3} in \cref{setting} below.
We also refer to \cref{figure1} for a graphical illustration of the architecture of the \ANNs\ formulated in \cref{setting}.

\begin{samepage}
\begin{setting}
	\label{setting}
	\begin{boxeddef}
	Let $\mathfrak{R}\colon(\cup_{k\in\N}\R^k)\to(\cup_{k\in\N}\R^k)$ and $\fabs{\cdot}\colon\pr
	{\cup_{k,l\in\N}\pr{\R^{k\times l}\times\R^k}}\to\R$
	satisfy for all $k,l\in\N$, $x=(x_1,\ldots,x_k)\in\R^k$, $W=(W_{i,j})_{(i,j)\in\{1,\ldots,k\}\times\{1,\ldots,l\}}\in\R^{k\times l}$ that 
	\begin{equation}
		\label{setting:eq1}
		\textstyle
		\mathfrak{R}(x)=(\max\{x_1,0\},\ldots,\max\{x_k,0\}) 
		\,\,\text{and}\,\,
		\fabs{(W,x)}=\max_{1\leq i\leq k}\max_{1\leq j \leq l}\max\{\vass{W_{i,j}},\vass{x_i}\},
	\end{equation} let $\ANNSET=\cup_{L\in\N}\cup_{l_0,l_1,\ldots,l_L\in\N}(\times_{k=1}^L(\R^{l_k\times l_{k-1}}\times\R^{l_k}))$, 
	let $\realisation\colon\ANNSET\to(\cup_{k,l\in\N} C(\R^k,\R^l))$, $\lengthANN\colon\ANNSET\to\N$, $\param\colon\ANNSET\to\N$, and $\intronorm{\cdot}\colon\ANNSET\to\R$
	satisfy for all $L\in\N$, $l_0,l_1,\ldots,l_L\in\N$,
	$
	\mathscr{f} 
	=((W_1,B_1),\allowbreak\ldots,\allowbreak(W_L,B_L))
	\in  \allowbreak
	\pr{ \times_{k = 1}^L\allowbreak(\R^{l_k \times l_{k-1}} \times \R^{l_k})}
	$,
	$x_0 \in \R^{l_0}$,
	$x_1 \in \R^{l_1}$, 
	$\ldots$, \allowbreak $x_{L-1} \in \R^{l_{L-1}}$ 
	with 
	$\forall \, k \in \N\cap(0,L) \colon\allowbreak x_k =\mathfrak{R}(W_k x_{k-1} + B_k)$  
	that
	\begin{equation}
		\label{setting:eq2}
		\realisation(\mathscr{f}) \in C(\R^{l_0},\R^{l_L}),
		\qquad
		( \realisation(\mathscr{f}) ) (x_0) = W_L x_{L-1} + B_L,
		\qquad
		\lengthANN(\mathscr{f})=L,
	\end{equation}
	\begin{equation}
		\label{setting:eq3}
		\textstyle
		\paramANN(\mathscr{f}) = \ssssum_{k = 1}^L l_k(l_{k-1} + 1),
		\qandq
		\intronorm{\f}=\max_{1\leq k\leq L}\fabs{(W_k,B_k)},
	\end{equation}
let $a\in\R$, $b\in(a,\infty)$,
and let $\Cost\colon\pr{\cup_{d\in\N}C(\R^d,\R)}\times[0,\infty]^3\to[0,\infty]$ satisfy for all $d\in\N$, $f\in C(\R^d,\R)$, $L,S,\eps\in[0,\infty]$ that
\begin{equation}
	\label{intro:Cost:def}
	\begin{split}
		\Cost(f,L,S,\eps)=
		&\min \pr*{ \pR*{ p \in \N \colon \PR*{\!\!
					\begin{array}{c}
						\exists\, \mathscr{f} \in \ANNSET \colon
						(\paramANN(\mathscr{f})\leq p)\land
						(\lengthANN(\mathscr{f})\leq L)
						\\\land{}
						(\intronorm{\mathscr{f}}\leq S)	\land{}
						(\realisation(\mathscr{f}) \in C(\R^d,\R))
						\\\land{}
						(\sup_{x\in[a,b]^d}\vass{(\realisation(\mathscr{f}))(x)-f(x)} \leq \varepsilon)\\
					\end{array} \!\! }
			}
			\cup\{\infty\}
		}.
	\end{split}
\end{equation}
\end{boxeddef}
\end{setting}
\end{samepage}

In \cref{setting} we also introduce the cost-functional which we employ to state our \ANN\ approximation results in \cref{Thm:gen}, \cref{Thm:sin}, and \cref{Thm:prod}.
Specifically, we note that \eqref{intro:Cost:def} in \cref{setting} ensures that for every dimension $d\in\N$, 
every continuous function $f\colon\R^d\to\R$, 
every upper bound $L\in[0,\infty]$ for the depths of the \ANNs,
every upper bound $S\in[0,\infty]$ for the size of the absolute values of the \ANN\ parameters, and 
every approximation precision $\eps\in[0,\infty]$ 
we have that $\Cost(f,L,S,\eps)$ represents the minimal number of \ANN\ parameters needed to approximate the target function $f\colon\R^d\to\R$ with the error tolerance $\eps$ among the class of \ANNs\ with at most $L$ affine linear transformations and the absolute values of the \ANN\ parameters being at most $S$.
\newpage

\usetikzlibrary{fadings,shapes,arrows.meta}
\def\layersep{4cm}
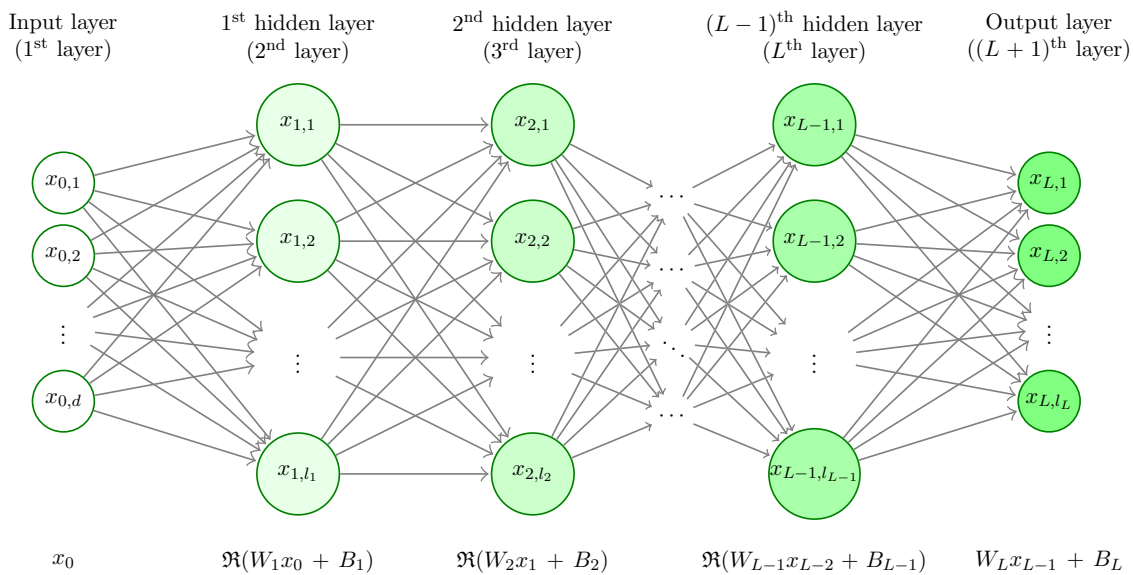
\begin{figure}
	\centering
	\begin{adjustbox}{width=\textwidth}
		\begin{tikzpicture}[shorten >=1pt,-latex,draw=black!100, node distance=\layersep,auto]
			\tikzstyle{every pin edge}=[<-,shorten <=1pt]
			\tikzstyle{neuron}=[circle,fill=black!25,draw=black!75,minimum size=40pt,inner sep=0pt,thick,font=\small]
			\tikzstyle{input neuron}=[neuron,
			draw={rgb:green,1;black,1}, 
			fill={rgb:white,5},minimum size=30pt];
			\tikzstyle{output neuron}=[neuron, draw={rgb:green,1;black,1}, fill={rgb:green,1;white,1},minimum size=30pt];
			\tikzstyle{hidden neuron1}=[neuron,
			draw={rgb:green,1;black,1}, 
			fill={rgb:green,1;white,10}];
			\tikzstyle{hidden neuron2}=[neuron,
			draw={rgb:green,1;black,1}, 
			fill={rgb:green,1;white,4}];
			\tikzstyle{hidden neuron3}=[neuron,
			draw={rgb:green,1;black,1}, 
			fill={rgb:green,1;white,2}];
			\tikzstyle{spacing neuron big}=[neuron,draw={rgb:white,5}, fill={rgb:white,40}];
			\tikzstyle{spacing neuron medium}=[neuron,draw={rgb:white,5}, fill={rgb:white,40},minimum size=30pt];
			\tikzstyle{spacing neuron small}=[neuron,draw={rgb:white,5}, fill={rgb:white,40},minimum size=20pt];
			\tikzstyle{annot} = [text width=9em, text centered,font=\small]
			\tikzstyle{annot2} = [text width=12em, text centered,font=\small]
			
			\foreach \name / \y in {1,...,2}
			\node[input neuron] (I-\name) at (0,-1.25*\y+1) {$x_{0,\y}$};
			\node[spacing neuron medium] (I-3) at (0,-3.75+1) {$\vdots$};
			\node[input neuron] (I-4) at (0,-5+1) {$x_{0,d}$};
			\foreach \name / \y in {1,...,2}
			\node[hidden neuron1] (H1-\name) at (\layersep,-2*\y+2.75) {$x_{1,\y}$};
			\node[spacing neuron big] (H1-3) at (\layersep,-6+2.75) {$\vdots$};
			\node[hidden neuron1] (H1-4) at (\layersep,-8+2.75)
			{$x_{1,l_1}$}; 
			\foreach \name / \y in {1,...,2}
			\node[hidden neuron2] (H2-\name) at (2*\layersep,-2*\y+2.75) {$x_{2,\y}$};
			\node[spacing neuron big] (H2-3) at (2*\layersep,-6+2.75) {$\vdots$};
			\node[hidden neuron2] (H2-4) at (2*\layersep,-8+2.75) {$x_{2,l_2}$};
			\node[spacing neuron small] (H3-1) at (2.6*\layersep,-2+1.5) {$\cdots$};
			\node[spacing neuron small] (H3-2) at (2.6*\layersep,-3.25+1.5) {$\cdots$};
			\node[spacing neuron small] (H3-3) at (2.6*\layersep,-4.5+1.5) {$\ddots$};
			\node[spacing neuron small] (H3-4) at (2.6*\layersep,-5.75+1.5) {$\cdots$};
			\foreach \name / \y in {1,...,2}
			\node[hidden neuron3] (H4-\name) at (3.2*\layersep,-2*\y+2.75) {$x_{L-1,\y}$};
			\node[spacing neuron big] (H4-3) at (3.2*\layersep,-6+2.75) {$\vdots$};
			\node[hidden neuron3] (H4-4) at (3.2*\layersep,-8+2.75) {$x_{L-1,l_{L-1}}$};
			\foreach \name / \y in {1,...,2}
			\node[output neuron] (H5-\y) at (4.2*\layersep,-1.25*\y+1 ) {$x_{L,\y}$};
			\node[spacing neuron small] (H5-3) at (4.2*\layersep,-3.75+1 ) {$\vdots$};
			\node[output neuron] (H5-4) at (4.2*\layersep,-5+1 ) {$x_{L,l_L}$};
			
			\foreach \target in {1,...,4}
			\foreach \source in {1,...,4}
			\path (I-\source) edge [-{Classical TikZ Rightarrow[scale=1.5]},draw=black!50,thick] (H1-\target);
			\foreach \target in {1,...,4}
			\foreach \source in {1,...,4}
			\path (H1-\source) edge [-{Classical TikZ Rightarrow[scale=1.5]},draw=black!50,thick] (H2-\target);
			\foreach \target in {1,...,4}
			\foreach \source in {1,...,4}
			\draw [-Classical TikZ Rightarrow, white,thick, path fading=west,postaction={draw,black!50,path fading=east}] (H2-\source) -- (H3-\target);
			\foreach \target in {1,...,4}
			\foreach \source in {1,...,4}
			\draw [-Classical TikZ Rightarrow, white,thick, path fading=east,postaction={draw,black!50,path fading=west}] (H3-\source) -- (H4-\target);
			\foreach \target in {1,...,4}
			\foreach \source in {1,...,4}
			\path (H4-\source) edge [-{Classical TikZ Rightarrow[scale=1.5]},draw=black!50,thick] (H5-\target);
			
			\node[annot,above of=H1-1, node distance=1.5cm, align=center] (hl) {$1^{\text{st}}$ hidden layer\\($2^{\text{nd}}$ layer)};
			\node[annot,right of=hl, align=center] (hl2) {$2^{\text{nd}}$ hidden layer\\($3^{\text{rd}}$ layer)};
			\node[annot,right of=hl2, align=center] (hl3) {};
			\node[annot,above of=H4-1, node distance=1.5cm,  align=center] (hl4) {$(L-1)^{\text{th}}$ hidden layer\\($L^{\text{th}}$ layer)};
			\node[annot,right of=hl4, align=center] (ol) {Output layer\\($(L+1)^{\text{th}}$ layer)};
			\node[annot,left of=hl, align=center] {Input layer\\ ($1^{\text{st}}$ layer)};
			\node[annot2,below of=H1-4, node distance=1.5cm, align=center] (sl1) {$\mathfrak{R}(W_1 x_{0} + B_1)$};
			\node[annot2,left of=sl1, align=center] (sl0) {$x_0$};
			\node[annot2,right of=sl1, align=center] (sl2) {$\mathfrak{R}(W_2 x_{1} + B_2)$};
			\node[annot2,right of=sl2, align=center] (sl3) {};
			\node[annot2,below of=H4-4, node distance=1.5cm, align=center] (sl4) {$\mathfrak{R}(W_{L-1} x_{L-2} + B_{L-1})$};
			\node[annot2,right of=sl4, align=center] (sl5) {$W_{L} x_{L-1} + B_{L}$};
		\end{tikzpicture}
	\end{adjustbox}
	\caption{\label{figure }Graphical illustration of the architecture of the \ANNs\ used in \eqref{setting:eq1} and \eqref{setting:eq2} in \cref{setting}: \ANNs\ with $L+1$ layers (with $L$ affine linear transformations) with an $l_0$-dimensional input layer ($l_0$ neurons on the input layer), an $l_1$-dimensional $1^{\text{st}}$ hidden layer ($l_1$ neurons on the $1^{\text{st}}$ hidden layer), an $l_2$-dimensional $2^{\text{nd}}$ hidden layer ($l_2$ neurons on the $2^{\text{nd}}$ hidden layer), $\dots$, an $l_{L-1}$-dimensional $(L-1)^{\text{th}}$ hidden layer ($l_{L-1}$ neurons on the $(L-1)^{\text{th}}$ hidden layer), and an $l_{L}$-dimensional output layer ($l_{L}$ neurons on the output layer).
		The realization function in \eqref{setting:eq2} in \cref{setting} 
		assigns the $l_0$-dimensional input vector $x_0=(x_{0,1},\ldots,x_{0,l_0})\in\R^{l_0}$ to the $l_1$-dimensional vector $x_1=(x_{1,1},\ldots,x_{1,l_1})\in\R^{l_1}$ with $x_1=\mathfrak{R}(W_1 x_{0} + B_1)$ on the $1^{\text{st}}$ hidden layer, 
		assigns the vector $x_1=(x_{1,1},\ldots,x_{1,l_1})\in\R^{l_1}$ on the $1^{\text{st}}$ hidden layer to the vector $x_2=(x_{2,1},\ldots,x_{2,l_2})\in\R^{l_2}$ with $x_2=\mathfrak{R}(W_2 x_{1} + B_2)$ on the $2^{\text{nd}}$ hidden layer, 
		$\dots$,
		assigns the vector $x_{L-2}=(x_{L-2,1},\ldots,x_{L-2,l_{L-2}})\in\R^{l_{L-2}}$ on the $(L-2)^{\text{th}}$ hidden layer to the vector $x_{L-1}=(x_{L-1,1},\ldots,x_{L-1,l_{L-1}})\in\R^{l_{L-1}}$ with $x_{L-1}=\mathfrak{R}(W_{L-1} x_{L-2} + B_{L-1})$ on the $(L-1)^{\text{th}}$ hidden layer, and
		assigns the vector $x_{L-1}=(x_{L-1,1},\ldots,x_{L-1,l_{L-1}})\in\R^{l_{L-1}}$ on the $(L-1)^{\text{th}}$ hidden layer to the vector $x_L=(x_{L,1},\ldots,x_{L,l_L})\in\R^{l_L}$ with $x_L=W_L x_{L-1} + B_L$ on the output layer.\\
	}
	\label{figure1}
\end{figure}

\subsection{The necessity of depth to approximate certain classes of functions}
\label{intro:bounded}
Using \cref{setting} we now formulate the above mentioned three representative key \ANN\ approximation results.
\begin{theorem}
	\label{Thm:gen}
			\begin{boxedthm}
	Assume \cref{setting}. Then there exist $\mathfrak{c}\in(0,\infty)$ and infinitely often differentiable $f_d\colon\R^d\to\R$, $d\in\N$,
	with compact support and $\sup_{d\in\N}\sup_{x\in\R^d}\vass{f_d(x)}\leq 1$
	such that for all $d,L\in\N$, $\eps\in(0,1)$ it holds that 
	\begin{equation}
		\label{Thm:gen:eq1}
		\Cost(f_d,L,\infty,\eps)\geq 2^{\nicefrac{d}{L}}
		\qandq
		\Cost(f_d,\mathfrak{c}d,1,\eps)\leq \mathfrak{c}d\eps^{-2}.
	\end{equation}
			\end{boxedthm}
\end{theorem}

\cref{Thm:gen} above is an immediate consequence of \cref{Cor:Thm6.1} in \cref{Upper and Lower Bounds} below and \cref{Cor:Thm6.1}, in turn, follows from \cref{Thm6.1} which is implied by \cref{Thm6} in \cref{Upper and Lower Bounds}.
Roughly speaking, \cref{Thm:gen} reveals that there exists a sequence of smooth and uniformly globally bounded functions $f_d\colon\R^d \to\R$ for $d\in\N$ with compact support which 
can neither be approximated without the curse of dimensionality by means of shallow \ANNs\ nor insufficiently deep \ANNs\ even if the \ANN\ parameters may be arbitrarily large (see the first inequality in \eqref{Thm:gen:eq1}) but which can be approximated without the curse of dimensionality by sufficiently deep \ANNs\ even if the absolute values of the \ANN\ parameters are assumed to be uniformly bounded by $1$ (see the second inequality in \eqref{Thm:gen:eq1}).
\Cref{Thm:gen} only asserts the existence of suitable smooth and uniformly globally bounded target functions which can be approximated without the curse of dimensionality by deep \ANNs\ but not by insufficiently deep \ANNs\ 
but it does not explicitly specify the employed target functions.
%
However, in the more general approximation result in \cref{Thm6} in \cref{Upper and Lower Bounds} we also explicitly specify a class of simple target functions which we use to prove \cref{Thm:gen}.
%
%
In addition, in this work we also reveal that the sine of the product functions serve as one possible sequence of smooth and uniformly globally bounded target functions for which the conclusion of \cref{Thm:gen} essentially applies 
and this is precisely the subject of our next representative key \ANN\  approximation result.

\begin{theorem}
		\label{Thm:sin}
		\begin{boxedthm}
	Assume \cref{setting}, assume $\max\{\vass{a},\vass{b}\}>1$, and for every $d\in\N$ let $f_d\colon\R^d\to\R$ satisfy for all $x=(x_1,\ldots,x_d)\in\R^d$ that $f_d(x)=\sin\prb{\sprod_{i=1}^d x_i}$.
	Then there exists $\mathfrak{c}\in(0,\infty)$ such that for all $d,L\in\N$, $\eps\in(0,1)$ it holds that 
	\begin{equation}
		\label{Thm:sin:eq1}
		\Cost(f_d,L,\infty,\eps)\geq 
		(1+\mathfrak{c}^{-1})^{\nicefrac{d}{L}}
		\qandq
		\Cost(f_d,\mathfrak{c}d^2\eps^{-1},1,\eps)\leq \mathfrak{c}d^3\eps^{-2}
		.
	\end{equation}
			\end{boxedthm}
\end{theorem}

\cref{Thm:sin} above is an immediate consequence of \cref{cor:sin:prod:oneplus} in \cref{Upper and Lower Bounds} below and \cref{cor:sin:prod:oneplus} is implied by \cref{Thm1} in \cref{Upper and Lower Bounds}.
\Cref{Thm:sin} shows that the sine of the product functions
can neither be approximated without the curse of dimensionality by means of shallow \ANNs\ nor insufficiently deep \ANNs\ even if the \ANN\ parameters may be arbitrarily large (see the first inequality in \eqref{Thm:sin:eq1}) but can be approximated without the curse of dimensionality by sufficiently deep \ANNs\ even if the absolute values of the \ANN\ parameters are assumed to be uniformly bounded by $1$ (see the second inequality in \eqref{Thm:sin:eq1}).
Actually, a bit modified variant of the conclusion of \cref{Thm:sin} applies also to the product functions themselves.
This is precisely the subject of our final representative key \ANN\ approximation result in this introductory section.

\begin{theorem}
	\label{Thm:prod}
			\begin{boxedthm}
	Assume \cref{setting}, assume $\max\{\vass{a},\vass{b}\}>1$, and for every $d\in\N$ let $f_d\colon\R^d\to\R$ satisfy for all $x=(x_1,\ldots,x_d)\in\R^d$ that $f_d(x)=\sprod_{i=1}^d x_i$.
	Then there exists $\mathfrak{c}\in(0,\infty)$ such that for all $d,L\in\N$, $\eps\in(0,1)$ it holds that
	\begin{equation}
		\label{eq:prod:Thm:1}
		\Cost(f_d,L,\infty,\eps)
		\geq 
		(1+\mathfrak c^{-1})^{\nicefrac{d}{L}}
		\qandq				
		\Cost(f_d,\mathfrak{c}d^2\eps^{-1},1,\eps)
		\leq  
		\mathfrak{c}d^3\eps^{-1}.
	\end{equation}
			\end{boxedthm}
\end{theorem}

\cref{Thm:prod} above is an immediate consequence of \cref{Cor:oneplus:for:prod} in \cref{Upper and Lower Bounds} below and \cref{Cor:oneplus:for:prod}, in turn, follows from \cref{Thm} in \cref{Upper and Lower Bounds}.
Loosely speaking, \cref{Thm:prod} proves that the plain vanilla product functions
can neither be approximated without the curse of dimensionality by means of shallow \ANNs\ nor insufficiently deep \ANNs\ even if the absolute values of the \ANN\ parameters are arbitrarily large (see the first inequality in \eqref{eq:prod:Thm:1}) but can be approximated without the curse of dimensionality by sufficiently deep \ANNs\ even if the absolute values of the \ANN\ parameters are assumed to be uniformly bounded by $1$ (see the second inequality in \eqref{eq:prod:Thm:1}).

\subsection{Literature review}
\label{intro:literature}

The purpose of this subsection is to discuss selected approximation
results which are related to the questions studied in this article.
In particular, we are interested in the dependence of approximation
bounds on the input dimension, the prescribed approximation precision,
the depth of the approximating \ANNs, and restrictions on their
parameters.
For broader discussions, we refer to the surveys
\cite{pinkus1999approximation,guehring2020expressivity}.

The existence of \ANN\ approximations with arbitrarily small errors
has been established under various assumptions on the activation
function and the target functions.
For shallow \ANNs\ with sigmoidal
activation functions, classical results can be found in
\cite{gallant1988there,carroll1989construction,Cybenko1989,
	funahashi1989approximate,hornik1989multilayer}; see also
\cite{irie1988capabilities,hornik1991approximation,leshno1993multilayer}
for related results and conditions on activation functions.
For deep \ANNs, \cite{Hanin2017,kidger2020universal} establish
universal approximation results with widths bounded independently of
the prescribed approximation precision, while the depth is allowed
to increase.
The existence of such approximations alone does not
determine whether their number of parameters can be bounded
polynomially in both the input dimension and the reciprocal of the
approximation precision.

Quantitative approximation results additionally describe how the
size of the approximating \ANNs\ depends on the regularity of the
target functions.
Lu et al.~\cite{LuShenYangZhang2020} and
Shen et al.~\cite{shen2020deep} establish estimates for smooth and
H\"older continuous target functions, respectively, which explicitly
account for the width and depth of the approximating \ReLU\ \ANNs.
Petersen \& Voigtlaender \cite{petersen2018optimal} study the
approximation of piecewise smooth functions and Voigtlaender \&
Petersen \cite{voigtlaender2019approximation} consider errors
measured in $L^p(\mu)$ for finite Borel measures $\mu$.
G\"uhring et al.~\cite{guhring2019error} establish upper and lower
bounds for Sobolev-regular functions with errors measured in weaker
Sobolev norms.
For target functions in Sobolev and Besov spaces,
Siegel \cite{Siegel2023Optimal} establishes optimal approximation
rates in $L^p$ under suitable compact embedding assumptions.
For the general function classes described above, these approximation results do not by themselves imply that \ReLU\ \ANNs\ overcome the curse of dimensionality.

For earlier quantitative results, including nearly geometric
convergence for target functions admitting analytic extensions to a
neighborhood of the approximation domain, we refer to Mhaskar
\cite{mhaskar1993}.
For analytic target functions in a fixed input dimension, E \& Wang
\cite{wang2018exponential} establish exponential convergence in a fractional power of network depth for deep \ANN\ approximations.
Moreover,
\cite{elbraechter2021deep,perekrestenko2018universal} provide
deep \ReLU\ \ANN\ constructions for polynomials and certain
oscillatory functions with depth growing at most polylogarithmically
in the reciprocal of the prescribed approximation accuracy.

Further approximation results make use of additional structure of
the target functions.
Lee et al.~\cite{lee2017ability} study
compositions of functions satisfying suitable Fourier conditions,
whereas Shaham et al.~\cite{ShahamCloningerCoifman2018} construct
\ANN\ approximations for functions defined on low-dimensional
manifolds.
Cheridito et al.~\cite{cheridito2021efficient} provide
a framework for approximating functions represented by suitable
catalog networks and obtain approximation results without the curse
of dimensionality.
These results identify assumptions under which
efficient high-dimensional \ANN\ approximations can be constructed.
For radially symmetric Lipschitz continuous target functions which
are constant outside a ball, McCane \& Szymanski
\cite{McCaneS17} construct deep \ReLU\ \ANN\ approximations and
bound the number of neurons sufficient to achieve a prescribed
approximation precision quadratically in the input dimension.

Efficient approximation is also possible with shallow \ANNs\ for
suitable classes of target functions.
In particular, the results in
\cite{barron1992neural,barron1993universal,barron1994approximation,
	jones1992simple} give estimates whose convergence exponents do not
worsen with the input dimension.
For shallow \ANNs\ with positively homogeneous activation functions,
including \ReLU, Bach \cite{bach2017breaking} establishes bounds
that exploit low-dimensional structure in the target functions.
Further estimates can be found, \eg, in
\cite{makovoz1996random,klusowski2018approximation}.

The choice of the activation function can also substantially affect
approximation complexity.
Maiorov \& Pinkus \cite{MaioPinkus1999}
construct an analytic sigmoidal activation function for which
\ANNs\ with two hidden layers and, for each fixed input dimension,
a fixed finite number of hidden neurons can approximate continuous
functions on compact domains with arbitrarily small errors; see also
\cite{GuliIsm2018a,GULIYEV2018296}.

Lower bounds provide complementary information on the complexity
of approximating prescribed classes of target functions.
Yarotsky \cite{yarotsky2017error} establishes lower bounds for
\ReLU\ \ANN\ approximation of Sobolev classes, while Petersen \&
Voigtlaender \cite{petersen2018optimal} establish optimality results
for piecewise smooth functions, including lower bounds under restrictions on the number of bits used to represent each parameter.
For approximation
in $L^2$, B\"olcskei et al.~\cite{BolcskeiGrohsKutyniokPetersen2019OptimalApproximation}
establish lower bounds for the number of non-zero weights and
memory requirements and show that these bounds can be attained
for suitable classes of target functions.

For target functions in Sobolev and Besov spaces, G\"uhring \&
Raslan \cite{GuehringRaslan2021} establish approximation bounds
under the requirement that each non-zero network parameter can be
represented using a number of bits which grows at most
logarithmically in the reciprocal of the approximation precision.
For Lipschitz continuous target functions on $[0,1]$ Ou et al.~\cite{OuSchenkelBolcskei2024} establish upper and lower
bounds for \ReLU\ \ANNs\ with finite-precision weights and
investigate how increasing depth can compensate for reduced
weight precision.

In this article we also investigate the dependence of memory
requirements on the input dimension.
More specifically,
\cref{cor:bit:both:bounds} in \cref{subsection:entropy_results} establishes for suitable 
functions exponential lower bounds for \ANNs\ with any fixed
number of hidden layers and polynomial upper bounds in both the
input dimension and the reciprocal of the approximation precision
for sufficiently deep \ANNs.

The piecewise affine structure of realizations of \ReLU\ \ANNs\
provides an important tool for establishing lower
bounds results.
For estimates concerning the number of affine pieces and
regions, we refer to
\cite{TakaiYuukiSannai20,Goujon22,Pascanu13,Petersen2020NeuralNT}.
We employ similar techniques among others in
\cref{S4_lower_bounds} to estimate the limits of the expressive power of \ReLU\ \ANNs\ with a limited depth.
Grohs \& Voigtlaender \cite{grohs2021proof} establish limitations on algorithms that construct approximations from finitely many evaluations of the target function, even for classes of functions that admit efficient \ANN\ approximations.

Of particular relevance to this article are results which compare
the number of parameters required by \ANNs\ of different depths
to approximate the same target functions.
Besides the results
discussed earlier in this section, Safran et al.~
\cite{SafranEldanShamir19} investigate how such comparisons depend
on the Lipschitz regularity, the approximation precision, and the
size of the approximation domain.
More recently, Safran et al.~\cite{SafranReichmanValiant2025} establish an exponential separation in the input dimension between shallow \ANNs\ and \ANNs\ with two hidden layers for approximation in $L^2$ with respect to probability measures supported in the unit ball, at a fixed approximation precision and under suitable exponential bounds on the absolute values of the parameters of the shallow \ANNs.

For approximation in $L^2$ with respect to the standard Gaussian
measure, Bortolotto et al.~\cite{BortolottoGanacimOliveira2026}
establish an exponential separation in the input dimension between
\ReLU\ \ANNs\ with one and two hidden layers for suitable radially
symmetric target functions, at a fixed approximation precision
and under subexponential bounds on the absolute values of the
parameters of the shallow \ANNs.

While two hidden layers are sufficient for efficient approximation in \cite{SafranReichmanValiant2025,BortolottoGanacimOliveira2026}, the target functions in \cref{Thm:gen} and \cref{Thm:sin} cannot be uniformly approximated efficiently on cubes by \ReLU\ \ANNs\ with any fixed depth, but admit efficient approximation by sufficiently deep \ReLU\ \ANNs.

Safran \cite{Safran2026EveryLayerCounts} establishes
exponential separations in the input dimension between
\ReLU\ \ANNs\ with any fixed positive number of hidden
layers and \ANNs\ with one additional hidden layer,
even if their parameters may be arbitrarily large.
These results concern approximation in $L^2$ at a fixed
approximation precision for globally bounded
$1$-Lipschitz continuous target functions and probability
measures supported at exponential distance from the origin.
The target functions are piecewise affine and depend on
the depths being compared.
In \cref{Thm:gen} and \cref{Thm:sin},
each of the considered families of infinitely often
differentiable target functions yields lower bounds
for every fixed number of hidden layers.

Concerning the dependence on depth and approximation precision,
Yarotsky \cite{Yarotsky2018Optimal} and Yarotsky \& Zhevnerchuk
\cite{YarotskyZhevnerchuk2020} identify approximation regimes in
which improved rates require increasing depth and discontinuous
parameter selection.
Our lower bounds in
\cref{Thm:gen}, \cref{Thm:sin}, and \cref{Thm:prod} also quantify the dependence on
the input dimension when the approximation precision is fixed.
Concerning parameter sizes, Vardi \& Shamir
\cite{VardiShamir20} show, under suitable assumptions on the target
functions and input distributions, that polynomial-size \ReLU\
\ANN\ approximations with unrestricted parameters and fixed depth
can be replaced by polynomial-size approximations of greater
constant depth with polynomially bounded parameters.
They also
explain why certain depth-separation results would resolve
longstanding open problems in computational complexity.

Interesting and simple target functions used repeatedly in the
literature are the product functions.
We refer to
\cite{Lin_2017,RolnickT17} for lower bounds for shallow \ANNs\
with smooth activation functions and to
\cite{cheridito2021efficient,beneventano21,Opschoor2021ExponentialRD}
for upper bounds for \ReLU\ \ANNs.
For uniform approximation on $[-k,k]^d$ with $k>1$,
Shapira \cite[Theorem~3]{Shapira2023} establishes exponential
lower bounds for the number of parameters in the input dimension for \ReLU\ \ANNs\ with
any fixed number of hidden layers, even if their parameters
may be arbitrarily large.
For $a,b\in\R$ with $a<b$ and
$\max\{\vass{a},\vass{b}\}>1$,
we establish in \cref{Thm} in \cref{Upper and Lower Bounds} below
exponential lower bounds for the number of non-zero parameters in the input dimension for
uniform approximations on $[a,b]^d$ by \ReLU\ \ANNs\ with any fixed
number of hidden layers, even if their parameters
may be arbitrarily large.
We employ \cref{Thm} to prove \cref{Thm:prod} above.
The same functions can be approximated without
the curse of dimensionality by sufficiently deep \ReLU\
\ANNs\ even if the absolute values of their parameters are
uniformly bounded by $1$.
On $[0,1]^d$, however, Blanchard \& Bennouna
\cite{BlanchardBennouna2022} construct uniform approximations
of the product functions using \ReLU\ \ANNs\ with two hidden
layers and a number of neurons bounded polynomially in both
the input dimension and the reciprocal of the approximation
precision.

Our upper bound constructions for the product functions adapt
the constructions in Beneventano et al.~\cite[Section~6]{beneventano21}.
In particular, \cref{lemma:2d_product} is a slightly extended
variant of \cite[Lemma~6.7]{beneventano21}.
Our constructions for
periodic functions also build on the properties of sawtooth
functions in Telgarsky \cite[Section~2.2]{Telgarsky15}.
The auxiliary approximation result for one-dimensional Lipschitz
continuous functions builds on Beneventano et al.~
\cite[Subsection~4.1]{beneventano21}.
The respective adaptations and additional estimates, as well as
the sources of further preparatory results, are specified in
\cref{S2_ANN_calculus,S5_upper_bounds}.

\subsection{Structure of this article}
\label{intro:structure}

%

The remainder of this article is organized as follows.
In \cref{S2_ANN_calculus} we present the concepts, operations, and elementary preparatory results regarding \ANNs\ that we frequently employ in \cref{S4_lower_bounds,S5_upper_bounds,Upper and Lower Bounds,subsection:entropy_results}.
In \cref{S4_lower_bounds} we establish suitable lower bounds for the minimal number of parameters of shallow or insufficiently deep \ANNs\ to approximate certain high-dimensional target functions.
In \cref{S5_upper_bounds} we establish suitable upper bounds for the minimal number of parameters of \ANNs\ to approximate the product functions and certain highly oscillating functions in the case where the absolute values of the parameters of the \ANNs\ are assumed to be uniformly bounded by $1$.
In \cref{Upper and Lower Bounds} we combine the main results from \cref{S4_lower_bounds,S5_upper_bounds} to obtain that the product functions and certain highly oscillating functions can be approximated without the curse of dimensionality by deep \ANNs\ but not by insufficiently deep \ANNs\ and, thereby, we prove our three representative key \ANN\ approximation results in \cref{Thm:gen}, \cref{Thm:sin}, and \cref{Thm:prod} above. 
In \cref{subsection:entropy_results} we utilize the main results from \cref{S5_upper_bounds,Upper and Lower Bounds} to show that the amount of memory to store \ANN\ approximations for certain simple functions on a computer
 must grow at least exponentially in the problem dimension for shallow or insufficiently deep \ANNs\
but
 grows at most polynomially in the problem dimension and the reciprocal of the prescribed approximation accuracy for sufficiently deep \ANNs.

\newpage

\section{Artificial neural network (ANN) calculus}
\label{S2_ANN_calculus}
The purpose of this section is to introduce the concepts, operations, and elementary preparatory results regarding \ANNs\ that we frequently employ in the later sections of this article.

In particular, in \cref{def:ANN} in \cref{subsection:SetofANNs} we recall the notion of the set of \ANNs\ $\ANNSET$ in the structured description as well as several basic functions acting on this set of \ANNs\ such as the parameter function 
$\paramANN \colon \ANNSET \to \N$ for \ANNs\
(counting the number of parameters of an \ANN), 
the length function
$\lengthANN \colon \ANNSET \to \N$ for \ANNs\
(measuring the number of affine linear transformations of an \ANN),
the input dimension function
$\inDimANN \colon \ANNSET \to \N$ for \ANNs\
(specifying the number of neurons on the input layer of an \ANN),
the output dimension function
$\outDimANN \colon \ANNSET \to \N$ for \ANNs\
(specifying the number of neurons on the output layer of an \ANN), 
the hidden layer function 
$\hidlengthANN\colon \ANNSET \to \N_0$ for \ANNs\
(counting the number of hidden layers of an \ANN),
the layer dimension vector function
$\dims\colon\ANNSET\to  \pr{\cup_{L\in\N} \N^{L}}$ for \ANNs\
(representing the numbers of neurons on the layers of an \ANN\ as a vector), and
the layer dimension functions
$\singledims_n \colon \ANNSET \to \N_0$, $n\in\N_0$, for \ANNs\
(counting the numbers of neurons on the layers of an \ANN).

In \cref{subsection:realization} we recall the concept of realization functions of \ANNs, 
in \cref{subsec:parallelization} we recall the concept of parallelizations of \ANNs,
in \cref{subsection:idANNs} we recall suitable \ANNs\ whose realization functions exactly coincide with the real identity functions,
in \cref{subsectio:compositions} we recall the concept of compositions of \ANNs,
and in \cref{subsection:size} we present elementary concepts and results regarding the sizes of the absolute values of the parameters of \ANNs.

Most of the concepts and results in this section are well-known and have appeared, often in a bit different form, in previous works in the literature (cf., \eg,
\cite{beneventano21,GrohsIbrgimovJentzen2021,petersen2018optimal,GrohsHornungJentzen2019,GrohsJentzenSalimova2019,%
	BeckJentzenKuckuck2019,ElbraechterSchwab2018,cheridito2021efficient,GononSchwab21,GononGrohsJentzenKoflerSiska22}).
In particular,
\cref{def:ANN} is a slightly extended version of, 
\eg, Grohs et al.~\cite[Definition~2.1]{GrohsHornungJentzen2019},
\cref{def:RELU} corresponds to, \eg, 
Grohs et al.~\cite[Definition~2.1]{GrohsIbrgimovJentzen2021},
\cref{def:ANNrealization} is a reformulated variant of, \eg, 
Grohs et al.~\cite[Definition~2.3]{GrohsHornungJentzen2019},
\cref{def:simpleParallelization} is a reformulated variant of, \eg, 
Grohs et al.~\cite[Definition~2.17]{GrohsHornungJentzen2019},
\cref{Lemma:PropertiesOfParallelizationEqualLength} is a slightly differently presented variant of, \eg, 
Grohs et al.~\cite[Lemma~2.18 and Proposition~2.19]{GrohsHornungJentzen2019},
\cref{def:ReLU_identity} corresponds to, \eg, 
Grohs et al.~\cite[Definition~2.13]{GrohsIbrgimovJentzen2021} and Beneventano et al.~\cite[Definition~2.12]{beneventano21},
\cref{Prop:identity_representation} corresponds to, 
\eg, Beneventano et al.~\cite[Proposition~2.13]{beneventano21},
\cref{def:ANNcomposition} is a reformulated variant of, \eg, 
Grohs et al.~\cite[Definition~2.5]{GrohsHornungJentzen2019},
\cref{Lemma:CompositionAssociative} corresponds to, \eg, 
Grohs et al.~\cite[Lemma~2.8]{GrohsHornungJentzen2019},
\cref{Lemma:PropertiesOfCompositions_n2} corresponds to, 
\eg, Beneventano et al.~\cite[Proposition~2.16]{beneventano21},
\eqref{eq:size:def} in \cref{def:size} is inspired by,
\eg, Grohs et al.~\cite[Definitions~2.21 and 2.22]{GrohsIbrgimovJentzen2021}, and
\cref{direct:comp:size} in \cref{lem:sizecomp} is a reformulated special case of, 
\eg, Grohs et al.~\cite[Lemma~2.23]{GrohsIbrgimovJentzen2021}.

\subsection{Set of ANNs}
\label{subsection:SetofANNs}

\begin{definition}[Set of ANNs]
	\label{def:ANN}
	We denote by $\ANNSET$ the set given by 
	\begin{equation}
	\begin{split}
	\ANNSET
	&=
	\textstyle
	\bigcup_{L \in \N}
	\bigcup_{ l_0,l_1,\ldots, l_L \in \N }
	\prb{
	\bigtimes_{k = 1}^L \pr{\R^{l_k \times l_{k-1}} \times \R^{l_k}}
	},
	\end{split}
	\end{equation}
	we denote by
	$\paramANN \colon \ANNSET \to \N$, 
	$\lengthANN \colon \ANNSET \to \N$,
	$\inDimANN \colon \ANNSET \to \N$,
	$\outDimANN \colon \ANNSET \to \N$, 
	$\hidlengthANN\colon \ANNSET \to \N_0$, and
	$\dims\colon\ANNSET\to  \pr{\cup_{L\in\N} \N^{L}}$
	the functions which satisfy
	for all $L\in\N$, $l_0,l_1,\ldots,\allowbreak l_L \in \N$, 
	$
	\mathscr{f} 
	\in  \allowbreak
	\pr{ \times_{k = 1}^L\allowbreak(\R^{l_k \times l_{k-1}} \times \R^{l_k})}$
	that
	\begin{equation}
		\label{eq:def:ANN:operators}
	\textstyle
	\paramANN(\mathscr{f}) =
	\sum_{k = 1}^L l_k(l_{k-1} + 1),
	\quad
	\lengthANN(\mathscr{f})=L,
	\quad 
	\inDimANN(\mathscr{f})=l_0,
	\quad
	\outDimANN(\mathscr{f})=l_L,
	\quad
	\hidlengthANN(\mathscr{f})=L-1,
	\end{equation} and
	$\dims(\mathscr{f})= (l_0,l_1,\ldots, l_L)$,
	for every $n \in \N_0$ we denote by $\singledims_n \colon \ANNSET \to \N_0$ the function which satisfies for all $L\in\N$, $l_0,l_1,\ldots,\allowbreak l_L \in \N$, 
	$
	\mathscr{f} 
	\in  \allowbreak
	\pr{ \times_{k = 1}^L\allowbreak(\R^{l_k \times l_{k-1}} \times \R^{l_k})}$ that
	\begin{equation}
	    \singledims_n(\mathscr{f})= \begin{cases}
	    l_n &\colon n\leq L\\
	    0   &\colon n>L,
	    \end{cases}
	\end{equation}
and for every $L\in\N$, $l_0,l_1,\ldots,\allowbreak l_L \in \N$, 
$
\mathscr{f}=((W_1, B_1),\allowbreak \ldots, (W_L,\allowbreak B_L))
\in  \allowbreak
\pr{ \times_{k = 1}^L\allowbreak(\R^{l_k \times l_{k-1}} \times \R^{l_k})}$ we denote by 
$\weight{(\cdot)}{\f}=(\weight{n}{\f})_{n\in\{1,2,\ldots,L\}}\colon\{1,2,\ldots,L\}\to\pr{\cup_{k,m\in\N}\R^{k\times m}}$ and 
$\bias{(\cdot)}{\f}=(\bias{n}{\f})_{n\in\{1,2,\ldots,L\}}\colon\{1,2,\ldots,L\}\to\pr{\cup_{k\in\N}\R^{k}}$ 
the functions which satisfy for all $n\in\{1,2,\ldots,L\}$ that $\weight{n}{\f}=W_n$ and $\bias{n}{\f}=B_n$.
\end{definition}
\cfclear

\subsection{Realization functions of ANNs}
\label{subsection:realization}

\begin{definition}[Multidimensional ReLU]
\label{def:RELU}
We denote by $\RELU\colon\pr{\cup_{d\in\N}\R^d}\to \pr{\cup_{d\in\N}\R^d}$ the function which satisfies for all $d\in\N$, $x=(x_1,\ldots,x_d)\in\R^d$ that
\begin{equation}
\RELU(x)=(\max\{x_1,0\},\max\{x_2,0\},\ldots,\max\{x_d,0\}).
\end{equation}
\end{definition}
\cfclear

\begin{definition}[Realization functions of ANNs]
	\label{def:ANNrealization}
	\cfconsiderloaded{def:ANNrealization}
	We denote by 
	$
	\realisation \colon\allowbreak \ANNSET \to\allowbreak \pr{\cup_{k,l\in\N}\,C(\R^k,\R^l)}
	$
	the function which satisfies
	for all
	$
	\mathscr{f}
	\in\ANNSET
	$,
	$x_0 \in \R^{\singledims_0(\f)}$, 
	$x_1 \in \R^{\singledims_1(\f)}$, $\ldots$, $x_{\hidlengthANN(\f)} \in \R^{\singledims_{\hidlengthANN(\f)}(\f)}$ 
	with
	$\forall \, k \in \N\cap[0,\hidlengthANN(\f)] \colon x_k =\RELU(\weight{k}{\f} x_{k-1} + \bias{k}{\f})$  
	that
	\begin{equation}
		\label{ANNrealization:ass2}
		\realisation(\mathscr{f}) \in C(\R^{\inDimANN(\f)},\R^{\outDimANN(\f)})\qandq
		( \realisation(\mathscr{f}) ) (x_0) =\weight{\lengthANN(\f)}{\f} x_{\hidlengthANN(\f)} + \bias{\lengthANN(\f)}{\f}
	\end{equation}
	\cfload.
\end{definition}
\cfclear

\subsection{Parallelizations of ANNs}
\label{subsec:parallelization}


\cfclear
\begin{definition}[Parallelization of ANNs]
	\label{def:simpleParallelization}
	\cfconsiderloaded{def:simpleParallelization}
	Let $n\in\N$. Then we denote by 
	\begin{equation}
	\parallelizationSpecial_{n}\colon \pRb{\f=(\mathscr{f}_1,\dots, \mathscr{f}_n)\in\ANNSET^n\colon \lengthANN(\mathscr{f}_1)= \lengthANN(\mathscr{f}_2)=\ldots =\lengthANN(\mathscr{f}_n) }\to \ANNSET
\end{equation}
	the function which satisfies for all
	 $\f=(\f_1,\ldots,\f_n)\in\ANNSET^n$, $k\in\{1,2,\ldots,\lengthANN(\f_1)\}$ with $\lengthANN(\mathscr{f}_1)= \lengthANN(\mathscr{f}_2)=\ldots =\lengthANN(\mathscr{f}_n)$ that 
		\begin{equation}\label{parallelisationSameLengthDef}\\
			\lengthANN(\parallelizationSpecial_{n}(\f))=\lengthANN(\mathscr{f}_1),
			\,\,\,
			\weight{k}{\parallelizationSpecial_{n}(\f)}=
			\begin{pmatrix}
				\weight{k}{\f_1}&\!\!\! 0&\!\!\!  \cdots&\!\!\! 0\\
					0&\!\!\! \weight{k}{\f_2}&\!\!\! \cdots&\!\!\! 0\\
					\vdots&\!\!\! \vdots&\!\!\! \ddots&\!\!\! \vdots\\
					0&\!\!\! 0&\!\!\! \cdots&\!\!\! \weight{k}{\f_n}
		    \end{pmatrix}\!,
	   \,\,\,\text{and}\,\,\,
	    \bias{k}{\parallelizationSpecial_{n}(\f)}=
	    \begin{pmatrix}
	    	\bias{k}{\f_1}\\
	    	\bias{k}{\f_2}\\
	    	\vdots\\
	    	\bias{k}{\f_n}
	    \end{pmatrix}  
	\end{equation}
	\cfout[.]
\end{definition}
\cfclear

\begin{athm}{prop}{Lemma:PropertiesOfParallelizationEqualLength}	
	Let 
	$n\in\N$, 
	$\mathscr{f}=(\mathscr{f}_1,\allowbreak\dots,\allowbreak \mathscr{f}_n)\in\ANNSET^n$
	satisfy  $\lengthANN(\mathscr{f}_1)= \lengthANN(\mathscr{f}_2)=\ldots =\lengthANN(\mathscr{f}_n)$
	\cfload.
	Then
	\begin{enumerate}[(i)]
		\item{\label{ParallelizationElementary:Display}it holds	for all $k\in\N_0$ that
	$
				\singledims_{k}\pr{\parallelizationSpecial_{n}(\mathscr{f})}=\sum_{j=1}^n\singledims_k(\mathscr{f}_j)
				$,
		}
		\item\label{PropertiesOfParallelizationEqualLength:ItemOne} it holds that 
			$
			\realisation(\parallelizationSpecial_{n}(\mathscr{f}))\in C\prb{\R^{[\sum_{j=1}^n \inDimANN(\mathscr{f}_j)]},\R^{[\sum_{j=1}^n \outDimANN(\mathscr{f}_j)]}}$,
		and
		\item\label{PropertiesOfParallelizationEqualLength:ItemTwo} it holds for all    $x_1\in\R^{\inDimANN(\mathscr{f}_1)},x_2\in\R^{\inDimANN(\mathscr{f}_2)},\dots, x_n\in\R^{\inDimANN(\mathscr{f}_n)}$ that 
		\begin{equation}\label{PropertiesOfParallelizationEqualLengthFunction}
			\begin{split}
				&\prb{\realisation\prb{\parallelizationSpecial_{n}(\mathscr{f}) }}(x_1,x_2,\dots, x_n) 
				=\prb{(\realisation(\mathscr{f}_1))(x_1), (\realisation(\mathscr{f}_2))(x_2),\dots,
					(\realisation(\mathscr{f}_n))(x_n) }
			\end{split}
		\end{equation}
	\end{enumerate}
	\cfout.
\end{athm}

\begin{proof}[Proof of \cref{Lemma:PropertiesOfParallelizationEqualLength}]
	\Nobs that, \eg, Grohs et al.~\cite[Lemma~2.18]{GrohsHornungJentzen2019}  \prove[ei]\ \cref{ParallelizationElementary:Display}. \Nobs that, \eg,
	Grohs et al.~\cite[Proposition~2.19]{GrohsHornungJentzen2019} \prove\ \cref{PropertiesOfParallelizationEqualLength:ItemOne,PropertiesOfParallelizationEqualLength:ItemTwo}.
	The proof of \cref{Lemma:PropertiesOfParallelizationEqualLength} is thus complete.
\end{proof}
\cfclear

\subsection{Identity ANNs}
\label{subsection:idANNs}

\begin{definition}
	[Identity ANNs]
	\label{def:ReLU_identity}
	\cfconsiderloaded{def:ReLU_identity}
	We denote by 
	$(\ReLUidANN{d})_{d \in \N} \subseteq \ANNSET$ 
	the \ANNs\ which satisfy 
	for all 
	$d \in \N\setminus\{1\}$ 
	that
	\begin{equation}
		\label{def:ReLU_identity_d_is_one}
		\begin{split}
			\ReLUidANN{1} = \pr*{ \!\pr*{\!\begin{pmatrix}
						1\\
						-1
					\end{pmatrix},
					\begin{pmatrix}
						0\\
						0
					\end{pmatrix}\! },
				\prbb{	\begin{pmatrix}
						1& -1
					\end{pmatrix}, 
					0 } \!
			}  \in \prb{(\R^{2 \times 1} \times \R^{2}) \times (\R^{1 \times 2} \times \R^1) }
		\end{split}
	\end{equation}
	and $\ReLUidANN{d} = \parallelizationSpecial_{d} (\ReLUidANN{1},\ReLUidANN{1},\dots, \ReLUidANN{1})$ \cfload.
\end{definition}
\cfclear

\begin{athm}{prop}{Prop:identity_representation}
	Let $d \in \N$. Then
	\begin{enumerate}[(i)]
		\item \label{identity_representation:2} it holds that $\realisation(\ReLUidANN{d}) = \id_{\R^d}$,
		\item \label{identity_representation:1} it holds that $\dims(\ReLUidANN{d})= (d, 2d, d)$, and
		\item \label{identity_representation:3} it holds that $\paramANN(\ReLUidANN{d}) = 4d^2+3d$
	\end{enumerate}
	\cfout.
\end{athm}

\begin{proof}
	[Proof of \cref{Prop:identity_representation}]
\Nobs that, \eg, Beneventano et al.~\cite[Proposition~2.13]{beneventano21}  \prove[ei]\ \cref{identity_representation:1,identity_representation:2,identity_representation:3}.
	The proof of \cref{Prop:identity_representation} is thus complete.
\end{proof}
\cfclear

\subsection{Compositions of ANNs}
\label{subsectio:compositions}

\begin{definition}[Compositions of ANNs]
	\label{def:ANNcomposition}
	\cfconsiderloaded{def:ANNcomposition}
	We denote by $\compANN{(\cdot)}{(\cdot)}\colon\allowbreak \{\f=(\mathscr{f}_1,\mathscr{f}_2)\allowbreak\in\ANNSET\times \ANNSET\colon \inDimANN(\mathscr{f}_1)=\outDimANN(\mathscr{f}_2)\}\allowbreak\to\ANNSET$ the function which satisfies for all 
	$ L,\mathfrak{L}\in\N$, $\f_1,\f_2\in\ANNSET$, $k\in\N\cap(0,L+\mathfrak{L})$ 
	with $\inDimANN(\mathscr{f}_1)=\outDimANN(\mathscr{f}_2)$, $\lengthANN(\f_2)=L$, and $\lengthANN(\f_1)=\mathfrak{L}$
	that $\lengthANN(\compANN{\mathscr{f}_1}{\mathscr{f}_2})=L+\mathfrak{L}-1$ and
		\begin{equation}\label{ANNoperations:Composition}
		\begin{split}
			&\prb{\weight{k}{\compANN{\mathscr{f}_1}{\mathscr{f}_2}},\bias{k}{\compANN{\mathscr{f}_1}{\mathscr{f}_2}}}=
			\begin{cases} 
				\prb{\weight{k}{\f_2},\bias{k}{\f_2}}&\colon k < L\\
			\prb{\weight{1}{\f_1}\weight{L}{\f_2},\weight{1}{\f_1}\bias{L}{\f_2}+\bias{1}{\f_1}}&\colon k = L\\
				\prb{\weight{k-L+1}{\f_1},\bias{k-L+1}{\f_1}}&\colon k > L
			\end{cases}
		\end{split}
	\end{equation}
	\cfload.
\end{definition}
\cfclear

\cfclear
\begin{athm}{lemma}{Lemma:CompositionAssociative}
	Let 
	$\mathscr{f}_1,\mathscr{f}_2,\mathscr{f}_3\in\ANNSET$
	satisfy
	$\inDimANN(\mathscr{f}_1)=\outDimANN(\mathscr{f}_2)$ and
	$\inDimANN(\mathscr{f}_2)=\outDimANN(\mathscr{f}_3)$ 
	\cfload.
	Then $
	\compANN{(\compANN{\mathscr{f}_1}{\mathscr{f}_2})}{\mathscr{f}_3}=\compANN{\mathscr{f}_1}{(\compANN{\mathscr{f}_2}{\mathscr{f}_3})} $
	\cfout.
\end{athm}
\begin{aproof}
    \Nobs that, \eg, Grohs et al.~\cite[Lemma~2.8]{GrohsHornungJentzen2019} \prove\
    $
	\compANN{(\compANN{\mathscr{f}_1}{\mathscr{f}_2})}{\mathscr{f}_3}=\compANN{\mathscr{f}_1}{(\compANN{\mathscr{f}_2}{\mathscr{f}_3})}$.
    \finishproofthus
\end{aproof}

\cfclear
\begin{athm}{prop}{Lemma:PropertiesOfCompositions_n2}
   Let $n \in \N\setminus\{1\}$,
	$\mathscr{f}_1, \f_2,  \ldots, \mathscr{f}_n \in \ANNSET$ 
	satisfy for all $k\in \N\cap(1,n]$ that
	$\inDimANN(\mathscr{f}_{k-1})=\outDimANN(\mathscr{f}_{k})$
	\cfload.
	Then
\begin{enumerate}[(i)]
	\item \label{PropertiesOfCompositions_n:Realization} 
	it holds that
	$
	\realisation({\mathscr{f}_1}\bullet{\mathscr{f}_2}\bullet \ldots \bullet \mathscr{f}_n)=[\realisation(\mathscr{f}_1)]\circ [\realisation(\mathscr{f}_2)]\circ \ldots \circ [\realisation(\mathscr{f}_n)]
	$ and
\item \label{PropertiesOfCompositions_n:Length}  it holds that
$\hidlengthANN(\compANN{\mathscr{f}_1}{{\mathscr{f}_2}}\bullet \ldots \bullet \mathscr{f}_n)=\sum_{k=1}^n\hidlengthANN(\mathscr{f}_k)$
\cfadd{Lemma:CompositionAssociative}
\end{enumerate}
\cfout.
\end{athm}

\begin{proof}
[Proof of \cref{Lemma:PropertiesOfCompositions_n2}]
Note that, \eg, Beneventano et al.~\cite[Proposition~2.16]{beneventano21} (see, \eg, also Grohs et al.~\cite[Proposition~2.6]{GrohsHornungJentzen2019})  \prove[ei]\ \cref{PropertiesOfCompositions_n:Length,PropertiesOfCompositions_n:Realization}.
The proof of \cref{Lemma:PropertiesOfCompositions_n2} is thus complete.
\end{proof}

\subsection{Sizes of parameters of ANNs}
\label{subsection:size}

\begin{definition}[Supremum norm]
	\label{def:matrixnorm}
	We denote by
	$\supn{\cdot}\colon\pr{\cup_{k,l\in\N}\R^{k\times l}}\to\R$ 
	the function which satisfies for all $k,l\in\N$, $W=(W_{i,j})_{(i,j)\in\{1,2,\ldots,k\}\times\{1,2,\ldots,l\}}\in\R^{k\times l}$ that 
	\begin{equation}
		\label{eq:supnorm}
		\supn{W}=\max_{i\in\{1,2,\ldots,k\}}\max_{j\in\{1,2,\ldots,l\}}\vass{W_{i,j}}.
	\end{equation}
\end{definition}

\begin{definition}[Sizes of parameters of ANNs]
	\label{def:size}
	\cfconsiderloaded{def:size}
	We denote by $\mathbb{S}_r\colon\ANNSET\to\R$, $r\in\{0,1\}$,
	the functions which satisfy
	for all $r\in\{0,1\}$,  
	$
	\mathscr{f} 
	\in  \ANNSET
	$ 
	that
	\begin{equation}
		\label{eq:size:def:IO}
		\begin{split}
			\mathbb{S}_{r}(\f)=\max\{\supn{\weight{r\hidlengthANN(\f)+1}{\f}},\supn{\bias{r\hidlengthANN(\f)+1}{\f}}\}
		\end{split}
	\end{equation}
and we denote by $\size\colon\ANNSET\to\R$ the function which satisfies for all $\f\in\ANNSET$ that
	\begin{equation}
	\label{eq:size:def}
	\begin{split}
		\size(\f)=\max_{k\in\{1,2,\ldots,\lengthANN(\f)\}}\max\{\supn{\weight{k}{\f}},\supn{\bias{k}{\f}}\}
	\end{split}
\end{equation}
	\cfload.
\end{definition}
\cfclear

\begin{athm}{lemma}{Lemma:ParallelizationElementary}[Sizes of ANN parameters of parallelizations]
	\cfconsiderloaded{Lemma:ParallelizationElementary}
	Let $n\in\N$, $\f_1,\f_2,\ldots,\mathscr{f}_n\in\ANNSET$ satisfy $\lengthANN(\mathscr{f}_1)= \lengthANN(\mathscr{f}_2)=\ldots =\lengthANN(\mathscr{f}_n)$
	\cfload[.]Then
	\begin{enumerate}[(i)]
		\item{\label{size:parallelization:overall}it holds that $\size\pr*{\parallelizationSpecial_{n}(\mathscr{f}_1,\f_2,\ldots,\mathscr{f}_n)}=\max\{\size(\f_1),\size(\f_2),\dots,\size(\f_n)\}$ and}
		\item{\label{size:parallelization:in}it holds for all $r\in\{0,1\}$ that $\rsize{r}(\parallelizationSpecial_{n}(\mathscr{f}_1,\f_2,\ldots,\mathscr{f}_n))=\max\{\rsize{r}(\f_1),\rsize{r}(\f_2),\dots,\rsize{r}(\f_n)\}$}
	\end{enumerate}
	\cfout[.]
\end{athm}

\begin{proof}[Proof of \cref{Lemma:ParallelizationElementary}]	
	\Nobs that \eqref{parallelisationSameLengthDef}  \proves[ei]\ \cref{size:parallelization:overall,size:parallelization:in}. 
	The proof of \cref{Lemma:ParallelizationElementary} is thus complete.
\end{proof}
\cfclear

\begin{athm}{cor}{Prop:identity_representationSIZE}[Sizes of identity ANNs]
	Let $d \in \N$. Then $\size(\ReLUidANN{d}) = \insize(\ReLUidANN{d})=\outsize(\ReLUidANN{d})=1$
	\cfout.
\end{athm}

\begin{proof}
	[Proof of \cref{Prop:identity_representationSIZE}]
	\Nobs that \eqref{def:ReLU_identity_d_is_one} and \cref{Lemma:ParallelizationElementary}  \prove[ei]\ $\size(\ReLUidANN{d}) = \insize(\ReLUidANN{d})=\outsize(\ReLUidANN{d})=1$.
	The proof of \cref{Prop:identity_representationSIZE} is thus complete.
\end{proof}
\cfclear

\cfclear
\begin{athm}{lemma}{lem:sizecomp}[Sizes of ANN parameters of compositions]
	Let $d\in\N$. Then 
	\begin{enumerate}[(i)]
		\item{
			\label{direct:comp:size}
			it holds for all $\f,\g\in\ANNSET$ with $d=\inDimANN(\f)=\outDimANN(\g)$ that 
			\begin{equation}
				\size(\compANN{\f}{\g})\leq \max\{\size(\f),\size(\g),\insize(\f)\outsize(\g)d+\insize(\f)\}
			\end{equation}
			and
		}
		\item{
			\label{direct:comp:insize}
			it holds for all $r\in\{0,1\}$, $\f_0,\f_1\in\ANNSET$ with $\inDimANN(\f_1)=\outDimANN(\f_0)$ and $\lengthANN(\f_r)> 1$ that
			\begin{equation}
				\rsize{r}(\compANN{\f_1}{\f_0})= \rsize{r}(\f_r)
			\end{equation}
		}
	\end{enumerate}
	\cfout.
\end{athm}
\begin{aproof}
	\Nobs that \eqref{eq:supnorm} \proves\ that for all $m,n\in\N$, $W\in\R^{m\times d}$, $B\in\R^{m}$, $\mathfrak{W}\in\R^{d\times n}$, $\mathfrak{B}\in\R^d$ it holds that
	\begin{equation}
		\supn{W\mathfrak{W}}\leq d\supn{W}\supn{\mathfrak{W}}
		\qandq
		\supn{W\mathfrak{B}+B}\leq d\supn{W}\supn{\mathfrak{B}}+\supn{B}
	\end{equation}
	\cfload. Combining this with \eqref{ANNoperations:Composition} and \eqref{eq:size:def} \proves\ that for all  $\f,\g\in\ANNSET$ with $d=\inDimANN(\f)=\outDimANN(\g)$ it holds that
	\begin{equation}
		\size(\compANN{\f}{\g})\leq \max\{\size(\f),\size(\g),\insize(\f)\outsize(\g)d+\insize(\f)\}
	\end{equation}
	\cfload. This  \proves[ei]\ \cref{direct:comp:size}.
	\Nobs that \eqref{ANNoperations:Composition} and \eqref{eq:size:def:IO} \prove\ that for all $r\in\{0,1\}$, $\f_0,\f_1\in\ANNSET$ with $\inDimANN(\f_1)=\outDimANN(\f_0)$ and $\lengthANN(\f_r)> 1$ it holds that
	\begin{equation}
		\rsize{r}(\compANN{\f_1}{\f_0})= \rsize{r}(\f_r).
	\end{equation}
	This  \proves[ei]\ \cref{direct:comp:insize}.
	\finishproofthus
\end{aproof}

\cfclear

\begin{athm}{prop}{Prop:identity_representation:prop}[Sizes of ANN parameters of compositions]
Let $d \in \N$, $\f,\g\in\ANNSET$ satisfy $\inDimANN(\f)=d=\outDimANN(\g)$ \cfload. Then
\begin{enumerate}[(i)]
    \item \label{identity_representation:prop:1} it holds that $\size(\compANN{\ReLUidANN{d}}{\g})= \max\{1,\size(\g)\}$,
    \item \label{identity_representation:prop:2} it holds that $\insize(\compANN{\ReLUidANN{d}}{\g})= \insize(\g)$,
    \item \label{identity_representation:prop:3} it holds that $\size(\compANN{\f}{\ReLUidANN{d}})= \max\{1,\size(\f)\}$,
		\item \label{identity_representation:prop:4} it holds that $\outsize(\compANN{\f}{\ReLUidANN{d}})= \outsize(\f)$,
		\item \label{identity_representation:prop:5} it holds that $\insize(\compANN{\f}{\ReLUidANN{d}})=\outsize(\compANN{\ReLUidANN{d}}{\g})=1$, and
		\item \label{identity_representation:prop:6} it holds that $\size(\compANN{\f}{\ReLUidANN{d}}\bullet\g)= \max\{\size(\f),\size(\g)\}$
\end{enumerate}
\cfout.
\end{athm}

\begin{proof}[Proof of \cref{Prop:identity_representation:prop}]
Throughout this proof 
let $A\in\R^{2d\times d}$ satisfy
\begin{equation}
\label{eq:Matrix:identity:Relu}
A=
\begin{pmatrix}
1 & 0 & \cdots & 0\\
-1 & 0 & \cdots & 0\\
0 & 1 & \cdots & 0\\
0 & -1 & \cdots & 0\\
\vdots & \vdots & \ddots & \vdots\\
0 & 0 & \cdots&1\\
0 & 0 & \cdots&-1
\end{pmatrix}
\ifnocf.
\end{equation}
\cfload[.]\Nobs that \eqref{eq:Matrix:identity:Relu} \proves\ that
\begin{equation}
\label{eq:left:comp:reluidann}
\supn{A\weight{\lengthANN(\g)}{\g}}=\supn{\weight{\lengthANN(\g)}{\g}}
\qandq
\supn{A\bias{\lengthANN(\g)}{\g}+\bias{1}{\ReLUidANN{d}}}=\supn{\bias{\lengthANN(\g)}{\g}}
\ifnocf.
\end{equation}
\cfload[.]\Moreover \eqref{def:ReLU_identity_d_is_one} and \eqref{eq:Matrix:identity:Relu} \prove\ that 
\begin{equation}
\label{eq:reluidann:explicit}
\ReLUidANN{d}=\prb{\prb{A,0},\prb{A^{*},0}}\in \prb{(\R^{2d \times d} \times \R^{2d}) \times (\R^{d \times 2d} \times \R^d) }.
\end{equation}
This and \eqref{ANNoperations:Composition} \prove\ that for all $k\in\{1,2,\ldots,\lengthANN(\g)+1\}$ it holds that
\begin{equation}
	\label{eq:weights:biases:Id:comp}
	\begin{split}
		&\prb{\weight{k}{(\compANN{\ReLUidANN{d}}{\g})},\bias{k}{(\compANN{\ReLUidANN{d}}{\g})}}=
		\begin{cases} 
			\prb{\weight{k}{\g},\bias{k}{\g}}&\colon k < \lengthANN(\g)\\
			\prb{A\weight{\lengthANN(\g)}{\g},A\bias{\lengthANN(\g)}{\g}}&\colon k = \lengthANN(\g)\\
			\prb{A^{*},0}&\colon k = \lengthANN(\g)+1.
		\end{cases}
	\end{split}
\end{equation}
Combining this and \eqref{eq:left:comp:reluidann} with \eqref{eq:size:def:IO} \proves[ei]\ that
\begin{equation}
	\insize(\compANN{\ReLUidANN{d}}{\g})
	=\max\{\supn{\weight{1}{(\compANN{\ReLUidANN{d}}{\g})}},\supn{\bias{1}{(\compANN{\ReLUidANN{d}}{\g})}}\}
	=\max\{\supn{\weight{1}{\g}},\supn{\bias{1}{\g}}\}
	=\insize(\g)
	.
\end{equation} This, 
\eqref{eq:Matrix:identity:Relu},
\eqref{eq:left:comp:reluidann}, and 
\eqref{eq:weights:biases:Id:comp} \prove[ei]\ \cref{identity_representation:prop:1,identity_representation:prop:2}.
\Nobs that \eqref{eq:Matrix:identity:Relu} and \eqref{eq:reluidann:explicit} \prove\ that
\begin{equation}
\label{eq:right:comp:reluidann}
\supn*{\weight{1}{\f}\pr*{A^{*}}}=\supn{\weight{1}{\f}}
\qandq
\supn{\weight{1}{\f}\bias{2}{\ReLUidANN{d}}+\bias{1}{\f}}=\supn{\bias{1}{\f}}
.
\end{equation}
\Moreover \eqref{eq:reluidann:explicit} and \eqref{ANNoperations:Composition} \prove\ that for all $k\in\{1,2,\ldots,\lengthANN(\f)+1\}$ it holds that
\begin{equation}
		\label{eq:weights:biases:comp:Id}
			\begin{split}
		&\prb{\weight{k}{(\compANN{\f}{\ReLUidANN{d}})},\bias{k}{(\compANN{\f}{\ReLUidANN{d}})}}=
		\begin{cases} 
			\prb{A,0}&\colon k=1\\
			\prb{\weight{1}{\f}\pr*{A^{*}},\bias{1}{\f}}&\colon k = 2\\
			\prb{\weight{k-1}{\f},\bias{k-1}{\f}}&\colon k > 2.
		\end{cases}
	\end{split}
\end{equation}
Combining this and \eqref{eq:right:comp:reluidann} with \eqref{eq:size:def:IO} \proves[ei]\ that
\begin{equation}
	\begin{split}
	\outsize(\compANN{\f}{\ReLUidANN{d}})
	&=\max\{\supn{\weight{\lengthANN(\compANN{\f}{\ReLUidANN{d}})}{(\compANN{\f}{\ReLUidANN{d}})}},\supn{\bias{\lengthANN(\compANN{\f}{\ReLUidANN{d}})}{(\compANN{\f}{\ReLUidANN{d}})}}\}
	\\&=\max\{\supn{\weight{\lengthANN(\f)}{\f}},\supn{\bias{\lengthANN(\f)}{\f}}\}
	=\outsize(\f).
\end{split}
\end{equation} 
This,
\eqref{eq:Matrix:identity:Relu},
\eqref{eq:right:comp:reluidann}, and
\eqref{eq:weights:biases:comp:Id} \prove[ei]\ \cref{identity_representation:prop:3,identity_representation:prop:4}.
\Nobs that 
\cref{lem:sizecomp},
\cref{Prop:identity_representationSIZE},
\eqref{eq:Matrix:identity:Relu}, and 
\eqref{eq:reluidann:explicit} \prove\ that
\begin{equation}
\insize(\compANN{\f}{\ReLUidANN{d}})=\insize(\ReLUidANN{d})=1=\outsize(\ReLUidANN{d})=\outsize(\compANN{\ReLUidANN{d}}{\g}).
\end{equation}
This  \proves[ei]\ \cref{identity_representation:prop:5}.
\Nobs that \eqref{eq:reluidann:explicit} and \eqref{ANNoperations:Composition} \prove[pis]\ that for all $k\in\{1,2,\ldots,\lengthANN(\f)+\lengthANN(\g)\}$ it holds that
\begin{equation}
	\label{eq:weights:biases:comp:Id:comp}
	\begin{split}
		&\prb{\weight{k}{(\compANN{\f}{\ReLUidANN{d}}\bullet\g)},\bias{k}{(\compANN{\f}{\ReLUidANN{d}}\bullet\g)}}=
			\begin{cases} 
	\prb{\weight{k}{\g},\bias{k}{\g}}&\colon k < \lengthANN(\g)\\
	\prb{A\weight{\lengthANN(\g)}{\g},A\bias{\lengthANN(\g)}{\g}}&\colon k = \lengthANN(\g)\\
		\prb{\weight{1}{\f}\pr*{A^{*}},\bias{1}{\f}}&\colon k = \lengthANN(\g)+1\\
	\prb{\weight{k-\lengthANN(\g)}{\f},\bias{k-\lengthANN(\g)}{\f}}&\colon k > \lengthANN(\g)+1.
\end{cases}
	\end{split}
\end{equation}
This,
 \eqref{eq:left:comp:reluidann},
  and \eqref{eq:right:comp:reluidann} \prove\ that
\begin{equation}
	\begin{split}
\size(\compANN{\f}{\ReLUidANN{d}}\bullet\g)
&=
\max_{k\in\{1,2,\ldots,\lengthANN(\f)+\lengthANN(\g)\}}\max\{\supn{\weight{k}{(\compANN{\f}{\ReLUidANN{d}}\bullet\g)}},\supn{\bias{k}{(\compANN{\f}{\ReLUidANN{d}}\bullet\g)}}\}
\\&= 
\max\{\size(\f),\size(\g)\}.
\end{split}
\end{equation}
This  \proves[ei]\ \cref{identity_representation:prop:6}.
The proof of \cref{Prop:identity_representation:prop} is thus complete.
\end{proof}

\cfclear
\begin{athm}{prop}{Prop:identity_representation:prop2}[Sizes of ANN parameters of iterated compositions]
Let $d\in\N$, $a_1,a_2,\ldots,\allowbreak a_{d} \in \N$, $\h_0,\h_1,\ldots,\h_d\in \ANNSET$ satisfy for all $i\in\{1,2,\ldots,d\}$ that $\inDimANN(\h_{i})=a_{i}=\outDimANN(\h_{i-1})$
and let $\f\in\ANNSET$ satisfy
    \begin{equation}
    \label{fdef:compositionswithidANNs}
        \f
        =
        \compANN{\mathscr{h}_d}{\ReLUidANN{a_{d}}}\bullet\compANN{\mathscr{h}_{d-1}}{\ReLUidANN{a_{d-1}}}\bullet\ldots\bullet\mathscr{h}_1\bullet\ReLUidANN{a_1}\bullet\mathscr{h}_0
        \ifnocf.
    \end{equation}
\cfload. Then $\size(\f)\leq\max\{\size(\h_0),\size(\h_1),\ldots,\size(\h_d)\}$ \cfout.
\end{athm}

\begin{proof}[Proof of \cref{Prop:identity_representation:prop2}]
\Nobs that \cref{Prop:identity_representation:prop} and induction \prove\ that 
\begin{equation}
\size(\f)\leq\max\{\size(\h_0),\size(\h_1),\ldots,\size(\h_d)\}
\end{equation}
\cfload. The proof of \cref{Prop:identity_representation:prop2} is thus complete.
\end{proof}

\cfclear

%

\renewcommand{\shft}{\lambda}

\section{Lower bounds for the minimal number of ANN parameters in the approximation of certain high-dimensional functions}
\label{S4_lower_bounds}

In this section we establish in 
\cref{cor:lower:bound:prod:sized}, \cref{prop:product:unrestricted}, 
\cref{LowerBound:prod}, and \cref{LowerBound:sum} below suitable lower bounds for the minimal number of parameters of shallow or insufficiently deep \ANNs\ to approximate certain high-dimensional target functions.

Our proof of \cref{cor:lower:bound:prod:sized} uses appropriate lower bounds for the minimal number of \ANN\ parameters to approximate the product functions in \cref{lemma:lower:bound:prod:sized}.
We derive \cref{lemma:lower:bound:prod:sized} from the well-known upper bounds for the absolute values of realization functions of \ANNs\ in \cref{lemma:upper:bound:real}.
\Cref{lemma:upper:bound:real} is a slightly modified variant of, \eg, Grohs et al.~\cite[Corollary~4.3]{GrohsIbrgimovJentzen2021}.

\Cref{prop:product:unrestricted} complements \cref{cor:lower:bound:prod:sized} by showing that the number of non-zero \ANN\ parameters needed to approximate the product functions at fixed depth and accuracy grows at least exponentially, even when the parameters may be arbitrarily large.
The proof of \cref{prop:product:unrestricted} employs
\cref{lemma:realisation:shallow} to bound the number of affine pieces of
the \ANN\ realization on a linear segment.
Our proofs of \cref{LowerBound:prod} and \cref{LowerBound:sum} employ the lower bound result for certain families of oscillating functions in \cref{gen:final:approximation:fail}. A result related to \cref{gen:final:approximation:fail} can be found, \eg, in Telgarsky \cite[Theorem~1.1]{Telgarsky15}.
Our proof of \cref{gen:final:approximation:fail} uses the essentially well-known upper bound result for the number of certain linear regions of realization functions of \ANNs\ in \cref{lemma:realisation:shallow}.
In the scientific literature results related to \cref{lemma:realisation:shallow} can be found, \eg, in Raghu et al.~\cite[Theorem~1]{raghu2017expressive}.
Our proof of \Cref{lemma:realisation:shallow}, in turn, utilizes the elementary \ANN\ representation result in \cref{lemma:deep:linear:structure}.
\Cref{lemma:deep:linear:structure} builds on the elementary concepts and results regarding intersections of half-spaces in \cref{subsection:geomertry}.

\newcommand{\proj}{\mathfrak{p}}

\subsection{Lower bounds for approximations of product functions}
\label{subsection:lower_bound_product}

\newcommand{\form}[2]{\mathcal{T}_{#1}^{#2}\cfadd{def:vect}}
\newcommand{\row}[2]{\mathscr{r}_{#1,#2}\cfadd{def:row}}
\newcommand{\col}[1]{\mathfrak{C}_{#1}}
\newcommand{\non}{Q}
\newcommand{\sngp}[1]{\form{#1}{\f}}
\newcommand{\lay}{\mathfrak{l}}

\cfclear
\begin{definition}[Parameters of an ANN]
	\label{def:vect}
	\cfconsiderloaded{def:vect}
	Let $\f\in\ANNSET$. Then we denote by 
	$\form{(\cdot)}{\f}=(\form{n}{\f})_{n\in\{1,2,\ldots,\param(\f)\}}\colon\allowbreak\{1,2,\ldots,\param(\f)\}\to\R$ the function 
	which satisfies for all $k\in\{1,2,\ldots,\lengthANN(\f)\}$, $i,j,s\in\N_0$ with $i=\singledims_{k-1}(\f)$, $j=\singledims_k(\f)$, and $s=\ssum_{l=1}^{k-1}\singledims_l(\f)(\singledims_{l-1}(\f)+1)$ that
	\begin{equation}
		\begin{split} 	
			&\weight{k}{\f}=
			\begin{pmatrix}
				\sngp{s+1}&\sngp{s+2}& \cdots&\sngp{s+i}\\
				\sngp{s+i+1}&\sngp{s+i+2}& \cdots&\sngp{s+2i}\\
				\sngp{s+2i+1}&\sngp{s+2i+2}& \cdots&\sngp{s+3i}\\
				\vdots&\!\!\! \vdots&\!\!\! \ddots&\!\!\! \vdots\\
				\sngp{s+(j-1)i+1}&\sngp{s+(j-1)i+2}& \cdots&\sngp{s+ji}\\
			\end{pmatrix}
			\qandq
			\bias{k}{\f}=
			\begin{pmatrix}
				\form{s+ji+1}{\f}\\
				\form{s+ji+2}{\f}\\
				\vdots\\
				\form{s+ji+j}{\f}
			\end{pmatrix}
		\end{split}
	\end{equation}
	\cfload.
\end{definition}

\cfclear
\begin{definition}[Rows of ANN weight matrices]
	\label{def:row}
	\cfconsiderloaded{def:row}
	Let $\f\in\ANNSET$, $n\in\N\cap[1,\sum_{k=1}^{\lengthANN(\f)}\singledims_k(\f)]$. Then we denote by 
	$\row{n}{\f}\in\pr{\cup_{d\in\N}\R^d}$ the vector which satisfies for all
	$i\in\{1,2,\ldots,\lengthANN(\f)\}$,
	$j\in\{1,2,\ldots,\singledims_i(\f)\}$,
	$s\in\N_0$ with 
	$j+\sum_{k=1}^{i-1}\singledims_k(\f)=n$ and
	$s=(j-1)\singledims_{i-1}(\f)+\ssum_{k=1}^{i-1}\singledims_k(\f)(1+\singledims_{k-1}(\f))$ that
	\begin{equation}
		\row{n}{\f}=(\sngp{s+1},\sngp{s+2},\ldots,\sngp{s+\singledims_{i-1}(\f)})
		\in\R^{\singledims_{i-1}(\f)}
	\end{equation}
	\cfload.
\end{definition}

\cfclear
\begin{definition}[Non-zero ANN parameters]
	\label{def:paramnonzero}
	\cfconsiderloaded{def:paramnonzero}
	We denote by $\paramnotzero\colon\ANNSET\to\N_0$ the function which satisfies for all $\f\in\ANNSET$ that
	\begin{equation}
		\textstyle
		\paramnotzero(\f)=\sum_{n=1}^{\param(\f)}\indicator{\R\setminus\{0\}}(\form{n}{\f})
	\end{equation}
	\cfload.
\end{definition}

\cfclear
\begin{athm}{lemma}{lem:zero_weights}
	Let $\f\in\ANNSET$, $j\in\{1,2,\ldots,\lengthANN(\f)\}$ satisfy $\weight{j}{\f}=0$ \cfload.
	Then it holds for all $x,y\in\R^{\inDimANN(\f)}$ that $(\realisation(\f))(x)=(\realisation(\f))(y)$
	\cfout.
\end{athm}

\begin{proof}[Proof of \cref{lem:zero_weights}]
	Throughout this proof assume without loss of generality that $j<\lengthANN(\f)$ and
	let 
	$f_i\colon\R^{\inDimANN(\f)}\to\R^{\singledims_i(\f)}$, $i\in\{0,1,\ldots,\lengthANN(\f)\}$, satisfy for all $x\in\R^{\inDimANN(\f)}$ that $f_0(x)=x$, $f_{\lengthANN(\f)}(x)=\weight{\lengthANN(\f)}{\f} f_{\hidlengthANN(\f)}(x) + \bias{\lengthANN(\f)}{\f}$, and
	\begin{equation}
		\label{eq:layerrealizationfunction}
		\fa{k}\N\cap(0,\lengthANN(\f))\colon f_k(x)=\RELU(\weight{k}{\f} f_{k-1}(x) + \bias{k}{\f})
	\end{equation}
	\cfload. Combining \eqref{eq:layerrealizationfunction} with induction \proves\ that for all $i\in\N_0\cap[0,\lengthANN(\f))$, $k\in\N\cap[i,\lengthANN(\f)]$, $a,b\in\R^{\inDimANN(\f)}$ with $\fa{x,y}\R^{\inDimANN(\f)}\colon f_i(x)=f_i(y)$ it holds that
	\begin{equation}
		\label{eq:constant:layer:functions}
		f_k(a)=f_k(b).
	\end{equation} 
	\Moreover \eqref{eq:layerrealizationfunction} and the assumption that $\weight{j}{\f}=0$ \prove\ that for all $x,y\in\R^{\inDimANN(\f)}$ it holds that
	\begin{equation}
		f_j(x)=\RELU(\weight{j}{\f} f_{j-1}(x) + \bias{j}{\f})=\RELU(\bias{j}{\f})=\RELU(\weight{j}{\f} f_{j-1}(y) + \bias{j}{\f})=f_j(y).
	\end{equation}
	Combining this, \eqref{eq:layerrealizationfunction}, and \eqref{eq:constant:layer:functions} with
	\eqref{ANNrealization:ass2} \proves\ that for all $x,y\in\R^{\inDimANN(\f)}$ it holds that
	\begin{equation}
		\pr{\realisation(\f)}(x)
		=f_{\lengthANN(\f)}(x)
		=f_{\lengthANN(\f)}(y)
		=\pr{\realisation(\f)}(y)
	\end{equation}
\cfload.
	\finishproofthus
\end{proof}

\cfclear
\begin{athm}{lemma}{lemma:upper:bound:real}
Let $a\in\R$, $b\in[a,\infty)$, $\f\in\ANNSET$ \cfload. Then it holds for all $x\in[a,b]^{\inDimANN(\f)}$ that 
\begin{equation}
\label{eq:upper:bound:real}
\supn{(\realisation(\f))(x)}
\leq \PR*{\paramnotzero(\f)\max\{\size(\f),1\}}^{\lengthANN(\f)} \max\{\vass{a},\vass{b},1\}
\end{equation}
\cfout.
\end{athm}

\begin{proof}[Proof of \cref{lemma:upper:bound:real}]
Throughout this proof 
assume without loss of generality that $\outDimANN(\f)=1$, $\paramnotzero(\f)> 0$, and $\supn{\weight{\lengthANN(\f)}{\f}}>0$,
let $d,L\in\N$, $l_0,l_1,\ldots,\allowbreak l_L\in\N$, $x_0\in[a,b]^{\inDimANN(\f)}$, $x_1\in\R^{l_1}$, $x_2\in\R^{l_2}$, $\dots$, $x_L\in\R^{l_L}$
satisfy for all $k\in\{1,2,\ldots,L\}$ that
\begin{equation}
\label{eq:setup:netw:limited}
\inDimANN(\f)=d,
\qquad
\lengthANN(\f)=L,
\qandq
x_k=\RELU(\weight{k}{\f} x_{k-1}+\bias{k}{\f}),
\end{equation}
let $\lay_k\in\N_0$, $k\in\N_0\cap[0,L)$, and $\col{k,s}\in\R^{l_k}$, $s\in\{1,2,\ldots,l_{k-1}\}$, $k\in\{1,2,\ldots,L\}$, satisfy for all $k\in\{1,2,\ldots,L\}$ that 
\begin{equation}
	\label{eq:setup:netw:limited:2}
\lay_{k-1}
=\textstyle
\sum_{n=1}^{l_{k-1}}\indicator{\R^{l_{k}}\setminus\{0\}}\pr{\col{k,n}}
\qandq
\weight{k}{\f}=	
\begin{pmatrix}
	\col{k,1} & \col{k,2} & \dots & \col{k,l_{k-1}}
\end{pmatrix}
\end{equation}
\cfload. \Nobs that \eqref{eq:setup:netw:limited} and \eqref{eq:setup:netw:limited:2} \prove\ that for all $y\in[a,b]^{\inDimANN(\f)}$ with $\paramnotzero(\f)= \ssum_{k=0}^{L-1}\lay_{k}$
it holds that
\begin{equation}
	\label{eq:worst:fit:non.zero:weights}
	\begin{split}
	\vass{(\realisation(\f))(y)}
	&\leq \PR*{\sprod_{k=0}^{L-1}\lay_k \size(\f)}\supn{y}
	\\&\leq \PR*{\max\{\lay_0,\lay_1,\ldots,\lay_{L-1}\}\size(\f)}^L\supn{y}
	\\&\leq \PR*{\paramnotzero(\f)\max\{\size(\f),1\}}^{\lengthANN(\f)} \max\{\vass{a},\vass{b},1\}.
	\end{split}
\end{equation}
\Moreover \eqref{eq:setup:netw:limited}, and \eqref{eq:setup:netw:limited:2} \prove\ that for all $k\in\{1,2,\ldots,L\}$ it holds that
\begin{equation}
\label{eq:size:step}
\begin{split}
\supn{x_k}
=\supn{\RELU(\weight{k}{\f} x_{k-1}+\bias{k}{\f})}
&\leq\supn{\weight{k}{\f} x_{k-1}}+\supn{\bias{k}{\f}}
\\&\leq \lay_{k-1}\supn{\weight{k}{\f}}\supn{x_{k-1}}+\supn{\bias{k}{\f}}
\\&\leq \lay_{k-1}\size(\f)\supn{x_{k-1}}+\size(\f)
\\&\leq \lay_{k-1}\size(\f)\max\{\supn{x_{k-1}},1\}+\size(\f)
\\&= \size(\f)\pr{\lay_{k-1}\max\{\supn{x_{k-1}},1\}+1}
\\&\leq \pr{\lay_{k-1}+1} \size(\f)\max\{\supn{x_{k-1}},1\}
\end{split}
\end{equation}
\cfload. Combining this and \eqref{eq:setup:netw:limited} with induction \proves\ that
\begin{equation}
	\label{eq:last:layer:estim}
\begin{split}
\vass{(\realisation(\f))(x_0)}
=\vass{\weight{L}{\f} x_{L-1} + \bias{L}{\f}}
&\leq \pr{\lay_{L-1}+1}\size(\f)\max\{\supn{x_{L-1}},1\}
\\&\leq \pr*{\sprod_{k=0}^{L-1}\pr{\lay_{k}+1}}\max\{\size(\f),1\}^L  \max\{\supn{x_{0}},1\}
\\&\leq \pr*{\sprod_{k=0}^{L-1}\pr{\lay_{k}+1}}\max\{\size(\f),1\}^L  \max\{\vass{a},\vass{b},1\}.
\end{split}
\end{equation}
This and the inequality of arithmetic and geometric means \prove\ that for all $y\in[a,b]^{\inDimANN(\f)}$ with $\ssum_{k=0}^{L-1}\lay_{k}<\paramnotzero(\f)$ and $y=x_0$ it holds that
\begin{equation}
\begin{split}
\vass{(\realisation(\f))(y)}
&\leq \pr*{\sprod_{k=0}^{L-1}\pr{\lay_{k}+1}}\max\{\size(\f),1\}^L \max\{\vass{a},\vass{b},1\}
\\&\leq \pr*{L^{-1}\ssum_{k=0}^{L-1}\pr{\lay_{k}+1}}^L\max\{\size(\f),1\}^L \max\{\vass{a},\vass{b},1\}
\\&\leq \pr*{1+\ssum_{k=0}^{L-1}\lay_{k}}^L\max\{\size(\f),1\}^L \max\{\vass{a},\vass{b},1\}
\\&\leq \pr*{\paramnotzero(\f)}^L\max\{\size(\f),1\}^L \max\{\vass{a},\vass{b},1\}
\\&=\PR*{\paramnotzero(\f)\max\{\size(\f),1\}}^L \max\{\vass{a},\vass{b},1\}
.
\end{split}
\end{equation}
This and \eqref{eq:worst:fit:non.zero:weights} \prove[ei]\ \eqref{eq:upper:bound:real}.
\finishproofthus
\end{proof}

\cfclear
\begin{athm}{lemma}{lemma:lower:bound:prod:gen}
	Let $a\in\R$, $b\in[a,\infty)$, $\f\in\ANNSET$, $d,L\in\N$, $\eps\in(0,\infty)$ satisfy
	\begin{equation}
		\label{eq:assumptions:prod:approx:lem}
		\textstyle
		\lengthANN(\f)\leq L,
		\qquad
		\realisation(\f)\in C(\R^d,\R),
		\qquad
		\sup_{x=(x_1,\ldots,x_d)\in[a,b]^d}\vass{(\realisation(\f))(x)-\prod_{i=1}^d x_i}\leq\eps,
	\end{equation} and $\max\{\vass{a},\vass{b}\}^d>\eps$ \cfload.
	Then
	\begin{equation}
		\label{eq:lower:bound:prod:gen}
		\paramnotzero(\f)\max\{\size(\f),1\}
		\geq
	\pr*{\frac{\max\{\vass{a},\vass{b}\}^{d}-\eps}{\max\{\vass{a},\vass{b},1\}}}^{\!\!\nicefrac{1}{L}}
	\end{equation}
	\cfout.
\end{athm}

\begin{proof}[Proof of \cref{lemma:lower:bound:prod:gen}]
	Throughout this proof
	let $x=(x_1,\ldots,x_d)\in[a,b]^d$ satisfy $\vass{\prod_{i=1}^d x_i}=\allowbreak \max\{\vass{a},\vass{b}\}^d$. \Nobs that \eqref{eq:assumptions:prod:approx:lem}, the fact that $\paramnotzero(\f) \max\{\size(\f),1\}\in\{0\}\cup[1,\infty)$, and \cref{lemma:upper:bound:real} (applied with
	$a \curvearrowleft a$,
	$b \curvearrowleft b$,
	$\f \curvearrowleft \f$
	in the notation of \cref{lemma:upper:bound:real}) \prove\ that
	\begin{equation}
		\begin{split}
			\max\{\vass{a},\vass{b}\}^d
			=\textstyle \vass{\prod_{i=1}^d x_i}
			&\leq \vass{(\realisation(\f))(x)}+\eps
			\\&\leq  \PR{\paramnotzero(\f) \max\{\size(\f),1\}}^{\lengthANN(\f)} \max\{\vass{a},\vass{b},1\}+\eps
			\\&\leq  \PR{\paramnotzero(\f) \max\{\size(\f),1\}}^{L} \max\{\vass{a},\vass{b},1\}+\eps
		\end{split}
	\end{equation}
	\cfload. This and the fact that $\max\{\vass{a},\vass{b}\}^d>\eps$ \prove\ that 
	\begin{equation}
		\paramnotzero(\f)\max\{\size(\f),1\}\geq
		\pr*{\frac{\max\{\vass{a},\vass{b}\}^{d}-\eps}{\max\{\vass{a},\vass{b},1\}}}^{\!\!\nicefrac{1}{L}}.
	\end{equation}
This \proves\ \eqref{eq:lower:bound:prod:gen}.
	\finishproofthus
	\cfload
\end{proof}

\cfclear
\begin{athm}{lemma}{lemma:lower:bound:prod:sized}
Let $a\in\R$, $b\in[a,\infty)$, $c\in[1,\infty)$, $d,L\in\N$, $\eps\in(0,\infty)$ satisfy 
$ \max\{\vass{a},\vass{b}\}>\max\{1,\eps\} $,
let $f\colon[a,b]^d\to\R$ satisfy for all $x=(x_1,\ldots,x_d)\in[a,b]^d$ that $f(x)=\prod_{i=1}^d x_i$,  
 and let $\f\in\ANNSET$ satisfy for all $x\in[a,b]^d$ that
 \begin{equation}
 	\textstyle
 	\size(\f)\leq c d^c,
 	\quad
 	\lengthANN(\f)\leq L,
 	\quad
 	\realisation(\f)\in C(\R^d,\R),
 	\quad\text{and}\quad 
 	\vass{(\realisation(\f))(x)-f(x)}\leq\eps
 \end{equation}\cfload.
Then  
\begin{equation}
\label{eq:lower:bound:prod:sized}
\begin{split}
&\paramnotzero(\f)
\geq
\pr*{\max\{\vass{a},\vass{b}\}^{d}-\eps}^{\!\!\nicefrac{1}{L}}\pr{\max\{\vass{a},\vass{b}\} cd^{c}}^{-1}
\end{split}
\end{equation}
\cfout.
\end{athm}
\begin{proof}[Proof of \cref{lemma:lower:bound:prod:sized}]
	Throughout this proof let $\dom\in(1,\infty)$ satisfy $\dom=\max\{\vass{a},\vass{b}\}$.
Combining this with \cref{lemma:lower:bound:prod:gen} (applied with
$a \curvearrowleft a$,
$b \curvearrowleft b$,
$\f \curvearrowleft \f$,
$d \curvearrowleft d$,
$L \curvearrowleft L$,
$\eps \curvearrowleft \eps$ in the notation of \cref{lemma:lower:bound:prod:gen}) \proves\ that 
\begin{equation}
\begin{split}
	\paramnotzero(\f)
	\geq 
	\pr*{\frac{\dom^{d}-\eps}{\dom}}^{\!\!\nicefrac{1}{L}}\max\{\size(\f),1\}^{-1}
	\geq 
	\pr*{\dom^{d}-\eps}^{\!\!\nicefrac{1}{L}}\pr{\dom cd^{c}}^{-1}
\end{split}
\end{equation}
\cfload. This \proves[ei]\ \eqref{eq:lower:bound:prod:sized}.
\finishproofthus
\cfload
\end{proof}

\cfclear

\begin{athm}{cor}{cor:lower:bound:prod:sized}
Let $a\in\R$, $b\in[a,\infty)$, $c\in[3,\infty)$, $\eps,R\in(0,\infty)$, $d,L\in\N$ satisfy $  \dom=\max\{\vass{a},\vass{b}\}>\max\{1,\eps\}$ and
let $f\colon[a,b]^d\to\R$ satisfy for all $x=(x_1,\ldots,x_d)\in[a,b]^d$ that $f(x)=\prod_{i=1}^d x_i$ \cfload.
Then  
	\begin{multline}
			\label{eq:lower:bound:prod:sized:nice}
		\min \pr*{ \pR*{ p \in \N \colon \PR*{\!\!
					\begin{array}{c}
						\exists\, \mathscr{f} \in \ANNSET \colon
						(\paramnotzero(\mathscr{f})\leq p)\land
						(\lengthANN(\mathscr{f})\leq L)
						\land{}\\
						(\size(\mathscr{f})\leq c d^c)\land{}
						(\realisation(\mathscr{f}) \in C(\R^d,\R))\land{}\\
						(\sup_{x\in[a,b]^d}\vass{(\realisation(\mathscr{f}))(x)-f(x)} \leq \varepsilon)\\
					\end{array} \!\! }
			}
			\cup\{\infty\}
		}
	\\\geq (\dom^{\nicefrac{d}{2}}-\eps^{\nicefrac{1}{2}})^{\nicefrac{1}{L}}\tfrac{\ln(\dom)^c}{\dom c^{c+1} L^c}
	\end{multline}
\cfout.
\end{athm}

\begin{proof}[Proof of \cref{cor:lower:bound:prod:sized}]
Throughout this proof 
let $g\colon(0,\infty)\to\R$ satisfy for all $x\in(0,\infty)$ that
\begin{equation}
	\label{lower:bound:param:function}
	\tgt(x)= \dom^{\frac{x}{2L}} x^{-c}.
\end{equation}
\Nobs that \eqref{lower:bound:param:function} \proves\ that for all $x\in(0,\infty)$ it holds that
\begin{equation}
	\begin{split}
	g'(x)
	=\ln\prb{\dom}(2L)^{-1}\dom^{\frac{x}{2L}} x^{-c}+\dom^{\frac{x}{2L}} (-c)x^{-c-1}
	= \dom^{\frac{x}{2L}}x^{-c-1} \prb{\tfrac{x\ln(\dom)}{2L}-c}.
	\end{split}
\end{equation}
This and the assumption that $\dom> 1$ \prove\ that for all $x\in\prb{0,\frac{2cL}{\ln(\dom)}}$, $y\in\prb{\frac{2cL}{\ln(\dom)},\infty}$ it holds that
\begin{equation}
	g'(x)<0,
	\qquad
	g'\prb{\tfrac{2cL}{\ln(\dom)}}=0,
	\qandq
	g'(y)>0.
\end{equation}
Combining this and \eqref{lower:bound:param:function} with the fact that $2\leq e$ \proves\ that
\begin{equation}
	\label{inf:of:param:function}
	\begin{split}
	&	\min_{x\in(0,\infty)} \tgt(x)
	= g\prb{\tfrac{2cL}{\ln(\dom)}}
	= \dom^{\frac{c}{\ln(\dom)}} \prb{\tfrac{2cL}{\ln(\dom)}}^{-c}
	= \prb{\tfrac{e\ln(\dom)}{2cL}}^{c}
	\geq \prb{\tfrac{\ln(\dom)}{cL}}^{c}.
	\end{split}
\end{equation}
This, the fact that $\dom^d-\eps>0$, and \eqref{lower:bound:param:function} \prove\ that 
\begin{equation}
	\begin{split}
		\pr*{\dom^{d}-\eps}^{\!\!\frac{1}{L}}\pr{\dom c d^{c}}^{-1}
&= (\dom^{\frac{d}{2}}-\eps^{\frac{1}{2}})^{\frac{1}{L}}\pr{\dom c }^{-1}(\dom^{\frac{d}{2}}+\eps^{\frac{1}{2}})^{\frac{1}{L}}d^{-c}
		\\&\geq
		(\dom^{\frac{d}{2}}-\eps^{\frac{1}{2}})^{\frac{1}{L}}\pr{\dom c }^{-1}\dom^{\frac{d}{2L}}d^{-c}
		\\&=(\dom^{\frac{d}{2}}-\eps^{\frac{1}{2}})^{\frac{1}{L}}\pr{\dom c }^{-1}\tgt(d)
		\\&\geq (\dom^{\frac{d}{2}}-\eps^{\frac{1}{2}})^{\frac{1}{L}}\pr{\dom c }^{-1}\prb{\tfrac{\ln(\dom)}{cL}}^{c}
		\\&\geq (\dom^{\frac{d}{2}}-\eps^{\frac{1}{2}})^{\frac{1}{L}}\tfrac{\ln(\dom)^c}{\dom c^{c+1} L^c}
		.
	\end{split}
\end{equation}
\Nobs that \cref{lemma:lower:bound:prod:sized} (applied with
$a \curvearrowleft a$,
$b \curvearrowleft b$,
$c \curvearrowleft c$,
$\eps \curvearrowleft \eps$,
$d \curvearrowleft d$,
$L \curvearrowleft L$,
$f \curvearrowleft f$ in the notation of \cref{lemma:lower:bound:prod:sized}) \hence \proves\ that
\begin{multline}
	\min \pr*{ \pR*{ p \in \N \colon \PR*{\!\!
				\begin{array}{c}
					\exists\, \mathscr{f} \in \ANNSET \colon
					(\paramnotzero(\mathscr{f})\leq p)\land
					(\lengthANN(\mathscr{f})\leq L)
					\land{}\\
					(\size(\f)\leq c d^c)\land{}
					(\realisation(\mathscr{f}) \in C(\R^d,\R))\land{}\\
					(\sup_{x\in[a,b]^d}\vass{(\realisation(\mathscr{f}))(x)-f(x)} \leq \varepsilon)\\
				\end{array} \!\! }
		}
		\cup\{\infty\}
	}
\\\geq 
\pr*{\dom^{d}-\eps}^{\!\!\nicefrac{1}{L}}\pr{\dom c d^{c}}^{-1}
\geq (\dom^{\nicefrac{d}{2}}-\eps^{\nicefrac{1}{2}})^{\nicefrac{1}{L}}\tfrac{\ln(\dom)^c}{\dom c^{c+1} L^c}
\end{multline}
\cfload.
This  \proves[ei]\ \eqref{eq:lower:bound:prod:sized:nice}.
	\finishproofthus
\end{proof}

\subsection{Intersections of half-spaces}
\label{subsection:geomertry}


\cfclear
\begin{definition}[Spaces of affine linear functions]
	\label{def:Lin}
	Let $d\in\N$ and let $D\subseteq\R^d$ be a set. Then we denote by $\Lin(D)\subseteq C(D,\R)$ the set given by
	\begin{equation}
		\label{eq:def:Lin}
		\begin{split}
			\Lin(D)=
			\pR*{ f\in C(D,\R) \colon \!\!\PR*{\!\!
					\begin{array}{c}
						\exists\,a_0,a_1,\ldots,a_d\in\R \\
						\forall\,x=(x_1,x_2,\ldots,x_d)\in D\colon \\
						f(x)= a_0+\ssum_{j=1}^d a_j x_j
					\end{array} \!\! }
			}
			.
		\end{split}
	\end{equation}
We say that $g\in\Lin(D)$ is constant if and only if it holds for all $x,y\in D$ that $\tgt(x)=\tgt(y)$.
\end{definition}
\cfclear

\begin{athm}{lemma}{cor:linear:realisations}
	Let $\f\in\ANNSET$ satisfy $\lengthANN(\f)=\outDimANN(\mathscr{f})=1$ \cfload.
	Then it holds that $\realisation(\f)\in\Lin(\R^{\inDimANN(\f)})$ \cfout.
\end{athm}

\begin{proof}[Proof of \cref{cor:linear:realisations}]
	\Nobs that \eqref{ANNrealization:ass2} and \eqref{eq:def:Lin} \prove\ that $\realisation(\f)\in\Lin\pr{\R^{\inDimANN(\f)}}$ \cfload.
	The proof of \cref{cor:linear:realisations} is thus complete.
\end{proof}

\cfclear

\begin{definition}[Intersections of half-spaces]
	\label{def:hyper2}
	\cfconsiderloaded{def:hyper2}
	Let $d,k\in\N$, $h=(h_1,\ldots,h_k)\in(\Lin(\R^d))^k$, $i=(i_1,\ldots,i_k)\in\{0,1\}^k$ \cfload. Then we denote by $\Hyp(h,i)\subseteq\R^d$ the set given by
	\begin{equation}
		\begin{split}
			\label{Hyperpl:eq1}
			\Hyp(h,i)=\cap_{j=1}^{k}\{x\in\R^d\colon (-1)^{i_j} h_j(x)\leq 0\}.
		\end{split}
	\end{equation}
\end{definition}

\cfclear
\begin{athm}{cor}{cor:index:relu:connection}
	Let $d,k\in\N$, $h=(h_1,\ldots,h_k)\in(\Lin(\R^d))^k$, $i=(i_1,\ldots,i_k)\in\{0,1\}^k$ and let $x\in\Hyp(h,i)$ \cfload. Then 
	\begin{equation}
		\label{eq:cor:index:relu:connection}
		\RELU\pr{\pr{h_1(x),h_2(x),\ldots,h_k(x)}}=
		\pr{i_1 h_1(x),i_2 h_2(x),\ldots,i_k h_k(x)}
	\end{equation}
	\cfout.
\end{athm}

\begin{proof}[Proof of \cref{cor:index:relu:connection}]
	\Nobs that 
	\eqref{Hyperpl:eq1} \proves\ that for all $j\in\{1,2,\ldots,k\}$ it holds that
	\begin{equation}
		\RELU(h_j(x))=i_j(h_j(x))
	\end{equation}
	\cfload. 
	This \proves[ei]\ \eqref{eq:cor:index:relu:connection}.
	The proof of \cref{cor:index:relu:connection} is thus complete.
\end{proof}

\cfclear
\begin{athm}{cor}{cor:hyperplanes:cutting:lines}
	Let $d,k\in\N$, $h=(h_1,\ldots,h_k)\in(\Lin(\R^d))^k$, $v\in\R^d\backslash\{0\}$, $w\in\R^d$, $\strch=\cup_{\lambda\in\R}\{w+\lambda v\}$ \cfload. Then there exist $i_0,i_1,\ldots,i_k\in\{0,1\}^{k}$ such that
	\begin{equation}
		\label{eq:hyperplanes:cutting:lines:statement}
		\strch\subseteq \prb{\cup_{j=0}^k \Hyp(h,i_j)}
	\end{equation}
	\cfout.
\end{athm}

\begin{proof}[Proof of \cref{cor:hyperplanes:cutting:lines}]
	Throughout this proof assume without loss of generality that there exist $\lambda_1,\lambda_2,\ldots,\lambda_k\in\R$ which satisfy for all $j\in\{1,2,\ldots,k\}$ that 
	$\{w+\lambda_j v\}=\{x\in\R^d\colon h_j(x)=0\}\cap\strch$ (\nobs\ that for each $j\in\{1,2,\dots,k\}$ with $\{x\in\R^d\colon h_j(x)=0\}\cap\strch\in\{\emptyset,\strch\}$ there exists $\mathfrak{i}\in\{0,1\}$ with $\strch\subseteq\Hyp(h_{j},\mathfrak{i})$)
	and
	$\lambda_{j}\leq\lambda_{\min\{j+1,k\}}$,
	let $\lambda_1,\lambda_2,\ldots,\lambda_{k+1}\in(-\infty,\infty]$ satisfy for all $j\in\{1,2,\ldots,k\}$ that
	\begin{equation}
		\label{eq:setup:intersection:points}
		\lambda_{j}\leq\lambda_{\min\{j+1,k\}}<\lambda_{k+1}=\infty
		\qandq
		\{w+\lambda_j v\}=\{x\in\R^d\colon h_j(x)=0\}\cap\strch
		\ifnocf,
	\end{equation}
	and let $\substrch_0,\substrch_1,\ldots,\substrch_k\subseteq\strch$ satisfy for all $j\in\{1,2,\ldots,k\}$ that
	\begin{equation}
		\label{eq:setup:intersection:sets}
		\substrch_0=\{w+\mu v\in\R^d\colon \mu \in(-\infty,\lambda_{1}) \}
		\qandq
		\substrch_j=\{w+\mu v\in\R^d\colon \mu \in[\lambda_j,\lambda_{j+1}) \}.
	\end{equation}
	\Nobs that \eqref{eq:setup:intersection:points} and the fact that for all $j\in\{1,2,\ldots,k\}$ it holds that $\{x\in\R^d\colon h_j(x)=0\}=\Hyp(h_j,0)\cap\Hyp(h_j,1)$ \prove\ that for all $j\in\{0,1,\ldots,k\}$, $\mathfrak{j}\in\{1,2,\ldots,k\}$ there exists $\mathfrak{i}\in\{0,1\}$ such that 
	\begin{equation}
		\substrch_j\subseteq \Hyp(h_{\mathfrak{j}},\mathfrak{i})
		\ifnocf.
	\end{equation}
	\cfload[.]This \proves\ that for all $j\in\{0,1,\ldots,k\}$ there exists $\mathfrak{i}\in\{0,1\}^k$ such that 
	\begin{equation}
		\substrch_j\subseteq \Hyp(h,\mathfrak{i}).
	\end{equation}
	Combining this with the fact that $\strch=\cup_{j=0}^k \substrch_j$ \proves\ \eqref{eq:hyperplanes:cutting:lines:statement}.
	The proof of \cref{cor:hyperplanes:cutting:lines} is thus complete.
	\cfload
\end{proof}

\cfclear
\begin{athm}{lemma}{lemma:Hyp:reduced}
	Let $d,k\in\N$, $h=(h_1,\ldots,h_k)\in(\Lin(\R^d))^k$, $i=(i_1,\ldots,i_k)\in\{0,1\}^k$ satisfy $\Hyp(h,i)\neq\emptyset$ and let $n\in\N$,
	$j_1,j_2,\ldots,j_n\in\{1,2,\dots,k\}$
	satisfy
	\begin{equation}
		j_1<j_2<\ldots<j_n
		\quad\text{and}\quad
		\pR{l\in\{1,2,\ldots,k\}\colon h_{l}\text{ is not constant }}
		\subseteq
		\{j_1,j_2,\ldots,j_n\}
	\end{equation} 
	\cfload. Then 
	\begin{equation}
		\Hyp(h,i)
		=\Hyp((h_{j_1},h_{j_2},\ldots,h_{j_n}),(i_{j_1},i_{j_2},\ldots,i_{j_n})).
	\end{equation}
\end{athm}

\begin{proof}[Proof of \cref{lemma:Hyp:reduced}]
	\Nobs that the fact that $\Hyp(h,i)=\cap_{s=1}^{k} \Hyp(h_s,i_s)$ and the assumption that $\Hyp(h,i)\neq\emptyset$ \prove\ that for all $s\in\{1,2,\ldots,k\}$ it holds that
	\begin{equation}
		 \Hyp(h_s,i_s)\neq\emptyset.
	\end{equation}
This and the fact that for all $f\in\Lin(\R^d)$ with $\fa{x,y}\R^d\colon f(x)=f(y)$ it holds that $\Hyp(f,0),\Hyp(f,1)\in\{\emptyset,\R^d\}$ \prove\ that
\begin{equation}
	\begin{split}
	\Hyp(h,i)=\cap_{s=1}^{k} \Hyp(h_s,i_s)
	&=\cap_{s=1}^{n} \Hyp(h_{j_s},i_{j_s})
	\\&=\Hyp((h_{j_1},h_{j_2},\ldots,h_{j_n}),(i_{j_1},i_{j_2},\ldots,i_{j_n}))
	.
	\end{split}
\end{equation}
\finishproofthus
\cfload
\end{proof}

\cfclear
\begin{definition}[Convex sets]
	\label{def:convex}
	Let $d\in\N$ and let $A\subseteq\R^d$ be a set. Then we denote by $\convex{d}(A)$ the set given by
	\begin{equation}
		\label{eq:def:convex}
		\convex{d}(A)=\{C\subseteq A\colon\,\forall\,x,y\in C,\, \lambda\in[0,1]\colon x+\lambda(y-x)\in C\}.
	\end{equation}
\end{definition}

\begin{athm}{cor}{cor:convex:sets:properties}
	Let $d,k\in\N$ and let $A\subseteq \R^d$ be a set.
	Then
	\begin{enumerate}[(i)]
		\item{\label{it:convex:basis:1}it holds for all $C_1,C_2,\ldots,C_k\in\convex{d}(A)$ that $\pr{\cap_{i=1}^k C_i}\in\convex{d}(A)$,}
		\item{\label{it:convex:basis:2}it holds for all $B\in\convex{d}(A)$, $C\in\convex{d}(\R^d)$ that $(B\cap C)\in\convex{d}(A)$, and}
		\item{\label{it:convex:basis:3}it holds for all $h\in\Lin(\R^d)$, $i\in\{0,1\}$ that $\Hyp(h,i)\in\convex{d}(\R^d)$}
		
	\end{enumerate}
	\cfout.
\end{athm}

\begin{proof}[Proof of \cref{cor:convex:sets:properties}]
	\Nobs that \eqref{eq:def:convex} \proves\ \cref{it:convex:basis:1,it:convex:basis:2}. 
	\Nobs that \eqref{Hyperpl:eq1} and \eqref{eq:def:convex}  \prove[ei]\ \cref{it:convex:basis:3}.
	The proof of \cref{cor:convex:sets:properties} is thus complete.
\end{proof}

\subsection{Lower bounds for approximations of certain classes of oscillating functions}
\label{subsection:lower_bound_oszillating}

\cfclear
\begin{athm}{lemma}{lem:inequ:nonzeroparam}
	Let $\f\in\ANNSET$ and
	let $L,l_0,l_1,\ldots,l_L\in\N$ satisfy 
	$\dims(\f)=(l_0,l_1,\ldots,l_L)$ and
	for every $s\in\{1,2,\ldots,L\}$ let $k_s,\lay_s\in\N_0$ satisfy 
	\begin{equation}
		\label{setup:amount:constant:polys0}
		\textstyle
		k_s=\sum_{i=1}^{s-1} l_i(l_{i-1}+1)
		\qandq
		\lay_s
		=\sum_{n=1}^{l_{s}}\indicator{\R^{l_{s-1}}\setminus\{0\}}\pr{\row{(\sum_{i=1}^{s-1}l_i)+n}{\f}}
	\end{equation}
\cfload. Then 
$
	\ssum_{s=1}^{L-1}\lay_s
	\leq \paramnotzero(\f)
$
\cfout.
\end{athm}

\begin{proof}[Proof of \cref{lem:inequ:nonzeroparam}]
	\Nobs that \eqref{setup:amount:constant:polys0} and the fact that $\param(\f)=\sum_{s=1}^{L}l_s(l_{s-1}+1)$ \prove\ that
	\begin{equation}
		\begin{split}
\paramnotzero(\f)
&=\param(\f)-\sum_{n=1}^{\param(\f)}\indicator{\{0\}}(\form{n}{\f})
		\\&=\sum_{n=1}^{\param(\f)}\indicator{\R\setminus\{0\}}(\form{n}{\f})
		\\&=\sum_{s=1}^{L}\sum_{n=1}^{l_s}\PRbb{\indicator{\R\setminus\{0\}}(\form{k_s+l_sl_{s-1}+n}{\f})+\sum_{u=1}^{l_{s-1}}\indicator{\R\setminus\{0\}}(\form{k_s+(n-1)l_{s-1}+u}{\f})}
		\\&\geq\sum_{s=1}^{L}\sum_{n=1}^{l_s}\PRbb{\indicator{\R\setminus\{0\}}(\form{k_s+l_sl_{s-1}+n}{\f})+\indicator{\R^{l_{s-1}}\setminus\{0\}}\pr{\row{(\sum_{i=1}^{s-1}l_i)+n}{\f}}}
		\\&\geq\sum_{s=1}^{L}\sum_{n=1}^{l_s}\indicator{\R^{l_{s-1}}\setminus\{0\}}\pr{\row{(\sum_{i=1}^{s-1}l_i)+n}{\f}}
		=\sum_{s=1}^{L}\lay_s\geq \sum_{s=1}^{L-1}\lay_s
\end{split}
	\end{equation}
\cfload. \finishproofthus
\end{proof}

\cfclear
\begin{athm}{lemma}{lemma:deep:linear:structure}
Let $d,L,l_0,l_1,\ldots,l_{L}\in\N$, $k_0,k_1,\ldots,k_{L}\in\N_0$ satisfy for all $s\in\N_0\cap[0,L]$ that $L\geq 2$, $l_0=d$, $l_L=1$, $k_s=\sum_{j=1}^s l_j$,
let $\f\in\pr{ \times_{\mathscr{k} = 1}^L\allowbreak(\R^{l_\mathscr{k} \times l_{\mathscr{k}-1}} \times \R^{l_\mathscr{k}})}\subseteq\ANNSET$,
for every $j\in\{1,2,\ldots,L\}$ let $\proj_j\colon\pr{\cup_{s=j}^L \R^{k_s}}\to\R^{k_j}$ satisfy for all $s\in\N\cap[j,L]$, $i=(i_1,\ldots,i_{k_s})\in\R^{k_s}$ that $\proj_j(i)=(i_1,i_2,\ldots,i_{k_j})$,
for every $\mathscr{k}\in\{1,2,\ldots,L\}$ let
\begin{equation}
	\label{eq:matrices:deep:linear}
\begin{split}
(W_{\mathscr{k},i,j})_{(i,j)\in\{1,2,\ldots,l_\mathscr{k}\}\times\{1,2,\ldots,l_{\mathscr{k}-1}\}}
=\weight{\mathscr{k}}{\f}
\qquad\text{and}\qquad
(B_{\mathscr{k},i})_{i\in\{1,2,\ldots,l_\mathscr{k}\}}
=\bias{\mathscr{k}}{\f},
\end{split}
\end{equation}
and let $G^i_s=(g^i_1,\ldots,g^i_{k_{s}})\in\pr{\Lin(\R^d)}^{k_{s}}$, $s\in\{1,2,\ldots,L\}$, $i\in\{0,1\}^{k_{L-1}}$, satisfy
 for all $i\in\{0,1\}^{k_{L-1}}$, $j\in\{1,2,\ldots,l_{1}\}$, $x=(x_1,\ldots,x_d)\in\R^d$ that
\begin{equation}
\label{item1:cuttinghyperplanes}
g^i_{j}(x)=B_{1,j}+\sum_{p=1}^{d}W_{1,j,p}x_p
\end{equation}
and assume for all $i=(i_1,\ldots,i_{k_{L-1}})\in\{0,1\}^{k_{L-1}}$, $s\in\{1,2,\ldots,L-1\}$, $j\in\{1,2,\ldots,l_{s+1}\}$, $x\in\R^d$ that
\begin{equation}
\label{item2:cuttinghyperplanes}
g^i_{k_s+j}(x)=B_{s+1,j}+\sum_{p=1}^{l_s}W_{s+1,j,p}i_{k_{s-1}+p}\prb{g^i_{k_{s-1}+p}(x)}
\end{equation}
\cfload. Then 
\begin{enumerate}[(i)]
\item{
\label{lemma:deep:linear:structure0}
it holds for all $i,j\in\{0,1\}^{k_{L-1}}$ that 
$
G_{1}^i=G_{1}^j
$,
}
\item{
\label{lemma:deep:linear:structure1}
it holds for all $s\in\{1,2,\ldots,L-1\}$, $i,j\in\{0,1\}^{k_{L-1}}$ with $\proj_s(i)=\proj_s(j)$ that 
\begin{equation}
G_{s+1}^i=G_{s+1}^j,
\end{equation}
}
\item{
\label{lemma:deep:linear:structure2}
it holds for all $x\in\R^d$ that there exists $i\in\{0,1\}^{k_{L-1}}$ such that for all $j\in\{0,1\}^{k_{L-1}}$ with $(L-2)\pr{\proj_{\max\{1,L-2\}}(i)-\proj_{\max\{1,L-2\}}(j)}=0$ it holds that 
\begin{equation}
x\in\Hyp\prb{G_{L-1}^j,i}=\Hyp\prb{G_{L-1}^i,i},
\cfadd{def:hyper2}
\end{equation}
}
\item{
\label{lemma:deep:linear:structure3}
it holds for all $i\in\{0,1\}^{k_{L-1}}$ that
\begin{equation}
\realisation(\f)|_{\Hyp(G_{L-1}^i,i)}=g^i_{k_{L}}|_{\Hyp(G_{L-1}^i,i)}\in\Lin\prb{\Hyp\prb{G_{L-1}^i,i}},
\end{equation}
and
}
\item{
	\label{lemma:deep:linear:structure:non_zero}
	it holds for all
	$s\in\{0,1,\ldots,L-1\}$, $j\in\{1,2,\ldots,l_{s+1}\}$, $i\in\{0,1\}^{k_{L-1}}$, $x\in\R^d$ with $\row{k_{s}+j}{\f}=0$ that $g_{k_{s}+j}^i(x)=B_{s+1,j}$.
}
\end{enumerate}
\cfout.
\end{athm}

\begin{proof}[Proof of \cref{lemma:deep:linear:structure}]
Throughout this proof let $x\in\R^d$.
\Nobs that \eqref{item1:cuttinghyperplanes} and the assumption that $l_1=k_1$ \prove\ that for all $i,j\in\{0,1\}^{k_{L-1}}$ it holds that 
\begin{equation}
\label{base:case:linear:structure:1}
G_{1}^i=G_{1}^j.
\end{equation}
This  \proves[ei]\ \cref{lemma:deep:linear:structure0}. Combining \eqref{item2:cuttinghyperplanes} and \eqref{base:case:linear:structure:1} \proves\ that for all $i,j\in\{0,1\}^{k_{L-1}}$ with $\proj_1(i)=\proj_1(j)$ it holds that 
\begin{equation}
\label{base:case:linear:structure:2}
G_{2}^i=G_{2}^j.
\end{equation}
\Moreover \eqref{item2:cuttinghyperplanes} \proves\ that 
for all $s\in\N\cap(0,L-1)$ 
with $\,\forall\,i,j\in\{0,1\}^{k_{L-1}}\colon\allowbreak[\proj_s(i)=\proj_s(j)]\Rightarrow[G_{s+1}^i=G_{s+1}^j]$ 
it holds that for all $i,j\in\{0,1\}^{k_{L-1}}$ 
with $\proj_{s+1}(i)=\proj_{s+1}(j)$ 
it holds that $G_{s+2}^i=G_{s+2}^j$.
Combining this
 and \eqref{base:case:linear:structure:2} 
with induction \proves\ that for all $s\in\{1,2,\ldots,L-1\}$, $i,j\in\{0,1\}^{k_{L-1}}$ with $\proj_s(i)=\proj_s(j)$ it holds that 
\begin{equation}
\label{base:case:linear:structure:equality}
G_{s+1}^i=G_{s+1}^j.
\end{equation}
This  \proves[ei]\ \cref{lemma:deep:linear:structure1}.
\Nobs that \eqref{base:case:linear:structure:1} and the fact that $l_1=k_1$ \prove\ that there exists $i\in\{0,1\}^{l_1}$ such that for all $j\in\{0,1\}^{k_{L-1}}$ it holds that
\begin{equation}
\label{base:case:linear:structure:1.1}
x\in \Hyp\prb{G^j_1,i}
\cfadd{def:hyper2}
\end{equation}
\cfload. This, \eqref{base:case:linear:structure:2}, and the fact that $l_1+l_2=k_1+l_2=k_2$ \prove\ that there exists $i\in\{0,1\}^{k_2}$ such that for all $j\in\{0,1\}^{k_{L-1}}$ with $\proj_{1}(i)=\proj_{1}(j)$ it holds that
\begin{equation}
\label{base:case:linear:structure:1.2}
x\in 
\Hyp\prb{G^j_2,i}.
\end{equation}
\Moreover \eqref{item2:cuttinghyperplanes} and \eqref{base:case:linear:structure:equality} \prove\ that for all $s\in\N\cap(1,L-1)$, $i\in\{0,1\}^{k_s}$ with
$\,\forall\,j\in\{0,1\}^{k_{L-1}}\colon[\proj_{s-1}(i)=\proj_{s-1}(j)]\Rightarrow[x\in \Hyp\prb{G_s^j,i}]$
there exists $\mathfrak{i}\in\{0,1\}^{k_{s+1}}$ with $\proj_{s}(\mathfrak{i})=i$ such that for all $j\in\{0,1\}^{k_{L-1}}$ with $\proj_{s}(\mathfrak{i})=i=\proj_{s}(j)$ it holds that
\begin{equation}
x\in \Hyp\prb{G^j_{s+1},\mathfrak{i}}.
\end{equation}
Combining this, 
\eqref{base:case:linear:structure:equality}, 
\eqref{base:case:linear:structure:1.1}, and 
\eqref{base:case:linear:structure:1.2} with 
\eqref{item1:cuttinghyperplanes} and induction \proves\ that there exists $i\in\{0,1\}^{k_{L-1}}$ such that for all $j\in\{0,1\}^{k_{L-1}}$ with $(L-2)\pr{\proj_{\max\{1,L-2\}}(i)-\proj_{\max\{1,L-2\}}(j)}=0$ it holds that
\begin{equation}
x\in \Hyp(G_{L-1}^j,i)=\Hyp(G_{L-1}^i,i)
\ifnocf.
\end{equation}
\cfload[.]This  \proves[ei]\ \cref{lemma:deep:linear:structure2}.
Next, let $i=(i_1,\ldots,i_{k_{L-1}})\in\{0,1\}^{k_{L-1}}$
and $y_s \in \R^{l_s}$, $s\in\{0,1,\ldots,L\}$, satisfy for all $s \in \{1,2,
    \ldots,L\}$ that 
		\begin{equation}
		\label{eq:seq:setup:for:ANN:calc}
		y_0\in \Hyp(G_{L-1}^i,i)
		\qandq
		y_s =\RELU(\weight{s}{\f} y_{s-1} + \bias{s}{\f})
	  %
		\end{equation}
		\cfload[.]
\Nobs that 
\eqref{eq:matrices:deep:linear},
\eqref{item1:cuttinghyperplanes}, 
\eqref{eq:seq:setup:for:ANN:calc},
and
\cref{cor:index:relu:connection} \prove\ that
\begin{equation}
\label{base:case:linear:at:point}
\begin{split}
y_{1} =\RELU(\weight{1}{\f} y_{0} + \bias{1}{\f})&=\RELU(g^i_1(y_{0}),g^i_2(y_{0}),\ldots,g^i_{l_1}(y_{0}))
\\&=(i_1(g^i_1(y_{0})),i_2(g^i_2(y_{0})),\ldots,i_{l_1}(g^i_{l_1}(y_{0}))).
\end{split}
\end{equation}
\Moreover 
\eqref{item2:cuttinghyperplanes},
\eqref{eq:seq:setup:for:ANN:calc},
and
\cref{cor:index:relu:connection} \prove\ that for all $s\in\{1,2,\ldots,L-2\}$ with
$y_{s} =(i_{k_{s-1}+1}(g^i_{k_{s-1}+1}(y_{0})),\allowbreak i_{k_{s-1}+2}(g^i_{k_{s-1}+2}(y_{0})),\ldots,i_{k_{s}}(g^i_{k_s}(y_{0})))$ it holds that
\begin{equation}
\begin{split}
y_{s+1} =\RELU(\weight{s+1}{\f} y_{s} + \bias{s+1}{\f})
&=\RELU(g^i_{k_s+1}(y_{0}),g^i_{k_s+2}(y_{0}),\ldots,g^i_{{k_{s+1}}}(y_{0}))
\\&=(i_{k_s+1}(g^i_{k_s+1}(y_{0})),i_{k_s+2}(g^i_{k_s+2}(y_{0})),\ldots,i_{k_{s+1}}(g^i_{{k_{s+1}}}(y_{0}))).
\end{split}
\end{equation}
Combining this and \eqref{base:case:linear:at:point} with induction \proves\ that 
\begin{equation}
\begin{split}
y_{L-1}=(i_{k_{L-2}+1}(g^i_{k_{L-2}+1}(y_{0})),i_{k_{L-2}+2}(g^i_{k_{L-2}+2}(y_{0})),\ldots,i_{k_{L-1}}(g^i_{{k_{L-1}}}(y_{0}))).
\end{split}
\end{equation}
This, 
\eqref{item2:cuttinghyperplanes}, 
\eqref{eq:seq:setup:for:ANN:calc},
and the fact that $( \realisation(\mathscr{f}) ) (y_0) = \weight{L}{\f} y_{L-1} + \bias{L}{\f}$ \prove\ that
\begin{equation}
\begin{split}
g^i_{k_{L}}(y_{0})
=g^i_{k_{L-1}+1}(y_{0})
&=B_{L,1}+\sum_{p=1}^{l_{L-1}}W_{L,1,p}i_{k_{L-2}+p}\prb{g^i_{k_{L-2}+p}(y_0)}
\\&=\weight{L}{\f} y_{L-1} + \bias{L}{\f}=(\realisation(\mathscr{f}))(y_0)
\end{split}
\end{equation}
\cfload.
This  \proves[ei]\ \cref{lemma:deep:linear:structure3}. 
\Nobs that \eqref{item1:cuttinghyperplanes},  \eqref{item2:cuttinghyperplanes}, and the fact that
$\forallDist s\in\{0,1,\ldots,L-1\}$, $j\in\{1,2,\ldots,l_{s+1}\}\colon\row{k_{s}+j}{\f}=(W_{s+1,j,1},W_{s+1,j,2},\ldots,W_{s+1,j,l_s})$ \prove\ that for all
$s\in\{0,1,\ldots,L-1\}$, $j\in\{1,2,\ldots,l_{s+1}\}$, $i\in\{0,1\}^{k_{L-1}}$ with $\row{k_{s}+j}{\f}=0$ it holds that 
\begin{equation}
	g_{k_s+j}^i=B_{s+1,j}
\end{equation}
\cfload.
This  \proves[ei]\ \cref{lemma:deep:linear:structure:non_zero}. 
The proof of \cref{lemma:deep:linear:structure} is thus complete.\cfload
\end{proof}

\cfclear

\begin{athm}{prop}{lemma:realisation:shallow}
Let $a\in[-\infty,\infty)$, $b\in[a,\infty]$, $d\in\N$, $\f\in\ANNSET$, $\mu\in\R^d$, $\nu\in\R^d\setminus\{0\}$ satisfy $\realisation(\f)\in C(\R^d,\R)$ and
let $\strch=[a,b]^d\cap\pr{\cup_{\lambda \in\R}\{\mu+\lambda\nu\}}\neq\emptyset$ \cfload. Then 
\begin{multline}
\label{hyperplanes:line:realisation}
	\min \pr*{ \pR*{ k \in \N \colon \PR*{\!\!
	\begin{array}{c}
	    \exists\, \substrch_1,\substrch_2,\ldots,\substrch_k\in\convex{d}(\strch) \colon
	    \PRb{(\strch=\cup_{j=1}^k \substrch_j)\land{}\\
        (\fa{j}\{1,2,\ldots,k\}\colon \realisation(\f)|_{\substrch_j}\in\Lin(\substrch_j))}\\
    \end{array} \!\! }
} \cup \{\infty\} }
\\\leq \PR*{1+\frac{\paramnotzero(\f)}{\max\{1,\hidlengthANN(\f)\}}}^{\hidlengthANN(\f)}
\end{multline}
\cfout.
\end{athm}

\begin{proof}[Proof of \cref{lemma:realisation:shallow}]
Throughout this proof 
assume without loss of generality that $\lengthANN(\f)-1=\hidlengthANN(\f)> 0$ \cfadd{cor:linear:realisations}\cfload\
and $\fa{j}\{1,2,\ldots,\lengthANN(\f)\}\colon\weight{j}{\f}\neq0$ \cfadd{lem:zero_weights}\cfload,
let $L,l_0,l_1,\ldots,l_L\in\N$, $k_0,k_1,\ldots,k_L\in\N_0$ satisfy for all $s\in\N_0\cap[0,L]$ that 
$\dims(\f)=(l_0,l_1,\ldots,l_L)$ and
$k_s=\sum_{j=1}^s l_j$,
for every $j\in\{1,2,\ldots,L\}$ let $\lay_j\in\N$ satisfy 
\begin{equation}
	\label{setup:amount:constant:polys}
 	\textstyle
\lay_j
=\sum_{n=1}^{l_{j}}\indicator{\R^{l_{j-1}}\setminus\{0\}}\pr{\row{k_{j-1}+n}{\f}}
=l_j-\sum_{n=1}^{l_{j}}\indicator{\{0\}}\pr{\row{k_{j-1}+n}{\f}}
,
\end{equation}
for every $s\in\N_0\cap[0,L)$ let $\tau_s\colon\{1,2,\dots,\lay_{s+1}\}\to\{t\in\{1,2,\dots,l_{s+1}\}\colon \row{k_{s}+t}{\f}\neq 0\}$ satisfy for all $t_1,t_2\in\{1,2,\dots,\lay_{s+1}\}$ with $t_1<t_2$ that $\tau_s(t_1)<\tau_s(t_2)$, 
let $P,Q\in\N$ satisfy $P=\prod_{n=1}^{L-1}(l_n+1)$ and $Q=\prod_{n=1}^{L-1}(\lay_n+1)$, and
for every $j\in\{1,2,\ldots,L\}$ let $\proj_j\colon\pr{\cup_{s=j}^L \R^{k_s}}\to\R^{k_j}$ satisfy for all $s\in\N\cap [j,L]$, $i=(i_1,\ldots,i_{k_s})\in\R^{k_s}$ that $\proj_j(i)=(i_1,i_2,\ldots,i_{k_j})$,
and let
$G^i_s=(g^i_1,\ldots,g^i_{k_{s}})\in\pr{\Lin(\R^d)}^{k_{s}}$, $s\in\{1,2,\ldots,L\}$, $i\in\{0,1\}^{k_{L-1}}$,
 satisfy that
\begin{enumerate}[(I)]
\item{
\label{item:deep:linear:structure0}
it holds for all $i,j\in\{0,1\}^{k_{L-1}}$ that 
$
G_{1}^i=G_{1}^j
$,
}
\item{
\label{item:deep:linear:structure1}
it holds for all $s\in\{1,2,\ldots,L-1\}$, $i,j\in\{0,1\}^{k_{L-1}}$ with $\proj_s(i)=\proj_s(j)$ that 
\begin{equation}
G_{s+1}^i=G_{s+1}^j,
\end{equation}
}
\item{
\label{item:deep:linear:structure3}
it holds for all $i\in\{0,1\}^{k_{L-1}}$ that
\begin{equation}
\realisation(\f)|_{\Hyp(G_{L-1}^i,i)}=g^i_{k_{L}}|_{\Hyp(G_{L-1}^i,i)}\in\Lin\prb{\Hyp\prb{G_{L-1}^i,i}},
\end{equation}
}
and
\item{
	\label{lemma:deep:linear:structure4}
	it holds for all
	$s\in\{0,1,\ldots,L-1\}$, $j\in\{1,2,\ldots,l_{s+1}\}$, $i\in\{0,1\}^{k_{L-1}}$, $x\in\R^d$ with $\row{k_{s}+j}{\f}=0$ that $g_{k_{s}+j}^i$ is constant.
}
\end{enumerate}
\cfadd{lemma:deep:linear:structure}
\cfload. 
\Nobs that \cref{item:deep:linear:structure0} and \cref{cor:hyperplanes:cutting:lines} (applied with
$d \curvearrowleft d$,
$k \curvearrowleft l_1$,
$h \curvearrowleft G_1^i$,
$v \curvearrowleft \nu$,
$w \curvearrowleft \mu$
for $i\in\{0,1\}^{k_{L-1}}$ in the notation of \cref{cor:hyperplanes:cutting:lines}) \prove\ that there exist $j_0,j_1,\ldots,j_{l_1}\in\{0,1\}^{l_1}$ such that for all $i\in\{0,1\}^{k_{L-1}}$ it holds that
\begin{equation}
\label{base:case:covering:Phi}
\strch\subseteq\cup_{m=0}^{l_{1}}\Hyp\prb{G_1^i,j_m}
=\cup_{m=0}^{l_{1}}\Hyp\prb{(g_{1}^i,g_2^i,\ldots,g_{l_1}^i),j_{m}}.
\end{equation}
This, \cref{item:deep:linear:structure0}, \cref{lemma:deep:linear:structure4}, and
\cref{lemma:Hyp:reduced} (applied 
with
$d \curvearrowleft d$,
$k \curvearrowleft l_{1}$,
$h \curvearrowleft (g_{1}^i,g_{2}^i,\ldots,g_{l_1}^i)$,
$i \curvearrowleft j_m$,
$n \curvearrowleft \lay_1$
for 
$i\in\{0,1\}^{k_{L-1}}$,
$m\in\N_0\cap[0,l_{1}]$, 
$j_m\in\{0,1\}^{l_1}$ with $\Hyp\prb{(g_{1}^i,g_{2}^i,\ldots,g_{l_1}^i),j_m}\neq\emptyset$
in the notation of \cref{lemma:Hyp:reduced})
\prove\ that there exist
$j_m=(j_{m,1},j_{m,2},\ldots,j_{m,{l_1}})\in\{0,1\}^{l_1}$, $m\in\N_0\cap[0,l_1]$, such that for all $i\in\{0,1\}^{k_{L-1}}$, $x\in\strch$ there exists $m\in\N_0\cap[0,l_1]$ such that 
\begin{equation}
	\label{base:case:covering:Phi:red}
	\begin{split}
x&\in \Hyp\prb{(g_{1}^i,g_2^i,\ldots,g_{l_1}^i),j_{m}}
\\&=
\Hyp\prb{(g_{\tau_0(1)}^i,g_{\tau_0(2)}^i,\ldots,g_{\tau_0(\lay_1)}^i),(j_{m,\tau_0(1)},j_{m,\tau_0(2)},\ldots,j_{m,\tau_0(\lay_1)})}.
\end{split}
\end{equation}
\Nobs that \cref{item:deep:linear:structure0} and \cref{cor:hyperplanes:cutting:lines} (applied with
$d \curvearrowleft d$,
$k \curvearrowleft \lay_1$,
$h \curvearrowleft (g_{\tau_0(1)}^i,g_{\tau_0(2)}^i,\ldots,g_{\tau_0(\lay_1)}^i)$,
$v \curvearrowleft \nu$,
$w \curvearrowleft \mu$
for $i\in\{0,1\}^{k_{L-1}}$ in the notation of \cref{cor:hyperplanes:cutting:lines}) hence \prove\ that there exist 
$j_m=(j_{m,1},j_{m,2},\ldots,j_{m,{\lay_1}})\in\{0,1\}^{\lay_1}$, $m\in\N_0\cap[0,\lay_1]$, such that for all $i\in\{0,1\}^{k_{L-1}}$,
$x\in\strch$
there exists 
$m\in\N_0\cap[0,\lay_1]$
such that
\begin{equation}
	\label{base:case:covering:Phi:red:PLUS}
	x\in
	\Hyp\prb{(g_{\tau_0(1)}^i,g_{\tau_0(2)}^i,\ldots,g_{\tau_0(\lay_1)}^i),j_m}.
\end{equation}
This and \eqref{base:case:covering:Phi:red} \prove\ that there exist 
$j_m=(j_{m,1},j_{m,2},\ldots,j_{m,{l_1}})\in\{0,1\}^{l_1}$, $m\in\N_0\cap[0,\lay_1]$, such that for all $i\in\{0,1\}^{k_{L-1}}$ it holds that
\begin{equation}
	\label{base:case:covering:Phi:red:PLUS:union}
	\begin{split}
	\strch &\subseteq
	\cup_{m=0}^{\lay_1}\Hyp\prb{(g_{1}^i,g_2^i,\ldots,g_{l_1}^i),j_{m}}
	\\&=
	\cup_{m=0}^{\lay_1}
	\Hyp\prb{(g_{\tau_0(1)}^i,g_{\tau_0(2)}^i,\ldots,g_{\tau_0(\lay_1)}^i),(j_{m,\tau_0(1)},j_{m,\tau_0(2)},\ldots,j_{m,\tau_0(\lay_1)})}.
	\end{split}
\end{equation}
\Moreover \cref{item:deep:linear:structure1} and \cref{cor:hyperplanes:cutting:lines} (applied with
$d \curvearrowleft d$,
$k \curvearrowleft l_{s+1}$,
$h \curvearrowleft (g_{k_s+1}^i,g_{k_s+2}^i,\ldots,g_{k_{s+1}}^i)$,
$v \curvearrowleft \nu$,
$w \curvearrowleft \mu$
for $s\in\N\cap(0,L-1)$, $i\in\{0,1\}^{k_{L-1}}$ in the notation of \cref{cor:hyperplanes:cutting:lines}) \prove\ that for all $s\in\N\cap(0,L-1)$, $\mathfrak{i}\in\{0,1\}^{k_s}$ there exist $j_0,j_1,\ldots,j_{l_{s+1}}\in\{0,1\}^{l_{s+1}}$ such that for all $i\in\{0,1\}^{k_{L-1}}$ with $\proj_s(i)=\mathfrak{i}$ it holds that
\begin{equation}
\label{IS:case:covering:Phi}
\strch\subseteq
\cup_{m=0}^{l_{s+1}}\Hyp\prb{(g_{k_s+1}^i,g_{k_s+2}^i,\ldots,g_{k_{s+1}}^i),j_{m}}
.
\end{equation}
This, \cref{item:deep:linear:structure1},
\cref{lemma:deep:linear:structure4}, and
\cref{lemma:Hyp:reduced} (applied with
$d \curvearrowleft d$,
$k \curvearrowleft l_{s+1}$,
$h \curvearrowleft (g_{k_s+1}^i,g_{k_s+2}^i,\allowbreak\ldots,g_{k_{s+1}}^i)$,
$i \curvearrowleft j_m$,
$n \curvearrowleft \lay_{s+1}$
for 
$i\in\{0,1\}^{k_{L-1}}$,
$s\in\N\cap(0,L-1)$,
$m\in\N_0\cap[0,l_{s+1}]$, 
$j_m\in\{0,1\}^{l_{s+1}}$ with $\Hyp\prb{(g_{k_s+1}^i,g_{k_s+2}^i,\ldots,g_{k_{s+1}}^i),j_m}\neq\emptyset$
in the notation of \cref{lemma:Hyp:reduced})
\prove\ that for all
$s\in\N\cap(0,L-1)$, $\mathfrak{i}\in\{0,1\}^{k_s}$ there exist
$j_m=(j_{m,1},j_{m,2},\ldots,j_{m,{l_{s+1}}})\in\{0,1\}^{l_{s+1}}$, $m\in\N_0\cap[0,l_{s+1}]$
 such that for all $i\in\{0,1\}^{k_{L-1}}$, $x\in\strch$
with $\proj_s(i)=\mathfrak{i}$
there exists $m\in\N_0\cap[0,l_{s+1}]$ such that  
\begin{equation}
	\begin{split}
	\label{base:case:covering:Phi:red:step}
x&\in \Hyp\prb{(g_{k_s+1}^i,g_{k_s+2}^i,\ldots,g_{k_{s+1}}^i),j_{m}}
\\&=
\Hyp\prb{(g_{k_s+\tau_{s}(1)}^i,g_{k_s+\tau_{s}(2)}^i,\ldots,g_{k_s+\tau_{s}(\lay_{s+1})}^i),(j_{m,\tau_{s}(1)},j_{m,\tau_{s}(2)},\ldots,j_{m,\tau_{s}(\lay_{s+1})})}
		.
\end{split}
\end{equation}
\Nobs that \cref{item:deep:linear:structure1} and \cref{cor:hyperplanes:cutting:lines} (applied with
$d \curvearrowleft d$,
$k \curvearrowleft \lay_{s+1}$,
$v \curvearrowleft \nu$,
$w \curvearrowleft \mu$,
$h \curvearrowleft (g_{k_s+\tau_{s}(1)}^i,g_{k_s+\tau_{s}(2)}^i,\ldots,g_{k_s+\tau_{s}(\lay_{s+1})}^i)$
for $i\in\{0,1\}^{k_{L-1}}$, $s\in\N\cap(0,L-1)$ in the notation of \cref{cor:hyperplanes:cutting:lines}) hence \prove\ that 
for all $s\in\N\cap(0,L-1)$, $\mathfrak{i}\in\{0,1\}^{k_s}$
there exist 
$j_m=(j_{m,1},j_{m,2},\ldots,j_{m,{\lay_{s+1}}})\in\{0,1\}^{\lay_{s+1}}$, $m\in\N_0\cap[0,\lay_{s+1}]$, 
such that
 for all $i\in\{0,1\}^{k_{L-1}}$, $x\in\strch$
with $\proj_s(i)=\mathfrak{i}$
there exists $m\in\N_0\cap[0,\lay_{s+1}]$ such that
\begin{equation}
	\label{base:case:covering:Phi:red:PLUS:IS}
	x\in
	\Hyp\prb{(g_{k_s+\tau_{s}(1)}^i,g_{k_s+\tau_{s}(2)}^i,\ldots,g_{k_s+\tau_{s}(\lay_{s+1})}^i),j_m}.
\end{equation}
\Nobs that \eqref{base:case:covering:Phi:red:step} hence \proves\ that for all
$s\in\N\cap(0,L-1)$, $\mathfrak{i}\in\{0,1\}^{k_s}$
there exist 
$j_m=(j_{m,1},j_{m,2},\ldots,j_{m,{l_{s+1}}})\in\{0,1\}^{l_{s+1}}$, $m\in\N_0\cap[0,\lay_{s+1}]$, such that for all $i\in\{0,1\}^{k_{L-1}}$ with $\proj_s(i)=\mathfrak{i}$ it holds that
\begin{equation}
	\label{base:case:covering:Phi:red:PLUS:IS:union}
	\begin{split}
		\strch &\subseteq
		\cup_{m=0}^{\lay_{s+1}}\Hyp\prb{(g_{k_s+1}^i,g_{k_s+2}^i,\ldots,g_{k_{s+1}}^i),j_{m}}
		\\&=
		\cup_{m=0}^{\lay_{s+1}}
		\Hyp\prb{(g_{k_s+\tau_{s}(1)}^i,g_{k_s+\tau_{s}(2)}^i,\ldots\allowbreak,g_{k_s+\tau_{s}(\lay_{s+1})}^i),(j_{m,\tau_{s}(1)},j_{m,\tau_{s}(2)},\ldots,j_{m,\tau_{s}(\lay_{s+1})})}.
	\end{split}
\end{equation}
Combining this,
\eqref{base:case:covering:Phi:red:PLUS:union}, and induction with the assumption that $P=\prod_{n=1}^{L-1}(l_n+1)$ and $Q=\prod_{n=1}^{L-1}(\lay_n+1)$ \proves\ that there exist
 $i_p=(i_{p,1},\ldots,i_{p,k_{L-1}})\in\{0,1\}^{k_{L-1}}$, $p\in\{1,2,\ldots,P\}$, and 
$\mathfrak{i}_q=(\mathfrak{i}_{q,1},\ldots,\mathfrak{i}_{q,k_{L-1}})\in\{i_1,i_2,\ldots,i_P\}$, $q\in\{1,2,\ldots,Q\}$, such that for all $x\in\strch$ there exist $p\in\{1,2,\ldots,P\}$ and $q\in\{1,2,\ldots,Q\}$ such that
\begin{align}
		\label{eq:reduced:index}
		x&\in\cap_{s=0}^{L-2}\Hyp\prb{(g_{k_{s}+1}^{i_p},g_{k_{s}+2}^{i_p},\ldots,g_{k_{s+1}}^{i_{p}}),(i_{p,k_{s}+1},i_{p,k_{s}+2},\ldots,i_{p,k_{s+1}})}
		\\&\nonumber=\cap_{s=0}^{L-2}\Hyp\prb{(g_{k_s+\tau_s(1)}^{\mathfrak{i}_q},g_{k_s+\tau_s(2)}^{\mathfrak{i}_q},\ldots,g_{k_s+\tau_s(\lay_{s+1})}^{\mathfrak{i}_{q}}),(\mathfrak{i}_{q,k_s+\tau_s(1)},\mathfrak{i}_{q,k_s+\tau_s(2)},\ldots,\mathfrak{i}_{q,k_s+\tau_s(\lay_{s+1})})}
.
	\end{align}
\Moreover for all $i=(i_1,\ldots,i_{k_{L-1}})\in\{0,1\}^{k_{L-1}}$ it holds that
\begin{equation}
	\begin{split}
\cap_{s=0}^{L-2}\Hyp\prb{(g_{k_{s}+1}^{i},g_{k_{s}+2}^{i},\ldots,g_{k_{s+1}}^{i}),(i_{k_{s}+1},i_{k_{s}+2},\ldots,i_{k_{s+1}})}
=\Hyp\prb{G_{L-1}^{i},i}
		.
	\end{split}
\end{equation}
Combining this and
\eqref{eq:reduced:index} with
\cref{item:deep:linear:structure3} \proves\ that there exist $i_1,i_2,\ldots,i_Q\in\{0,1\}^{k_{L-1}}$ such that for all $t\in\{1,2,\ldots,Q\}$ it holds that
\begin{equation}
	\strch=\cup_{q=1}^Q \PRb{\strch\cap\Hyp\prb{G_{L-1}^{i_q},i_q}}
	\qandq 
	\realisation(\f)|_{\strch\cap\Hyp\pr{G_{L-1}^{i_t},i_t}}\in\Lin(\strch\cap\Hyp\prb{G_{L-1}^{i_t},i_t}).
\end{equation}
This,
\cref{cor:convex:sets:properties}, and
the fact that $\strch\in\convex{d}(\strch)$ \prove\ that
\begin{equation}
\label{result:1}
\begin{split}
	\min \pr*{ \pR*{ k \in \N \colon \PR*{\!\!
	\begin{array}{c}
	    \exists\, \substrch_1,\substrch_2,\ldots,\substrch_k\in\convex{d}(\strch) \colon
	    \PRb{(\strch=\cup_{b=1}^k \substrch_b)\land{}\\
        (\fa{b}\{1,2,\ldots,k\}\colon \realisation(\f)|_{\substrch_b}\in\Lin(\substrch_b))}\\
    \end{array} \!\! }
} \cup \{\infty\} }
\leq Q.
\end{split}
\end{equation}
\Moreover \cref{lem:inequ:nonzeroparam} and the inequality of arithmetic and geometric means \prove\ that
\begin{equation}
	\begin{split}
		Q=\sprod_{k=1}^{L-1}(\lay_k+1)
		\leq \PR*{\frac{\ssum_{k=1}^{L-1}(\lay_k+1)}{L-1}}^{L-1} 
		&\leq \PR*{\frac{L-1+\paramnotzero(\f)}{L-1}}^{L-1}
		= \PR*{1+\frac{\paramnotzero(\f)}{L-1}}^{L-1}.
	\end{split}
\end{equation}
\cfload[.]This and \eqref{result:1}  \prove[ei]\ \eqref{hyperplanes:line:realisation}.
The proof of \cref{lemma:realisation:shallow} is thus complete.\cfload
\end{proof}

\cfclear
\begin{athm}{cor}{prop:product:unrestricted}
	Let $a\in\R$, $b\in(a,\infty)$, $\dom\in(1,\infty)$,
	$d\in\N\cap[2,\infty)$, $\eps\in(0,\infty)$ satisfy
	$\dom=\max\{\vass{a},\vass{b}\}$, let $f\colon[a,b]^d\to\R$
	satisfy for all $x=(x_1,\ldots,x_d)\in[a,b]^d$ that
	$f(x)=\sprod_{i=1}^d x_i$, and for every $L\in\N$ let
	$p_L\in\N\cup\{\infty\}$ satisfy
	\begin{equation}
		\llabel{eq:minimum}
		p_L=\min\pr*{\pR*{p\in\N\colon\PR*{\!\!
					\begin{array}{c}
						\exists\,\f\in\ANNSET\colon
						(\paramnotzero(\f)\leq p)\land(\lengthANN(\f)\leq L)\\
						\land\,(\realisation(\f)\in C(\R^d,\R))\\
						\land\,(\sup_{x\in[a,b]^d}\vass{(\realisation(\f))(x)-f(x)}
						\leq\eps)
					\end{array}\!\!}} \cup \{\infty\}}
	\end{equation}
	\cfload. Then
	\begin{enumerate}[(i)]
		\item\llabel{it:1} it holds for all $L\in\N\cap[2,\infty)$ that
		\begin{equation}
			p_L\geq
			\max\pR*{1,(L-1)\pr*{
					\pr*{\frac{(b-a)\dom^{(d-2)/2}}{\sqrt{8\eps}}}^{\!1/(L-1)}-1}}
		\end{equation}
		and
		\item\llabel{it:2} it holds for all $n\in\N$ with $p_1=n$ that $8\eps\geq(b-a)^2\dom^{d-2}$.
	\end{enumerate} 
\end{athm}

\begin{proof}[Proof of \cref{prop:product:unrestricted}]
	Throughout this proof assume without loss of generality that $\infty>\inf_{L\in\N}p_L$,
	let $L\in\N\cap[2,\infty)$ satisfy $p_L<\infty$,
	let $\f\in\ANNSET$, $H\in\N_0$ satisfy
	$H=\hidlengthANN(\f)$, $\paramnotzero(\f)\leq p_L$,
	$\lengthANN(\f)\leq L$, $\realisation(\f)\in C(\R^d,\R)$, and
	\begin{equation}
		\label{eq:product:unrestricted:approx}
		\textstyle
		\sup_{x\in[a,b]^d}\vass{(\realisation(\f))(x)-f(x)}\leq\eps,
	\end{equation}
	let $r\in\{a,b\}$ satisfy $\vass{r}=\dom$,
	let $\xi=\pr{\xi_1,\dots,\xi_d}\colon\R\to\R^d$,
	$q\colon\R\to\R$, and
	$g\colon\R\to\R$ satisfy for all $i\in\{1,2,\dots,d\}$, $t\in\R$ that
	\begin{equation}
		\llabel{eq:setup}
		\xi_i(t)=\indicator{[1,2]}(i)t+\indicator{(2,\infty)}(i)r,
		\quad
		q(t)=r^{d-2}t^2,
		\qandq
		g(t)=(\realisation(\f))(\xi(t)),
	\end{equation}
	let $\strch\subseteq \R^d$ satisfy $\strch=\xi([a,b])$,
	and let $K\in\N\cup\{\infty\}$ satisfy
	\begin{equation}
		\llabel{eq:lin:seg}
		K=\min \pr*{ \pR*{ k \in \N \colon \PR*{\!\!
					\begin{array}{c}
						\exists\, \substrch_1,\substrch_2,\ldots,\substrch_k\in\convex{d}(\strch) \colon
						\PRb{(\strch=\cup_{j=1}^k \substrch_j)\land{}\\
							(\fa{j}\{1,2,\ldots,k\}\colon \realisation(\f)|_{\substrch_j}\in\Lin(\substrch_j))}\\
					\end{array} \!\! }
			} \cup \{\infty\} }.
	\end{equation}
	\Nobs that \eqref{eq:product:unrestricted:approx} and the fact that
	for all $t\in[a,b]$ it holds that $\xi(t)\in[a,b]^d$ and
	$f(\xi(t))=q(t)$ \prove\ that
	\begin{equation}
		\label{eq:product:unrestricted:lineerror}
		\textstyle
		\sup_{t\in[a,b]}\vass{g(t)-q(t)}
		\leq \sup_{x\in[a,b]^d}\vass{(\realisation(\f))(x)-f(x)}
		\leq\eps.
	\end{equation}
	\Moreover 
	\lref{eq:lin:seg},
	the fact that
	$[a,b]^d\cap\xi(\R)=\strch$, and
	\cref{lemma:realisation:shallow} (applied with
	$\mu\curvearrowleft\xi(0)$, $\nu\curvearrowleft\xi(1)-\xi(0)$,
	$\strch\curvearrowleft\strch$ in the notation of
	\cref{lemma:realisation:shallow}) \prove\ that $K\in\N$ and that there exist nonempty
	$\substrch_1,\substrch_2,\ldots,\substrch_K\in\convex{d}(\strch)$
	such that for all $i\in\{1,2,\ldots,K\}$ it holds that
	\begin{equation}
		\llabel{eq:cover:number}
		\strch=\cup_{j=1}^K \substrch_j,
		\quad
		\realisation(\f)|_{\substrch_i}\in\Lin(\substrch_i),
		\qandq
		K\leq\pr*{1+\frac{\paramnotzero(\f)}{\max\{1,H\}}}^{H}.
	\end{equation}
\Nobs that \lref{eq:setup} hence \proves\ that there exist  $u_1,u_2,\dots,u_K,v_1,v_2,\dots,v_K\in[a,b]$, $I_1,I_2,\ldots,\allowbreak I_K\in\convex{1}([a,b])$ which satisfy
for all $i\in\{1,2,\ldots,K\}$ that
\begin{equation}
	\label{eq:product:unrestricted:cover}
	[a,b]=\cup_{j=1}^K I_j,
	\quad
	g|_{I_i}\in\Lin(I_i),
	\quad
	u_i=\inf(I_i),
	\qandq
	v_i=\sup(I_i).
\end{equation}
	\cfload. This and the fact that $g\in C(\R,\R)$ \prove\ that for all
	$i\in\{1,2,\ldots,K\}$ it holds that 
	\begin{equation}
		\llabel{eq:loc:lin}
		g|_{[u_i,v_i]}\in\Lin([u_i,v_i]).
	\end{equation}
	Combining this 
	and
	\eqref{eq:product:unrestricted:cover} with 
	\eqref{eq:product:unrestricted:lineerror} hence \proves\ that
	for all $k\in\{1,2,\dots,K\}$ with $v_k-u_k= \max_{j\in\{1,2,\dots,K\}}(v_j-u_j)$ it holds that
	\begin{equation}
		\llabel{eq:eps:via:lin:seg}
		\begin{split}
			&\dom^{d-2}(b-a)^2K^{-2}
			\\&\textstyle\leq \dom^{d-2}\prb{\sum_{j=1}^K (v_j-u_j)}^2K^{-2}
			\\&\leq \dom^{d-2}(v_k-u_k)^2
			\\&= \dom^{d-2}\vass{(u_k)^2+(v_k)^2-2 u_k v_k}
			\\&= \vass{2r^{d-2}(u_k)^2+2r^{d-2}(v_k)^2-r^{d-2}(u_k+v_k)^2}
			\\&= 2\vass{q(u_k)+q(v_k)-2q((u_k+v_k)2^{-1})}
			\\&= 2\vass{q(u_k)-g(u_k)+q(v_k)-g(v_k)+2g((u_k+v_k)2^{-1})-2q((u_k+v_k)2^{-1})}
			\\&\leq 
			2\vass{q(u_k)-g(u_k)}
			+2\vass{q(v_k)-g(v_k)}
			+4\vass{q((u_k+v_k)2^{-1})-g((u_k+v_k)2^{-1})}
			\leq 8\eps.
		\end{split}
	\end{equation}
	This and \lref{eq:cover:number} \prove[pe]\ \lref{it:2}.
	\Nobs that \lref{eq:eps:via:lin:seg} and \lref{eq:cover:number} with the inequality of arithmetic and geometric means
	\proves\ that for all $h\in\N\cap[H,\infty)$ it holds that
	\begin{equation}
		\label{eq:product:unrestricted:pieces}
		\begin{split}
		\frac{(b-a)\dom^{(d-2)/2}}{\sqrt{8\eps}}\leq K
		&\leq\pr*{1+\frac{\paramnotzero(\f)}{\max\{1,H\}}}^{H}
		\\&= \pr*{1+\frac{\paramnotzero(\f)}{\max\{1,H\}}}^{H} 1^{h-H} 
		\\
		&\leq \pr*{
			\frac{1}{h}\PR*{H\pr*{1+\frac{\paramnotzero(\f)}{\max\{1,H\}}}+(h-H)}
		}^{h} \\
		&= \pr*{1+\frac{H}{\max\{1,H\}}\frac{\paramnotzero(\f)}{h}}^{h}
		\leq \pr*{1+\frac{\paramnotzero(\f)}{h}}^{h}.
		\end{split}
	\end{equation}
	\Nobs that \eqref{eq:product:unrestricted:pieces} \proves\ that
	it holds that
	\begin{equation}
		p_L\geq\max\{1,\paramnotzero(\f)\}
		\geq \max\pR*{1,(L-1)\pr*{
				\pr*{\frac{(b-a)\dom^{(d-2)/2}}{\sqrt{8\eps}}}^{\!1/(L-1)}-1}}.
	\end{equation}
	This \proves[pe]\ \lref{it:1}.
	\finishproofthus
\end{proof}

\cfclear
\begin{definition}
	[Euclidean norm]
	\label{def:Euclidean_norm}
	We denote by $\norm{\cdot} \colon \pr{\cup_{d \in \N} \R^d} \rightarrow \R$ the function which satisfies for all $d \in \N$, $x = (x_1,\dots,x_d) \in \R^d$ that
	$
	\norm{x} = \PRb{ \textstyle{\sum}_{j=1}^{d} \vass{x_j}^2 }^{\nicefrac{1}{2}}
	$.
\end{definition}

\cfclear
\begin{athm}{lemma}{lem:1:dim:line:split}
	Let $S\in[1,\infty)$, $N\in\N$,
	let $C_1,C_2,\dots,C_N\in\convex{1}([0,S])$ satisfy $[0,\floor{S}-1]\subseteq \cup_{j=1}^NC_j $, and
	let $G\in C([0,S],\R)$ satisfy for all $k\in\{1,2,\dots,\floor{S}\}$, $j\in\{1,2,\dots,N\}$ that $(-1)^{k}G(k-1)>0$ and $G|_{C_j}\in\Lin(C_j)$ \cfload. Then $N\geq \floor{S}-1$ \cfout.
\end{athm}

\begin{proof}
	Throughout this proof assume without loss of generality that $S\geq 3$.
	\Nobs that the assumption that for all $k\in\{1,2,\dots,\floor{S}\}$ it holds that $G\in C([0,S],\R)$ and $(-1)^{k}G(k-1)>0$ \proves\ that for every $k\in\{1,2,\dots,\floor{S}-1\}$ there exists $z_k\in(k-1,k)$ which satisfies that
	\begin{equation}
		G(z_k)=0.
	\end{equation}
	\Moreover the fact that for all $k\in\{1,2,\dots,\floor{S}-1\}$ it holds that $(-1)^{k}G(k-1)>0=G(z_k)$ and $k-1<z_k<k$ \proves\ that for all $k,K\in\{1,2,\dots,\floor{S}-1\}$ with $k<K$ it holds that
	\begin{equation}
		G|_{[z_k,z_K]}\notin\Lin([z_k,z_K]).
	\end{equation}
	Combining this with the assumption that for all $j\in\{1,2,\dots,N\}$ it holds that $C_j\in\convex{1}([0,S])$ and
	$G|_{C_j}\in\Lin(C_j)$ \proves\ that
	for all $j\in\{1,2,\dots,N\}$, $k,K\in\{1,2,\dots,\floor{S}-1\}$ with $z_k,z_K\in C_j$ it holds that
	\begin{equation}
		k=K.
	\end{equation}
	Combining this with the fact that $[0,\floor{S}-1]\subseteq\cup_{j=1}^NC_j $ \proves\ that for all 
	$k\in\{1,2,\dots,\floor{S}-1\}$ there exists
	$j\in\{1,2,\dots,N\}$ which satisfies
	\begin{equation}
		z_k\in C_j.
	\end{equation}
	This \proves\ $N\geq \floor{S} -1$. \finishproofthus
\end{proof}

\cfclear
\begin{athm}{lemma}{lemma:sin:hyperplanes}
Let $a\in\R$, $b\in[a,\infty)$, $d\in\N$, $\kappa\in(0,\infty)$, $S\in[1,\infty)$,
let $(v_k)_{k\in\N}\subseteq \R^d$ satisfy for all $k\in\N$ that $v_{k+1}-v_{k}=v_2-v_1$,
let $\strch=[a,b]^d\cap\pr{\cup_{\lambda \in\R}\{\lambda v_1+(1-\lambda)v_2\}}$,
assume $\fa{k}\N\cap(0,S]\colon v_k\in A$,
let $f\colon\R^d\to\R$ satisfy for all $k\in\N\cap[3
,\infty)\cap(0,S]$ that $f(v_{k-1})-f(v_{k-2})=f(v_{k-1})-f(v_{k})\in\{-2\kappa,2\kappa\}$,
and let $g\in C(\R^d,\R)$ satisfy for all $x\in [a,b]^d$ that $\vass{f(x)-\tgt(x)}<\kappa$ \cfload.
Then 
\begin{equation}
\begin{split}
\label{hyperplanes:line}
	\min \pr*{ \pR*{ k \in \N \colon \PR*{\!\!
	\begin{array}{c}
	    \exists\, \substrch_1,\substrch_2,\ldots,\substrch_k\in\convex{d}(\strch) \colon
	    \PRb{(\strch=\cup_{i=1}^k \substrch_i)\land{}\\
        (\fa{i}\{1,2,\ldots,k\}\colon g|_{\substrch_i}\in\Lin(\substrch_i))}\\
    \end{array} \!\! }
} \cup \{\infty\} }
\geq \floor{S}-1
\end{split}
\end{equation}
\cfout.
\end{athm}

\begin{proof}[Proof of \cref{lemma:sin:hyperplanes}]
Throughout this proof assume without loss of generality for all $k\in\N\cap(0,S]$ that $S\geq3$, $f(v_k)=(-1)^k\kappa$,
 and
\begin{equation}
\begin{split}
\label{hyperplanes:line.1}
	\min \pr*{ \pR*{ k \in \N \colon \PR*{\!\!
	\begin{array}{c}
	    \exists\, \substrch_1,\substrch_2,\ldots,\substrch_k\in\convex{d}(\strch) \colon
	    \PRb{(\strch=\cup_{i=1}^k \substrch_i)\land{}\\
        (\fa{i}\{1,2,\ldots,k\}\colon g|_{\substrch_i}\in\Lin(\substrch_i))}\\
    \end{array} \!\! }
} \cup \{\infty\} }
<\infty,
\end{split}
\end{equation}
let $N\in\N$, $\substrch_1,\substrch_2,\ldots,\substrch_N\in\convex{d}(\strch)$,
$C_1,C_2,\dots,C_N\in\convex{1}([0,S])$ satisfy for all $j\in\{1,2,\ldots,N\}$ that
$\strch=\cup_{i=1}^N \substrch_i$,
$g|_{\substrch_j}\in\Lin(\substrch_j)$,
and $C_j=\{t\in[0,S]\colon v_1+t(v_2-v_1)\in B_j\}$,
and let 
$G\colon[0,S]\to\R$ and
$\operatorname{sign}\colon\R\to\{-1,0,1\}$ satisfy for all $t\in[0,S]$, $x\in\R$ that
\begin{equation}
	\label{eq:easy:fcts:for:hlp}
	G(t)=\tgt(v_1+t(v_2-v_1))
	\qquad\text{and}\qquad
	\operatorname{sign}(x)=
	\begin{cases}
		\frac{x}{\vass{x}} &\colon x\neq0\\
		0 &\colon x=0\\
	\end{cases}
\end{equation} \cfload.
\Nobs that \eqref{eq:easy:fcts:for:hlp} and the assumption that for all $k\in\N$ it holds that $v_{k+1}-v_{k}=v_2-v_1$ \prove\ that for all $k\in\N$ it holds that
\begin{equation}
	\label{eq:proj:eq0}
	G(0)=\tgt(v_1)
	\quad\text{and}\quad
	v_1+k(v_2-v_1)=v_1+\sum_{j=1}^k(v_2-v_1)=v_1+\sum_{j=1}^k(v_{j+1}-v_j)=v_{k+1}.
\end{equation}
Combining this with induction \proves\ that for all $k\in\N_0\cap[0,S]$ it holds that
\begin{equation}
	G(k)=\tgt(v_{k+1}).
\end{equation}
This and
the assumption that $\fa{x}[a,b]^d\colon\vass{f(x)-\tgt(x)}<\kappa$ 
and the assumption that
$\fa{k}\N\cap(0,S]\colon f(v_k)=(-1)^k\kappa$ \prove\ that for all $k\in\N\cap(0,S]$ it holds that
\begin{equation}
	\label{eq:proj:eq1}
	\operatorname{sign}\pr{G(k-1)}
	=\operatorname{sign}\pr{\tgt(v_{k})}
	=\operatorname{sign}\pr{f(v_{k})}
	=\operatorname{sign}\pr{(-1)^{k}\kappa}
	=(-1)^{k}.
\end{equation}
\Moreover \eqref{eq:easy:fcts:for:hlp} and the assumption that for all $j\in\{1,2,\ldots,N\}$ it holds that
$g|_{\substrch_j}\in\Lin(\substrch_j)$,
and $C_j=\{t\in[0,S]\colon v_1+t(v_2-v_1)\in B_j\}$ \prove\ that for all $j\in\{1,2,\ldots,N\}$ it holds that
\begin{equation}
	\label{eq:proj:eq2}
	G|_{C_j}\in\Lin(C_j).
\end{equation}
\Moreover 
\eqref{eq:proj:eq0},
the assumption that 
$v_1,v_{\floor{S}}\in\strch$,
the fact that $v_{\floor{S}}=v_1+(\floor{S}-1)(v_2-v_1)$, and 
the fact that $A$ is convex \prove\ that
for all $t\in[0,\floor{S}-1]$ it holds that 
\begin{equation}
	v_1+t(v_2-v_1)\in A.
\end{equation}
Combining this and the assumption that $\strch=\cup_{j=1}^N \substrch_j$ \proves\ that 
\begin{equation}
	\label{eq:proj:eq3}
	\begin{split}
	[0,\floor{S}-1]
	&\subseteq
	\{t\in[0,S]\colon v_1+t(v_2-v_1)\in A\}
	\\&=\cup_{j=1}^N \{t\in[0,S]\colon v_1+t(v_2-v_1)\in B_j\}
	=\cup_{j=1}^N C_j
	.
	\end{split}
\end{equation}
Combining this,
\eqref{eq:proj:eq1},
\eqref{eq:proj:eq2},
and \cref{lem:1:dim:line:split} with
the fact that $G\in C([0,S],\R)$ \proves\ that
\begin{equation}
	N\geq  \floor{S}-1.
\end{equation}
The proof of \cref{lemma:sin:hyperplanes} is thus complete \cfload.
\end{proof}

\cfclear
\begin{athm}{prop}{prop:low:dim:gen:low:len:cases}
 For every $k\in\{1,2\}$ let $\f_k\in\ANNSET$ satisfy $\realisation(\f_k) \in C(\R^k,\R)$,
let $a\in\R$, $b\in[a+7,\infty)$, $\varepsilon \in (0,1)$ and
let $f\colon\R^2\to\R$ satisfy for all $x,y\in[a,b]$ that
\begin{equation}
	\label{eq:f:properties:low:dim:LB}
	 f(x,y)=\sin(xy)
	\qandq
	\vass{f(x,y) - (\realisation\pr{\f_2}) (x,y)} \leq \varepsilon
\end{equation}
\cfload. Then it holds for all $k\in\{1,2\}$ that 
\begin{equation}
	\label{eq:prop:low:dim:gen:low:len:statement}
	\param(\f_k)\geq 2^{\nicefrac{k}{\max\{1,\hidlengthANN(\f_k)\}}}.
\end{equation}
\end{athm}

\begin{proof}[Proof of \cref{prop:low:dim:gen:low:len:cases}]
	Throughout this proof let $x,y_1,y_2,y_3\in[a,b]$ satisfy $\vass{x}\geq 3$, $y_1\leq a+\vass{x}^{-1}\pi$, $\vass{\sin(x y_1)}=1$, $y_2=y_1+\vass{x}^{-1}\pi$, and $y_3=y_2+\vass{x}^{-1}\pi$.
\Nobs that \eqref{eq:def:ANN:operators} \proves\ that
\begin{equation}
\label{eq:prop:low:dim:gen:low:len:cases:d=1}
\param(\f_1)=\sum_{k=1}^{\lengthANN(\f_1)}\singledims_k(\f_1)(\singledims_{k-1}(\f_1)+1)
\geq 
2 \geq 2^{\nicefrac{1}{\max\{1,\hidlengthANN(\f_1)\}}}.
\end{equation}
\Moreover \eqref{eq:f:properties:low:dim:LB} \proves\ that for all $g\in\Lin(\R^2)$ with $\vass{\tgt(x,y_1)-f(x,y_1)}\leq\eps$ and $\vass{\tgt(x,y_2)-f(x,y_2)}\leq\eps$ it holds that
\begin{equation}
	\begin{split}
		\vass{\tgt(x,y_3)-f(x,y_3)}
		&=\vass*{\tgt( x,y_2)+\pr*{\frac{\tgt( x,y_2)-\tgt( x,y_1)}{\norm{(x,y_2)-(x,y_1)}}}\norm{(x,y_3)-(x,y_2)}-f(x,y_1)}
		\\&=\vass{2\tgt( x,y_2)-\tgt( x,y_1)-f(x,y_1)}
		\\&\geq \vass{2f(x,y_2)-f(x,y_1)-f(x,y_1)}-3\eps
		\\&= 2\vass{f(x,y_2)-f(x,y_1)}-3\eps
		\\&= 4-3\eps>\eps
	\end{split}
\end{equation}
\cfload.
Combining this with \cref{cor:linear:realisations} \proves\ that for all $\g\in\ANNSET$ with $\hidlengthANN(\g)=0$ and $\realisation(\g)\in C(\R^2,\R)$ there exists $x\in[a,b]^2$ such that
\begin{equation}
	\label{eq:prop:low:dim:gen:low:len:cases:d=2:L=1}
	\vass{(\realisation(\g))(x)-f(x)}>\eps.
\end{equation}
\Moreover for all $\g\in\ANNSET$ with $\realisation(\g)\in C(\R^2,\R)$ and $\hidlengthANN(\g)\geq 1$ it holds that
\begin{equation}
	\label{eq:prop:low:dim:gen:low:len:cases:d=2:L=2}
	\param(\g)=\sum_{k=1}^{\lengthANN(\g)}\singledims_k(\g)(\singledims_{k-1}(\g)+1)\geq 2\lengthANN(\g)=2(\hidlengthANN(\g)+1)\geq 4\geq 2^{\frac{2}{\max\{1,\hidlengthANN(\g)\}}}.
\end{equation}
Combining this with \eqref{eq:prop:low:dim:gen:low:len:cases:d=1}
and
\eqref{eq:prop:low:dim:gen:low:len:cases:d=2:L=1}  \proves[ei]\ \eqref{eq:prop:low:dim:gen:low:len:statement}.
The proof of \cref{prop:low:dim:gen:low:len:cases} is thus complete.\cfload
\end{proof}

\cfclear
\begin{athm}{prop}{gen:final:approximation:fail}
    Let $a\in\R$, $b\in[a,\infty)$, $d,H\in \N$, $\nu\in\R^d\setminus\{0\}$, $\kappa\in(0,\infty)$, $S\in[2,\infty)$, $\sigma\in\{-2\kappa,2\kappa\}$, $\varepsilon \in [0,\kappa)$, 
let $v_k\in[a,b]^d$, $k\in\N\cap(0,S]$,
let $f\colon\R^d\to\R$ satisfy for all $k\in\N\cap(1,S]$ that $v_k=v_{k-1}+\nu$ and $f(v_{k})-f(v_{k-1})=\sigma(-1)^k$, and let	$\mathscr{f} \in \ANNSET$ satisfy for all $x\in[a,b]^d$ that
\begin{equation}
	\label{eq:fin:ap:net:prop:fail}
	\realisation(\mathscr{f}) \in C(\R^d,\R),
	\qquad
	\vass{f(x) - (\realisation\pr{\mathscr{f}}) (x)} \leq \varepsilon,
	\qandq
	H=\max\{1,\hidlengthANN(\f)\}	
\end{equation}
    \cfload. 
		Then $\paramnotzero(\f)\geq
		\max\pRb{\tfrac{H}{2}\pr{\floor{S}-1}^{\nicefrac{1}{H}},\hidlengthANN(\f)+1}$ \cfout.
\end{athm}

\begin{proof}[Proof of \cref{gen:final:approximation:fail}]
Throughout this proof let $v_k\in\R^d$, $k\in\N\cap(S,\infty)$, and
let $\strch=[a,b]^d\cap\pr{\cup_{\lambda\in\R}\{\lambda v_1+(1-\lambda)v_2\}}$.
\Nobs that 
\cref{lem:zero_weights},
the assumption that $S\geq 2$, and the fact that for all $\g\in\ANNSET$ with $\paramnotzero(\g)<\lengthANN(\g)$ it holds that $\fa{x,y}\R^{\inDimANN(\g)}\colon(\realisation(\g))(x)=(\realisation(\g))(y)$ \prove\ that for all $\g\in\ANNSET$ with $\paramnotzero(\g)<\lengthANN(\g)$ and $\realisation(\g)\in C(\R^d,\R)$ it holds that
\begin{equation}
	\begin{split}
	\max_{i\in\{1,2\}}\{\vass{f(v_i)-(\realisation\pr{\mathscr{g}}) (v_i)}\}
	&\geq
	\tfrac{1}{2}\vass{f(v_1)-(\realisation\pr{\mathscr{g}}) (v_1)}
	+
	\tfrac{1}{2}\vass{f(v_2)-(\realisation\pr{\mathscr{g}}) (v_2)}
	\\&=
	\tfrac{1}{2}\vass{f(v_1)-(\realisation\pr{\mathscr{g}}) (v_1)}
	+
	\tfrac{1}{2}\vass{f(v_1)+\sigma-(\realisation\pr{\mathscr{g}}) (v_1)}
	\\&\geq
	\tfrac{\vass{\sigma}}{2}
	=
	\kappa
	>
	\eps
	.
	\end{split}
\end{equation}
\cfload.
Combining this with the fact that $\paramnotzero(\f)\in\N_0$,
 the assumption that $v_1,v_2\in[a,b]^d$,
 and \eqref{eq:fin:ap:net:prop:fail} \proves\ that
\begin{equation}
		\label{eq:zero:nets:fail:approx}
	\paramnotzero(\f)\geq\lengthANN(\f)=\hidlengthANN(\f)+1\geq H.
\end{equation}
\Moreover \cref{eq:fin:ap:net:prop:fail} and \cref{lemma:sin:hyperplanes} (applied with
		$a \curvearrowleft  a$,
		$b \curvearrowleft  b$,
		$d \curvearrowleft  d$,
		$\kappa \curvearrowleft \kappa$,
		$S \curvearrowleft S$,
		$\strch \curvearrowleft  \strch$,
		$(v_k)_{k\in\N} \curvearrowleft 
		\pr{\N\ni k\mapsto v_1+(k-1)\nu\in\R^d}$,
		$g \curvearrowleft \realisation(\f)$,
		$f \curvearrowleft f$
in the notation of \cref{lemma:sin:hyperplanes}) \prove\ that
\begin{equation}
\begin{split}
\label{eq:gen:sat:hyperplanes:line}
	\min \pr*{ \pR*{ k \in \N \colon \PR*{\!\!
	\begin{array}{c}
	    \exists\, \substrch_1,\substrch_2,\ldots,\substrch_k\in\convex{d}(\strch) \colon
	    \PRb{(\strch=\cup_{i=1}^k \substrch_i)\land{}\\
        (\fa{i}\{1,2,\ldots,k\}\colon \realisation(\f)|_{\substrch_i}\in\Lin(\substrch_i))}\\
    \end{array} \!\! }
} \cup \{\infty\} }
\geq \floor{S}-1
\end{split}
\end{equation}
\cfload.
Combining this and \eqref{eq:zero:nets:fail:approx} with \cref{lemma:realisation:shallow} (applied with
		$a \curvearrowleft  a$,
		$b \curvearrowleft  b$,
		$d \curvearrowleft  d$,
		$\f \curvearrowleft  \f$,
		$\strch \curvearrowleft  \strch$
in the notation of \cref{lemma:realisation:shallow}) \proves\ that
\begin{align}
\label{eq:gen:sat:hyperplanes:line2}
	&\PR*{\frac{2\paramnotzero(\f)}{H}}^{H}
	\geq
	\PR*{\frac{2\paramnotzero(\f)}{H}}^{\hidlengthANN(\f)}
	\geq
	\PR*{1+\frac{\paramnotzero(\f)}{H}}^{\hidlengthANN(\f)}
	\\&\nonumber\geq
	\min \pr*{ \pR*{ k \in \N \colon \PR*{\!\!
	\begin{array}{c}
	    \exists\, \substrch_1,\substrch_2,\ldots,\substrch_k\in\convex{d}(\strch) \colon
	    \PRb{(\strch=\cup_{i=1}^k \substrch_i)\land{}\\
        (\fa{i}\{1,2,\ldots,k\}\colon \realisation(\f)|_{\substrch_i}\in\Lin(\substrch_i))}\\
    \end{array} \!\! }
} \cup \{\infty\} }
\geq \floor{S}-1.
\end{align}
This and \eqref{eq:zero:nets:fail:approx} \prove[e]\ $\paramnotzero(\f)\geq \max\pRb{\tfrac{H}{2}\pr{\floor{S}-1}^{\nicefrac{1}{H}},\hidlengthANN(\f)+1}$. 
The proof of \cref{gen:final:approximation:fail} is thus complete.
\end{proof}

\subsection{Oscillation properties of certain families of functions}
\label{subsection:lower_bound_properties}
\newcommand{\scl}{\beta}

\cfclear
\begin{athm}{cor}{sin:cos:prod:targetfunctions}
Let $\varphi\in\R$, $\kappa,\dom\in(0,\infty)$, $\gamma\in(0,1]$, $\scl\in[1,\infty)$, $a\in\R$, $b\in(a,\infty)$, $g\in C(\R,\R)$ satisfy for all $x\in\R$ that $\tgt(x)=\kappa\sin(x+\varphi)$ and $\dom=\max\{\vass{a},\vass{b}\}>\scl^{-1}$,
let $f_d\in C(\R^d,\R)$, $d\in\N$, satisfy for all $d\in\N$, $x=(x_1,\ldots,x_d)\in[a,b]^d$ that $f_d(x)=g\prb{\gamma\scl^d\prod_{i=1}^d x_i}$, and
let $S\colon\N\to[0,\infty)$ satisfy for all $d\in\N$ that
		\begin{equation}
			\label{eq:iteration_counter_prod}
		S(d)=\floor{(b-a)\gamma(\dom\scl)^{d-1}\pi^{-1}}.
		\end{equation}
Then there exist $(\nu_d,\sigma_d)\in\prb{\R^d\setminus\{0\}}\times\{-2\kappa,2\kappa\}$, $d\in\N$, and $v_{k,d}\in[a,b]^d$, $(k,d)\in\N^2$, such that
\begin{enumerate}[(i)]
\item{
\label{prop:sin:cos:recursion}
it holds for all $d\in\N$, $k\in\N\cap(1,\infty)\cap[0,S(d)]$ that $v_{k,d}=v_{k-1,d}+\nu_d$,
}

\item{
\label{prop:sin:cos:alternation}
it holds for all $d\in\N$, $k\in\N\cap(1,\infty)\cap[0,S(d)]$ that $f_d(v_{k,d})-f_d(v_{k-1,d})=(-1)^k\sigma_d$,
}

\item{
\label{prop:sin:cos:periodicity}
it holds for all $x\in\R$, $k\in\Z$ that $\tgt(x+2k\pi)=\tgt(x)\in[-\kappa,\kappa]$, and
}

\item{
\label{prop:sin:cos:lipschitz}
it holds for all $x,y\in\R$ that $\vass{\tgt(x)-\tgt(y)}\leq \kappa\vass{x-y}$.
}
\end{enumerate}
\end{athm}

\begin{proof}[Proof of \cref{sin:cos:prod:targetfunctions}]
	Throughout this proof 
	let $\mathscr{N}=\{d\in\N\colon S(d)>1\}$
	and
	let $c\in[a,b]$ satisfy $\vass{c}=\dom$.
	\Nobs that for all $d\in\N\setminus\{1\}$, $v=(\alpha,c,c,\ldots,c)\in[a,b]^d$ it holds that
	\begin{equation}
		\label{eq:prod:function:on:line}
		f_1(\alpha)=\kappa\sin(\alpha \gamma\scl +\varphi)
		\qandq
		f_d(v)=g\prb{\gamma\scl^d\alpha\textstyle\prod_{i=2}^d c}=\kappa\sin(\alpha \gamma\scl^d c^{d-1} +\varphi).
	\end{equation}
	\Moreover the fact that for all $d\in\N$ it holds that $\vass{(a+\pi\gamma^{-1}\scl^{-d}\vass{c}^{1-d})\gamma\scl^d c^{d-1}- a\gamma\scl^d c^{d-1}}= \pi$ \proves\ that for all $d\in\N$ there exists $\alpha\in[a,a+\pi\gamma^{-1}\scl^{-d}\vass{c}^{1-d})$ such that
	\begin{equation}
			\label{eq:sequence:start}
		\vass{\sin(\alpha \gamma\scl^d c^{d-1}+\varphi)}=1.
	\end{equation}
	\Moreover for all $d\in\N$, $\alpha\in\R$, $k\in\Z$ with $\vass{\sin(\alpha\gamma\scl^d c^{d-1}+\varphi)}=1$ it holds that
	\begin{equation}
		\label{eq:sequence:recursion}
		\begin{split}
			\sin((\alpha+k\pi\gamma^{-1}\scl^{-d}\vass{c}^{1-d}) \gamma\scl^d c^{d-1}+\varphi)
			&=\sin(\alpha \gamma\scl^d c^{d-1}+\varphi+k\vass{c}^{1-d}c^{d-1}\pi)
			\\&=(-1)^k \sin(\alpha \gamma\scl^d c^{d-1}+\varphi).
		\end{split}
	\end{equation}
	\Moreover the assumption that $\scl\geq 1$ and \eqref{eq:iteration_counter_prod}
	\prove\ that for all $d\in \N$, $k\in\N\cap(0,S(d)]$ it holds that
	\begin{equation}
		\begin{split}
		a
		\leq
		a+k\pi\gamma^{-1}\scl^{-d}\vass{c}^{1-d}
		&\leq
		a+S(d)\pi\gamma^{-1}\scl^{-d}\vass{c}^{1-d}
		\\&\leq
		a+(b-a)\gamma(\dom\scl)^{d-1}\pi^{-1}\pi\gamma^{-1}\scl^{-d}\vass{c}^{1-d}
		\\&=
		a+(b-a)\scl^{-1}
		\leq
		 b.
		\end{split}
	\end{equation}
This, 
\eqref{eq:prod:function:on:line}, 
\eqref{eq:sequence:start}, and
\eqref{eq:sequence:recursion}
\prove\ that there exist 
$v_{k,d}\in[a,b]^d$, $(k,d)\in\N^2$,
which satisfy that
\begin{enumerate}[(I)]
\item{
\label{it:vect:sequence:start:0}
it holds for all $d\in\N\cap(0,S(1))$ that $v_{1,1}\in[a,a+\pi\gamma^{-1}\scl^{-1})\subseteq[a,b]$,
}
\item{
	\label{it:vect:sequence:start}
	it holds for all $d\in\mathscr{N}\setminus\{1\}$ that  $v_{1,d}\in[a,a+\pi\gamma^{-1}\scl^{-d}\vass{c}^{1-d})\times\{c\}^{d-1}\subseteq[a,b]^d$,
}
\item{
	\label{it:vect:sequence:recursion:d=1}
	it holds for all $k\in\N\cap(1,\infty)\cap(0,S(1)]$ that $
		v_{k,1}=v_{k-1,1}+\pi\gamma^{-1}\scl^{-1}\in[a,b]$,
}
\item{
\label{it:vect:sequence:recursion}
it holds for all $d\in\N\setminus\{1\}$, $k\in\N\cap(1,\infty)\cap(0,S(d)]$ that 
\begin{equation}
	v_{k,d}=v_{k-1,d}+(\pi\gamma^{-1}\scl^{-d}\vass{c}^{1-d},0,0,\allowbreak\ldots,0)\in[a,b]\times\{c\}^{d-1}\subseteq[a,b]^d,
\end{equation}and
}
\item{
\label{it:vect:sequence:alternating}
it holds for all $d\in\N$ that there exists $s\in\{-\kappa,\kappa\}$ such that for all $k\in\N\cap(0,S(d)]$ it holds that $f_d(v_{k,d})=s (-1)^k$.
}
\end{enumerate}
Combining \cref{it:vect:sequence:start:0},
\cref{it:vect:sequence:start},
\cref{it:vect:sequence:recursion:d=1},
\cref{it:vect:sequence:recursion}, and
\cref{it:vect:sequence:alternating} with the fact that
for all $x,y\in\R$, $k\in\Z$ it holds that $\sin(x+2k\pi)=\sin(x)$ and $\vass{\sin(x)-\sin(y)}\leq \vass{x-y}$
 \proves[ei]\ \cref{prop:sin:cos:recursion,prop:sin:cos:periodicity,prop:sin:cos:lipschitz,prop:sin:cos:alternation}.
The proof of \cref{sin:cos:prod:targetfunctions} is thus complete.
\end{proof}

\cfclear
\begin{athm}{cor}{sin:cos:sum:targetfunctions}
Let $\varphi\in\R$, $\gamma,\kappa\in(0,\infty)$, $a\in\R$, $b\in[a+\pi\gamma^{-1},\infty)$,  $g\in C(\R,\R)$ satisfy for all $x\in\R$ that $\tgt(x)=\kappa\sin(x+\varphi)$,
let $f_d\in C(\R^d,\R)$, $d\in\N$, satisfy for all $d\in\N$, $x=(x_1,\ldots,x_d)\in[a,b]^d$ that $f_d(x)=g\prb{\gamma2^d\sum_{i=1}^d x_i}$,
let $S\colon\N\to\N$ satisfy for all $d\in\N$ that
		\begin{equation}
		S(d)=
		2^d d.
		\end{equation}
Then there exist $(\nu_d,\sigma_d)\in\prb{\R^d\setminus\{0\}}\times\{-2\kappa,2\kappa\}$, $d\in\N$, and $v_{k,d}\in[a,b]^d$, $(k,d)\in\N^2$, such that
\begin{enumerate}[(i)]
\item{
\label{prop:sum:sin:cos:recursion}
it holds for all $d\in\N\setminus\{1\}$, $k\in\N\cap(1,S(d)]$ that $v_{k,d}=v_{k-1,d}+\nu_d$,
}

\item{
\label{prop:sum:sin:cos:alternation}
it holds for all $d\in\N\setminus\{1\}$, $k\in\N\cap(1,S(d)]$ that $f_d(v_{k,d})-f_d(v_{k-1,d})=(-1)^k\sigma_d$,
}

\item{
\label{prop:sum:sin:cos:periodicity}
it holds for all $x\in\R$, $k\in\Z$ that $\tgt(x+2k\pi)=\tgt(x)\in[-\kappa,\kappa]$, and
}

\item{
\label{prop:sum:sin:cos:lipschitz}
it holds for all $x,y\in\R$ that $\vass{\tgt(x)-\tgt(y)}\leq \kappa\vass{x-y}$.
}
\end{enumerate}
\end{athm}

\begin{proof}[Proof of \cref{sin:cos:sum:targetfunctions}]
\Nobs that for all $d\in\N$, $v=(\alpha,\alpha,\ldots,\alpha)\in[a,b]^d$ it holds that
\begin{equation}
\label{eq:sum:function:on:line}
f_d(v)=g\prb{\gamma2^d\textstyle\sum_{i=1}^d \alpha}=\kappa\sin(\gamma 2^d d\alpha +\varphi).
\end{equation}
\Moreover the fact that for all $d\in\N$ it holds that $\pi=\vass{\gamma 2^d da-\gamma 2^d d(a+\pi\gamma^{-1}2^{-d}d^{-1})}$ \proves\ that for all $d\in\N$ there exists $\alpha\in[a,a+\pi\gamma^{-1}2^{-d}d^{-1})$ such that
\begin{equation}
\label{eq:sum:sequence:start}
\vass{\sin(\gamma 2^d d\alpha +\varphi)}=1.
\end{equation}
\Moreover for all $d\in\N$, $\alpha\in\R$, $k\in\Z$ with $\vass{\sin(\gamma 2^d d\alpha+\varphi)}=1$ it holds that
\begin{equation}
\label{eq:sum:sequence:recursion}
\sin(\gamma 2^d d(\alpha+k\pi\gamma^{-1}2^{-d}d^{-1})+\varphi)
=\sin(\gamma 2^d d\alpha+\varphi+k\pi)
=(-1)^k \sin(\gamma 2^d d\alpha+\varphi).
\end{equation}
\Moreover the assumption that for all $d\in\N$ it holds that $S(d)= 2^d d$
\proves\ that for all $d\in \N$, $k\in\{1,2,\ldots,S(d)\}$ it holds that
\begin{equation}
\label{eq:sum:sequence:location}
a
\leq
a+k\pi\gamma^{-1}2^{-d}d^{-1}\leq a+\pi\gamma^{-1}\leq b.
\end{equation}
This, 
\eqref{eq:sum:function:on:line},
\eqref{eq:sum:sequence:start}, and
\eqref{eq:sum:sequence:recursion}
\prove\ that there exist $v_{k,d}\in\R^d$, $(k,d)\in\N^2$, which satisfy that
\begin{enumerate}[(I)]
\item{
\label{it:vect:sum:sequence:start}
it holds for all $d\in\N$ that $v_{1,d}\in[a,a+\pi\gamma^{-1}2^{-d}d^{-1})^{d}\subseteq[a,b]^d$,
}
\item{
\label{it:vect:sum:sequence:recursion}
it holds for all $d\in\N\setminus\{1\}$, $k\in\N\cap(1,S(d)]$ that 
\begin{equation}
v_{k,d}=v_{k-1,d}+(\pi\gamma^{-1}2^{-d}d^{-1},\pi\gamma^{-1}2^{-d}d^{-1},\ldots,\pi\gamma^{-1}2^{-d}d^{-1})\in[a,b]^d,
\end{equation}
and
}
\item{
\label{it:vect:sum:sequence:alternating}
it holds for all $d\in\N$ that there exists $\sigma\in\{-\kappa,\kappa\}$ such that for all $k\in\{1,2,\ldots,S(d)\}$ it holds that $f_d(v_{k,d})=\sigma (-1)^k$.
}
\end{enumerate}
Combining 
\cref{it:vect:sum:sequence:start,it:vect:sum:sequence:recursion,it:vect:sum:sequence:alternating} with the fact that
for all $x,y\in\R$, $k\in\Z$ it holds that $\sin(x+2k\pi)=\sin(x)$ and $\vass{\sin(x)-\sin(y)}\leq \vass{x-y}$
 \proves[ei]\ \cref{prop:sum:sin:cos:recursion,prop:sum:sin:cos:periodicity,prop:sum:sin:cos:lipschitz,prop:sum:sin:cos:alternation}.
The proof of \cref{sin:cos:sum:targetfunctions} is thus complete.
\end{proof}

\subsection{Lower bounds for approximations of specific families of oscillating functions}
\label{subsection:lower_bound_results}

\cfclear
\begin{athm}{lemma}{final:approximation:fail:length:inequality}
	Let $a\in\R$, $b\in[a,\infty)$, $d \in \N$, $\varepsilon \in (0,\infty)$, $P\in[1,\infty)$, let $h\colon\N\to\N$ satisfy for all $l\in\N$ that $h(l)=\max\{1,l-1\}$, and let $f\colon\R^d\to\R$ satisfy for all $L\in\N$ that
	\begin{equation}
		\label{eq:approx:fail:length:equality}
		\min \pr*{ \pR*{ p \in \N \colon \PR*{\!\!
					\begin{array}{c}
						\exists\, \mathscr{f} \in \ANNSET \colon
						(\paramnotzero(\mathscr{f})=p)\land
						(\lengthANN(\mathscr{f})= L)
						\land{}\\
						(\realisation(\mathscr{f}) \in C(\R^d,\R))\land{}\\
						(\sup_{x\in[a,b]^d}\vass{(\realisation(\mathscr{f}))(x)-f(x)} \leq \eps)\\
					\end{array} \!\! }
			}
			\cup\{\infty\}
		}
		\geq \tfrac{h(L)}{2}\pr{2P}^{\nicefrac{1}{h(L)}}
	\end{equation}
	\cfload. Then it holds for all $H\in\N_0$ that
	\begin{equation}
		\label{eq:lower:inequality:bound}
		\min \pr*{ \pR*{ p \in \N \colon \PR*{\!\!
					\begin{array}{c}
						\exists\, \mathscr{f} \in \ANNSET \colon
						(\paramnotzero(\mathscr{f})\leq p)\land
						(\hidlengthANN(\mathscr{f})\leq H)
						\land{}\\
						(\realisation(\mathscr{f}) \in C(\R^d,\R))\land{}\\
						(\sup_{x\in[a,b]^d}\vass{(\realisation(\mathscr{f}))(x)-f(x)} \leq \eps)\\
					\end{array} \!\! }
			}
			\cup\{\infty\}
		}
		\geq P^{\nicefrac{1}{\max\{1,H\}}}
		.
	\end{equation}
\end{athm}

\begin{proof}[Proof of \cref{final:approximation:fail:length:inequality}]
	Throughout this proof assume without loss of generality that $P>1$.
	The assumption that $P>1$,
	the fact that there exists $\f\in\ANNSET$ with $\paramnotzero(\mathscr{f})=1$, $\lengthANN(\mathscr{f})=2$, $\realisation(\mathscr{f}) \in C(\R^d,\R)$, and $\forall\,x\in\R^d\colon(\realisation(\mathscr{f}))(x)=0$, and
	 \eqref{eq:approx:fail:length:equality} \prove\ that for all $\f\in\ANNSET$ with $\paramnotzero(\mathscr{f})=0$ and $\realisation(\mathscr{f}) \in C(\R^d,\R)$ it holds that
	\begin{equation}
				\llabel{eq:triv:case}
				\sup_{x\in[a,b]^d}\vass{(\realisation(\mathscr{f}))(x)-f(x)}
				=
				\sup_{x\in[a,b]^d}\vass{f(x)}
				 > \eps.
	\end{equation}
	This and
	 the fact that for all $\f\in\ANNSET$ it holds that $\paramnotzero(\mathscr{f})\in\N_0$ \prove\ that for all $H\in\N_0$ it holds that
	\begin{equation}
		\label{eq:lower:inequality:bound:1}
				\begin{split}
			&\min \pr*{ \pR*{ p \in \N \colon \PR*{\!\!
						\begin{array}{c}
							\exists\, \mathscr{f} \in \ANNSET \colon
							(\paramnotzero(\mathscr{f})\leq p)\land
							(\hidlengthANN(\mathscr{f})\leq H)
							\land{}\\
							(\realisation(\mathscr{f}) \in C(\R^d,\R))\land{}\\
							(\sup_{x\in[a,b]^d}\vass{(\realisation(\mathscr{f}))(x)-f(x)} \leq \varepsilon)\\
						\end{array} \!\! }
				}
				\cup\{\infty\}
			}
		\\&\geq
		\min \pr*{ \pR*{ p \in \N_0 \colon \PR*{\!\!
					\begin{array}{c}
						\exists\, \mathscr{f} \in \ANNSET \colon
						(\paramnotzero(\mathscr{f})= p)\land
						(\hidlengthANN(\mathscr{f})\leq H)
						\land{}\\
						(\realisation(\mathscr{f}) \in C(\R^d,\R))\land{}\\
						(\sup_{x\in[a,b]^d}\vass{(\realisation(\mathscr{f}))(x)-f(x)} \leq \varepsilon)\\
					\end{array} \!\! }
			}
			\cup\{\infty\}
		}
		\\&=
		\min \pr*{ \pR*{ p \in \N \colon \PR*{\!\!
					\begin{array}{c}
						\exists\, \mathscr{f} \in \ANNSET \colon
						(\paramnotzero(\mathscr{f})= p)\land
						(\hidlengthANN(\mathscr{f})\leq H)
						\land{}\\
						(\realisation(\mathscr{f}) \in C(\R^d,\R))\land{}\\
						(\sup_{x\in[a,b]^d}\vass{(\realisation(\mathscr{f}))(x)-f(x)} \leq \varepsilon)\\
					\end{array} \!\! }
			}
			\cup\{\infty\}
		}.
		\end{split}
	\end{equation}
	This,
	\eqref{eq:def:ANN:operators},
	 \eqref{eq:approx:fail:length:equality}, and the fact that for all $L\in\N$ it holds that $h(L)\leq h(L+1)$ \prove\ that for all $H,L\in\N_0$ with $H=L-1$ it holds that
	\begin{equation}
		\begin{split}
			&\min \pr*{ \pR*{ p \in \N \colon \PR*{\!\!
						\begin{array}{c}
							\exists\, \mathscr{f} \in \ANNSET \colon
							(\paramnotzero(\mathscr{f})\leq p)\land
							(\hidlengthANN(\mathscr{f})\leq H)
							\land{}\\
							(\realisation(\mathscr{f}) \in C(\R^d,\R))\land{}\\
							(\sup_{x\in[a,b]^d}\vass{(\realisation(\mathscr{f}))(x)-f(x)} \leq \varepsilon)\\
						\end{array} \!\! }
				}
				\cup\{\infty\}
			}
			\\
			&\geq\min \pr*{ \pR*{ p \in \N \colon \PR*{\!\!
						\begin{array}{c}
							\exists\, \mathscr{f} \in \ANNSET \colon
							(\paramnotzero(\mathscr{f})= p)\land
							(\hidlengthANN(\mathscr{f})\leq H)
							\land{}\\
							(\realisation(\mathscr{f}) \in C(\R^d,\R))\land{}\\
							(\sup_{x\in[a,b]^d}\vass{(\realisation(\mathscr{f}))(x)-f(x)} \leq \varepsilon)\\
						\end{array} \!\! }
				}
				\cup\{\infty\}
			}
			\\&
			=\min \pr*{ \pR*{ p \in \N \colon \PR*{\!\!
						\begin{array}{c}
							\exists\, \mathscr{f} \in \ANNSET \colon
							(\paramnotzero(\mathscr{f})= p)\land
							(\lengthANN(\mathscr{f})\leq L)
							\land{}\\
							(\realisation(\mathscr{f}) \in C(\R^d,\R))\land{}\\
							(\sup_{x\in[a,b]^d}\vass{(\realisation(\mathscr{f}))(x)-f(x)} \leq \varepsilon)\\
						\end{array} \!\! }
				}
				\cup\{\infty\}
			}
			\\&
			=\min_{l\in\{1,2,\ldots,L\}}
			\min \pr*{ \pR*{ p \in \N \colon \PR*{\!\!
						\begin{array}{c}
							\exists\, \mathscr{f} \in \ANNSET \colon
							(\paramnotzero(\mathscr{f})= p)\land
							(\lengthANN(\mathscr{f})=l)
							\land{}\\
							(\realisation(\mathscr{f}) \in C(\R^d,\R))\land{}\\
							(\sup_{x\in[a,b]^d}\vass{(\realisation(\mathscr{f}))(x)-f(x)} \leq \eps)\\
						\end{array} \!\! }
				}
				\cup\{\infty\}
			}
			\\&
			\geq \min_{l\in\{1,2,\ldots,L\}} \tfrac{h(l)}{2}\pr{2P}^{\nicefrac{1}{h(l)}}
			\geq \min_{l\in\{1,2,\ldots,L\}} P^{\nicefrac{1}{h(l)}}
			=P^{\nicefrac{1}{h(L)}}
			=P^{\nicefrac{1}{\max\{1,H\}}}.
		\end{split}
	\end{equation}
This and \eqref{eq:lower:inequality:bound:1} \prove\ \eqref{eq:lower:inequality:bound}.
	The proof of \cref{final:approximation:fail:length:inequality} is thus complete.
\end{proof}

\cfclear
\begin{athm}{prop}{LowerBound:prod}
	Let $\varphi\in\R$, $\gamma\in(0,1]$, $\beta\in[1,\infty)$, $a\in\R$, $b\in[a,\infty)$, $\dom,\kappa\in(0,\infty)$ satisfy $\dom=\max\{\vass{a},\vass{b}\}>\scl^{-1}$ and
	for every $d\in\N$ let
	$f_d\colon\R^d\to\R$ satisfy for all $x=(x_1,\ldots,x_d)\in\R^d$ that $f_d(x)=\kappa\sin\prb{\gamma\beta^d\prb{\sprod_{i = 1}^d x_i}+\varphi}$. Then\label{it:LowerBound:prod:eq1}
			it holds for all $d\in\N$, $H\in\N_0$, $\eps\in(0,\kappa)$ that 
				\begin{multline}
					\label{eq:statement:lower:bound:nonzero:param}
					\min \pr*{ \pR*{ p \in \N \colon \PR*{\!\!
								\begin{array}{c}
									\exists\, \mathscr{f} \in \ANNSET \colon
									(\paramnotzero(\mathscr{f})\leq p)\land
									(\hidlengthANN(\mathscr{f})\leq H)
									\land{}\\
									(\realisation(\mathscr{f}) \in C(\R^d,\R))\land{}\\
									(\sup_{x\in[a,b]^d}\vass{(\realisation(\mathscr{f}))(x)-f_d(x)} \leq \varepsilon)\\
								\end{array} \!\! }
						}
						\cup\{\infty\}
					}
				\\\geq 	\vass[\Big]{\frac{\floor{(b-a)\gamma(\dom\scl)^{d-1}\pi^{-1}-1}}{2}}^{\nicefrac{1}{\max\{1,H\}}}
				\end{multline}
	\cfout.
\end{athm}

\begin{proof}[Proof of \cref{LowerBound:prod}]
	Throughout this proof 
	assume without loss of generality that $b>a$
	 and
	let $h\colon\N\to\N$ and $S\colon\N\to[0,\infty)$ satisfy for all $d\in\N$ that
	\begin{equation}
		\label{eq:setup:S:and:h:fct}
		h(d)=\max\{1,d-1\}
		\qandq
		S(d)=\floor{(b-a)\gamma(\dom\scl)^{d-1}\pi^{-1}}.
	\end{equation}
	\Nobs that \eqref{eq:setup:S:and:h:fct} \proves\ that for all
	$d\in\N$, $H\in\N_0$, $\eps\in(0,\kappa)$ with $S(d)\leq3$ it holds that 
		\begin{multline}
			\label{eq:for:tiny:S:cover}
			\min \pr*{ \pR*{ p \in \N \colon \PR*{\!\!
						\begin{array}{c}
							\exists\, \mathscr{f} \in \ANNSET \colon
							(\paramnotzero(\mathscr{f})\leq p)\land
							(\hidlengthANN(\mathscr{f})\leq H)
							\land{}\\
							(\realisation(\mathscr{f}) \in C(\R^d,\R))\land{}\\
							(\sup_{x\in[a,b]^d}\vass{(\realisation(\mathscr{f}))(x)-f_d(x)} \leq \varepsilon)\\
						\end{array} \!\! }
				}
				\cup\{\infty\}
			}
			\\\geq
			1
			\geq  \vass[\Big]{\frac{S(d)-1}{2}}^{\nicefrac{1}{\max\{1,H\}}}.
	\end{multline}
	\Moreover \cref{sin:cos:prod:targetfunctions} (applied with
	$\varphi \curvearrowleft  \varphi$,
	$\kappa \curvearrowleft  \kappa$,
	$\dom \curvearrowleft  \dom$,
	$\gamma \curvearrowleft  \gamma$,
	$\scl \curvearrowleft  \scl$,
	$a \curvearrowleft  a$,
	$b \curvearrowleft  b$,
	$g\curvearrowleft \pr{\R\ni x\mapsto \kappa\sin(x+\varphi)\in\R}$,
	$(f_d)_{d\in\N} \curvearrowleft (f_d)_{d\in\N}$,
	$S \curvearrowleft  S$
	in the notation of \cref{sin:cos:prod:targetfunctions}) \proves\ that there exist $(\nu_d,\sigma_d)\in\prb{\R^d\setminus\{0\}}\times\{-2\kappa,2\kappa\}$, $d\in\N$, and $v_{k,d}\in[a,b]^d$, $(k,d)\in\N^2$, such that
	\begin{enumerate}[(I)]
		\item{
			\label{it:prop:sin:cos:recursion}
			it holds for all $d\in\N$, $k\in\N\cap(1,\infty)\cap(0,S(d)]$ that $v_{k,d}=v_{k-1,d}+\nu_d$ and
		}
		
		\item{
			\label{it:prop:sin:cos:alternation}
			it holds for all $d\in\N$, $k\in\N\cap(1,\infty)\cap(0,S(d)]$ that $f_d(v_{k,d})-f_d(v_{k-1,d})=(-1)^k\sigma_d$.
		}
%
%
	\end{enumerate}
	\Nobs that \cref{gen:final:approximation:fail} (applied with
	$a \curvearrowleft  a$,
	$b \curvearrowleft  b$,
	$d \curvearrowleft  d$,
	$H \curvearrowleft  h(L)$,
	$\nu \curvearrowleft  \nu_d$,
	$\kappa \curvearrowleft  \kappa$,
	$\sigma \curvearrowleft  \sigma_d$,
	$\eps \curvearrowleft  \eps$,
	$S \curvearrowleft  S(d)$,
	$(v_k)_{k\in\N\cap(0,S(d)]} \curvearrowleft (v_{k,d})_{k\in\N\cap(0,S(d)]}$,
	$f \curvearrowleft  f_d$
	for $d,L\in\N$, $\eps\in(0,\kappa)$ in the notation of \cref{gen:final:approximation:fail}) \hence \proves\ that for all $d,L\in\N$, $\eps\in(0,\kappa)$ with $S(d)>3$ it holds that 
	\begin{equation}
		\begin{split}
			&\min \pr*{ \pR*{ p \in \N \colon \PR*{\!\!
						\begin{array}{c}
							\exists\, \mathscr{f} \in \ANNSET \colon
							(\paramnotzero(\mathscr{f})=p)\land
							(\lengthANN(\mathscr{f})= L)
							\land{}\\
							(\realisation(\mathscr{f}) \in C(\R^d,\R))\land{}\\
							(\sup_{x\in[a,b]^d}\vass{(\realisation(\mathscr{f}))(x)-f_d(x)} \leq \varepsilon)\\
						\end{array} \!\! }
				}
				\cup\{\infty\}
			}
		\geq \tfrac{h(L)}{2}\pr{S(d)-1}^{\nicefrac{1}{h(L)}}
		\end{split}
	\end{equation}
	\cfload. 
	This and \cref{final:approximation:fail:length:inequality} 
	(applied with
	$a \curvearrowleft  a$,
	$b \curvearrowleft  b$,
	$d \curvearrowleft  d$,
	$\eps \curvearrowleft  \eps$,
	$P \curvearrowleft  \frac{S(d)-1}{2}$,
	$h \curvearrowleft  h$,
	$f \curvearrowleft  f_d$
	for $d\in\N$, $\eps\in(0,\kappa)$ in the notation of \cref{final:approximation:fail:length:inequality}) \prove\ that for all $d\in\N$, $H\in\N_0$, $\eps\in(0,\kappa)$ with $S(d)\geq 3$ it holds that 
	\begin{equation}
		\begin{split}
			&\min \pr*{ \pR*{ p \in \N \colon \PR*{\!\!
						\begin{array}{c}
							\exists\, \mathscr{f} \in \ANNSET \colon
							(\paramnotzero(\mathscr{f})\leq p)\land
							(\hidlengthANN(\mathscr{f})\leq H)
							\land{}\\
							(\realisation(\mathscr{f}) \in C(\R^d,\R))\land{}\\
							(\sup_{x\in[a,b]^d}\vass{(\realisation(\mathscr{f}))(x)-f_d(x)} \leq \varepsilon)\\
						\end{array} \!\! }
				}
				\cup\{\infty\}
			}\geq  \prb{\tfrac{S(d)-1}{2}}^{\nicefrac{1}{\max\{1,H\}}}.
		\end{split}
	\end{equation}
Combining this and
\eqref{eq:setup:S:and:h:fct} with \eqref{eq:for:tiny:S:cover} 
\proves\ 
	that for all $d\in\N$, $H\in\N_0$, $\eps\in(0,\kappa)$ it holds that 
\begin{multline}
	\min \pr*{ \pR*{ p \in \N \colon \PR*{\!\!
				\begin{array}{c}
					\exists\, \mathscr{f} \in \ANNSET \colon
					(\paramnotzero(\mathscr{f})\leq p)\land
					(\hidlengthANN(\mathscr{f})\leq H)
					\land{}\\
					(\realisation(\mathscr{f}) \in C(\R^d,\R))\land{}\\
					(\sup_{x\in[a,b]^d}\vass{(\realisation(\mathscr{f}))(x)-f_d(x)} \leq \varepsilon)\\
				\end{array} \!\! }
		}
		\cup\{\infty\}
	}
	\\\geq  \vass[\Big]{\frac{S(d)-1}{2}}^{\nicefrac{1}{\max\{1,H\}}}
	=
		 \vass[\Big]{\frac{\floor{(b-a)\gamma(\dom\scl)^{d-1}\pi^{-1}-1}}{2}}^{\nicefrac{1}{\max\{1,H\}}}
	.
\end{multline}
	The proof of \cref{LowerBound:prod} is thus complete.
\end{proof}

\cfclear
\begin{athm}{prop}{LowerBound:sum}
	Let $\varphi\in\R$, $\gamma,\kappa\in(0,\infty)$, $a\in\R$, $b\in[a+\pi\gamma^{-1},\infty)$ and
	for every $d\in\N$ let
	$f_d\colon\R^d\to\R$ satisfy for all $x=(x_1,\ldots,x_d)\in\R^d$ that $f_d(x)=\kappa\sin\prb{\gamma 2^d\prb{\ssum_{i = 1}^d x_i}+\varphi}$. Then\label{it:LowerBound:sum:eq1}
	it holds for all $d\in\N$, $H\in\N_0$, $\eps\in(0,\kappa)$ that 
	\begin{equation}
		\begin{split}
			&\min \pr*{ \pR*{ p \in \N \colon \PR*{\!\!
						\begin{array}{c}
							\exists\, \mathscr{f} \in \ANNSET \colon
							(\paramnotzero(\mathscr{f})\leq p)\land
							(\hidlengthANN(\mathscr{f})\leq H)
							\land{}\\
							(\realisation(\mathscr{f}) \in C(\R^d,\R))\land{}\\
							(\sup_{x\in[a,b]^d}\vass{(\realisation(\mathscr{f}))(x)-f_d(x)} \leq \varepsilon)\\
						\end{array} \!\! }
				}
				\cup\{\infty\}
			}
		\geq  \prb{d 2^{d-1}}^{\nicefrac{1}{\max\{1,H\}}}
		\end{split}
	\end{equation}
	\cfout.
\end{athm}

\begin{proof}[Proof of \cref{LowerBound:sum}]
	Throughout this proof 
	let $h\colon\N\to\N$ satisfy for all $L\in\N$ that $h(L)=\max\{1,L-1\}$ and let $S_d\in\N$, $d\in\N$, satisfy for all $d\in\N$ that $S_d=2^d d$.
	\Nobs that \cref{sin:cos:sum:targetfunctions} (applied with
	$\varphi \curvearrowleft  \varphi$,
	$\kappa \curvearrowleft  \kappa$,
	$\gamma \curvearrowleft  \gamma$,
	$a \curvearrowleft  a$,
	$b \curvearrowleft  b$,
	$f_d \curvearrowleft f_d$,
	$S(d) \curvearrowleft  S_d$
	for $d\in\N$ in the notation of \cref{sin:cos:sum:targetfunctions}) \proves\ that there exist $(\nu_d,\sigma_d)\in\prb{\R^d\setminus\{0\}}\times\{-2\kappa,2\kappa\}$, $d\in\N$, and $v_{k,d}\in[a,b]^d$, $k\in\{1,2,\ldots,S_d\}$, $d\in\N$, such that
	\begin{enumerate}[(I)]
		\item{
			\label{it2:prop:sin:cos:recursion}
			it holds for all $d\in\N\setminus\{1\}$, $k\in\N\cap(1,S_d]$ that $v_{k,d}=v_{k-1,d}+\nu_d$ and
		}
		\item{
			\label{it2:prop:sin:cos:alternation}
			it holds for all $d\in\N\setminus\{1\}$, $k\in\N\cap(1,S_d]$ that $f_d(v_{k,d})-f_d(v_{k-1,d})=(-1)^k\sigma_d$.
		}
	\end{enumerate}
	\Nobs that \cref{gen:final:approximation:fail} (applied with
	$a \curvearrowleft  a$,
	$b \curvearrowleft  b$,
	$d \curvearrowleft  d$,
	$H \curvearrowleft  h(L)$,
	$\nu \curvearrowleft  \nu_d$,
	$\kappa \curvearrowleft  \kappa$,
	$\sigma \curvearrowleft  \sigma_d$,
	$\eps \curvearrowleft  \eps$,
	$S \curvearrowleft  S_d$,
	$(v_1,v_2,\ldots,v_{S_d}) \curvearrowleft (v_{1,d},v_{2,d},\ldots,v_{S_d,d})$,
	$f \curvearrowleft  f_d$
	for $d,L\in\N$, $\eps\in(0,\kappa)$ in the notation of \cref{gen:final:approximation:fail}) \proves\ that for all $d\in\N\setminus\{1\}$, $L\in\N$, $\eps\in(0,\kappa)$ it holds that 
	\begin{equation}
		\begin{split}
			&\min \pr*{ \pR*{ p \in \N \colon \PR*{\!\!
						\begin{array}{c}
							\exists\, \mathscr{f} \in \ANNSET \colon
							(\paramnotzero(\mathscr{f})=p)\land
							(\lengthANN(\mathscr{f})= L)
							\land{}\\
							(\realisation(\mathscr{f}) \in C(\R^d,\R))\land{}\\
							(\sup_{x\in[a,b]^d}\vass{(\realisation(\mathscr{f}))(x)-f_d(x)} \leq \varepsilon)\\
						\end{array} \!\! }
				}
				\cup\{\infty\}
			} \geq \tfrac{h(L)}{2}\pr{S_d-1}^{\nicefrac{1}{h(L)}}
		\end{split}
	\end{equation}
	\cfload. 
	This and \cref{final:approximation:fail:length:inequality} (applied with
	$a \curvearrowleft  a$,
	$b \curvearrowleft  b$,
	$d \curvearrowleft  d$,
	$\eps \curvearrowleft  \eps$,
	$P \curvearrowleft  \frac{S_d-1}{2}$,
	$h \curvearrowleft  h$,
	$f \curvearrowleft  f_d$
	for $d\in\N\setminus\{1\}$, $\eps\in(0,\kappa)$ in the notation of \cref{final:approximation:fail:length:inequality}) \prove\ that for all $d\in\N\setminus\{1\}$, $H\in\N_0$, $\eps\in(0,\kappa)$ it holds that 
	\begin{multline}
			\min \pr*{ \pR*{ p \in \N \colon \PR*{\!\!
						\begin{array}{c}
							\exists\, \mathscr{f} \in \ANNSET \colon
							(\paramnotzero(\mathscr{f})\leq p)\land
							(\hidlengthANN(\mathscr{f})\leq H)
							\land{}\\
							(\realisation(\mathscr{f}) \in C(\R^d,\R))\land{}\\
							(\sup_{x\in[a,b]^d}\vass{(\realisation(\mathscr{f}))(x)-f_d(x)} \leq \varepsilon)\\
						\end{array} \!\! }
				}
				\cup\{\infty\}
			}
		 \\\geq \prbb{\frac{S_d-1}{2}}^{\nicefrac{1}{\max\{1,H\}}}
		 	=
		 \prbb{\frac{2^d d-1}{2}}^{\nicefrac{1}{\max\{1,H\}}}
		 \geq
		 \prbb{2^{d-1} d -\frac{1}{2}}^{\nicefrac{1}{\max\{1,H\}}}
		.
		\end{multline}
		Combining this with the fact that for all $d,h,p\in\N$ with $p^h\geq d 2^{d-1}-\frac{1}{2}$ it holds that $p^h\geq d 
		2^{d-1}$ \proves\ that for all $d\in\N$, $H\in\N_0$, $\eps\in(0,\kappa)$ it holds that 
		\begin{multline}
			\min \pr*{ \pR*{ p \in \N \colon \PR*{\!\!
						\begin{array}{c}
							\exists\, \mathscr{f} \in \ANNSET \colon
							(\paramnotzero(\mathscr{f})\leq p)\land
							(\hidlengthANN(\mathscr{f})\leq H)
							\land{}\\
							(\realisation(\mathscr{f}) \in C(\R^d,\R))\land{}\\
							(\sup_{x\in[a,b]^d}\vass{(\realisation(\mathscr{f}))(x)-f_d(x)} \leq \varepsilon)\\
						\end{array} \!\! }
				}
				\cup\{\infty\}
			}
			\geq
			\prb{d2^{d-1}}^{\nicefrac{1}{\max\{1,H\}}}
			.
		\end{multline}
	The proof of \cref{LowerBound:sum} is thus complete.
\end{proof}

\section[Upper bounds for the minimal number of ANN parameters in the approximation of certain high-dimensional functions]{%
	Upper bounds for the minimal number of ANN parameters in the
	approximation of certain high-dimen\-sional functions}
\label{S5_upper_bounds}

In this section we establish in 
\cref{cor:downsized:products}, 
\cref{final:approximation:imp}, and \cref{cor:rescaled:final:approximation3} 
below suitable upper bounds for the minimal number of parameters of \ANNs\ to approximate the product functions (\cref{cor:downsized:products}) and certain highly oscillating functions (\cref{final:approximation:imp} and \cref{cor:rescaled:final:approximation3}) in the case where the absolute values of the parameters of the \ANNs\ are assumed to be uniformly bounded by $1$.
%

Our proof of \cref{cor:downsized:products} employs
the elementary result regarding the reduction of the absolute value of the size of the parameters of an \ANN\ without changing its realization function in \cref{lemma:network:halfed:paramsize} 
and the essentially well-known upper bounds for the minimal number of parameters of \ANNs\ to approximate certain scaled product functions in \cref{lemma:d_product}.
\cref{lemma:d_product} is an extended variant of, \eg, Beneventano et al.~\cite[Proposition~6.8]{beneventano21}.
Our proof of \cref{lemma:d_product} utilizes
the elementary result regarding suitable deep \ANNs\ whose realization functions agree with appropriate one-dimensional scaling functions in \cref{gen:scaling:networks} 
and the essentially well-known upper bound result for the minimal number of parameters of \ANNs\ approximating the product functions in \cref{lemma:2d_product}.
\cref{lemma:2d_product} is a slightly extended variant of
Beneventano et al.~\cite[Lemma~6.7]{beneventano21}.
For its proof and the auxiliary results in \cref{high:dim:prod},
we adapt the constructions and proofs in
Beneventano et al.~\cite[Section~6]{beneventano21}.
The respective adaptations and additional estimates are
specified before the individual results.
%

Our proof of \cref{final:approximation:imp} employs
the elementary result regarding the reduction of the absolute value of the size of the parameters of an \ANN\ without changing its realization function in \cref{lemma:network:halfed:paramsize} 
and the upper bounds for the minimal number of parameters of \ANNs\ approximating compositions of certain periodic functions and certain scaled product functions in \cref{final:approximation:streched}.
Our proof of \cref{final:approximation:streched}, in turn, combines \cref{gen:scaling:networks} 
and \cref{lemma:d_product}
with the essentially well-known upper bound result for the minimal number of parameters of \ANNs\ approximating certain periodic functions in \cref{sin:approximation}.
Our proof of \cref{sin:approximation} 
 employs the essentially well-known \ANN\ approximation result for certain one-dimensional Lipschitz continuous functions in \cref{propo:sin:approximation} and
builds on the essentially well-known properties of sawtooth functions (suitable one-dimensional piecewise linear functions with compact support) in \cref{edgy:sin:prop} and \cref{lem:edgy:sin}. 
The results in \cref{edgy:sin:prop} and \cref{lem:edgy:sin} are extensions of, \eg, Telgarsky \cite[Section~2.2]{Telgarsky15}
and \cref{propo:sin:approximation} is inspired by Beneventano et al.~\cite[Subsection~4.1]{beneventano21}.
%

Our proof of \cref{cor:rescaled:final:approximation3} employs \cref{lemma:network:halfed:paramsize} 
as well as the upper bounds for the minimal number of parameters of \ANNs\ approximating compositions of certain periodic functions and scaled sum functions in \cref{cor:rescaled:final:approximation2}.
Our proof of \cref{cor:rescaled:final:approximation2}, in turn, utilizes \cref{gen:scaling:networks} 
and \cref{sin:approximation}.

%

\subsection{One-dimensional scaling ANNs}
\label{subsec:multiplication}

\cfclear

\begin{athm}{lemma}{scaling:networks}
	Let $\beta\in(0,\infty)$, $\wdt,n\in\N$. Then there exists $\f\in\ANNSET$ such that
	\begin{enumerate}[(i)]
		\item{
			\label{scaling:realisation}
			it holds for all $x\in\R$ that $(\realisation(\f))(x)=(\wdt\beta)^n x$,
		}
		\item{
			\label{scaling:dims}
			it holds that $\dims(\f)=(1,2\wdt,2\wdt,\ldots,2\wdt,1)\in\N^{n+2}$,
		}
		\item{
			\label{scaling:size}
			it holds that $\insize(\f)=1$, $\outsize(\f)=\beta$, and $\size(\f)=\max\{1,\beta\}$, and
		}
		\item{
			\label{scaling:param}
			it holds that $\param(\f)=(4n-4)\wdt^2+(2n+4)\wdt+1$
		}
	\end{enumerate}
	\cfout.
\end{athm}

\begin{proof}[Proof of \cref{scaling:networks}]
	Throughout this proof let $W_1\in\R^{2\wdt \times 1}$, $W_2\in\R^{1\times 2\wdt}$, $W_3\in\R^{2\wdt\times 2\wdt}$ satisfy
	\begin{equation}
		\label{def:scaling:ANNs}
		\begin{split}
			W_1=
			\begin{pmatrix}
				1\\
				-1\\
				1\\
				-1\\
				\vdots\\
				1\\
				-1
			\end{pmatrix}
			,
			\qquad
			W_2=\begin{pmatrix}
				\beta &  -\beta & \beta &  -\beta &\cdots & \beta &  -\beta
			\end{pmatrix}
			,
			\qandq
			W_3=W_1W_2,
					\end{split}
				\end{equation}
				and let $\f_k\in\ANNSET$, $k\in\N$, satisfy for all $k\in\N$ that
				$\f_1=\pr{\pr{W_1,0},\pr{W_2,0}}\in\ANNSET$
				, and $\f_{k+1}=\compANN{\f_{k}}{\f_1}$ \cfload. \Nobs that \eqref{def:scaling:ANNs} \proves\ that for all $x\in\R$ it holds that
				\begin{equation}
					\label{eq:scaling:basecase}
					(\realisation(\f_1))(x)=\wdt(\beta(\RELU(x))-\beta(\RELU(-x)))=\wdt\beta x.
				\end{equation}
				This and \cref{Lemma:PropertiesOfCompositions_n2} \prove\ that for all $x\in\R$, $k\in\N$ with $\fa{y}\R\colon(\realisation(\f_k))(y)=(\wdt\beta)^k y$ it holds that
				\begin{equation}
					\label{eq:scaling:inductionstep}
					\begin{split}
					(\realisation(\f_{k+1}))(x)=(\realisation(\compANN{\f_{k}}{\f_1}))(x)
					=(\realisation(\f_k))\prb{(\realisation(\f_1))(x)}
					&=(\realisation(\f_k))\pr{\wdt\beta x}
					\\&=(\wdt\beta)^{k+1}x.
					\end{split}
				\end{equation}
				Combining this and \eqref{eq:scaling:basecase} with induction \proves\ that for all $x\in\R$ it holds that
				\begin{equation}
					\label{eq:scaling:proving:it:1}
					(\realisation(\f_n))(x)=(\wdt\beta)^{n} x.
				\end{equation}
				\Moreover \eqref{ANNoperations:Composition} and \eqref{def:scaling:ANNs} \prove\ that
				\begin{equation}
					\f_2=\f_1\bullet\f_1=
					\prb{ \pr{W_1,0},\pr{W_1 W_2,W_1\cdot 0+0},\pr{W_2,0}	 }=\prb{\pr{W_1,0},\pr{W_3,0},\pr{W_2,0}}.
				\end{equation}
				Combining this, \eqref{ANNoperations:Composition}, and \eqref{def:scaling:ANNs} with induction \proves\ that
				\begin{equation}
					\label{eq:structured:form:of:scaling:network}
					\begin{split}
						\f_n
						&=\prb{\pr{W_1,0},\pr{W_3,0},\pr{W_3,0},\dots,\pr{W_3,0},\pr{W_2,0}}
						\\&\in
						\prb{(\R^{2\wdt \times 1} \times \R^{2\wdt}) \times\prb{\times_{k=1}^{n-1}\pr{\R^{2\wdt \times 2\wdt} \times \R^{2\wdt}}}\times (\R^{1 \times 2\wdt} \times \R^1) }.
					\end{split}
				\end{equation}
				This and \eqref{def:scaling:ANNs} \prove\ that
				\begin{equation}
					\label{eq:scaling:proving:it:2}
					\insize(\f_n)=\insize(\f_1)=1,
					\quad
					\outsize(\f_n)=\outsize(\f_1)=\beta,
					\quad\text{and}\quad
					\size(\f_n)=\size(\f_1)=\max\{1,\beta\}
				\end{equation}
				\Moreover \eqref{eq:structured:form:of:scaling:network} \proves\ that
				\begin{equation}
					\label{eq:scaling:proving:it:3}
					\dims(\f_n)=(1,2\wdt,2\wdt,\ldots,2\wdt,1)\in\N^{n+2}.
				\end{equation}
				This \proves\ that
				\begin{equation}
					\begin{split}
						\param(\f_n)
						=\textstyle\sum_{k=1}^{n+1}\singledims_k(\f_n)(\singledims_{k-1}(\f_n)+1)
						&=2\wdt(1+1)+(n-1)(2\wdt(2\wdt+1))+1(2\wdt+1)
						\\&=4(n-1)\wdt^2+(4+2(n-1)+2)\wdt+1
						\\&=(4n-4)\wdt^2+(2n+4)\wdt+1.
					\end{split}
				\end{equation}
				\cfload. Combining this, \eqref{eq:scaling:proving:it:1}, 
				\eqref{eq:scaling:proving:it:2},
				\eqref{eq:scaling:proving:it:3}, and
				\eqref{eq:structured:form:of:scaling:network}
				 \proves[ei]\ \cref{scaling:realisation,scaling:dims,scaling:size,scaling:param}.
				The proof of \cref{scaling:networks} is thus complete.
			\end{proof}

			\cfclear
			\begin{athm}{cor}{cor:extended:network}
				Let $\f\in\ANNSET$, $L\in\N$ satisfy $\outDimANN(\f)=1$ and $L>\lengthANN(\f)$. Then there exists $\g\in\ANNSET$ which satisfies that
				\begin{enumerate}[(i)]
					\item{
						\label{extension:function}
						it holds that $\realisation(\f)=\realisation(\g)$,
					}
					\item{
						\label{extension:dims}
						it holds for all $k\in\N_0\cap[0,L]$ that 
						\begin{equation}
							\singledims_k(\g)=
							\begin{cases}
								\singledims_k(\f) & \colon k\in\N_0\cap[0,\lengthANN(\f))\\
								2 & \colon k\in\N\cap[\lengthANN(\f),L)\\
								1 & \colon k=L,\\
							\end{cases}
						\end{equation}
					}
					\item{
						\label{extension:length}
						it holds that $\lengthANN(\g)=L$,
						and
					}
					\item{
						\label{extension:size}
						it holds that $\insize(\g)=\insize(\f)$, $\outsize(\g)=1$, and $\size(\g)=\max\{1,\size(\f)\}$
					}
				\end{enumerate}
				\cfout.
			\end{athm}
			
			\begin{proof}[Proof of \cref{cor:extended:network}]
				Throughout this proof let $\h_1,\h_2,\ldots,\h_{L}\in\ANNSET$ satisfy for all $k\in\N\cap(1,L]$ that
				\begin{equation}
					\label{eq:extending:networks}
					\h_1=\ReLUidANN{1}
					\qandq
					\h_{k}=\compANN{\ReLUidANN{1}}{\h_{k-1}}
				\end{equation}
				\cfload. Combining 
				\eqref{eq:extending:networks},
				\dimsbene, and
				\cref{Prop:identity_representation} with induction \proves\ that
				\begin{equation}
					\lengthANN(\h_{L-\lengthANN(\f)})=L-\lengthANN(\f)+1
					\qandq
					\dims(\h_{L-\lengthANN(\f)})=(1,2,2,\ldots,2,1)\in\N^{L-\lengthANN(\f)+2}.
				\end{equation}
				This, \cref{Lemma:PropertiesOfCompositions_n2}, and \dimsbene\ \prove\ that for all $k\in\N_0\cap[0,L]$ it holds that
				\begin{equation}
					\label{eq:extended:dims}
					\lengthANN(\compANN{\h_{L-\lengthANN(\f)}}{\f})=L
					\qandq
					\singledims_k(\compANN{\h_{L-\lengthANN(\f)}}{\f})=
					\begin{cases}
						\singledims_k(\f) & \colon k\in\N_0\cap[0,\lengthANN(\f))\\
						2 & \colon k\in\N\cap[\lengthANN(\f),L)\\
						1 & \colon k=L.\\
					\end{cases}
				\end{equation}
				\Moreover \cref{Lemma:PropertiesOfCompositions_n2} and
				\cref{Prop:identity_representation} \prove\ that for all $k\in\N\cap(0,L-\lengthANN(\f))$ it holds that
				\begin{equation}
					\realisation(\h_{k+1})
					=\realisation(\compANN{\ReLUidANN{1}}{\h_k})
					=(\realisation(\ReLUidANN{1}))\circ(\realisation(\h_{k}))
					=\realisation(\h_{k}).
				\end{equation}
				This, \cref{Prop:identity_representation}, and induction \prove\ that for all $x\in\R$ it holds that
				\begin{equation}
					(\realisation(\h_{L-\lengthANN(\f)}))(x)=(\realisation(\h_{1}))(x)=(\realisation(\ReLUidANN{1}))(x)=x.
				\end{equation}
				Combining this and the assumption that $\outDimANN(\f)=1$ with \cref{Lemma:PropertiesOfCompositions_n2} \proves\ that
				\begin{equation}
					\label{eq:extended:function}
					(\realisation(\compANN{\h_{L-\lengthANN(\f)}}{\f}))
					=(\realisation(\h_{L-\lengthANN(\f)}))\circ(\realisation(\f))
					=\realisation(\f).
				\end{equation}
				\Moreover \eqref{eq:extending:networks},
				\cref{Prop:identity_representation:prop}, and induction \prove\ that
				\begin{equation}
					\outsize(\compANN{\h_{L-\lengthANN(\f)}}{\f})=1,
					\quad\quad
					\insize(\compANN{\h_{L-\lengthANN(\f)}}{\f})=\insize(\f),
				\end{equation}
				\begin{equation}
					\qandq
					\size(\compANN{\h_{L-\lengthANN(\f)}}{\f})=\max\{1,\size(\f)\}.
				\end{equation}
				Combining this, \eqref{eq:extended:dims}, and \eqref{eq:extended:function}  \proves[ei]\
				\cref{extension:length,extension:dims,extension:function,extension:size}.
				The proof of \cref{cor:extended:network} is thus complete.
			\end{proof}

			\cfclear
			\begin{athm}{cor}{gen:scaling:networks}
				Let $\beta\in\R\backslash\{0\}$, $L\in\N_0$ satisfy $L\geq \log_2(\abs{\beta})$ \cfload. Then there exists $\f\in\ANNSET$ such that
				\begin{enumerate}[(i)]
					\item{
						\label{gen:scaling:realisation}
						it holds for all $x\in\R$ that $(\realisation(\f))(x)=\beta x$,
					}
					\item{
						\label{gen:scaling:dims}
						it holds that $\dims(\f)=(1,2,2,\ldots,2,1)\in\N^{L+2}$,
					}
					\item{
						\label{gen:scaling:size}
						it holds that $\insize(\f)\leq 1$, $\outsize(\f)\leq 2$, and $\size(\f)\leq 2$, and
					}
					\item{
						\label{gen:scaling:param}
						it holds that $\param(\f)\leq 6\max\{L,1\}+1$
					}
				\end{enumerate}
				\cfout.
			\end{athm}
			
			\begin{proof}[Proof of \cref{gen:scaling:networks}]
				Throughout this proof assume without loss of generality that 
				$L= \min\pr{\N_0\cap[\log_2(\abs{\beta}),\infty)}$ \cfadd{cor:extended:network}\cfload{} and 
				$\vass{\beta}> 1$ (otherwise consider $((\beta),0)\in(\R^{1 \times 1} \times \R)\subseteq\ANNSET$),
				and let $\g_2\in\ANNSET$ satisfy
				\begin{equation}
					\label{eq:sgn:netw}
					\g_2=
					\pr*{\!\begin{pmatrix}
							\tfrac{\beta}{\abs{\beta}}
						\end{pmatrix},
						0}   \in \prb{\R^{1 \times 1} \times \R^{1} }
				\end{equation}
				\cfload. \Nobs that \cref{scaling:networks} (applied with 
				$\beta \curvearrowleft \abs{\beta}^{\nicefrac{1}{L}}$,
				$n \curvearrowleft L$,
				$\wdt \curvearrowleft 1$
				in the notation of \cref{scaling:networks}) \proves\ that there exists $\g_1\in\ANNSET$ which satisfies that
				\begin{enumerate}[(I)]
					\item{
						\label{for:gen:scaling:realisation}
						it holds for all $x\in\R$ that $(\realisation(\g_1))(x)=\vass{\beta} x$,
					}
					\item{
						\label{for:gen:scaling:dims}
						it holds that $\dims(\g_1)=(1,2,2,\ldots,2,1)\in\N^{L+2}$,
					}
					\item{
						\label{for:gen:scaling:size}
						it holds that $\insize(\g_1)=1$, $\outsize(\g_1)=\abs{\beta}^{\nicefrac{1}{L}}$, and $\size(\g_1)=\abs{\beta}^{\nicefrac{1}{L}}$, and
					}
					\item{
						\label{for:gen:scaling:param}
						it holds that $\param(\g_1)=6L+1$
					}
				\end{enumerate}
				\cfload. \Nobs that \eqref{eq:sgn:netw}, \cref{for:gen:scaling:dims}, and \dimsbene\ \prove\ that
				\begin{equation}
					\label{eq:dims:gen:scaling}
					\dims(\compANN{\g_2}{\g_1})=\dims(\g_1)=(1,2,2,\ldots,2,1)\in\N^{L+2}
					\quad\text{and}\quad
					\param(\compANN{\g_2}{\g_1})=\param(\g_1)=6L+1
				\end{equation}
				\cfload. \Moreover \eqref{eq:sgn:netw}, \cref{for:gen:scaling:realisation}, and \cref{Lemma:PropertiesOfCompositions_n2} \prove\ that for all $x\in\R$ it holds that
				\begin{equation}
					\label{eq:function:gen:scaling}
					(\realisation(\compANN{\g_2}{\g_1}))(x)=(\realisation(\g_2))(\realisation(\g_1)(x))=\tfrac{\beta}{\abs{\beta}}(\vass{\beta} x)=\beta x.
				\end{equation}
				\Moreover 
				\eqref{ANNoperations:Composition},
				\eqref{eq:sgn:netw}, 
				\cref{for:gen:scaling:size}, 
				\cref{lem:sizecomp},
				and the fact that $\abs{\beta}^{\nicefrac{1}{L}}\leq 2$ \prove\ that for all $x\in\R$ it holds that $\insize(\compANN{\g_2}{\g_1})=\insize(\g_1)=1$,
				\begin{equation}
					\label{eq:size:gen:scaling}
					\outsize(\compANN{\g_2}{\g_1})=\vass[\big]{\tfrac{\beta}{\abs{\beta}}}\outsize(\g_1)\leq 2,
					\quad\text{and}\quad
					\size(\compANN{\g_2}{\g_1})\leq \max\{\size(\g_1),\outsize(\compANN{\g_2}{\g_1})\}\leq 2.
				\end{equation}
				Combining this and \eqref{eq:dims:gen:scaling} with \eqref{eq:function:gen:scaling}  \proves[ei]\
				\cref{gen:scaling:realisation,gen:scaling:dims,gen:scaling:size,gen:scaling:param}.
				The proof of \cref{gen:scaling:networks} is thus complete.
			\end{proof}

\subsection{Trade-off between the number and the size of ANN parameters}
\label{subsec:tradeoff}

\newcommand{\dnsc}{\lambda}

\cfclear
\begin{athm}{cor}{cor:ANN:downscaled:param}
	Let $\f\in\ANNSET$, $L,d\in\N$, $\dnsc\in(0,1]$ satisfy $\realisation(\f)\in C(\R^d,\R)$ and $\lengthANN(\f)=L$ \cfload. Then there exists $\g\in\ANNSET$ such that
	\begin{enumerate}[(i)]
		\item{\label{ANN:downscaled:param:function}it holds for all $x\in\R^{d}$ that $\pr{\realisation(\g)}(x)=\dnsc^{L}\pr{\realisation(\f)}(x)$,}
		\item{\label{ANN:downscaled:param:length}it holds that $\lengthANN(\g)=L$,}
		\item{\label{ANN:downscaled:param:size}it holds that $\size(\g)\leq \dnsc\size(\f)$, and}
		\item{\label{ANN:downscaled:param:dims}it holds that $\dims(\g)=\dims(\f)$}
	\end{enumerate}
	\cfout.
\end{athm}

\begin{proof}[Proof of \cref{cor:ANN:downscaled:param}]
	Throughout this proof let $l_0,l_1,\ldots,l_L\in\N$ satisfy $\dims(\f)=(l_0,l_1,\ldots,l_L)$,
	let $\g\in\ANNSET$ satisfy for all $k\in\{1,2,\ldots,L\}$ that
	\begin{equation}
		\label{eq:setup:downscaled:params}
		\lengthANN(\g)=L,
		\qquad
		\weight{k}{\g}=\dnsc\weight{k}{\f},
		\qandq
		\bias{k}{\g}=\dnsc^{k}\bias{k}{\f},
	\end{equation}
	and let $x_0, y_0 \in \R^{l_0},\, x_1,y_1 \in \R^{l_1}, 
	\ldots,\, \allowbreak x_{L},y_L \in \R^{l_{L}}$ satisfy for all $k\in\{1,2,\ldots,L\}$ that 
	\begin{equation}
		\label{eq:setup:downscaled:seq}
		x_0=y_0,
		\qquad
		x_k=\RELU(\weight{k}{\f} x_{k-1}+\bias{k}{\f}),
		\qandq
		y_k=\RELU(\weight{k}{\g} y_{k-1}+\bias{k}{\g}).
	\end{equation}
	\Nobs that \eqref{eq:setup:downscaled:params} \proves\ that
	\begin{equation}
		\label{eq:size:and:dims:for:downscaled:ANN}
		\size(\g)\leq \dnsc\size(\f)
		\qandq
		\dims(\g)=\dims(\f).
	\end{equation}
	\Moreover \eqref{eq:setup:downscaled:seq} \proves\ that for all $k\in\N\cap(0,L)$ with $y_k=\dnsc^{k}x_k$ it holds that
	\begin{equation}
		\begin{split}
			y_{k+1}
			=\RELU(\weight{k+1}{\g} y_{k}+\bias{k+1}{\g})
			&=\RELU(\weight{k+1}{\g} (\dnsc^{k}x_k)+\bias{k+1}{\g})
			\\&=\RELU(\dnsc^{k+1}(\weight{k+1}{\f} x_k+\bias{k+1}{\f}))
			\\&=\dnsc^{k+1}x_{k+1}.
		\end{split}
	\end{equation}
	Combining this and \eqref{eq:setup:downscaled:seq} with induction \proves\ that $y_{L-1}=\dnsc^{L-1}x_{L-1}$. This,
	\eqref{ANNrealization:ass2},
	 and \eqref{eq:setup:downscaled:seq} \prove\ that
	\begin{equation}
		\begin{split}
			\label{eq:scaled:realisation:at:x0}
			\pr{\realisation(\g)}(x_0)
			=\weight{L}{\g} y_{L-1}+\bias{L}{\g}
			&=\weight{L}{\g}(\dnsc^{L-1}x_{L-1})+\bias{L}{\g}
			\\&=\dnsc^{L}(\weight{L}{\f} x_{L-1}+\bias{L}{\f})
			\\&=\dnsc^{L}\pr{\realisation(\f)}(x_0)
			.
		\end{split}
	\end{equation}
	Combining this and \eqref{eq:size:and:dims:for:downscaled:ANN}  \proves[ei]\ \cref{ANN:downscaled:param:function,ANN:downscaled:param:length,ANN:downscaled:param:size,ANN:downscaled:param:dims}.
	The proof of \cref{cor:ANN:downscaled:param} is thus complete.
\end{proof}

\begin{definition}[Ceiling and floor of real numbers]
	\label{def:ceiling}
	We denote by $\ceil{\cdot} \colon \R \to \Z$ and $\floor{\cdot} \colon \R \to \Z$ the functions which satisfy for all $x \in \R$ that $\ceil{x} = \min(\Z \cap [x, \infty))$ and $\floor{x}= \max(\Z \cap (-\infty,x])$.
\end{definition}

\cfclear

\cfclear
\begin{athm}{lemma}{lemma:network:halfed:paramsize}
	Let $\f\in\ANNSET$, $d\in\N$, $\dnsc\in(0,1)$ satisfy $\realisation(\f)\in C(\R^d,\R)$ \cfload. Then there exists $\g\in\ANNSET$ such that
	\begin{enumerate}[(i)]
		\item{\label{network:halfed:paramsize:function}
			it holds that $\realisation(\g)=\realisation(\f)$,
		}
		\item{\label{network:halfed:paramsize:length}
			it holds that $\lengthANN(\g)\leq \ceil[\big]{1-\tfrac{\ln(\dnsc)}{\ln(2)}}\lengthANN(\f)+1$,
		}
		\item{\label{network:halfed:paramsize:dims}
			it holds for all $k\in\N_0\cap[0,\lengthANN(\g)]$ that
			\begin{equation}
				\singledims_k(\g)=
				\begin{cases}
					\singledims_k(\f) &\colon k\in\N_0\cap[0,\lengthANN(\f))\\
					2 &\colon k=\lengthANN(\f)\\
					4 &\colon k\in\N\cap(\lengthANN(\f),\lengthANN(\g))\\
					1 &\colon k=\lengthANN(\g),
				\end{cases}
			\end{equation}
		}
		\item{\label{network:halfed:paramsize:param}
			it holds that $\param(\g)\leq\param(\f)+\singledims_{\hidlengthANN(\f)}(\f)+20\ceil[\big]{\tfrac{\vass{\ln(\dnsc)}}{\ln(2)}}\lengthANN(\f)\leq 2\param(\f)+20\ceil[\big]{\tfrac{\vass{\ln(\dnsc)}}{\ln(2)}}\lengthANN(\f)$, and
		}
		\item{\label{network:halfed:paramsize:size}
			it holds that $\size(\g)\leq\max\pR{1,\dnsc\size(\f)}$
		}
	\end{enumerate}
	\cfout.
\end{athm}

\begin{proof}[Proof of \cref{lemma:network:halfed:paramsize}]
	Throughout this proof 
	let $n,L\in\N$, $\beta\in\R$ satisfy $L=\lengthANN(\f)$, $n=\ceil{L\ln(\dnsc^{-1})\ln(2)^{-1}}$, and $\beta=2^{-1}\dnsc^{-\frac{L}{n}}$ \cfload.
	\Nobs that \cref{cor:ANN:downscaled:param} (applied with 
	$\f \curvearrowleft \f$,
	$L \curvearrowleft L$,
	$d \curvearrowleft d$,
	$\dnsc \curvearrowleft \dnsc$
	in the notation of \cref{cor:ANN:downscaled:param})
	\proves\ that there exists $\g_1\in\ANNSET$ which satisfies that
	\begin{enumerate}[(I)]
		\item{\label{it:ANN:downscaled:param:function}it holds for all $x\in\R^{d}$ that $\pr{\realisation(\g_1)}(x)=\dnsc^L\pr{\realisation(\f)}(x)$,}
		\item{\label{it:ANN:downscaled:param:length}it holds that $\lengthANN(\g_1)=L$,}
		\item{\label{it:ANN:downscaled:param:size}it holds that $\size(\g_1)\leq \dnsc\size(\f)$, and}
		\item{\label{it:ANN:downscaled:param:dims}it holds that $\dims(\g_1)=\dims(\f)$}
	\end{enumerate}
	\cfload. \Nobs that the fact that $\beta\leq 1$ and \cref{scaling:networks} (applied with 
	$\beta \curvearrowleft \beta$,
	$\wdt \curvearrowleft 2$,
	$n \curvearrowleft n$
	in the notation of \cref{scaling:networks})
	\prove\ that there exists $\g_2\in\ANNSET$ which satisfies that
	\begin{enumerate}[(A)]
		\item{
			\label{it:cor:scaling:realisation}
			it holds for all $x\in\R$ that $(\realisation(\g_2))(x)=(2\beta)^n x$,
		}
		\item{
			\label{it:cor:scaling:dims}
			it holds that $\dims(\g_2)=(1,4,4,\ldots,4,1)\in\N^{n+2}$,
		}
		\item{
			\label{it:cor:scaling:size}
			it holds that $\size(\g_2)=1$, and
		}
		\item{
			\label{it:cor:scaling:param}
			it holds that $\param(\g_2)=(4n-4)2^2+(2n+4)2+1=20n-7\leq 20n$.
		}
	\end{enumerate}
	\Nobs that \cref{it:ANN:downscaled:param:function}, \cref{it:cor:scaling:realisation}, \cref{Lemma:PropertiesOfCompositions_n2}, \cref{Prop:identity_representation},
	and the fact that $(2\beta)^n=\dnsc^{-L}$ \prove\ that for all $x\in\R^d$ it holds that
	\begin{equation}
		\label{eq:network:halfed:paramsize:function}
		\begin{split}
		(\realisation(\compANN{\g_2}{\ReLUidANN{1}}\bullet\g_1))(x)
		=\PR{\realisation(\g_2)\circ\realisation(\g_1)}(x)
		=(2\beta)^n\prb{\dnsc^{L}(\realisation(\f))(x)}
		=(\realisation(\f))(x)
		\end{split}
	\end{equation}
	\cfload. \Moreover \cref{it:ANN:downscaled:param:length}, \cref{it:cor:scaling:dims}, \cref{Lemma:PropertiesOfCompositions_n2}, and \cref{Prop:identity_representation} \prove\ that
	\begin{equation}
		\label{eq:network:halfed:paramsize:length}
		\begin{split}
		\lengthANN(\compANN{\g_2}{\ReLUidANN{1}}\bullet\g_1)
		=\lengthANN(\g_2)+\lengthANN(\ReLUidANN{1})+\lengthANN(\g_1)-2
		&=(n+1)+2+L-2
		\\&=n+L+1
		\\&=\ceil{L\ln(\dnsc^{-1})\ln(2)^{-1}}+L+1
		\\&\leq \ceil{\ln(\dnsc^{-1})\ln(2)^{-1}+1}L+1
		.
		\end{split}
	\end{equation}
	Combining this, 
	\cref{it:ANN:downscaled:param:dims}, 
	\cref{it:cor:scaling:dims}, 
	\dimsbene,
	and \cref{Prop:identity_representation} \proves\ that for all $k\in\N_0\cap[0,\lengthANN(\compANN{\g_2}{\ReLUidANN{1}}\bullet\g_1)]$ it holds that
	\begin{equation}
		\label{eq:network:halfed:paramsize:dims}
		\singledims_k(\compANN{\g_2}{\ReLUidANN{1}}\bullet\g_1)=
		\begin{cases}
			\singledims_k(\f) &\colon k\in\N_0\cap[0,L)\\
			2 &\colon k=L\\
			4 &\colon k\in\N\cap(L,\lengthANN(\compANN{\g_2}{\ReLUidANN{1}}\bullet\g_1))\\
			1 &\colon k=\lengthANN(\compANN{\g_2}{\ReLUidANN{1}}\bullet\g_1).
		\end{cases}
	\end{equation}
	This, \eqref{eq:network:halfed:paramsize:length}, and the fact that $\singledims_{L}(\f)=1$ \prove\ that
	\begin{equation}
		\label{eq:network:halfed:paramsize:param}
		\begin{split}
			\param(\compANN{\g_2}{\ReLUidANN{1}}\bullet\g_1)
			&=\sum_{k=1}^{n+L+1}\singledims_k(\compANN{\g_2}{\ReLUidANN{1}}\bullet\g_1)(\singledims_{k-1}(\compANN{\g_2}{\ReLUidANN{1}}\bullet\g_1)+1)
			\\&=\PR*{\sum_{k=1}^{L-1}\singledims_k(\f)(\singledims_{k-1}(\f)+1)}+2(\singledims_{L-1}(\f)+1)+4(2+1)
			\\&\quad+(n-1)(4(4+1))+1(4+1)
			\\&= \PR*{\sum_{k=1}^{L}\singledims_k(\f)(\singledims_{k-1}(\f)+1)}+\singledims_{L-1}(\f)+1+12+20n-20+5
			\\&\leq\param(\f)+\singledims_{L-1}(\f)+20n
			\\&\leq\param(\f)+\singledims_{L-1}(\f)+20\ceil[\big]{\tfrac{\vass{\ln(\dnsc)}
				}{\ln(2)}}L
			.
		\end{split}
	\end{equation}
	\Moreover \cref{it:ANN:downscaled:param:size}, \cref{it:cor:scaling:size}, and \cref{Prop:identity_representation:prop} \prove\ that
	\begin{equation}
		\label{eq:network:halfed:paramsize:size}
		\size(\compANN{\g_2}{\ReLUidANN{1}}\bullet\g_1)
		=\max\{\size(\g_2),\size(\g_1)\}
		\leq\max\pR{1,\dnsc\size(\f)}.
	\end{equation}
	Combining this,
	\eqref{eq:network:halfed:paramsize:function},
	\eqref{eq:network:halfed:paramsize:length},
	\eqref{eq:network:halfed:paramsize:dims},
	and
	\eqref{eq:network:halfed:paramsize:param}
	 \proves[ei]\
	\cref{network:halfed:paramsize:function,network:halfed:paramsize:length,network:halfed:paramsize:dims,network:halfed:paramsize:size,network:halfed:paramsize:param}.
	The proof of \cref{lemma:network:halfed:paramsize} is thus complete.
\end{proof}

\subsection{Upper bounds for approximations of product functions}
\label{high:dim:prod}

\cfclear
\begin{athm}{lemma}{lemma:base:square}
Let $N\in\N$. Then there exists $\f\in\ANNSET$ such that
\begin{enumerate}[(i)]
    \item \label{base_continuous} it holds 
    that $\realisation(\mathscr{f}) \in C(\R,\R)$,

    \item\label{base_approx} it holds for all $x\in[0,1]$ that
    $0\leq  \pr*{\realisation\pr*{\mathscr{f}}} (x) -x^2\leq 4^{-N-1}$,
		
		\item\label{base_approx_outside} it holds for all 
    $x \in \R\setminus[0,1]$ that 
    $\pr*{\realisation\pr*{\mathscr{f}}} (x)=\RELU(x)$,
		
		\item\label{base_lipschitz} it holds for all 
    $x,y \in \R$ that 
    $\vass[\big]{\pr*{\realisation\pr*{\mathscr{f}}} (x) - \pr*{\realisation\pr*{\mathscr{f}}} (y)} \leq 2\vass{x-y}$,

    \item\label{base_dims} it holds 
    that
    $\dims(\mathscr{f}) =(1,4,4,\ldots,4,1)\in\N^{N+2}
    $, and
		
		\item\label{base_size} it holds 
    that
    $\size(\mathscr{f}) \leq 4
    $
\end{enumerate}
\cfout.
\end{athm}

\begin{proof}[Proof of \cref{lemma:base:square}]
\Nobs that Lemma 5.1 and Lemma 5.2 in Grohs et al.~\cite{GrohsIbrgimovJentzen2021} \prove\ that there exists $\f\in\ANNSET$ which satisfies that
\begin{enumerate}[(I)]
    \item \label{origin:base_continuous} it holds 
    that $\realisation(\mathscr{f}) \in C(\R,\R)$,

    \item\label{origin:base_approx}  it holds that
    $\sup_{x \in [0,1]}\vass[\big]{x^2 - \pr*{\realisation\pr*{\mathscr{f}}} (x)} \leq 4^{-N-1}$,
		
		\item\label{origin:base_approx_outside} it holds for all 
    $x \in \R\setminus[0,1]$ that 
    $\pr*{\realisation\pr*{\mathscr{f}}} (x)=\RELU(x)$,
		
		\item\label{origin:base_lipschitz} it holds for all 
    $k\in\N_0\cap[0,2^N)$, $x \in \big[\tfrac{k}{2^N},\tfrac{k+1}{2^N}\big)$ that 
    $\pr*{\realisation\pr*{\mathscr{f}}} (x) =\prb{\frac{2k+1}{2^N}}x-\frac{k^2+k}{2^{2N}}$,

    \item\label{origin:base_dims} it holds 
    that
    $\dims(\mathscr{f}) =(1,4,4,\ldots,4,1)\in\N^{N+2}
    $, and
		
		\item\label{origin:base_size} it holds 
    that
    $\size(\mathscr{f}) \leq 4
    $
\end{enumerate}
\cfload[.]\Nobs that \cref{origin:base_lipschitz} \proves\ that for all $k\in\N_0\cap[0,2^N)$, $x \in \big[\tfrac{k}{2^N},\tfrac{k+1}{2^N}\big)$ it holds that
\begin{equation}
	\begin{split}
	0
	\leq
	\prb{\tfrac{k+1}{2^N}-x}\prb{x-\tfrac{k}{2^N}}
	=
	\prb{\tfrac{2k+1}{2^N}}x-\tfrac{k^2+k}{2^{2N}}-x^2
	=
	\pr*{\realisation\pr*{\mathscr{f}}} (x) -x^2.
	\end{split}
\end{equation}
This, \cref{origin:base_continuous}, and \cref{origin:base_lipschitz} \prove\ \cref{base_continuous,base_approx}.
\Nobs that \cref{origin:base_continuous}, \cref{origin:base_approx_outside}, \cref{origin:base_lipschitz}, and the triangle inequality \prove\ that for all $x,y\in\R$ it holds that
\begin{equation}
\vass[\big]{\pr*{\realisation\pr*{\mathscr{f}}} (x) - \pr*{\realisation\pr*{\mathscr{f}}} (y)}\leq\max\pRb{1,\tfrac{2^{N+1}-1}{2^N}}\vass{x-y}\leq 2\vass{x-y}.
\end{equation}
Combining this with
\cref{origin:base_approx_outside,origin:base_dims,origin:base_size}
\proves[ie]\
 \cref{base_approx_outside,base_dims,base_size,base_lipschitz}.
The proof of \cref{lemma:base:square} is thus complete.
\end{proof}

In the following result, \cref{lemma:general:square} below,
we adapt the rescaling construction in
Beneventano et al.~\cite[Lemma~6.2]{beneventano21}.
The additional scaling layers allow us to bound the absolute
values of the parameters independently of the radius $R$.

\cfclear
\begin{athm}{lemma}{lemma:general:square}
Let $N\in\N$, $R\in(1,\infty)$. Then there exists $\f\in\ANNSET$ such that
\begin{enumerate}[(i)]
    \item \label{general:square_continuous} it holds 
    that $\realisation(\mathscr{f}) \in C(\R,\R)$,

    \item\label{general:square_approx} it holds for all $x\in[-\radius,\radius]$ that
    $0\leq  \pr*{\realisation\pr*{\mathscr{f}}} (x) -x^2\leq R^2 4^{-N-1}$,

		\item\label{general:square_lipschitz} it holds for all 
    $x,y \in \R$ that 
    $\vass[\big]{\pr*{\realisation\pr*{\mathscr{f}}} (x) - \pr*{\realisation\pr*{\mathscr{f}}} (y)} \leq 2R\vass{x-y}$,

    \item\label{general:square_length} it holds 
    that
    $\lengthANN(\mathscr{f}) =N+\ceil{\log_2(R)}+4
    $,

    \item\label{general:square_dims} it holds 
		for all $k\in\N_0\cap[0,\lengthANN(\f)]$ that
    \begin{equation} 
        \singledims_k(\mathscr{f})=
\begin{cases}
1 &\colon k=0\\
2 &\colon k\in\N\cap(0,2]\\
4 &\colon k\in\N\cap(2,N+2]\\
2 &\colon k\in\N\cap(N+2, N+\ceil{\log_2(R)}+4)\\
1 &\colon k=N+\ceil{\log_2(R)}+4,
\end{cases}
\end{equation}
and
		
		\item\label{general:square_size} it holds 
    that
    $\size(\mathscr{f}) \leq 4
    $
\end{enumerate}
\cfout.
\end{athm}

\begin{proof}[Proof of \cref{lemma:general:square}]
Throughout this proof let $n\in\N$ satisfy $n=\ceil{\log_2(R)}$ and
        let 
            $\g_1\in\ANNSET$
        satisfy
        \begin{equation}
        \label{Lemma:square_r:networks}
            \g_1
            =
            \pr*{\!\pr*{\!
            \begin{pmatrix}
                    \radius^{-1}\\
                    -\radius^{-1}
            \end{pmatrix}
            ,
            \begin{pmatrix}
                    0\\
                    0
            \end{pmatrix}
            \!}
            ,
            \pr*{
            \begin{pmatrix}
                    1 & 1
            \end{pmatrix}
            ,
            0
            }\!}\in\pr*{\pr*{\R^{2\times 1}\times\R^2}\times\pr*{\R^{1\times 2}\times\R}}
        \end{equation}
    \cfload.
		 \Moreover
        \cref{lemma:base:square}  (applied with
		$N \curvearrowleft  N$
in the notation of \cref{lemma:base:square})
    \proves\ that there exists 
        $\g_2 \in \ANNSET$ 
    which satisfies that
    \begin{enumerate}[(I)]
    \item \label{it:base_continuous} it holds 
    that $\realisation(\g_2) \in C(\R,\R)$,

    \item\label{it:base_approx}
    it holds for all $x\in[0,1]$ that
    $0\leq  \pr*{\realisation\pr*{\g_2}} (x) -x^2\leq 4^{-N-1}$,
		
		\item\label{it:it:base_approx_outside} it holds for all 
    $x \in \R\setminus[0,1]$ that 
    $\pr*{\realisation\pr*{\g_2}} (x)=\RELU(x)$,
		
		\item\label{it:base_lipschitz} it holds for all 
    $x,y \in \R$ that 
    $\vass[\big]{\pr*{\realisation\pr*{\g_2}} (x) - \pr*{\realisation\pr*{\g_2}} (y)} \leq 2\vass{x-y}$,

    \item\label{it:base_dims} it holds 
    that
    $\dims(\g_2) =(1,4,4,\ldots,4,1)\in\N^{N+2}
    $, and

		
		\item\label{it:base_size} it holds 
    that
    $\size(\g_2) \leq 4
    $\ifnocf.
    \end{enumerate}
    \cfload[.]%
		\Nobs that \cref{scaling:networks} (applied with
		$\beta \curvearrowleft  R^{\frac{2}{n}}$,
		$B \curvearrowleft  1$,
		$n \curvearrowleft  n$
in the notation of \cref{scaling:networks}) \proves\ that there exists $\g_3\in\ANNSET$ which satisfies that
		\begin{enumerate}[(A)]
\item{
\label{general:square:scaling:realisation}
it holds for all $x\in\R$ that $(\realisation(\g_3))(x)=R^2 x$,
}
\item{
\label{general:square:scaling:dims}
it holds that $\dims(\g_3)=(1,2,2,\ldots,2,1)\in\N^{n+2}$, and
}
\item{
\label{general:square:scaling:size}
it holds that $\insize(\g_3)=1$, $\outsize(\g_3)=R^{\frac{2}{n}}$, and $\size(\g_3)=R^{\frac{2}{n}}$.
}
\end{enumerate}
    Next let $\mathscr{f} \in \ANNSET$ satisfy
    \begin{equation}
    \label{Square_r:phi}
        \mathscr{f} 
        = 
        \compANN{\compANN{\g_3\bullet\ReLUidANN{1}}{\g_2\bullet\ReLUidANN{1}}}{\g_1}
    \end{equation}
    \cfload.
    \Nobs that 
        \eqref{Lemma:square_r:networks}, 
        \cref{it:base_dims}, 
        \cref{general:square:scaling:dims}, 
		\cref{Lemma:PropertiesOfCompositions_n2},
        \dimsbene, and
		\cref{Prop:identity_representation}
    \prove\ that for all $k\in\N_0\cap[0,\lengthANN(\f)]$ it holds that
    \begin{equation} 
    \label{Square_r:eq_pr_-1}
		\lengthANN(\f)=N+n+4
		\qandq
        \singledims_k(\mathscr{f})=
\begin{cases}
	1 &\colon k=0\\
	2 &\colon k\in\N\cap(0,2]\\
	4 &\colon k\in\N\cap(2,N+2]\\
	2 &\colon k\in\N\cap(N+2, N+n+4)\\
	1 &\colon k=N+n+4.
\end{cases}
\end{equation}
    \Moreover 
        \eqref{Lemma:square_r:networks}, 
        \eqref{Square_r:phi}, 
				\cref{general:square:scaling:realisation},
        \cref{Lemma:PropertiesOfCompositions_n2},
				and
				\cref{Prop:identity_representation}
    \prove\ that for all 
        $x \in \R$ 
    it holds that
    \begin{equation}
    \label{Square_r:phi_reduced}
        \realisation(\mathscr{f}) \in C(\R,\R)
        \qandq
        (\realisation(\mathscr{f}))(x) 
        = 
        \radius^2 \PRb{(\realisation(\g_2))\prb{\tfrac{\vass{x}}{\radius}}}
        \ifnocf.
    \end{equation}
    \cfload[.]%
    Combining 
        this 
    with 
        \cref{it:base_approx} 
    \proves\ that for all $x\in[-R,R]$ it holds that
    \begin{equation}
    \begin{split}
    \label{Square_r:eq_pr_0}
        &
        (\realisation(\mathscr{f}))(x) -x^2 
        = \radius^2 \PRb{(\realisation(\g_2))\prb{\tfrac{\vass{x}}{\radius}}
        - \prb{\tfrac{\vass{x}}{\radius}}^2}
        \in
        [0,R^2 4^{-N-1}]
        .
    \end{split}
    \end{equation}
    \Moreover 
        \eqref{Square_r:phi_reduced} and 
        \cref{it:base_lipschitz}  
    \prove\ that for all 
        $x,y\in\R$ 
    it holds that
    \begin{equation}
    \begin{split}
    \label{Square_r:eq_pr_2}
        \vass*{(\realisation(\mathscr{f}))(x)-(\realisation(\mathscr{f}))(y) }
        &=
        \vass[\big]{\radius^2 \PRb{(\realisation(\g_2))\prb{\tfrac{\vass{x}}{\radius}}}-\radius^2 \PRb{(\realisation(\g_2))\prb{\tfrac{\vass{y}}{\radius}}}}
        \\&\leq 
        2\radius^2\vass[\big]{\tfrac{\vass{x}}{\radius}-\tfrac{\vass{y}}{\radius}}
        \leq 
        2\radius\vass{x-y}.
    \end{split}
    \end{equation}
		\Moreover \eqref{Lemma:square_r:networks}, 
		\eqref{Square_r:phi}, 
		\cref{it:base_size}, 
		\cref{general:square:scaling:size}, 
		\cref{Prop:identity_representation:prop2}, and the fact that $1< R\leq 2^n$ \prove\ that
		\begin{equation}
		\begin{split}
		\size(\f)&\leq\max\{\size(\g_3),\size(\g_2),\size(\g_1)\}
		\leq\max\pRb{R^{\frac{2}{n}},4,1}= 4
		.
		\end{split}
		\end{equation}
    Combining 
        this 
    with 
        \eqref{Square_r:eq_pr_-1}, 
        \eqref{Square_r:phi_reduced}, 
        \eqref{Square_r:eq_pr_0}, and
        \eqref{Square_r:eq_pr_2}
    \proves[ie]\ 
        \cref{general:square_size,general:square_approx,general:square_dims,general:square_length,general:square_continuous,general:square_lipschitz}. The proof of \cref{lemma:general:square} is thus complete.
\end{proof}

In the following result, \cref{lemma:product} below, we adapt
the construction and proof of \benetwoprod\ using the
square approximation networks in \cref{lemma:general:square}.
We additionally establish the uniform bound on the absolute
values of the parameters of the approximating \ANN\ in \cref{product_size}.

\cfclear
\begin{athm}{lemma}{lemma:product}
Let $N\in\N$, $R\in(1,\infty)$. Then there exists $\f\in\ANNSET$ such that
\begin{enumerate}[(i)]
    \item \label{product_continuous} it holds 
    that $\realisation(\mathscr{f}) \in C(\R^2,\R)$,

    \item\label{product_approx} it holds that
    $\sup_{x,y \in [-R,R]}\vass[\big]{xy - \pr*{\realisation\pr*{\mathscr{f}}} (x,y)} \leq R^2 2^{-2N}$,

		\item\label{product_lipschitz} it holds for all 
    $x,y \in \R^2$ that 
    $\vass[\big]{\pr*{\realisation\pr*{\mathscr{f}}} (x) - \pr*{\realisation\pr*{\mathscr{f}}} (y)} \leq \sqrt{32}\radius \norm{x-y}$,

    \item\label{product_length} it holds 
    that
    $\lengthANN(\mathscr{f}) =N+\ceil{\log_2(R)}+7
    $,

    \item\label{product_dims} it holds 
		for all $k\in\N_0\cap[0,\lengthANN(\f)]$ that
    \begin{equation} 
        \singledims_k(\mathscr{f})=
\begin{cases}
2 &\colon k=0\\
6 &\colon k\in\N\cap(0,3]\\
12 &\colon k\in\N\cap(3,N+3]\\
6 &\colon k\in\N\cap(N+3, N+\ceil{\log_2(R)}+7)\\
1 &\colon k=N+\ceil{\log_2(R)}+7,
\end{cases}
\end{equation}
and

		
		\item\label{product_size} it holds 
    that
    $\size(\mathscr{f}) \leq 4
    $
\end{enumerate}
\cfout.
\end{athm}

\begin{proof}
[Proof of \cref{lemma:product}.]
Throughout this proof let $\g_1 \in\pr*{\R^{3\times 2}\times\R^3}\subseteq \ANNSET$, $\g_3 \in\pr*{\R^{1\times 3}\times\R}\subseteq \ANNSET$ satisfy
\begin{equation}
\label{Lemma:prod_r:eq:networks}
    \g_1 = \pr*{ \!\begin{pmatrix} 1 & 1 \\ 1 & 0 \\ 0 & 1   \end{pmatrix}, \begin{pmatrix} 0 \\ 0 \\ 0    \end{pmatrix}\!}\qandq
		\g_3 = \pr*{ \begin{pmatrix} \tfrac{1}{2} & -\tfrac{1}{2} & -\tfrac{1}{2}    \end{pmatrix}, 0}
\end{equation}
\cfload. \Nobs that \eqref{Lemma:prod_r:eq:networks} \proves\ that
\begin{equation}
\label{Lemma:prod_r:eq:varphi_new}
\dims(\g_1) = (2, 3),\;\;
\dims(\g_3) = (3,1),\;\;
\realisation(\g_1) \in C(\R^2,\R^3),
\;\;\text{and}\;\;
\realisation(\g_3) \in C(\R^3,\R)
\ifnocf. 
\end{equation}
\cfload[.]%
\Moreover \eqref{Lemma:prod_r:eq:networks} \proves\ that for all $x,y,z \in \R$ it holds that
\begin{equation}
\label{Lemma:prod_r:eq:varphi}
(\realisation\pr{\g_1})(x,y) = (x+y,x,y)\qandq(\realisation\pr{\g_3})(x,y,z) = \frac{x-y-z}{2}
\ifnocf. 
\end{equation}
\cfload[.]\Nobs that \cref{lemma:general:square} (applied with
		$N \curvearrowleft N$,
		$R \curvearrowleft 2R$
in the notation of \cref{lemma:general:square}) \proves\ that there exists $\g_2 \in \ANNSET$ which satisfies that
\begin{enumerate}[(I)]
    \item \label{it:general:square_continuous} it holds 
    that $\realisation(\g_2) \in C(\R,\R)$,

    \item\label{it:general:square_approx} it holds for all $x\in[-2\radius,2\radius]$ that
    $0\leq  \pr*{\realisation\pr*{\g_2}} (x) -x^2\leq R^2 4^{-N}$,

		\item\label{it:general:square_lipschitz} it holds for all 
    $x,y \in \R$ that 
    $\vass[\big]{\pr*{\realisation\pr*{\g_2}} (x) - \pr*{\realisation\pr*{\g_2}} (y)} \leq 4R\vass{x-y}$,

    \item\label{it:general:square_length} it holds 
    that
    $\lengthANN(\g_2) =N+\ceil{\log_2(2R)}+4
    $,

    \item\label{it:general:square_dims} it holds 
		for all $k\in\N_0\cap[0,\lengthANN(\g_2)]$ that
    \begin{equation} 
        \singledims_k(\g_2)=
\begin{cases}
1 &\colon k=0\\
2 &\colon k\in\N\cap(0,2]\\
4 &\colon k\in\N\cap(2,N+2]\\
2 &\colon k\in\N\cap(N+2, N+\ceil{\log_2(2R)}+4)\\
1 &\colon k=N+\ceil{\log_2(2R)}+4,\\
\end{cases}
\end{equation}
and

		
		\item\label{it:general:square_size} it holds 
    that
    $\size(\g_2) \leq 4
    $
\end{enumerate}
\cfload. Next let $\mathscr{f} \in \ANNSET$ satisfy
\begin{equation}
\label{prod_r:phi}
	\mathscr{f} = \compANN{\compANN{\g_3\bullet\ReLUidANN{3}}{(\parallelizationSpecial_{3}(\g_2,\g_2,\g_2))}}{\ReLUidANN{3}\bullet\g_1}
\end{equation}
\cfload. 
\Nobs that 
    \eqref{Lemma:prod_r:eq:varphi_new}, 
    \eqref{prod_r:phi},
    \cref{it:general:square_length},
		\cref{Lemma:PropertiesOfCompositions_n2},
		\cref{Lemma:PropertiesOfParallelizationEqualLength},
		and
		\cref{Prop:identity_representation}
\prove\ that 
\begin{equation}
\begin{split}
\label{prod_r:phi:length}
    \lengthANN(\mathscr{f}) 
		&=\lengthANN(\g_3)+\lengthANN(\ReLUidANN{3})+\lengthANN(\parallelizationSpecial_{3}(\g_2,\g_2,\g_2))+\lengthANN(\ReLUidANN{3})+\lengthANN(\g_1)-4
		\\&=1+2+\lengthANN(\g_2)+2+1-4
		\\&= N+\ceil{\log_2(2R)}+6
		= N+\ceil{\log_2(R)}+7
    .
		\end{split}
		\end{equation}
This, 
    \eqref{Lemma:prod_r:eq:varphi_new}, 
    \eqref{prod_r:phi},
    \cref{it:general:square_dims},
		\dimsbene,
		\cref{Lemma:PropertiesOfParallelizationEqualLength},
		and
		\cref{Prop:identity_representation}
\prove\ that for all $k\in\N_0\cap[0,\lengthANN(\f)]$ it holds that
\begin{equation}
\label{prod_r:phi:dims}
    \singledims_k(\f)=
\begin{cases}
2 &\colon k=0\\
6 &\colon k\in\N\cap(0,3]\\
12 &\colon  k\in\N\cap(3,N+3]\\
6 &\colon  k\in\N\cap(N+3,N+\ceil{\log_2(R)}+7)\\
1 &\colon k=N+\ceil{\log_2(R)}+7,
\end{cases}
\end{equation}
Next \nobs that 
\eqref{Lemma:prod_r:eq:varphi}, 
\eqref{prod_r:phi}, 
\cref{Lemma:PropertiesOfCompositions_n2},
	and
\cref{Prop:identity_representation} 
\prove\ that for all $x,y\in\R$ it holds that $\realisation(\mathscr{f}) \in C(\R^2,\R)$ and
\begin{equation}
\label{prod_r:phi:prep}
    \pr{\realisation\pr{\mathscr{f}}}(x,y)
    =
    \tfrac{1}{2}\PRb{\pr{\realisation\pr{\g_2}}(x+y)
    -\pr{\realisation\pr{\g_2}}(x)
    -\pr{\realisation\pr{\g_2}}(y)}.
\end{equation}
This, \cref{it:general:square_approx}, and the fact that $\fa{x,y}[-R,R]\colon x+y\in[-2R,2R]$ \prove\ that for all $x,y \in [-\radius, \radius]$ it holds that
\begin{equation}
\begin{split}
\label{prod_r:phi:eq2}
	&\vass{xy - \pr{\realisation\pr{\mathscr{f}}}(x,y)}
	\\&= \tfrac{1}{2}\abs[\big]{(x+y)^2  - \pr{\realisation\pr{\g_2}}(x+y) + \pr{\realisation\pr{\g_2}}(x)- x^2 + \pr{\realisation\pr{\g_2}}(y) - y^2} 
	\\&\leq\tfrac{1}{2}\max\pRb{ \pr{\realisation\pr{\g_2}}(x+y)-(x+y)^2  ,\pr{\realisation\pr{\g_2}}(x)- x^2 + \pr{\realisation\pr{\g_2}}(y) - y^2}
	\\&\leq\tfrac{1}{2}\max\pRb{ R^2 4^{-N} ,2R^2 4^{-N}}
	=R^2 4^{-N}
	.
\end{split}
\end{equation}
\Moreover 
    \eqref{prod_r:phi:prep} and 
    \cref{it:general:square_lipschitz} 
\prove\ that for all 
    $x_1,x_2,y_1,y_2\in \R$ 
it holds that
\begin{equation}
\begin{split}
\label{prod_r:phi:eq3}
	&\vass{\pr{\realisation\pr{\mathscr{f}}}(x_1,x_2) - \pr{\realisation\pr{\mathscr{f}}}(y_1,y_2)}
	\\&\leq
    \tfrac{1}{2}\prb{ \vass{\pr{\realisation\pr{\g_2}}(x_1+x_2) - \pr{\realisation\pr{\g_2}}(y_1+y_2)} \\
	&\quad+ \vass{\pr{\realisation\pr{\g_2}}(x_1) - \pr{\realisation\pr{\g_2}}(y_1)} + \vass{\pr{\realisation\pr{\g_2}}(x_2) - \pr{\realisation\pr{\g_2}}(y_2)} } \\
	&\leq
    2\radius \prb{ \vass{(x_1+x_2)-(y_1+y_2)} + \vass{x_1-y_1} + \vass{x_2-y_2} } \\
	&\leq
    4\radius \pr*{ \vass{x_1-y_1} + \vass{x_2-y_2} } 
    \leq
    \sqrt{32}\radius \norm{(x_1-y_1,x_2-y_2)}
\end{split}
\end{equation}
\cfload. 
\Moreover 
\eqref{Lemma:prod_r:eq:networks},
\eqref{prod_r:phi}, 
\cref{it:general:square_size}, 
\cref{Lemma:ParallelizationElementary}, 
and
\cref{Prop:identity_representation:prop2}
\prove\ that
\begin{equation}
\begin{split}
{}\size(\f)
{}= \max\pR*{\size(\g_3),\size(\parallelizationSpecial_{3}(\g_2,\g_2,\g_2)),\size(\g_1)}
&{}= \max\pR*{\size(\g_3),\size(\g_2),\size(\g_1)}
\\&{}\leq \max\pR*{\tfrac{1}{2},4,1}=4.
\end{split}
\end{equation}
    This, 
    \eqref{prod_r:phi:dims}, 
    \eqref{prod_r:phi:prep}, 
    \eqref{prod_r:phi:eq2}, and
    \eqref{prod_r:phi:eq3}
\prove[ie]\ 
    \cref{product_lipschitz,product_size,product_dims,product_approx,product_continuous,product_length}.
The proof of \cref{lemma:product} is thus complete\cfload.
\end{proof}

In the following result, \cref{lemma:2d_d_product} below, we adapt
the construction and proof of \benepairprod\ using the
multiplication networks in \cref{lemma:product}.
In particular, we establish the uniform bound on the absolute
values of the parameters in \cref{2d_d_product_size}
and record the resulting bounds on the network architecture.

\cfclear
\begin{athm}{lemma}{lemma:2d_d_product}
Let $d,N\in\N$, $R\in(1,\infty)$. Then there exists $\f\in\ANNSET$ such that
\begin{enumerate}[(i)]
    \item \label{2d_d_product_continuous} it holds 
    that $\realisation(\mathscr{f}) \in C(\R^{2d},\R^d)$,

    \item\label{2d_d_product_approx} it holds for all 
                $x = (x_1, \ldots, x_{2d}) \in [-\radius, \radius]^{2d}$ 
            that
            \begin{equation}
                \norm*{\pr*{x_1x_2,x_3x_4,\ldots,x_{2d-1}x_{2d}} - (\realisation(\f)) (x)} 
                \leq 
                R^2 d^{\nicefrac{1}{2}} 2^{-2N},
            \end{equation}

		\item\label{2d_d_product_lipschitz} it holds for all 
    $x,y \in \R^{2d}$ that 
    $\norm[\big]{\pr*{\realisation\pr*{\mathscr{f}}} (x) - \pr*{\realisation\pr*{\mathscr{f}}} (y)} \leq \sqrt{32}\radius \norm{x-y}$,

    \item\label{2d_d_product_length} it holds 
    that
    $\lengthANN(\mathscr{f}) = N+\ceil{\log_2(R)}+7
    $,

    \item\label{2d_d_product_dims} it holds 
		for all $k\in\N_0\cap[0,\lengthANN(\f)]$ that
    \begin{equation} 
        \singledims_k(\mathscr{f})=
\begin{cases}
	2d &\colon k=0\\
	6d &\colon k\in\N\cap(0,3]\\
	12d &\colon  k\in\N\cap(3,N+3]\\
	6d &\colon  k\in\N\cap(N+3,N+\ceil{\log_2(R)}+7)\\
	d &\colon k=N+\ceil{\log_2(R)}+7,
\end{cases}
\end{equation}

    \item\label{2d_d_product_cost} it holds 
    that
    $\paramANN(\mathscr{f}) = 162d^2+37d+N(144d^2+12d)+\ceil{\log_2(R)}(36d^2+6d)
    $, and
		
		\item\label{2d_d_product_size} it holds 
    that
    $\size(\mathscr{f}) \leq 4
    $
\end{enumerate}
\cfout.
\end{athm}

\begin{proof}[Proof of \cref{lemma:2d_d_product}]
    \Nobs that 
        \cref{lemma:product}
				(applied with
				$N \curvearrowleft N$,
		$\radius \curvearrowleft \radius$
in the notation of \cref{lemma:product})
    \proves\ that there exists
        $\g\in \ANNSET$ 
    which satisfies that
    \begin{enumerate}[(I)]
    \item \llabel{it:single_realization} 
        it holds that 
            $\realisation({\g}) \in C(\R^2,\R)$,
    \item \llabel{it:single_lipschitz} 
        it holds for all 
            $x,y\in\R^{2}$ 
        that 
        $
            \vass{(\realisation\pr{\g})(x) - (\realisation\pr{\g})(y)}
            \leq 
            \sqrt{32} R \norm{x-y}
        $,
    \item \llabel{it:single_approx} 
        it holds that 
        $
            \sup_{x,y \in [-R, R] } \vass{xy - (\realisation\pr{\g})(x,y)} 
            \leq 
            R^2 2^{-2N}
        $,
    \item \llabel{it:single_length} 
        it holds that 
        $
            \lengthANN(\g) 
            = 
            N+\ceil{\log_2(R)}+7
        $,   
    \item \llabel{it:single_dims} 
        it holds for all $k\in\N_0\cap[0,\lengthANN(\g)]$ that
    \begin{equation} 
        \singledims_k(\mathscr{g})=
\begin{cases}
	2 &\colon k=0\\
	6 &\colon k\in\N\cap(0,3]\\
	12 &\colon  k\in\N\cap(3,N+3]\\
	6 &\colon  k\in\N\cap(N+3,N+\ceil{\log_2(R)}+7)\\
	1 &\colon k=N+\ceil{\log_2(R)}+7,
\end{cases}
\end{equation}
and
    \item \llabel{it:single_size} 
        it holds that
        $
            \size(\g) 
            \leq 
           4
        $
    \end{enumerate}
    \cfload. 
    Next let 
        $\mathscr{f} \in \ANNSET$ 
    satisfy 
    \begin{equation}
    \begin{split}
    \label{eq:fdef}
        \mathscr{f}
        &=
        \parallelizationSpecial_{d} \pr*{\g,\g,\ldots,\g}
        \ifnocf.
    \end{split}
    \end{equation}
		\cfload[.]%
    \Nobs that 
        \eqref{eq:fdef}, 
        \lref{it:single_length}, 
        \lref{it:single_dims}, and 
        \cref{Lemma:PropertiesOfParallelizationEqualLength} 
    \prove\ that for all $k\in\N_0\cap[0,\lengthANN(\f)]$ it holds that
\begin{equation}
\label{2d_d_product:dims}
    \lengthANN(\mathscr{f}) = N+\ceil{\log_2(R)}+7
    \,\,\text{and}\,\,
    \singledims_k(\f)=
\begin{cases}
	2d &\colon k=0\\
6d &\colon k\in\N\cap(0,3]\\
12d &\colon  k\in\N\cap(3,N+3]\\
6d &\colon  k\in\N\cap(N+3,N+\ceil{\log_2(R)}+7)\\
d &\colon k=N+\ceil{\log_2(R)}+7.
\end{cases}
\end{equation}
This and \eqref{eq:def:ANN:operators} \prove\ that
    \begin{equation}
    \begin{split}
    \label{eq:params}
        &\paramANN(\mathscr{f})
				\\&=\sum_{k=1}^{\lengthANN(\f)}\singledims_k(\f)\pr*{\singledims_{k-1}(\f)+1}
        \\&=
				\singledims_1(\f)\pr*{\singledims_{0}(\f)+1}
				+\singledims_2(\f)\pr*{\singledims_{1}(\f)+1}
				+\singledims_3(\f)\pr*{\singledims_{2}(\f)+1}
				+\singledims_4(\f)\pr*{\singledims_{3}(\f)+1}
				\\&\quad
				+\PR*{\sum_{k=5}^{N+3}\singledims_k(\f)\pr*{\singledims_{k-1}(\f)+1}}
				+\singledims_{N+4}(\f)\pr*{\singledims_{N+3}(\f)+1}
				\\&\quad+\PR*{\sum_{k=N+5}^{\lengthANN(\f)-1}\singledims_k(\f)\pr*{\singledims_{k-1}(\f)+1}}
				+\singledims_{\lengthANN(\f)}(\f)\pr*{\singledims_{\lengthANN(\f)-1}(\f)+1}
				\\&=6d\pr*{2d+1}+2\pr*{6d\pr*{6d+1}}+12d\pr*{6d+1}
				+\PR*{\sum_{k=1}^{N-1}12d\pr*{12d+1}}+6d\pr*{12d+1}
				\\&\quad
				+\PR*{\sum_{k=1}^{\lengthANN(\f)-N-5}6d\pr*{6d+1}}
				+d\pr*{6d+1}
				\\&=(12+72+72+72+6)d^2+(6+12+12+6+1)d+(N-1)(144d^2+12d)
				\\&\quad+(\lengthANN(\f)-N-5)(36d^2+6d)
				\\&=234d^2+37d+(N-1)(144d^2+12d)+\pr*{\ceil{\log_2(R)}+2}(36d^2+6d)
				\\&=162d^2+37d+N(144d^2+12d)+\ceil{\log_2(R)}(36d^2+6d)
        .
    \end{split}
    \end{equation}
    \Moreover 
        \eqref{eq:fdef},
        \lref{it:single_approx}, and
        \cref{Lemma:PropertiesOfParallelizationEqualLength}
    \prove\ that for all 
        $x=(x_1,\ldots,\allowbreak x_{2d})\in[-R,R]^{2d}$ 
    it holds that
        $\realisation(\mathscr{f}) \in C(\R^{2d},\R^d)$
        and
    \begin{equation}
    \begin{split}
    \label{eq:approx}
        &
        \norm*{\pr*{x_1x_2,x_3x_4,\ldots,\allowbreak x_{2d-1}x_{2d}} - (\realisation(\f)) (x)}
        \\&=
        \PR*{\ssum_{i=1}^d\vass*{x_{2i-1}x_{2i} - (\realisation(\g)) (x_{2i-1},x_{2i})}^2}^{\nicefrac{1}{2}}
        \leq
        \PR*{\ssum_{i=1}^d R^4 2^{-4N}}^{\nicefrac{1}{2}}
        = 
        R^2 d^{\nicefrac{1}{2}} 2^{-2N}
        .
    \end{split}
    \end{equation}
		\Moreover \lref{it:single_size} and \cref{Lemma:ParallelizationElementary} \prove\ that
		\begin{equation}
		\label{eq:2d_d_product:size}
		\size(\f)=\size(\g)\leq 4.
		\end{equation}
    Combining
     \eqref{eq:fdef}, 
        \lref{it:single_lipschitz}, and 
        \cref{Lemma:PropertiesOfParallelizationEqualLength}
    with \benecorlip\ (applied with
                $L\is \sqrt{32}R$,
                $d\is d$,
                $(g_1,g_2,\dots,g_d)\is(\realisation({\g}),\realisation({\g}),\dots,\realisation({\g}))$,
                $f\is \realisation(\mathscr f)$
            in the notation of \benecorlip)
    \proves\ that for all 
        $x,y\in\R^{2d}$ 
    it holds that
    \begin{equation}
        \norm{(\realisation(\f))(x) - (\realisation(\f))(y)}
        \leq 
        \sqrt{32} R \norm{x-y}
        .
    \end{equation}
        This,
				\eqref{2d_d_product:dims},
				\eqref{eq:params},
        \eqref{eq:approx}, and
        \eqref{eq:2d_d_product:size}      
    \prove[ie]\
        \cref{2d_d_product_continuous,2d_d_product_approx,2d_d_product_lipschitz,2d_d_product_length,2d_d_product_dims,2d_d_product_size,2d_d_product_cost}. 
The proof of \cref{lemma:2d_d_product} is thus complete.
\end{proof}

In the following result, \cref{lemma:2d_product} below, we adapt \benepotprod\ using the
multiplication networks constructed in
\cref{lemma:2d_d_product}.
We recall the composition argument of the original proof
and additionally establish the parameter-size, length,
and layer-dimension bounds in
\cref{2d_product_size,2d_product_length,2d_product_dims} in \cref{lemma:2d_product}.

\cfclear
\begin{athm}{lemma}{lemma:2d_product}
Let $d\in\N$, $\eps\in(0,1)$, $R\in(1,\infty)$. Then there exists $\f\in\ANNSET$ such that
\begin{enumerate}[(i)]
    \item \label{2d_product_continuous} it holds 
    that $\realisation(\mathscr{f}) \in C(\R^{(2^d)},\R)$,

    \item\label{2d_product_approx} it holds for all 
                $x = (x_1, \ldots, x_{2^d}) \in [-\radius, \radius]^{(2^d)}$ 
            that
            $
                \vass[\big]{\sprod_{i=1}^{2^d}x_i - (\realisation(\f)) (x)} 
                \leq 
                \eps,
           $

		\item\label{2d_product_lipschitz} it holds for all 
    $x,y \in \R^{(2^d)}$ that 
    $\vass[\big]{\pr*{\realisation\pr*{\mathscr{f}}} (x) - \pr*{\realisation\pr*{\mathscr{f}}} (y)} \leq 2^{\frac{5d}{2}}\! R^{(2^{d}-1)} \norm{x-y}$,

    \item\label{2d_product_length} it holds that
    $\lengthANN(\mathscr{f}) \leq d2^{d+2}+d2^{d}\ceil{\log_2(R)}-\tfrac{d\log_2(\eps)}{2}
    $,
		
		\item\label{2d_product_dims} it holds that
				$\singledims_1(\f)=2^d\,3
		$ and
				$\singledims_{\hidlengthANN(\mathscr{f})}(\f)=6
		$,

    \item\label{2d_product_cost} it holds 
    that
    $\paramANN(\mathscr{f}) \leq 2^{3d+9}+2^{3d+8}\ceil{\log_2(R)}-2^{2d+7}\log_2(\eps)
    $, and
		
		\item\label{2d_product_size} it holds 
    that
    $\size(\mathscr{f}) \leq 4
    $
\end{enumerate}
\cfout.
\end{athm}

\begin{proof}[Proof of \cref{lemma:2d_product}]
	Throughout this proof let $c=\ceil{\log_2(R)}$ and
	$t=\log_2(\eps^{-1})$,
	for every $i\in\{1,2,\ldots,d\}$ let $q_i,n_i,N_i,\ell_i\in\N$, $r_i\in(1,\infty)$ satisfy
	\begin{equation}
		\label{eq:2d_product:stage:parameters}
		\begin{split}
			q_i=2^{d-i},
			\qquad 
			r_i=R^{2^{i-1}},
			\qquad
			D_i=[-r_i,r_i]^{2q_i},
			\qquad
			n_i=\ceil{2^{i-1}\log_2(R)},\\
			N_i=\ceil{\tfrac{3(d-i)}{2}+\tfrac{\log_2(d)}{2}
				+2^{d-1}\log_2(R)+\tfrac{t}{2}},
			\qandq \ell_i=N_i+n_i+7,
		\end{split}
	\end{equation}
	and for every $i\in\{1,2,\ldots,d\}$ 
				let $p_i\colon D_i\to\R^{q_i}$
        satisfy for all 
            $x=(x_1,x_2,\dots,x_{2q_i})\in D_i$ 
        that 
        \begin{equation}
				\label{eq:def:partial:products}
            p_i(x)
            =
            (x_1x_2,x_3x_4,\ldots,x_{2q_i-1}x_{2q_i}).
        \end{equation}
				\Nobs that the fact that for all $i\in\{1,2,\ldots,d\}$ it holds that $2N_i\geq 3(d-i)+\log_2(d)
				+2^{d}\log_2(R)+t$ \proves\ that for all $i\in\{1,2,\ldots,d\}$ it holds that
				\begin{equation}
					\label{approx:prep:N_i:prop}
					\begin{split}
					2^{3(d-i)-2N_i}R^{(2^{d})}
					\leq 2^{3(d-i)-3(d-i)-\log_2(d)-2^d\log_2(R)-t}R^{(2^{d})}
					= d^{-1}\eps
					.
				\end{split}
				\end{equation}
    \Nobs that
        \cref{lemma:2d_d_product}
			(applied for every $i\in\{1,2,\ldots,d\}$ with
                $d \curvearrowleft q_i$,
								$N \curvearrowleft N_i$,
                $R \curvearrowleft r_i$
            in the notation of \cref{lemma:2d_d_product}) 				
    \proves\ that for every $i\in\{1,2,\dots,d\}$ there exists
        $\h_i\in\ANNSET$ 
    which satisfies that
    \begin{enumerate}[(I)]
    \item \label{it:2d_d_product_continuous} it holds 
    that $\realisation(\h_i) \in C\prb{\R^{2q_i},\R^{q_i}}$,

    \item\label{it:2d_d_product_approx} it holds that
            $
                \sup_{x \in D_i}\norm*{p_i(x) - (\realisation(\h_i)) (x)} 
                \leq
                (r_i)^2\sqrt{q_i}\,2^{-2N_i}
            $,

		\item\label{it:2d_d_product_lipschitz} it holds for all 
    $x,y \in \R^{2q_i}$ that 
    $\norm[\big]{\pr*{\realisation\pr*{\h_i}} (x) - \pr*{\realisation\pr*{\h_i}} (y)} \leq \sqrt{32}\,r_i \norm{x-y}$,

    \item\label{it:2d_d_product_length} it holds 
    that
    $\lengthANN(\h_i) = \ell_i
    $,

    \item\label{it:2d_d_product_dims} it holds 
		for all $k\in\N_0\cap[0,\lengthANN(\h_i)]$ that
    \begin{equation} 
        \singledims_k(\h_i)=
\begin{cases}
2q_i &\colon k=0\\
6q_i  &\colon k\in\N\cap(0,3]\\
12q_i &\colon k\in\N\cap(3, N_i +3]\\
6q_i  &\colon k\in\N\cap(N_i +3, \ell_i)\\
q_i  &\colon k=\ell_i,\\
\end{cases}
\end{equation}

    \item\label{it:2d_d_product_cost} it holds 
    that
    $\paramANN(\mathscr{\h_i}) 
    =162(q_i)^2+37q_i+N_i(144(q_i)^2+12q_i)+n_i(36(q_i)^2+6q_i)
    $, and
		
		\item\label{it:2d_d_product_size} it holds 
    that
    $\size(\mathscr{\h_i}) \leq 4
    $
\end{enumerate}
    \cfload[.]%
    Next let 
        $\f\in\ANNSET$ and $\g_n\in\ANNSET$, $n\in\{1,2,\dots,d\}$,
    satisfy that
        \begin{equation}
    	\label{fdef:2d_product}
    	\g_1=\h_1,
    	\qquad
    	\prb{\forall\,k\in\{1,2,\dots,d-1\}\colon \g_{k+1}=\mathscr{h}_{k+1}\bullet\ReLUidANN{2^{d-k}}\bullet\g_{k}},
    	\qandq
    	\f
    	=
    	\g_d
    	\ifnocf.
    \end{equation}
    \cfload[.]%
    \Nobs that 
        \eqref{fdef:2d_product}, 
        \cref{it:2d_d_product_continuous},
				\cref{Lemma:PropertiesOfCompositions_n2}, and
        \cref{Prop:identity_representation}
    \prove\ that
    \begin{equation}
    \label{eq:d:dim:product:functioncheck}
        \realisation(\f)
        =
        \PR*{\realisation(\mathscr{h}_d)}\circ\PR*{\realisation(\mathscr{h}_{d-1})}\circ\ldots\circ\PR*{\realisation(\mathscr{h}_1)}
        \in C\prb{\R^{(2^d)},\R}
        \ifnocf.
    \end{equation}
    \cfload[.]%
    Combining this and
        \cref{it:2d_d_product_lipschitz} with
        induction 
    \proves\ that for all 
        $x,y\in\R^{(2^d)}$ 
    it holds that
    \begin{equation}
    \label{eq:d:dim:product:lipschitz}
    \begin{split}
        \vass[\big]{
            (\realisation(\f))(x)
            -
            (\realisation(\f))(y)
        }
        \leq 
        \textstyle\prb{\prod_{i=1}^d\sqrt{32}r_i}\norm{x-y}
        &=
        \prb{\sqrt{32}}^{d} R^{\PR{\sum_{i=1}^d2^{i-1}}}\norm{x-y}
        \\&= 
        2^{\frac{5d}{2}}\! R^{(2^{d}-1)}\norm{x-y}
        .
    \end{split}
    \end{equation}
    Next \nobs that
        the fact that
            for all
                $i\in\{1,2,\dots,d\}$,
                $x,y\in[-R^{(2^{i-1})},R^{(2^{i-1})}]$
            it holds that
                $xy\in[-R^{(2^{i})},R^{(2^{i})}]$
    \proves\ that for all
        $i\in\{1,2,\dots,d-1\}$
    it holds that
        $p_i(D_i)\subseteq D_{i+1}$.
    Combining
        this,
				\eqref{eq:def:partial:products},
				\eqref{approx:prep:N_i:prop},
				\eqref{fdef:2d_product},
        \eqref{eq:d:dim:product:functioncheck},
        \cref{it:2d_d_product_approx}, and
        \cref{it:2d_d_product_lipschitz}
    with
    \benecorcomp\
            (applied with
                $n \is d$,
                $(d_0,\allowbreak d_1,\dots,d_n)\is(2^d,\allowbreak 2^{d-1},\dots,2^0)$,
                $(L_i)_{i\in\{1,2,\ldots,n\}} \is \prb{\sqrt{32}\, r_i}_{i\in\{1,2,\ldots,d\}}$,
                $(\varepsilon_i)_{i\in\{1,2,\ldots,n\}} \is \prb{(r_i)^2\sqrt{q_i}\,2^{-2N_i}}_{i\in\{1,2,\ldots,d\}}$,
                $(D_1,\allowbreak D_2,\dots,D_n)\is(D_1,\allowbreak D_2,\dots\allowbreak,D_d)$,
                $(f_1,\allowbreak f_2,\dots,f_n)\is(p_1,\allowbreak p_2,\dots,p_d)$,
                $(g_1,\allowbreak g_2,\dots,g_n)\is \pr{\realisation(\mathscr{h}_1),\allowbreak \realisation(\mathscr{h}_2),\dots,\realisation(\mathscr{h}_d)}$
            in the notation of \benecorcomp)
    \proves\ that for all
        $x=(x_1,\dots,x_{2^d})\in[-R,R]^{(2^d)}$
    it holds that
    \begin{equation}
    \begin{split}
    \label{eq:d:dim:product:approxcheck}
        &\abs[\bigg]{\PRbbb{\sprod_{i=1}^{{2^d}} x_i} - (\realisation(\f))(x)}
        \\&=
        \abs[\big]{
            \pr*{p_d\circ p_{d-1}\circ \ldots\circ p_1}(x)
            -
            \prb{\PR*{\realisation(\mathscr{h}_d)}\circ\PR*{\realisation(\mathscr{h}_{d-1})}\circ\ldots\circ\PR*{\realisation(\mathscr{h}_1)}}(x)
        }
        \\&\leq
        \ssum_{i=1}^d\PRbbb{\prbbb{\sprod_{j=i+1}^d \sqrt{32} r_j} (r_i)^2\sqrt{q_i}\,2^{-2N_i}}
        \\&=
        \PRbbb{\ssum_{i=1}^{d} 2^{\frac{5d-5i}{2}} R^{\PR{\sum_{j=i+1}^{d}2^{j-1}}}R^{(2^{i})} 2^{\frac{d-i}{2}} 2^{-2N_i}}
        \\&=
        \PRbbb{\ssum_{i=1}^{d} 2^{3(d-i)} R^{(2^{d}-2^i)}R^{(2^{i})} 2^{-2N_i}}
        =
        \PRbbb{\ssum_{i=1}^{d} 2^{3(d-i)-2N_i}R^{(2^{d})} } 
        \leq \PRbbb{\ssum_{i=1}^{d} d^{-1}\eps}=
        \eps
        .
    \end{split}
    \end{equation}\cfload
		\Moreover \eqref{fdef:2d_product}, \cref{it:2d_d_product_dims}, and \dimsbene\ \prove\ that 
		\begin{equation}
		\label{eq:d_product:hidden:dims}
		\singledims_{\lengthANN(\f)-1}(\f)=\singledims_{\lengthANN(\h_d)-1}(\h_d)=6q_d=6,
		\qandq
		\singledims_1(\f)=\singledims_1(\h_1)=6q_1=2^d\,3.
		\end{equation}\cfload
		\Moreover \eqref{fdef:2d_product}, \cref{it:2d_d_product_size}, and \cref{Prop:identity_representation:prop2} \prove\ that 
		\begin{equation}
		\label{eq:d_product:size}
		\size(\f)\leq\max\{\size(\h_1),\size(\h_2),\ldots,\size(\h_d)\}\leq 4.
		\end{equation}\cfload
    \Moreover 
    the fact that $ \log_2(d)\leq d -1$ and
        the assumption that $d\in\N$, $\eps\in(0,1)$, $R\in[1,\infty)$ \prove\ that
            for all
                $i\in\{1,2,\dots,d\}$ 
            it holds that
    \begin{equation}
		\label{eq:est:of:N_is}
    \begin{split}
N_i=\ceil{\tfrac{3(d-i)}{2}+\tfrac{\log_2(d)}{2}
	+2^{d-1}\log_2(R)+\tfrac{t}{2}}
	&\leq\tfrac{4(d-1)}{2}
		+2^{d-1}c+\tfrac{t}{2}+1
	\\&=2d+2^{d-1}c+\tfrac{t}{2}-1
\\&\leq2^d+2^{d-1}c+\tfrac{t}{2}-1
.
    \end{split}
    \end{equation}\cfload
    Thus, the fact that for all $i\in\{1,2,\dots,d\}$ it holds that $n_i\leq2^{i-1}c\leq2^{d-1}c$ and $q_i\geq 1$, and \cref{it:2d_d_product_cost} \prove\ that for all
        $i\in\{1,2,\dots,d\}$ 
    it holds that
    \begin{equation}
    	\begin{split}
    		\llabel{param_h_i}
    		\paramANN(\h_i)
    		&=162(q_i)^2+37q_i+N_i(144(q_i)^2+12q_i)+n_i(36(q_i)^2+6q_i)
    		\\&\leq (q_i)^2(199+156N_i+42n_i)
    		\\&\leq (q_i)^2(199+2^d 156+2^d 78c+78t+2^d 21c)
    		\leq (q_i)^2(2^{d+9}+2^{d+8}c+2^7 t)
.
    	\end{split}
    \end{equation}\cfload
		Combining this and \cite[Proposition~2.19]{beneventano21}
		with the fact that $\sum_{i=0}^{d-1}4^{i}= \frac{4^d-1}{3}$ \proves\ that
		\begin{equation}
    \begin{split}
    \label{eq:d:dim:product:paramcheck}
        \paramANN(\f)
				\leq
					3\PR*{ \ssum_{i=1}^d\paramANN(\h_i) }-\paramANN(\h_1)-\paramANN(\h_d)
					&\leq
					3\PR*{ \ssum_{i=1}^d (q_i)^2(2^{d+9}+2^{d+8}c+2^7 t) }
					\\&=
					3(2^{d+9}+2^{d+8}c+2^7 t)\PR*{ \ssum_{i=1}^d 4^{d-i} }
					\\&=
					(2^{d+9}+2^{d+8}c+2^7 t)\pr{4^d-1}
					\\&\leq
					2^{3d+9}+2^{3d+8}c+2^{2d+7} t
        .
    \end{split}
    \end{equation}
		\Moreover  
		\eqref{fdef:2d_product}, 
		\eqref{eq:est:of:N_is}, 
		\cref{it:2d_d_product_length},
		\cref{Prop:identity_representation}, and \cref{Lemma:PropertiesOfCompositions_n2} \prove\ that
		\begin{equation}
		\begin{split}
		\lengthANN(\f)&=\hidlengthANN(\f)+1
		\\&=\PR*{\ssum_{i=1}^d\hidlengthANN(\h_i)}
		+\PR*{\ssum_{i=1}^{d-1}\hidlengthANN(\ReLUidANN{2^{i}})}+1
		\\&=\PR*{\ssum_{i=1}^d\pr{\ell_i-1}}
		+\PR*{\ssum_{i=1}^{d-1}1}+1
		\\&=\ssum_{i=1}^d \pr{N_i+n_i+7}
		\\&\leq\ssum_{i=1}^d \pr{2d+2^{d-1}c+\tfrac{t}{2}+2^{d-1}c+6}
		=d\pr{2d+2^{d}c+\tfrac{t}{2}+6}
		\leq d2^{d+2}+ d 2^d c + \tfrac{td}{2}
		.
		\end{split}
		\end{equation}
    Combining 
        this 
    with 
        \eqref{eq:d:dim:product:functioncheck}, 
        \eqref{eq:d:dim:product:lipschitz},
        \eqref{eq:d:dim:product:approxcheck},
				\eqref{eq:d_product:hidden:dims},
				\eqref{eq:d_product:size},				
				and
				\eqref{eq:d:dim:product:paramcheck}
    \proves[ie]\ 
        \cref{2d_product_continuous,2d_product_approx,2d_product_lipschitz,2d_product_cost,2d_product_length,2d_product_size,2d_product_dims}\cfload. 
    The proof of \cref{lemma:2d_product} is thus complete.
    \end{proof}

In the following result, \cref{lemma:d_product} below, we adapt
the construction and proof of \beneprod\ to the scaled
product functions considered here.
We also establish the bounds on the absolute values of the
parameters, the length, and the layer dimensions of the approximating \ANN\ in
\cref{d_product_size,d_product_length,d_product_dims}.

\newcommand{\D}{D}

\cfclear
\begin{athm}{lemma}{lemma:d_product}
Let $d\in\N$, $\eps\in(0,1)$, $R\in(1,\infty)$, $\gamma\in(0,1]$, $\scl\in[1,\infty)$. Then there exists $\f\in\ANNSET$ such that
\begin{enumerate}[(i)]
    \item \label{d_product_continuous} it holds 
    that $\realisation(\mathscr{f}) \in C(\R^{d},\R)$,

    \item\label{d_product_approx} it holds that
            $
                \sup_{x = (x_1,\ldots, x_{d}) \in [-\radius, \radius]^{d}}\vass[\big]{\gamma\scl^d\sprod_{i=1}^{d}x_i - (\realisation(\f)) (x)} 
                \leq 
                \eps
            $,

		\item\label{d_product_lipschitz} it holds for all 
    $x,y \in \R^{d}$ that 
    $\vass[\big]{\pr*{\realisation\pr*{\mathscr{f}}} (x) - \pr*{\realisation\pr*{\mathscr{f}}} (y)} \leq \sqrt{32}d^{\frac{5}{2}}\beta^d R^{2d-1} \norm{x-y}$,

    \item\label{d_product_length} it holds 
    that
    $\lengthANN(\mathscr{f}) \leq 8d^2+2d^2\ceil{\log_2(R)}+d\log_2(\eps^{-1})+d^2\ceil{\log_2(\scl)}+2
    $,
		
		\item\label{d_product_dims} it holds that
				$\singledims_1(\f)\leq 4d
		$ and
				$\singledims_{\hidlengthANN(\f)}(\f)=2
		$,

    \item\label{d_product_cost} it holds 
    that
    $\paramANN(\mathscr{f}) \leq 4107 d^3+2048 d^3\ceil{\log_2(R)}-512 d^2 \log_2(\eps)+514d^3\log_2(\scl)
    $, and
		
		\item\label{d_product_size} it holds 
    that
    $\size(\mathscr{f}) \leq 4
    $
\end{enumerate}
\cfout.
\end{athm}

\begin{proof}[Proof of \cref{lemma:d_product}]
Throughout this proof 
assume without loss of generality that 
    $d>1$ \cfadd{gen:scaling:networks}\cfload,
let 
    $b,\D\in\N$ 
satisfy 
$b=\max\{1,\ceil{\log_2(\beta)}\}$ and
    $\D=2^{\ceil{\log_2(d)}}$, 
let 
    $A\in\R^{\D\times d}$, $B\in\R^{\D}$
satisfy for all 
    $x=(x_1,\ldots,x_d)\in\R^d$ 
that
\begin{equation}
\label{affinedef}
    Ax+B=(\gamma x_1,x_2,\dots,x_d,1,1,\dots,1),
\end{equation}
 and let $\g_1\in\pr{\R^{D\times d}\times\R^D}\subseteq\ANNSET$ satisfy $\g_1=(A,B)$ \cfload.
\Nobs that \cref{lemma:2d_product} (applied with 
$d \curvearrowleft \ceil{\log_2(d)}$,
$R \curvearrowleft R$,
$\eps \curvearrowleft \eps\beta^{-d}$
 in the notation of \cref{lemma:2d_product}) \proves\ that there exists $\g_2\in\ANNSET$ which satisfies that
\begin{enumerate}[(I)]
    \item \label{it:2d_product_continuous} it holds 
    that $\realisation(\g_2) \in C(\R^{D},\R)$,

    \item\label{it:2d_product_approx} it holds that $
                \sup_{x = (x_1, \dots, x_{D}) \in [-\radius, \radius]^{D}}\vass[\big]{\sprod_{i=1}^{D}x_i - (\realisation(\g_2)) (x)} 
                \leq 
                \eps\beta^{-d}$,

		\item\label{it:2d_product_lipschitz} it holds for all 
    $x,y \in \R^{D}$ that 
    $\vass[\big]{\pr*{\realisation\pr*{\g_2}} (x) - \pr*{\realisation\pr*{\g_2}} (y)} \leq D^{\frac{5}{2}}\! R^{D-1} \norm{x-y}$,

    \item\label{it:2d_product_length} it holds 
    that
    $\lengthANN(\g_2) \leq 4\ceil{\log_2(d)}D+\ceil{\log_2(d)}D\ceil{\log_2(R)}+\tfrac{\ceil{\log_2(d)}(d\log_2(\scl)- \log_2(\eps))}{2}
    $,

		\item\label{it:2d_product_dims} it holds that
				$\singledims_1(\g_2)=3D
		$ and
				$\singledims_{\hidlengthANN(\g_2)}(\g_2)=6
		$,

    \item\label{it:2d_product_cost} it holds 
    that
    $\paramANN(\g_2) \leq 2^{9} D^3+2^8 D^3\ceil{\log_2(R)}+2^7 D^2(d\log_2(\scl)- \log_2(\eps))
    $, and
		
		\item\label{it:2d_product_size} it holds 
    that
    $\size(\g_2) \leq 4
    $
\end{enumerate}
\cfload. 
\Nobs that \cref{gen:scaling:networks} (applied with
$\beta \curvearrowleft  \beta^d$,
$L \curvearrowleft  d\ceil{\log_2(\beta)}$
in the notation of \cref{gen:scaling:networks}) \proves\ that there exists $\g_3\in\ANNSET$ which satisfies that
\begin{enumerate}[(A)]
	\item{
		\label{domainup:scaling:realisation}
		it holds for all $x\in\R$ that $(\realisation(\g_3))(x)=\beta^d x$,
	}
	\item{
		\label{domainup:scaling:dims}
		it holds that $\dims(\g_3)=(1,2,2,\ldots,2,1)\in\N^{d\ceil{\log_2(\beta)}+2}$,
	}
	\item{
		\label{domainup:scaling:size}
		it holds that $\insize(\g_3)\leq 1$, $\outsize(\g_3)\leq 2$, and $\size(\g_3)\leq 2$, and
	}
	\item{
		\label{domainup:scaling:param}
		it holds that $\param(\g_3)\leq 6db+1$
	}
\end{enumerate}
\Nobs that 
    \cref{Lemma:PropertiesOfCompositions_n2}, 
    \cref{it:2d_product_continuous}, 
    \cref{domainup:scaling:dims}, 
    and the fact that 
        $\realisation(\g_1)\in C(\R^d,\R^\D)$ 
\prove\ that 
\begin{equation}
\label{functioncheck}
\realisation(\g_3\bullet\ReLUidANN{1}\bullet\g_2\bullet\compANN{\ReLUidANN{D}}{\g_1})=\PR{\realisation(\mathscr{g}_3)}\circ\PR{\realisation(\mathscr{g}_2)}\circ\PR{\realisation(\g_1)}\in C(\R^{d},\R)\ifnocf.
\end{equation}
\cfload[.]%
Combining this,
    \cref{it:2d_product_lipschitz},
    \cref{domainup:scaling:realisation}, 
    \eqref{affinedef} and
    the fact that 
        $D\leq 2d$
	with the assumption that $R> 1$
\proves\ that for all 
    $x,y\in\R^d$ 
it holds that
\begin{equation}
\label{lipschitzcheck}
\begin{split}
    &\vass[\big]{(\realisation(\g_3\bullet\ReLUidANN{1}\bullet\g_2\bullet\compANN{\ReLUidANN{D}}{\g_1}))(x)-(\realisation(\g_3\bullet\ReLUidANN{1}\bullet\g_2\bullet\compANN{\ReLUidANN{D}}{\g_1}))(y)}
    \\&=
    \vass*{\prb{\PR{\realisation(\mathscr{g}_3)}\circ\PR{\realisation(\mathscr{g}_2)}\circ\PR{\realisation(\g_1)}}(x)-\prb{\PR{\realisation(\mathscr{g}_3)}\circ\PR{\realisation(\mathscr{g}_2)}\circ\PR{\realisation(\g_1)}}(y)}
    \\&\leq 
    \scl^d D^{\frac{5}{2}}\! R^{\D-1}\norm{\pr{\realisation(\g_1)}(x)-\pr{\realisation(\g_1)}(y)}
    \\&=
    \scl^d D^{\frac{5}{2}}\! R^{\D-1}\norm{Ax+B-(Ay+B)}
    \\&\leq
    D^{\frac{5}{2}}\scl^d R^{\D-1}\norm{x-y}
    \leq 
    \sqrt{32}d^{\frac{5}{2}}\scl^d R^{2d-1}\norm{x-y}
    .
\end{split}
\end{equation}
\Moreover
    \eqref{affinedef} and 
    the assumption that 
        $R> 1\geq \gamma$
\prove\ that for all $x\in[-R,R]^d$ it holds that $Ax+B\in[-R,R]^D$.
This,
    \cref{it:2d_product_approx}, 
    \cref{domainup:scaling:realisation},
    \eqref{affinedef}, and
    \eqref{functioncheck}
\prove\ that for all 
    $x=(x_1,\dots,x_d)\in[-R,R]^d$ 
it holds that
\begin{equation}
\label{approxcheck}
\begin{split}
    &\vass*{\PRbbb{\gamma\scl^d\sprod_{i=1}^{d} x_i} - (\realisation(\g_3\bullet\ReLUidANN{1}\bullet\g_2\bullet\compANN{\ReLUidANN{D}}{\g_1}))(x)}
    \\&=
    \vass*{\PRbbb{\gamma\scl^d\sprod_{i=1}^{d} x_i}- \prb{\PR{\realisation(\mathscr{g}_3)}\circ\PR{\realisation(\mathscr{g}_2)}\circ\PR{\realisation(\g_1)}}(x)}
    \\&=\scl^d
    \vass*{\PRbbb{1^{D-d}\gamma\sprod_{i=1}^{d} x_i} - \pr{\realisation(\mathscr{g}_2)}(\gamma x_1,x_2,\dots,x_d,1,1,\dots,1)}\leq\scl^d\varepsilon\scl^{-d}=\eps
    .
\end{split}
\end{equation}
\Moreover \eqref{affinedef}, 
\cref{it:2d_product_size}, 
\cref{domainup:scaling:size}, and \cref{Prop:identity_representation:prop} \prove\ that
\begin{equation}
\label{eq:size:d_product}
\begin{split}
\size(\g_3\bullet\ReLUidANN{1}\bullet\g_2\bullet\compANN{\ReLUidANN{D}}{\g_1})&= \max\{\size(\g_3),\size(\g_2),\size(\g_1)\}\leq \max\{2,4,1\}=4.
\end{split}
\end{equation}
\Moreover 
\eqref{affinedef},
\cref{it:2d_product_length}, 
\cref{domainup:scaling:dims}, 
\cref{Lemma:PropertiesOfCompositions_n2}, 
\cref{Prop:identity_representation}, 
and the fact that $\max\{D,2\ceil{\log_2(d)}\}\leq 2d$ \prove\ that
\begin{align}
	\llabel{loc:length}
		\nonumber&\lengthANN(\g_3\bullet\ReLUidANN{1}\bullet\g_2\bullet\compANN{\ReLUidANN{D}}{\g_1})
		\\\nonumber&=\lengthANN(\g_3)+\lengthANN(\g_2)+\lengthANN(\g_1)+4-4
		\\&=\lengthANN(\g_3)+\lengthANN(\g_2)+1
		\\\nonumber&\leq
		(d\ceil{\log_2(\scl)}+1)+ \prb{4\ceil{\log_2(d)}D+\ceil{\log_2(d)}D\ceil{\log_2(R)}+\tfrac{\ceil{\log_2(d)}(d\log_2(\scl)- \log_2(\eps))}{2}}+1
		\\\nonumber&\leq 8d^2+2d^2\ceil{\log_2(R)}+\log_2(\eps^{-1})d+d^2\ceil{\log_2(\scl)}+2
\end{align}
This,
    \dimsbene, 
		\cref{it:2d_product_dims},
		\cref{domainup:scaling:dims}, and
		\eqref{affinedef}
\prove\ that  for all 
    $k\in\N_0\cap[0,\lengthANN(\g_3)+\lengthANN(\g_2)+1]$ 
it holds that
\begin{equation}
\label{eq:d_product:dims}
    \singledims_k(\g_3\bullet\ReLUidANN{1}\bullet\g_2\bullet\compANN{\ReLUidANN{D}}{\g_1})=
		\begin{cases}
		d & \colon k=0\\
		2D& \colon k=1\\
		\singledims_{k-1}(\g_2)&\colon k\in\N\cap(1,\lengthANN(\g_2)]\\
		2& \colon k\in\N\cap(\lengthANN(\g_2),\lengthANN(\g_3)+\lengthANN(\g_2)+1)\\
		1& \colon k=\lengthANN(\g_3)+\lengthANN(\g_2)+1
		.
		\end{cases}
\end{equation}
Combining this, 
\cref{it:2d_product_dims}, 
\cref{it:2d_product_cost},
\cref{domainup:scaling:dims}, and
\cref{domainup:scaling:param} with the fact that $\singledims_0(\g_2)=D$, $D=2^{\ceil{\log_2(d)}}\leq 2(d-1)$, and $\max\{4,D\}\leq 2d$ \proves\ that
\begin{equation}
\begin{split}
\label{eq:d_product:param}
&\paramANN(\g_3\bullet\ReLUidANN{1}\bullet\g_2\bullet\compANN{\ReLUidANN{D}}{\g_1})
\\&=\ssum_{k=1}^{\lengthANN(\g_3)+\lengthANN(\g_2)+1}\singledims_k(\g_3\bullet\ReLUidANN{1}\bullet\g_2\bullet\compANN{\ReLUidANN{D}}{\g_1})(\singledims_{k-1}(\g_3\bullet\ReLUidANN{1}\bullet\g_2\bullet\compANN{\ReLUidANN{D}}{\g_1})+1)
\\&=2D(d+1)
+\singledims_1(\g_2)(2\singledims_0(\g_2)+1)
+\PR*{\ssum_{k=2}^{\lengthANN(\g_2)-1}\singledims_k(\g_2)(\singledims_{k-1}(\g_2)+1)}
\\&\quad+2\singledims_0(\g_3)(\singledims_{\lengthANN(\g_2)-1}(\g_2)+1)
+\singledims_1(\g_3)(2\singledims_0(\g_3)+1)
\\&\quad+\PR*{\ssum_{k=2}^{\lengthANN(\g_3)}\singledims_k(\g_3)(\singledims_{k-1}(\g_3)+1)}
\\&=2D(d+1)
+\singledims_1(\g_2)\singledims_0(\g_2)
+\PR*{\ssum_{k=1}^{\lengthANN(\g_2)}\singledims_k(\g_2)(\singledims_{k-1}(\g_2)+1)}
\\&\quad+(\singledims_{\hidlengthANN(\g_2)}(\g_2)+1)
+\singledims_1(\g_3)
+\PR*{\ssum_{k=1}^{\lengthANN(\g_3)}\singledims_k(\g_3)(\singledims_{k-1}(\g_3)+1)}
\\&= 2D(d+1)+3D^2+\param(\g_2)+7+\singledims_1(\g_3)+\param(\g_3)
\\&\leq 2D(d+1)+3D^2+2^{9} D^3+2^8 D^3\ceil{\log_2(R)}+2^7 D^2(d\log_2(\scl)- \log_2(\eps))
\\&\quad+6db+10
\\&\leq 4(d-1)(d+1)+12(d-1)^2+2^{12} (d-1)^3+2^{11} (d-1)^3\ceil{\log_2(R)}
\\&\quad-2^9 (d-1)^2 \log_2(\eps)+(2^9 d (d-1)^2+6d) \log_2(\scl)+6d+10
\\&\leq (2+6+2^{12}+3) d^3+2^{11} d^3\ceil{\log_2(R)}-2^9 d^2 \log_2(\eps)+(2^9+2) d^3\log_2(\scl)
\\&= 4107 d^3+2048 d^3\ceil{\log_2(R)}-512 d^2 \log_2(\eps)+514d^3\log_2(\scl).
\end{split}
\end{equation}
Combining this, \eqref{functioncheck}, 
\eqref{lipschitzcheck}, 
\eqref{approxcheck}, 
\eqref{eq:size:d_product},
\lref{loc:length},
and
\eqref{eq:d_product:dims}
with
\eqref{eq:d_product:param} 
\proves[ei]\ \cref{d_product_continuous,d_product_approx,d_product_lipschitz,d_product_length,d_product_dims,d_product_cost,d_product_size}. The proof of \cref{lemma:d_product} is thus complete\cfload.
\end{proof}

\cfclear
\begin{athm}{cor}{cor:downsized:products}
Let $d\in\N$, $\eps\in(0,\infty)$, $a\in\R$, $b\in(a,\infty)$, $\gamma\in(0,1]$, $\scl\in[1,\infty)$. Then there exists $\f\in\ANNSET$ such that
\begin{enumerate}[(i)]
    \item \label{d_downsized:product_continuous} it holds 
    that $\realisation(\mathscr{f}) \in C(\R^{d},\R)$,

    \item\label{d_downsized:product_approx} it holds that
    $
    \sup_{x = (x_1,\ldots, x_{d}) \in [a,b]^{d}}\vass[\big]{\gamma\scl^d\sprod_{i=1}^{d}x_i - (\realisation(\f)) (x)} 
                \leq 
                \eps$,

    \item\label{d_downsized:product_length} it holds 
    that
    $\lengthANN(\mathscr{f}) \leq 43 d^2 \max\pR{1,\ceil{\log_2(\max\{\vass{a},\vass{b}\})},\log_2(\eps^{-1}),\ceil{\log_2(\scl)}}
    $,

    \item\label{d_downsized:product_cost} it holds 
    that
    $\paramANN(\mathscr{f}) \leq 7743 d^3 \max\pR{1,\ceil{\log_2(\max\{\vass{a},\vass{b}\})},\log_2(\eps^{-1}),\ceil{\log_2(\scl)}}
    $, and
		
		\item\label{d_downsized:product_size} it holds 
    that
    $\size(\mathscr{f}) \leq 1
    $
\end{enumerate}
\cfout.
\end{athm}

\begin{proof}[Proof of \cref{cor:downsized:products}]
Throughout this proof 
assume without loss of generality that
$\eps<1<\max\{\vass{a},\vass{b}\}$ and
let $R\in(1,\infty)$ satisfy $R=\max\{\vass{a},\vass{b}\}$.
\Nobs that \cref{lemma:d_product} (applied with 
$d \curvearrowleft d$,
$\eps \curvearrowleft \eps$,
$R \curvearrowleft R$,
$\gamma \curvearrowleft \gamma$,
$\scl \curvearrowleft \scl$
in the notation of \cref{lemma:d_product}) \proves\ that there exists $\g\in\ANNSET$ which satisfies that
\begin{enumerate}[(I)]
    \item \label{it3:d_downsized:product_continuous} it holds 
    that $\realisation(\g) \in C(\R^{d},\R)$,

    \item\label{it3:d_downsized:product_approx} it holds that
    $
    \sup_{x = (x_1,\ldots, x_{d}) \in [-\radius, \radius]^{d}}\vass[\big]{\gamma\scl^d\sprod_{i=1}^{d}x_i - (\realisation(\g)) (x)} 
                \leq 
                \eps$,


    \item\label{it3:d_downsized:product_length} it holds 
    that
    $\lengthANN(\g) \leq 8d^2+2d^2\ceil{\log_2(R)}+\log_2(\eps^{-1})d+d^2\ceil{\log_2(\scl)}+2
    $,
		
		\item\label{it3:d_downsized:product_dims} it holds that
				$\singledims_{\hidlengthANN(\g)}(\g)=2
		$,

    \item\label{it3:d_downsized:product_cost} it holds 
    that
    $\paramANN(\g) \leq 4107 d^3+2048 d^3\ceil{\log_2(R)}-512 d^2 \log_2(\eps)+514d^3\log_2(\scl)
    $, and
		
		\item\label{it3:d_downsized:product_size} it holds 
    that
    $\size(\g) \leq 4
    $
\end{enumerate}
\Nobs that \cref{it3:d_downsized:product_length} \proves\ that
\begin{equation}
\label{eq:length:downsized}
\begin{split}
\lengthANN(\g) &\leq 8d^2+2d^2\ceil{\log_2(R)}+d\log_2(\eps^{-1})+d^2\ceil{\log_2(\scl)}+2
\\&\leq (8+2+1+1+2)d^2\max\pR{1,\ceil{\log_2(R)},\log_2(\eps^{-1}),\ceil{\log_2(\scl)}}
\\&= 14d^2\max\pR{\ceil{\log_2(R)},\log_2(\eps^{-1}),\ceil{\log_2(\scl)}}.
\end{split}
\end{equation}
\Moreover \cref{it3:d_downsized:product_cost} \proves\ that
\begin{equation}
\label{eq:param:downsized}
\begin{split}
\param(\g) &\leq 4107 d^3+2048 d^3\ceil{\log_2(R)}-512 d^2 \log_2(\eps)+514d^3\log_2(\scl)
\\&\leq (4107+2048+512+514)d^3 \max\pR{1,\ceil{\log_2(R)},\log_2(\eps^{-1}),\ceil{\log_2(\scl)}}
\\&=7181 d^3 \max\pR{\ceil{\log_2(R)},\log_2(\eps^{-1}),\ceil{\log_2(\scl)}}
.
\end{split}
\end{equation}
Combining this, 
\eqref{eq:length:downsized},
\cref{it3:d_downsized:product_continuous}, 
\cref{it3:d_downsized:product_dims},
and \cref{it3:d_downsized:product_size} with \cref{lemma:network:halfed:paramsize} (applied with 
$\f \curvearrowleft \g$,
$d \curvearrowleft d$,
$\dnsc \curvearrowleft 4^{-1}$
in the notation of \cref{lemma:network:halfed:paramsize}) \proves\ that there exists $\f\in\ANNSET$ which satisfies that
\begin{enumerate}[(A)]
    \item \label{it:d_downsized:product_continuous} it holds 
    that $\realisation(\mathscr{f})=\realisation(\mathscr{g}) \in C(\R^{d},\R)$,

    \item\label{it:d_downsized:product_length} it holds 
    that
    \begin{equation}
    	\begin{split}
    \lengthANN(\mathscr{f}) &\leq 3(14d^2\max\pR{\ceil{\log_2(R)},\log_2(\eps^{-1}),\ceil{\log_2(\scl)}})+1
    \\&\leq 43 d^2 \max\pR{\ceil{\log_2(R)},\log_2(\eps^{-1}),\ceil{\log_2(\scl)}}
    ,
    	\end{split}
	\end{equation}

    \item\label{it:d_downsized:product_cost} it holds 
    that
    \begin{equation}
		\begin{split}
		\paramANN(\mathscr{f}) 
		&\leq
		(7181 d^3+2+40(14d^2)) \max\pR{\ceil{\log_2(R)},\log_2(\eps^{-1}),\ceil{\log_2(\scl)}}
		\\&\leq
		7743 d^3 \max\pR{\ceil{\log_2(R)},\log_2(\eps^{-1}),\ceil{\log_2(\scl)}}
		,
		\end{split}
		\end{equation}
    and
		
		\item\label{it:d_downsized:product_size} it holds 
    that
    $\size(\mathscr{f}) \leq 1
    $.
\end{enumerate}
\Nobs that \cref{it3:d_downsized:product_approx} and \cref{it:d_downsized:product_continuous} \prove\ that for all $x=(x_1,\ldots,x_d)\in[a,b]^d\subseteq[-R,R]^d$ it holds that
\begin{equation}
	\textstyle
\vass[\big]{\gamma\scl^d\sprod_{i=1}^{d}x_i - (\realisation(\f)) (x)} 
                \leq 
                \eps.
\end{equation}
Combining this with
\cref{it:d_downsized:product_continuous,it:d_downsized:product_length,it:d_downsized:product_cost,it:d_downsized:product_size} \proves[ie]\
\cref{d_downsized:product_continuous,d_downsized:product_approx,d_downsized:product_length,d_downsized:product_cost,d_downsized:product_size}.
\finishproofthus
\end{proof}

\subsection{Upper bounds for approximations of periodic functions}
\label{subsection:upper_periodic}


\cfclear
\begin{athm}{lemma}{edgy:sin:prop}
	Let $\lambda\in(0,\infty)$ and let $\edgy\colon\R\to\R$ satisfy for all $k\in\Z$, $x\in[2k,2k+2)$ that $\edgy(x)=1-\vass{x-2k-1}$. Then it holds for all $x\in\R$ that 
	\begin{equation}
		\label{step:edgy:sin:eq}
		\edgy(2\lambda{}x)=2\edgy(\lambda{}x)-4\RELU\prb{\edgy(\lambda{}x)-\tfrac{1}{2}}
	\end{equation}
	\cfout.
\end{athm}

\begin{proof}[Proof of \cref{edgy:sin:prop}]
	\Nobs that 
	the fact that for all $k\in\Z$, $x\in[2k,2k+1)$ it holds that $\edgy(x)=x-2k$ \proves\ that
	for all $k\in\Z$, $x\in\R$ with $2\lambda{}x\in[4k,4k+1)$ it holds that
	\begin{equation}
		\label{step:edgy:sin:case1}
		\begin{split}
			2\edgy(\lambda{}x)-4\RELU\prb{\edgy(\lambda{}x)-\tfrac{1}{2}}
			&=2(\lambda{}x-2k)-4\RELU\prb{(\lambda{}x-2k)-\tfrac{1}{2}}
			\\&=(2\lambda{}x-4k)-2\RELU\prb{2\lambda{}x-4k-1}
			\\&=2\lambda{}x-4k
			\\&=\edgy(2\lambda{}x)
		\end{split}
	\end{equation}
	\cfload. \Moreover 
	the fact that for all $k\in\Z$, $x\in[2k,2k+1)$, $y\in[2k+1,2k+2)$ it holds that $\edgy(x)=x-2k$ and $\edgy(y)=-y+2k+2$ \proves\ that
	for all $k\in\Z$, $x\in\R$ with $2\lambda{}x\in[4k+1,4k+2)$ it holds that
	\begin{equation}
		\label{step:edgy:sin:case2}
		\begin{split}
			2\edgy(\lambda{}x)-4\RELU\prb{\edgy(\lambda{}x)-\tfrac{1}{2}}
			&=2(\lambda{}x-2k)-4\RELU\prb{(\lambda{}x-2k)-\tfrac{1}{2}}
			\\&=(2\lambda{}x-4k)-2\RELU\prb{2\lambda{}x-4k-1}
			\\&=-2\lambda{}x+(4k+2)
			\\&=\edgy(2\lambda{}x).
		\end{split}
	\end{equation}
	\Moreover 
	the fact that for all $k\in\Z$, $x\in[2k,2k+1)$, $y\in[2k+1,2k+2)$ it holds that $\edgy(x)=x-2k$ and $\edgy(y)=-y+2k+2$ \proves\ that
	for all $k\in\Z$, $x\in\R$ with $2\lambda{}x\in[4k+2,4k+3)$ it holds that
	\begin{equation}
		\label{step:edgy:sin:case3}
		\begin{split}
			&2\edgy(\lambda{}x)-4\RELU\prb{\edgy(\lambda{}x)-\tfrac{1}{2}}
			\\&=2(-\lambda{}x+2k+2)-4\RELU\prb{(-\lambda{}x+2k+2)-\tfrac{1}{2}}
			\\&=-2\lambda{}x+4k+4-2\RELU\prb{-2\lambda{}x+4k+3}
			\\&=2\lambda{}x-(4k+2)
			\\&=\edgy(2\lambda{}x).
		\end{split}
	\end{equation}
	\Moreover 
	the fact that for all $k\in\Z$, $y\in[2k+1,2k+2)$ it holds that $\edgy(y)=-y+2k+2$ \proves\ that
	for all $k\in\Z$, $x\in\R$ with $2\lambda{}x\in[4k+3,4k+4)$ it holds that
	\begin{equation}
		\label{step:edgy:sin:case4}
		\begin{split}
			&2\edgy(\lambda{}x)-4\RELU\prb{\edgy(\lambda{}x)-\tfrac{1}{2}}
			\\&=2(-\lambda{}x+2k+2)-4\RELU\prb{(-\lambda{}x+2k+2)-\tfrac{1}{2}}
			\\&=-2\lambda{}x+4k+4-2\RELU\prb{-2\lambda{}x+4k+3}
			\\&=-2\lambda{}x+(4k+4)
			\\&=\edgy(2\lambda{}x).
		\end{split}
	\end{equation}
	Combining this with \eqref{step:edgy:sin:case1}, \eqref{step:edgy:sin:case2}, and \eqref{step:edgy:sin:case3} \proves\ \eqref{step:edgy:sin:eq}.
	The proof of \cref{edgy:sin:prop} is thus complete.
\end{proof}

			\cfclear
			\begin{athm}{lemma}{lem:edgy:sin}
				Let $B\in(0,\infty)$, 
				let $\edgy\colon\R\to\R$ satisfy for all $k\in\Z$, $x\in[2k,2k+2)$ that $\edgy(x)=1-\vass{x-2k-1}$, 
				let $\g\in \prb{(\R^{4 \times 1} \times \R^{4}) \times (\R^{1 \times 4} \times \R^1) }\subseteq \ANNSET$ satisfy 
				\begin{equation}
					\label{lem:stepnet}
					\begin{split}
						\g &= \pr*{ \!\pr*{\!\begin{pmatrix}
									1\\
									1\\
									1\\
									1
								\end{pmatrix},
								\begin{pmatrix}
									0\\
									0\\
									\nicefrac{-1}{2}\\
									\nicefrac{-1}{2}
								\end{pmatrix}\! },
							\prbb{	\begin{pmatrix}
									1 & 1 &-2 & -2
								\end{pmatrix}, 
								0 } \!
						}  ,
					\end{split}
				\end{equation}
				and let $\f_n\in\ANNSET$, $n\in\N_0$, satisfy for all $n\in\N$ that
				\begin{equation}
					\label{lem:steprule}
					\begin{split}
						\f_0 &= \pr*{ \!\pr*{\!\begin{pmatrix}
									\nicefrac{2}{B}\\
									\nicefrac{2}{B}
								\end{pmatrix},
								\begin{pmatrix}
									0\\
									-1
								\end{pmatrix}\! },
							\prbb{	\begin{pmatrix}
									1& -2
								\end{pmatrix}, 
								0 } \!
						}
						\qandq
						\f_n = \compANN{\g}{\f_{n-1}}
					\end{split}
				\end{equation}
				\cfload. Then
				\begin{enumerate}[(i)]
					\item{
						\label{eq:edgy:sin1}
						it holds for all $n\in\N_0$, $x\in[0,B]$ that $(\realisation(\f_n))(x)=\edgy(2^{n+1}B^{-1}x)$,
					}
					\item{
						\label{eq:edgy:sin2}
						it holds for all $n\in\N$, $x\in\pr{\R\setminus[0,B]}$ that $(\realisation(\f_n))(x)=0$,
					}
					\item{
						\label{eq:edgy:sin4}
						it holds for all $n\in\N$ that
						$\dims(\f_0)=(1,2,1)$ and $\dims(\f_n)=(1,2,4,4,\ldots,4,1)\in\N^{\lengthANN(\f_n)+1}$,
					}
					\item{
						\label{eq:edgy:sin3}
						it holds for all $n\in\N_0$ that $\lengthANN(\f_n)=n+2$, and
					}
					
					\item{
						\label{eq:edgy:sin6}
						it holds for all $n\in\N_0$ that $\size(\f_n)\leq \max\{2B^{-1},2\}$, $\insize(\f_n)= \max\{2B^{-1},1\}$, and $\outsize(\f_n)\leq 2$
					}
				\end{enumerate}
				\cfout.
			\end{athm}
			
			\begin{proof}[Proof of \cref{lem:edgy:sin}]
				Throughout this proof let
				$W\in\R^{4\times 4}$, $\mathfrak{W}\in\R^{4\times 2}$
				satisfy 
				\begin{equation}
					\label{eq:matrices:networks:scaled:hat}
					W=\weight{1}{\g}\weight{2}{\g}
					\qandq
					\mathfrak{W}=\weight{1}{\g}\weight{2}{\f_0}
					.
				\end{equation}
				\Nobs that \eqref{lem:stepnet}, 
				\eqref{lem:steprule}, and 
				\cref{Lemma:PropertiesOfCompositions_n2} \prove\ that for all $n\in\N$ it holds that
				\begin{equation}
					\label{eq:ind:step:length:scaled:hat}
					\lengthANN(\f_0)=2
					\qandq
					\lengthANN(\f_n)=\lengthANN(\compANN{\g}{\f_{n-1}})=\lengthANN(\g)+\lengthANN(\f_{n-1})-1=\lengthANN(\f_{n-1})+1.
				\end{equation} 
				Combining this with induction \proves\ that for all $n\in\N_0$ it holds that
				\begin{equation}
					\label{eq:length:scaled:hat}
					\lengthANN(\f_n)=n+2.
				\end{equation}
				This  \proves[ei]\ \cref{eq:edgy:sin3}.
				\Moreover 
				\eqref{ANNoperations:Composition}, 
				\eqref{lem:stepnet}, 
				\eqref{lem:steprule}, 
				\eqref{eq:matrices:networks:scaled:hat} and the fact that
				$
				\weight{1}{\g}\bias{2}{\f_0}+\bias{1}{\g}=\bias{1}{\g}
				$
				\prove\ that
				\begin{equation}
					\begin{split}
						\label{eq:structured:form:of:scaled:hat:n1}
						\f_1=\compANN{\g}{\f_{0}}=\prb{(\weight{1}{\f_0},\bias{1}{\f_0}),(\mathfrak{W},\bias{1}{\g}),(\weight{2}{\g},\bias{2}{\g})}.
					\end{split}
				\end{equation}
				This, 
				\eqref{ANNoperations:Composition}, 
				\eqref{lem:stepnet}, 
				\eqref{lem:steprule}, 
				\eqref{eq:matrices:networks:scaled:hat} and the fact that
				$
				\weight{1}{\g} \bias{2}{\g}+\bias{1}{\g}=\bias{1}{\g}
				$
				\prove\ that
				\begin{equation}
					\begin{split}
						\label{eq:structured:form:of:scaled:hat:n2}
						\f_2=\compANN{\g}{\f_{1}}=\prb{(\weight{1}{\f_0},\bias{1}{\f_0}),(\mathfrak{W},\bias{1}{\g}),(W,\bias{1}{\g}),(\weight{2}{\g},\bias{2}{\g})}.
					\end{split}
				\end{equation}
				Combining this, 
				\eqref{ANNoperations:Composition}, 
				\eqref{eq:length:scaled:hat},
				and \eqref{eq:matrices:networks:scaled:hat} with induction \proves\ that for all $n\in\N\setminus\{1\}$ it holds that
				\begin{equation}
					\label{eq:structured:form:of:scaled:hat}
					\begin{split}
						\f_n
						&=((\weight{1}{\f_0},\bias{1}{\f_0}),(\mathfrak{W},\bias{1}{\g}),(W,\bias{1}{\g}),(W,\bias{1}{\g}),\ldots,(W,\bias{1}{\g}),(\weight{2}{\g},\bias{2}{\g}))
						\\&\in
						\prb{(\R^{2 \times 1} \times \R^{2}) \times
						(\R^{4 \times 2} \times \R^{4}) \times\prb{\times_{k=1}^{n-1}\pr{\R^{4 \times 4} \times \R^{4}}}\times (\R^{1 \times 4} \times \R^1) }.
					\end{split}
				\end{equation}
				This, 
				\eqref{lem:stepnet},
				\eqref{eq:length:scaled:hat},
				\eqref{eq:structured:form:of:scaled:hat:n1},
				\eqref{eq:structured:form:of:scaled:hat:n2}, and
				\eqref{eq:structured:form:of:scaled:hat} \prove\ that for all $n\in\N_0$, $k\in\N_0\cap[0,\lengthANN(\f_n)]$ it holds that
				\begin{equation}
					\singledims_k(\f_n)=
					\begin{cases}
						1 & \colon k\in\{0,\lengthANN(\f_n)\}\\
						2 & \colon k=1\\
						4 & \colon k\in\N\cap(1,\lengthANN(\f_n)).
					\end{cases}
				\end{equation}
				This  \proves[ei]\ \cref{eq:edgy:sin4}. 
				\Moreover \eqref{lem:steprule} 
				and the fact that for all $x\in\R$ it holds that 
				$\RELU(\vass{x})=\RELU(x)+\RELU(-x)$
				and
				$\RELU(\vass{x}-1)=\RELU(x-1)+\RELU(-x-1)$ 
				\prove\ that for all $x\in\R$ it holds that
				\begin{equation}
					\label{basecase:edgy:sin:pos}
					(\realisation(\f_0))(x)=\RELU(2B^{-1}x+0)-2\RELU(2B^{-1}x-1)=
					\begin{cases} 
						0 &\colon x\in(-\infty,0)\\
						2B^{-1}x &\colon x\in[0,B2^{-1}]\\
						2-2B^{-1}x &\colon x\in[B2^{-1},\infty)
					\end{cases}
				\end{equation}
				\cfload. This and the assumption that for all $x\in[0,1]$ it holds that $\edgy(x)=x=\edgy(2-x)$
				 \prove\ that for all $x\in[0,\tfrac{B}{2}]$, $y\in[\tfrac{B}{2},B]$, $z\in\pr{\R\setminus[0,B]}$ it holds that 
				\begin{equation}
					\label{base:edgy:sin00}
					(\realisation(\f_0))(x)=\edgy(2B^{-1}x),
					\qquad
					(\realisation(\f_0))(y)=\edgy(2B^{-1}y),
					\qquad\text{and}\qquad
					(\realisation(\f_0))(z)\leq 0.
				\end{equation}
				\Moreover \eqref{lem:stepnet}, \eqref{lem:steprule}, and \cref{Lemma:PropertiesOfCompositions_n2} \prove\ that for all $n\in\N$, $x\in\R$, $y\in(-\infty,0]$ it holds that 
				\begin{equation}
					\label{eq:inductionstep:functioncomp:g:at:zero}
					(\realisation(\f_n))(x)=(\realisation(\g\bullet\f_{n-1}))(x)=(\realisation(\g))\pr*{(\realisation(\f_{n-1}))(x)}
					\qandq
					(\realisation(\g))\pr*{y}=0.
				\end{equation}
				\Moreover
				\eqref{lem:stepnet}, 
				\eqref{eq:inductionstep:functioncomp:g:at:zero}, 
				\cref{Lemma:PropertiesOfCompositions_n2}, 
				and \cref{edgy:sin:prop} (applied with
				$\lambda \curvearrowleft  2^{n} B^{-1}$
				for $n\in\N$ in the notation of \cref{edgy:sin:prop}) \prove\ that for all $x\in[0,B]$, $n\in\N$ with $(\realisation(\f_{n-1}))(x)=\edgy(2^{n}B^{-1}x)$ it holds that
				\begin{equation}
					\label{step:edgy:sin}
					\begin{split}
						(\realisation(\f_{n}))(x)
						&=(\realisation(\g))\pr*{(\realisation(\f_{n-1}))(x)}
						\\&=(\realisation(\g))\pr*{\edgy(2^{n}B^{-1}x)}
						\\&=2\RELU(\edgy(2^{n}B^{-1}x))-4\RELU\prb{\edgy(2^{n}B^{-1}x)-\tfrac{1}{2}}
						\\&=2\edgy(2^{n}B^{-1}x)-4\RELU\prb{\edgy(2^{n}B^{-1}x)-\tfrac{1}{2}}
						=\edgy(2^{n+1}B^{-1}x).
					\end{split}
				\end{equation}
				Combining this and
				\eqref{base:edgy:sin00}
				with induction  \proves[ei]\ \cref{eq:edgy:sin1}.
				\Moreover \eqref{basecase:edgy:sin:pos}
				and \eqref{eq:inductionstep:functioncomp:g:at:zero} \prove\ that for all $x\in\pr{\R\setminus[0,B]}$ it holds that
				\begin{equation}
					\label{upper:area}
					\begin{split}
						(\realisation(\f_1))(x)=(\realisation(\g))\pr*{(\realisation(\f_{0}))(x)}
						=0.
					\end{split}
				\end{equation}
				Combining this and
				\eqref{eq:inductionstep:functioncomp:g:at:zero}
				with induction  \proves[ei]\ \cref{eq:edgy:sin2}.
				\Moreover 
				\eqref{lem:stepnet},
				\eqref{eq:matrices:networks:scaled:hat},
				\eqref{eq:structured:form:of:scaled:hat:n1},
				\eqref{eq:structured:form:of:scaled:hat:n2}, and
				\eqref{eq:structured:form:of:scaled:hat}
				\prove\ that for all $n\in\N_0$ it holds that
				\begin{equation}
					\label{eq:indims:outdims:edgysin:prep}
					\insize(\f_n)=\insize(\f_0)=\max\{2B^{-1},1\}
					\qandq
					\outsize(\f_n)\leq\max\{\outsize(\g),\outsize(\f_0)\}=2.
				\end{equation}
				\Moreover 
				\eqref{lem:stepnet},
				\eqref{eq:matrices:networks:scaled:hat},
				\eqref{eq:structured:form:of:scaled:hat:n1},
				\eqref{eq:structured:form:of:scaled:hat:n2}, and
				\eqref{eq:structured:form:of:scaled:hat}
				\prove\ that for all $n\in\N_0$ it holds that
				\begin{equation}
					\size(\f_n)\leq\max\{\size(\g),\size(\f_0),\supn{W},\supn{\mathfrak{W}}\}=\max\{2,2B^{-1},2,2\}=\max\{2B^{-1},2\}
				\end{equation}
				\cfload. 
				This and \eqref{eq:indims:outdims:edgysin:prep}  \prove[ei]\ \cref{eq:edgy:sin6}.
				The proof of \cref{lem:edgy:sin} is thus complete.
			\end{proof}

			\cfclear
			\begin{athm}{prop}{single:hat:network:prep}
				Let $n,N\in\N$, $a,b\in[0,2\pi]$, $c\in\R$ satisfy $b=a+\frac{4\pi}{N+1}$,
				let $\edgy\colon\R\to\R$ satisfy for all $k\in\Z$, $x\in[2k,2k+2)$ that $\edgy(x)=1-\vass{x-2k-1}$,
				let $f\colon\R\to\R$ satisfy for all $x\in[0,2\pi)$, $y\in(-\infty,-2^n\pi)\cup[2^n\pi,\infty)$, $k\in\Z\cap[-2^{n-1},2^{n-1})$ that
				\begin{equation}
					\label{eq:setup:hatfunctions0:prep}
					f(y)=0
					\qandq
					f(x+2k\pi)=
					f(x)=
					\begin{cases}
						c\,\edgy\prb{\tfrac{(x-a)(N+1)}{2\pi}} &\colon x\in [a,b]\\
						0 &\colon x\notin [a,b],
					\end{cases}
				\end{equation}
				let $\mathscr{a},\mathscr{b},\mathscr{c}\in\R$ satisfy $\mathscr{a}=a-\frac{(N-1)\pi}{N+1}$, $\mathscr{b}=\frac{N-1}{N+1}$, and $\mathscr{c}=\frac{(N+1)c}{2}$,
				and let $g\colon\R\to\R$ satisfy for all $x\in\R$ that
				\begin{equation}
					\label{eq:param:introd:scr:prep}
					\tgt(x)=
					\begin{cases}
						\mathscr{c}\RELU\prb{\edgy(\pi^{-1}(x-a)+\mathscr{b})-\mathscr{b}} &\colon x\in [-2^n\pi+\mathscr{a},2^n\pi+\mathscr{a}]\\
						0 &\colon x\notin [-2^n\pi+\mathscr{a},2^n\pi+\mathscr{a}]
					\end{cases}
				\end{equation}
				\cfload. Then $f=g$.
			\end{athm}
			
			\begin{proof}[Proof of \cref{single:hat:network:prep}]
				\Nobs that \eqref{eq:param:introd:scr:prep} and the fact that for all $x\in[0,\infty)$ it holds that $\edgy(x)\leq x$ \prove\ that for all $x\in[\mathscr{a},a]$ it holds that
				\begin{equation}
					\label{eq:realisation:composition:interval1}
					\begin{split}
						\tgt(x)
						&=\mathscr{c}\RELU\prb{\edgy(\pi^{-1}(x-a)+\mathscr{b})-\mathscr{b}}=0.
					\end{split}
				\end{equation}
				\Moreover \eqref{eq:param:introd:scr:prep} and the fact that for all $x\in[0,1]$ it holds that $\edgy(x)=x=\RELU(x)$ \prove\ that for all $x\in\PR*{a,a+\tfrac{2\pi}{N+1}}$ it holds that
				\begin{equation}
					\label{eq:realisation:composition:interval2}
					\begin{split}
						\tgt(x)
						&=\mathscr{c}\RELU\prb{\edgy(\pi^{-1}(x-a)+\mathscr{b})-\mathscr{b}}
						\\&=\mathscr{c}\RELU\prb{\edgy\pr*{\pi^{-1}(x-a)+\tfrac{N-1}{N+1}}-\tfrac{N-1}{N+1}}
						\\&=\mathscr{c}\RELU\prb{\pi^{-1}(x-a)+\tfrac{N-1}{N+1}-\tfrac{N-1}{N+1}}
						\\&=\prb{\tfrac{2\mathscr{c}}{N+1}}\RELU\prb{\tfrac{(x-a)(N+1)}{2\pi}}
						\\&=c\,\RELU\prb{\tfrac{(x-a)(N+1)}{2\pi}}
						\\&=c\,\edgy\prb{\tfrac{(x-a)(N+1)}{2\pi}}
						.
					\end{split}
				\end{equation}
				\Moreover \eqref{eq:param:introd:scr:prep} and the fact that for all $x\in[1,2]$ it holds that $\edgy(x)=\edgy(2-x)=2-x$ and the fact that for all $x\in\PR*{a+\frac{2\pi}{N+1},b}$ it holds that $\pi^{-1}(x-a)\in\PR*{\frac{2}{N+1},\frac{4}{N+1}}$ and $b\leq a+\frac{(N+3)\pi}{N+1}=\mathscr{a}+2\pi$ \prove\ that for all $x\in\PR*{a+\frac{2\pi}{N+1},b}$ it holds that
				\begin{equation}
					\label{eq:realisation:composition:interval3}
					\begin{split}
						\tgt(x)
						&=\mathscr{c}\RELU\prb{\edgy(\pi^{-1}(x-a)+\mathscr{b})-\mathscr{b}}
						\\&=\mathscr{c}\RELU\prb{\edgy\pr*{\pi^{-1}(x-a)+\tfrac{N-1}{N+1}}-\tfrac{N-1}{N+1}}
						\\&=\mathscr{c}\RELU\prb{2-\pi^{-1}(x-a)-\tfrac{N-1}{N+1}-\tfrac{N-1}{N+1}}
						\\&=\mathscr{c}\RELU\prb{\tfrac{4}{N+1}-\pi^{-1}(x-a)}
						\\&=\prb{\tfrac{2\mathscr{c}}{N+1}}\RELU\prb{2-\tfrac{(x-a)(N+1)}{2\pi}}
						\\&=c\,\RELU\prb{2-\tfrac{(x-a)(N+1)}{2\pi}}
						\\&=c\prb{2-\tfrac{(x-a)(N+1)}{2\pi}}
						\\&=c\,\edgy\prb{\tfrac{(x-a)(N+1)}{2\pi}}
						.
					\end{split}
				\end{equation}
				\Moreover \eqref{eq:param:introd:scr:prep} and the fact that for all $x\in\PRb{\frac{4}{N+1},\frac{N+3}{N+1}}$ it holds that $\edgy(x+\mathscr{b})\leq\edgy\pr*{\frac{N+3}{N+1}}= \mathscr{b}$ \prove\ that for all $x\in[b,2\pi+\mathscr{a}]=\PRb{a+\frac{4\pi}{N+1},a+\frac{(N+3)\pi}{N+1}}$ it holds that
				\begin{equation}
					\label{eq:realisation:on:interval4}
					\begin{split}
						\tgt(x)
						&=\mathscr{c}\RELU\prb{\edgy(\pi^{-1}(x-a)+\mathscr{b})-\mathscr{b}}
						=0.
					\end{split}
				\end{equation}
				Combining this, 
				\eqref{eq:setup:hatfunctions0:prep},
				\eqref{eq:param:introd:scr:prep}, 
				\eqref{eq:realisation:composition:interval1}, 
				\eqref{eq:realisation:composition:interval2}, and
				\eqref{eq:realisation:composition:interval3} with the fact that $\mathscr{a}\leq a$ \proves\ that for all 
				$x\in[a,b]$, $y\in[\mathscr{a},a)\cup(b,2\pi+\mathscr{a}]$
				it holds that
				\begin{equation}
					\label{eq:realisation:on:atopi}
					\tgt(x)=c\,\edgy\prb{\tfrac{(x-a)(N+1)}{2\pi}}=f(x)
					\qandq
					\tgt(y)=0=f(y).
				\end{equation}
				\Moreover \eqref{eq:param:introd:scr:prep} 
				and the fact that for all $x\in\R$, $k\in\Z$ it holds that $\edgy(x+2k)= \edgy(x)$ 
				\prove\ that for all $x\in[\mathscr{a},2\pi+\mathscr{a}]$, $k\in\Z\cap[-2^{n-1},2^{n-1})$ it holds that
				\begin{equation}
					\label{eq:realisation:on:otoB}
					\begin{split}
						\tgt(x+2k\pi)
						=\mathscr{c}\RELU\prb{\edgy(\pi^{-1}(x-\mathscr{a})+2k)-\mathscr{b}}
						=\mathscr{c}\RELU\prb{\edgy(\pi^{-1}(x-\mathscr{a}))-\mathscr{b}}
						=\tgt(x).
					\end{split}
				\end{equation}
				This and \eqref{eq:setup:hatfunctions0:prep} \prove\ that for all $x\in[\mathscr{a},\mathscr{a}+2\pi)$, $k\in\Z\cap[-2^{n-1},2^{n-1})$ with $x+2k\pi\in[\max\{\mathscr{a},0\}-2^n \pi,\min\{\mathscr{a},0\}+2^n\pi)$ it holds that
				\begin{equation}
					\label{eq:2pi:periodic:area:for:f:and:g}
					\tgt(x+2k\pi)=\tgt(x)
					\qandq
					f(x+2k\pi)=f(x).
				\end{equation}
				This and \eqref{eq:realisation:on:atopi}
				\prove\ that for all $x\in[\max\{\mathscr{a},0\}-2^n\pi,\min\{\mathscr{a},0\}+2^n\pi)$ it holds that
				\begin{equation}
					\label{eq:match:middle:large:interval}
					\tgt(x)=f(x).
				\end{equation}
				\Moreover 
				\eqref{eq:setup:hatfunctions0:prep},
				\eqref{eq:param:introd:scr:prep},
				and
				the fact that $2\pi\geq b$ and $\mathscr{a}\leq a$
				\prove\ that for all $x\in\R$ with $-2^n\pi\leq x \leq \mathscr{a}-2^n\pi$ it holds that
				\begin{equation}
					\label{eq:lower:part:matha:pos}
					f(x)=f(x+2^n\pi)=0=\tgt(x).
				\end{equation}
				\Moreover 
				\eqref{eq:setup:hatfunctions0:prep},
				\eqref{eq:realisation:composition:interval1},
				\eqref{eq:realisation:on:otoB}, and
				the fact that $\max\{\mathscr{a},0\}\leq a$
				\prove\ that for all $x\in\R$ with $\mathscr{a}-2^n\pi\leq x \leq -2^n\pi$ it holds that
				\begin{equation}
					\label{eq:lower:part:matha:neg}
					\tgt(x)=\tgt(x+2^n\pi)=0=f(x).
				\end{equation}
				\Moreover 
				\eqref{eq:setup:hatfunctions0:prep},
				\eqref{eq:param:introd:scr:prep},
				and
				the fact that $\mathscr{a}\geq b-2\pi$ and $0\leq a$
				\prove\ that for all $x\in\R$ with $\mathscr{a}+2^n\pi\leq x \leq 2^n\pi$ it holds that
				\begin{equation}
					\label{eq:upper:part:matha:pos}
					f(x)=f(x-2^n\pi)=0=\tgt(x).
				\end{equation}
				\Moreover 
				\eqref{eq:setup:hatfunctions0:prep},
				\eqref{eq:realisation:composition:interval1},
				\eqref{eq:realisation:on:atopi}, and
				the fact that $2\pi\geq b$
				\prove\ that for all $x\in\R$ with $2^n\pi\leq x \leq \mathscr{a}+2^n\pi$ it holds that
				\begin{equation}
					\label{eq:upper:part:matha:neg}
					\tgt(x)=\tgt(x-2(2^{n-1}-1)\pi)=0=f(x).
				\end{equation}
				Combining this,
				\eqref{eq:lower:part:matha:pos},
				\eqref{eq:lower:part:matha:neg}, and
				\eqref{eq:upper:part:matha:pos}
				with 
				\eqref{eq:setup:hatfunctions0:prep} and
				\eqref{eq:param:introd:scr:prep}
				\proves\ that for all $x\in(-\infty,\max\{\mathscr{a},0\}-2^n\pi]$, $y\in[\min\{\mathscr{a},0\}+2^n\pi,\infty)$ it holds that
				\begin{equation}
					f(x)=\tgt(x)=0=f(y)=\tgt(y).
				\end{equation}
				This and \eqref{eq:match:middle:large:interval}  \prove[ei]\ $f=g$.
				The proof of \cref{single:hat:network:prep} is thus complete.
			\end{proof}

			\cfclear

			\begin{athm}{lemma}{single:hat:network}
				Let $n,N\in\N$, $C\in[1,\infty)$, $a,b\in[0,2\pi]$, $c\in[-C,C]$ satisfy $b=a+\tfrac{4\pi}{N+1}$ and $n\geq 2$, 
				let $\edgy\colon\R\to\R$ satisfy for all $k\in\Z$, $x\in[2k,2k+2)$ that $\edgy(x)=1-\vass{x-2k-1}$,
				and let $f\colon\R\to\R$ satisfy for all $x\in[0,2\pi)$, $y\in(-\infty,-2^n\pi)\cup[2^n\pi,\infty)$, $k\in\Z\cap[-2^{n-1},2^{n-1})$ that
				\begin{equation}
					\label{eq:setup:hatfunctions0}
					f(y)=0
					\qandq
					f(x+2k\pi)=
					f(x)=
					\begin{cases}
						c\,\edgy\prb{\tfrac{(x-a)(N+1)}{2\pi}} &\colon x\in [a,b]\\
						0 &\colon x\notin [a,b].
					\end{cases}
				\end{equation}
				Then there exists $\f\in\ANNSET$ such that 
				\begin{enumerate}[(i)]
					\item{
						\label{eq:edgy:base1.2}
						it holds for all $x\in\R$ that
						$
						(\realisation(\f))(x)= \prb{\tfrac{2}{C(N+1)}}f(x)
						$,
					}
					\item{
						\label{eq:edgy:base2.1}
						it holds that $\lengthANN(\f)=n+5$,
					}
					\item{
						\label{eq:edgy:base3.1}
						it holds for all $k\in\N_0\cap[0,\lengthANN(\f)]$ that
						\begin{equation}
							\singledims_k(\f)=
							\begin{cases}
								1 & \colon k\in\{0,n+4,n+5\}\\
								2 & \colon k \in \{1,2,3,n+3\}\\
								4 & \colon k\in\N\cap(3,n+3),\\
							\end{cases}
						\end{equation}
						and
					}
					\item{
						\label{eq:edgy:base5.1}
						it holds that $\size(\f)\leq 2$
					}
				\end{enumerate}
				\cfout.
			\end{athm}

			\begin{proof}[Proof of \cref{single:hat:network}]
				Throughout this proof let $\mathscr{a}\in[-\pi,a]$, $\mathscr{b}\in[0,1]$, $\mathscr{c}\in\R$ satisfy
				\begin{equation}
					\label{eq:param:introd:scr}
					\mathscr{a}=a-\tfrac{(N-1)\pi}{N+1},
					\qquad 
					\mathscr{b}=\tfrac{N-1}{N+1},
					\qandq
					\mathscr{c}=\tfrac{c}{C}
				\end{equation}
				and let $\g_1\in \prb{(\R^{2 \times 1} \times \R^{2}) \times (\R^{1 \times 2} \times \R^1) }\subseteq\ANNSET$ and $\g_3\in \prb{(\R^{1 \times 1} \times \R^{1}) \times (\R^{1 \times 1} \times \R^1) }\subseteq\ANNSET$ satisfy
				\begin{equation}
					\label{pre:post:network:def}
					\begin{split}
						\g_1= \pr*{ \!\pr*{\!\begin{pmatrix}
									\nicefrac{1}{2}\\
									\nicefrac{-1}{2}
								\end{pmatrix},
								\begin{pmatrix}
									\nicefrac{-\mathscr{a}}{2}\\
									\nicefrac{\mathscr{a}}{2}
								\end{pmatrix}\! },
							\prbb{	\begin{pmatrix}
									2& 2
								\end{pmatrix}, 
								0 } \!
						}
						\qandq
						\g_3 = \prb{ \pr{
								1,
								-\mathscr{b}},
							\pr{
								\mathscr{c}, 
								0 } 
						}
					\end{split}
				\end{equation}
				\cfload. \Nobs that \cref{lem:edgy:sin} (applied with
				$B \curvearrowleft  2^{n}\pi$,
				$n \curvearrowleft  n-1$
				in the notation of \cref{lem:edgy:sin}) \proves\ that there exists $\g_2\in\ANNSET$ which satisfies that
				\begin{enumerate}[(I)]
					\item{
						\label{prop:eq:edgy:sin1}
						it holds for all $x\in[0,2^n\pi]$ that $(\realisation(\g_2))(x)=\edgy(2^{n}(2^n\pi)^{-1}x)=\edgy(\pi^{-1}x)$,
					}
					\item{
						\label{prop:eq:edgy:sin2}
						it holds for all $x\in\pr{\R\setminus[0,2^n\pi]}$ that $(\realisation(\g_2))(x)=0$,
					}
					\item{
						\label{prop:eq:edgy:sin3}
						it holds that $\lengthANN(\g_2)=(n-1)+2$,
					}
					\item{
						\label{prop:eq:edgy:sin4}
						it holds that $\dims(\g_2)=(1,2,4,\ldots,4,1)\in\N^{n+2}$, and
					}
					\item{
						\label{prop:eq:edgy:sin6}
						it holds that $\size(\g_2)\leq 2$, $\insize(\g_2)= 1$, and $\outsize(\g_2)\leq 2$
					}
				\end{enumerate}
				\cfload. 
				Next let $\f\in\ANNSET$ satisfy
				\begin{equation}
					\label{eq:edgy:fdef}
					\f=\compANN{\g_3}{\ReLUidANN{1}}\bullet\g_2\bullet\ReLUidANN{1}\bullet\g_1
				\end{equation}
				\cfload. \Nobs that \eqref{pre:post:network:def} and the fact that $\vass{\mathscr{a}}\leq \pi$ \prove\ that for all $x\in\R$ it holds that $\size(\g_1)\leq 2$, $\dims(\g_1)=(1,2,1)$, and
				\begin{equation}
					\label{eq:prop:pre:net:real}
					(\realisation(\g_1))(x)
					=2\RELU\pr*{\tfrac{x}{2}-\nicefrac{\mathscr{a}}{2}}+
					2\RELU\pr*{\tfrac{-x}{2}+\nicefrac{\mathscr{a}}{2}}
					=\vass{x-\mathscr{a}}.
				\end{equation}
				\Moreover \eqref{pre:post:network:def} \proves\ that for all $x\in\R$ it holds that 
				\begin{equation}
					\label{eq:prop:post:net:real}
					\size(\g_3)\leq 1,\qquad
					\dims(\g_3)=(1,1,1),\qandq
					(\realisation(\g_3))(x)=\mathscr{c}\RELU(x-\mathscr{b}).
				\end{equation}
				Combining this and \eqref{eq:prop:pre:net:real} with \cref{Lemma:PropertiesOfCompositions_n2} and \cref{Prop:identity_representation} \proves\ that for all $x\in\R$ it holds that
				\begin{equation}
					\label{eq:realisation:composition1}
					(\realisation(\f))(x)=\pr*{\PR{\realisation(\g_3)}\circ\PR{\realisation(\g_2)}\circ\PR{\realisation(\g_1)}}(x)
					=\mathscr{c}\RELU\prb{(\realisation(\g_2))(\vass{x-\mathscr{a}})-\mathscr{b}}.
				\end{equation}
				This, the fact that $\forall\,x\in\R\colon\edgy(x)=\edgy(-x)$, and \cref{prop:eq:edgy:sin1} \prove\ that for all $x\in[-2^n\pi+\mathscr{a},2^n\pi+\mathscr{a}]$ it holds that
				\begin{equation}
					\label{eq:realisation:composition2}
					\begin{split}
						(\realisation(\f))(x)
						&=\mathscr{c}\RELU\prb{(\realisation(\g_2))(\vass{x-\mathscr{a}})-\mathscr{b}}
						\\&=\mathscr{c}\RELU\prb{\edgy(\pi^{-1}\vass{x-\mathscr{a}})-\mathscr{b}}
						\\&=\mathscr{c}\RELU\prb{\edgy(\pi^{-1}x-\pi^{-1}\mathscr{a})-\mathscr{b}}
						\\&=\mathscr{c}\RELU\prb{\edgy(\pi^{-1}x-(\pi^{-1}a-\mathscr{b}))-\mathscr{b}}
						\\&=\mathscr{c}\RELU\prb{\edgy(\pi^{-1}(x-a)+\mathscr{b})-\mathscr{b}}
						\\&=\prb{\tfrac{2}{(N+1)C}}\prb{\tfrac{(N+1)c}{2}}\RELU\prb{\edgy(\pi^{-1}(x-a)+\mathscr{b})-\mathscr{b}}.
					\end{split}
				\end{equation}
				\Moreover \cref{prop:eq:edgy:sin2}, \eqref{eq:realisation:composition1}, and
				the fact that $\mathscr{b}\geq0$ \prove\ that for all $x\in\pr{\R\setminus[-2^n\pi+\mathscr{a},2^n\pi+\mathscr{a}]}$ it holds that
				\begin{equation}
					\label{eq:realisation:composition3}
					\begin{split}
						(\realisation(\f))(x)
						=\mathscr{c}\RELU\prb{(\realisation(\g_2))(\vass{x-\mathscr{a}})-\mathscr{b}}
						&=\mathscr{c}\RELU\prb{-\mathscr{b}}=0.
					\end{split}
				\end{equation}
				Combining this, \eqref{eq:setup:hatfunctions0}, and \eqref{eq:realisation:composition2} with \cref{single:hat:network:prep} (applied with
				$n \curvearrowleft  n$,
				$N \curvearrowleft  N$,
				$a \curvearrowleft  a$,
				$b \curvearrowleft  b$,
				$c \curvearrowleft  c$,
				$f \curvearrowleft  f$,
				$\mathscr{a} \curvearrowleft  \mathscr{a}$,
				$\mathscr{b} \curvearrowleft  \mathscr{b}$,
				$\mathscr{c} \curvearrowleft  \tfrac{(N+1)c}{2}$,
				$g \curvearrowleft  \prb{\tfrac{C(N+1)}{2}}\realisation(\f)$
				in the notation of \cref{single:hat:network:prep})
				\proves\ that for all $x\in\R$ it holds that
				\begin{equation}
					\label{eq:function:for:composed:netw:for:single:hat}
					(\realisation(\f))(x)
					=\prb{\tfrac{2}{C(N+1)}}\prb{\tfrac{C(N+1)}{2}}(\realisation(\f))(x)
					=\prb{\tfrac{2}{C(N+1)}}f(x).
				\end{equation}
				\Nobs that \cref{Lemma:PropertiesOfCompositions_n2}, 
				\cref{Prop:identity_representation}, 
				\eqref{pre:post:network:def}, and
				\cref{prop:eq:edgy:sin3} \prove\ that
				\begin{equation}
					\label{eq:length:for:composed:netw:for:single:hat}
					\lengthANN(\f)=\lengthANN(\g_3)+\lengthANN(\g_2)+\lengthANN(\g_1)+2\lengthANN(\ReLUidANN{1})-4=2+(n+1)+2=n+5.
				\end{equation}
				This, \eqref{pre:post:network:def}, 
				\eqref{eq:edgy:fdef}, 
				\eqref{eq:prop:pre:net:real}, 
				\eqref{eq:prop:post:net:real}, 
				\cref{prop:eq:edgy:sin4}, 
				\dimsbene, 
				and \cref{Prop:identity_representation} \prove\ that for all $k\in\N_0\cap[0,\lengthANN(\f)]$ it holds that
				\begin{equation}
					\label{eq:dims:for:composed:netw:for:single:hat}
					\singledims_k(\f)=
					\begin{cases}
						1 & \colon k\in\{0,n+4,n+5\}\\
						2 & \colon k \in \{1,2,3,n+3\}\\
						4 & \colon k\in\N\cap(3,n+3).\\
					\end{cases}
				\end{equation}
				\Moreover \eqref{eq:param:introd:scr}, 
				\eqref{pre:post:network:def}, 
				\eqref{eq:edgy:fdef}, 
				\eqref{eq:prop:post:net:real}, 
				\cref{prop:eq:edgy:sin6}, and \cref{Prop:identity_representation:prop2} \prove\ that
				\begin{equation}
					\begin{split}
						\size(\f)&\leq \max\{\size(\g_3),\size(\g_2),\size(\g_1)\}\leq \max\{1,2,2\}=2.
					\end{split}
				\end{equation}
				Combining this,
				\eqref{eq:function:for:composed:netw:for:single:hat},
				\eqref{eq:length:for:composed:netw:for:single:hat}, and
 				\eqref{eq:dims:for:composed:netw:for:single:hat} 
				 \proves[ei]\
				\cref{eq:edgy:base1.2,eq:edgy:base2.1,eq:edgy:base3.1,eq:edgy:base5.1}. 
				The proof of \cref{single:hat:network} is thus complete.
			\end{proof}

			\cfclear
			\begin{athm}{lemma}{sum:hat:networks:eps}
				Let $n,N\in\N\setminus\{1\}$, $C\in[1,\infty)$, $\shft\in[-2,2]$,
				let $\edgy\colon\R\to\R$ satisfy for all $k\in\Z$, $x\in[2k,2k+2)$ that $\edgy(x)=1-\vass{x-2k-1}$,
				for every $j\in\{1,2,\ldots,N\}$ let $a_j,b_j\in[0,2\pi]$, $c_j\in[-C,C]$ satisfy $b_j=a_j+\tfrac{4\pi}{N+1}$,
				and for every $j\in\{1,2,\ldots,N\}$ let $f_j\colon\R\to\R$ satisfy for all $x\in[0,2\pi)$, $y\in(-\infty,-2^n\pi)\cup[2^n\pi,\infty)$, $k\in\Z\cap[-2^{n-1},2^{n-1})$ that
				
				\begin{equation}
					\label{eq:setup:hatfunctions:eps}
					f_j(y)=0
					\qandq
					f_j(x+2k\pi)=f_j(x)=
					\begin{cases}
						c_j\edgy\pr*{\frac{(x-a_j)(N+1)}{2\pi}} &\colon x\in [a_j,b_j]\\
						0 &\colon x\notin [a_j,b_j].
					\end{cases}
				\end{equation}
				Then there exists $\f\in\ANNSET$ such that 
				\begin{enumerate}[(i)]
					\item{
						\label{sum:hat:networks1:eps}
						it holds for all $x\in\R$ that
						$
						(\realisation(\f))(x)=\shft+\sum_{j=1}^N f_j(x)$,
					}
					\item{
						\label{sum:hat:networks3:eps}
						it holds that $\lengthANN(\f)\leq n+\log_2(C)+9$,
					}
					\item{
						\label{sum:hat:networks4:eps}
						it holds that $\singledims_0(\f)=\singledims_{\lengthANN(\f)}(\f)=1$, $\singledims_1(\f)=2N$, and $\singledims_{\hidlengthANN(\f)}(\f)=2$,
					}
					\item{
						\label{sum:hat:networks5:eps}
						it holds that $\param(\f)\leq \pr*{14+18n+5\log_2(C)}N^2$, and
					}
					\item{
						\label{sum:hat:networks6:eps}
						it holds that $\size(\f)\leq 2$
					}
				\end{enumerate}
				\cfout.
			\end{athm}
			
			\begin{proof}[Proof of \cref{sum:hat:networks:eps}]
				Throughout this proof 
				let $\mathscr{n},\wdt\in\N$, $\beta\in(0,2]$ satisfy 
				\begin{equation}
					\label{def:summing:netw:post:eps0}
					\mathscr{n}=\max\{1,\ceil{\log_2(C)}\},
					\qquad
					\wdt=\ceil*{\pr*{\tfrac{N+1}{2}}^{\tfrac{1}{\mathscr{n}}}},
					\qandq
					\beta=\wdt^{-1}\pr*{\tfrac{C(N+1)}{2}}^{\tfrac{1}{\mathscr{n}}},
				\end{equation}
				and
				let $\g_1,\g_2,\g_4\in\ANNSET$ satisfy $\g_4=((1),\shft)\in(\R^{1 \times 1} \times \R^1)$,
				\begin{equation}
					\label{def:summing:netw:post:eps}
					\begin{split}
						&\g_1 =
						\pr*{\!
							\begin{pmatrix} 1\\1\\\vdots\\1 \end{pmatrix},\!
							\,0
							\!}
						\in (\R^{N \times 1} \times \R^N),
						\quad\text{and}\quad
						\g_2 = 
						\prb{
							\begin{pmatrix} 1 & 1 & \cdots & 1 \end{pmatrix},\!
							\,0
						}
						\in (\R^{1 \times N} \times \R^1)
					\end{split}
				\end{equation}
				%
				\cfload. \Nobs that \cref{single:hat:network} (applied with
				$n \curvearrowleft  n$,
				$N \curvearrowleft  N$,
				$C \curvearrowleft  C$,
				$a \curvearrowleft  a_j$,
				$b \curvearrowleft  b_j$,
				$c \curvearrowleft  c_j$,
				$f \curvearrowleft  f_j$ for $j\in\{1,2,\ldots,N\}$
				in the notation of \cref{single:hat:network}) \proves\ that there exist $\f_1,\f_2,\ldots,\f_N\in\ANNSET$ such that
				\begin{enumerate}[(I)]
					\item{
						\label{eq:edgy:base1.2:prop:eps}
						it holds for all $j\in\{1,2,\ldots,N\}$, $x\in\R$ that
						$
						(\realisation(\f_j))(x)= \prb{\tfrac{2}{C(N+1)}}f_j(x)
						$,
					}
					\item{
						\label{eq:edgy:base2.1:prop:eps}
						it holds for all $j\in\{1,2,\ldots,N\}$ that $\lengthANN(\f_j)=n+5$,
					}
					\item{
						\label{eq:edgy:base3.1:prop:eps}
						it holds for all $j\in\{1,2,\ldots,N\}$, $k\in\N_0\cap[0,\lengthANN(\f_j)]$ that
						\begin{equation}
							\singledims_k(\f_j)=
							\begin{cases}
								1 & \colon k\in\{0,n+4,n+5\}\\
								2 & \colon k \in \{1,2,3,n+3\}\\
								4 & \colon k\in\N\cap(3,n+3),\\
							\end{cases}
						\end{equation}
						and
					}
					\item{
						\label{eq:edgy:base5.1:prop:eps}
						it holds for all $j\in\{1,2,\ldots,N\}$ that $\size(\f_j)\leq 2$
					}
				\end{enumerate}
				\cfload. 
				\Moreover \cref{scaling:networks} (applied with
				$\beta \curvearrowleft  \beta$,
				$\wdt \curvearrowleft  \wdt$,
				$n \curvearrowleft  \mathscr{n}$
				in the notation of \cref{scaling:networks}) \proves\ that there exists $\g_3\in\ANNSET$ which satisfies that
				\begin{enumerate}[(A)]
					\item{
						\label{it:for:sum:gen:scaling:realisation:eps}
						it holds for all $x\in\R$ that $(\realisation(\g_3))(x)=\tfrac{C(N+1) x}{2}$,
					}
					\item{
						\label{it:for:sum:gen:scaling:dims:eps}
						it holds that $\dims(\g_3)=(1,2\wdt,2\wdt,\ldots,2\wdt,1)\in\N^{\mathscr{n}+2}$,
					}
					\item{
						\label{it:for:sum:gen:scaling:size:eps}
						it holds that $\insize(\g_3)\leq 1$, $\outsize(\g_3)\leq 2$, and $\size(\g_3)\leq 2$, and
					}
					\item{
						\label{it:for:sum:gen:scaling:param:eps}
						it holds that $\param(\g_3)=(4\mathscr{n}-4)\wdt^2+(2\mathscr{n}+4)\wdt+1$.
					}
				\end{enumerate}
				Next let $\h\in\ANNSET$ satisfy
				\begin{equation}
					\label{eq:edgy:composition:fdef:eps}
					\h=\g_4\bullet\ReLUidANN{1}\bullet\g_3\bullet\g_2\bullet\compANN{\ReLUidANN{N}}{\pr*{\parallelizationSpecial_{N}(\f_1,\f_2,\ldots,\mathscr{f}_N)}}\bullet\ReLUidANN{N}\bullet\g_1
				\end{equation}
				\cfload. \Nobs that 
				\eqref{def:summing:netw:post:eps}, 
				\eqref{eq:edgy:composition:fdef:eps}, 
				\cref{Lemma:PropertiesOfCompositions_n2}, and 
				\cref{Lemma:PropertiesOfParallelizationEqualLength} \prove\ that for all $x\in\R$ it holds that
				\begin{equation}
					\begin{split}
						\label{eq:edgy:composition:function:eps}
						&(\realisation(\h))(x)
						\\&=(\PR{\realisation(\g_4)}\circ\PR{\realisation(\g_3)}\circ\PR{\realisation(\g_2)}\circ\PR{\realisation(\parallelizationSpecial_{N}(\f_1,\f_2,\ldots,\mathscr{f}_N))})((\realisation(\g_1))(x))
						\\&=(\PR{\realisation(\g_4)}\circ\PR{\realisation(\g_3)}\circ\PR{\realisation(\g_2)}\circ\PR{\realisation(\parallelizationSpecial_{N}(\f_1,\f_2,\ldots,\mathscr{f}_N))})(x,x,\ldots,x)
						\\&=(\PR{\realisation(\g_4)}\circ\PR{\realisation(\g_3)}\circ\PR{\realisation(\g_2)})\prb{\pr*{\realisation\pr*{\parallelizationSpecial_{N}(\f_1,\f_2,\ldots,\mathscr{f}_N)}}(x,x,\ldots,x)}
						\\&=(\PR{\realisation(\g_4)}\circ\PR{\realisation(\g_3)}\circ\PR{\realisation(\g_2)})\prb{(\realisation(\f_1))(x),(\realisation(\f_2))(x),\ldots,(\realisation(\f_N))(x)}
						\\&=(\PR{\realisation(\g_4)}\circ\PR{\realisation(\g_3)})\pr*{\textstyle\sum_{j=1}^N(\realisation(\f_j))(x)}
						\\&=(\realisation(\g_4))\pr*{\prb{\tfrac{C(N+1)}{2}}\textstyle\sum_{j=1}^N \prb{\tfrac{2}{C(N+1)}}f_j(x)}
						=\shft+\textstyle\sum_{j=1}^N f_j(x)
						\ifnocf.
					\end{split}
				\end{equation}
				\cfload[.]\Nobs that 
				\eqref{def:summing:netw:post:eps}, 
				\eqref{eq:edgy:composition:fdef:eps}, 
				\cref{eq:edgy:base2.1:prop:eps}, 
				\cref{it:for:sum:gen:scaling:dims:eps}, 
				\cref{Lemma:PropertiesOfCompositions_n2}, 
				and \cref{Lemma:PropertiesOfParallelizationEqualLength} \prove\ that
				\begin{equation}
					\label{eq:edgy:composition:length:eps}
					\begin{split}
						\lengthANN\pr*{\h}
						&=\lengthANN(\g_4)+\lengthANN(\g_3)+\lengthANN(\g_2)+\lengthANN\pr*{\parallelizationSpecial_{N}(\f_1,\f_2,\ldots,\mathscr{f}_N)}+\lengthANN(\g_1)+3\lengthANN(\ReLUidANN{1})-7
						\\&=1+\lengthANN(\g_3)+1+\lengthANN\pr*{\f_1}+1+6-7
						\\&=(n+5)+(\mathscr{n}+1)+2
						\leq n+\log_2(C)+9.
					\end{split}
				\end{equation}
				\Nobs that \cref{eq:edgy:base2.1:prop:eps}, \cref{eq:edgy:base3.1:prop:eps}, and \cref{Lemma:PropertiesOfParallelizationEqualLength} \prove\ that for all $k\in\N_0\cap[0,n+5]$ it holds that
				\begin{equation}
					\singledims_k\pr*{\parallelizationSpecial_{N}(\f_1,\f_2,\ldots,\mathscr{f}_N)}=
					\begin{cases}
						N & \colon k\in\{0,n+4,n+5\}\\
						2N & \colon k \in \{1,2,3,n+3\}\\
						4N & \colon k\in\N\cap(3,n+3).\\
					\end{cases}
				\end{equation}
				Combining this,  
				\eqref{eq:edgy:composition:fdef:eps}, 
				\eqref{eq:edgy:composition:length:eps}, 
				\cref{it:for:sum:gen:scaling:dims:eps},
				\dimsbene, and
				\cref{Prop:identity_representation}
				with the fact that $\dims(\g_1)=(1,N)$, $\dims(\g_2)=(N,1)$, and $\dims(\g_4)=(1,1)$ \proves\ that for all $k\in\N_0\cap[0,n+\mathscr{n}+8]$ it holds that
				\begin{equation}
					\label{eq:edgy:composition:dims:eps}
					\singledims_k\pr*{\h}=
					\begin{cases}
						1 & \colon k\in\{0,n+\mathscr{n}+8\}\\
						2N & \colon k \in \{1,2,3,4,n+4,n+6\}\\
						4N & \colon k\in\N\cap(4,n+4)\\
						N & \colon k\in\{n+5\}\\
						2B & \colon k\in\N\cap(n+6,n+\mathscr{n}+7)\\
						2 & \colon k=n+\mathscr{n}+7
						.\\
					\end{cases}
				\end{equation}
				This and the fact that $B\leq \ceil*{\tfrac{N+1}{2}}\leq N\geq 2$ \prove\ that
				\begin{equation}
					\label{eq:edgy:composition:param:eps}
					\begin{split}
						\param\pr*{\h}
						&=\sum_{k = 1}^{n+\mathscr{n}+8} \singledims_k\pr*{\h}\pr*{\singledims_{k-1}\pr*{\h} + 1}
						\\&=2N(1+1)+3(2N(2N+1))+4N(2N+1)+\PR*{\sum_{k = 6}^{n+3} 4N\pr*{4N + 1}}
						\\&\quad +2N(4N+1)+N(2N+1)+2N(N+1)+2B(2N+1)
						\\&\quad+\PR*{\sum_{k = 3}^{\mathscr{n}+1} 2B\pr*{2B + 1}}
						+2(2B+1)+1(2+1)
						\\&=32N^2+19N+4BN+6B+5+(n-2)(16N^2+4N)+(\mathscr{n}-1)(4B^2+2B)
						\\&\leq36N^2+25N+5+(n-2)(16N^2+4N)+(\mathscr{n}-1)(4N^2+2N)
						\\&=15N+5+n(16N^2+4N)+\mathscr{n}(4N^2+2N)
						\\&\leq9N^2+18nN^2+5\mathscr{n}N^2
						\\&\leq \pr*{14+18n+5\log_2(C)}N^2
						.
					\end{split}
				\end{equation}
				\Moreover  
				\eqref{def:summing:netw:post:eps}, 
				\eqref{eq:edgy:composition:fdef:eps}, 
				\cref{eq:edgy:base5.1:prop:eps},
				\cref{it:for:sum:gen:scaling:size:eps},  
				\cref{lem:sizecomp}, 
				\cref{Lemma:ParallelizationElementary}, and 
				\cref{Prop:identity_representation:prop2} \prove\ that
				\begin{equation}
					\begin{split}
						&\size(\h)
						\\&\leq \max\{\size(\g_4),\size(\g_3\bullet\g_2),\size\pr*{\parallelizationSpecial_{N}(\f_1,\f_2,\ldots,\mathscr{f}_N)},\size\pr*{\g_1}\}
						\\&\leq \max\{\size(\g_4),\size(\g_3),\size(\g_2),\insize(\g_3)(\outsize(\g_2)+1),\size\pr*{\mathscr{f}_1},\size\pr*{\mathscr{f}_2},\ldots,\size\pr*{\mathscr{f}_N},\size\pr*{\g_1}\}
						\\&\leq\max\{\vass{\shft},2,1,2,2,2,\ldots,2,1\}=2.
					\end{split}
				\end{equation}
				Combining this,
				\eqref{eq:edgy:composition:function:eps},
				\eqref{eq:edgy:composition:length:eps},
				\eqref{eq:edgy:composition:dims:eps}, and
				\eqref{eq:edgy:composition:param:eps}
				 \proves[ei]\ 
				\cref{sum:hat:networks1:eps,sum:hat:networks3:eps,sum:hat:networks4:eps,sum:hat:networks5:eps,sum:hat:networks6:eps}. 
				The proof of \cref{sum:hat:networks:eps} is thus complete.
				\cfload
			\end{proof}
			
			\cfclear
			\begin{athm}{lemma}{sin:approximator}
				Let $\varepsilon\in(0,1)$, $n,N\in\N$ satisfy $\eps (N+1) \geq 2\pi$, 
				let $\edgy\colon\R\to\R$ satisfy for all $k\in\Z$, $x\in[2k,2k+2)$ that $\edgy(x)=1-\vass{x-2k-1}$,
				for every $j\in\N_0\cap[0,N+1]$
				let $c_j=\frac{2j\pi}{N+1}
				$,
				let $\tgt\colon\R\to\R$ satisfy for all $x,y\in\R$, $k\in\Z$ that \begin{equation}
					\vass{\tgt(0)}\leq2,
					\qquad
					\tgt(x+2k\pi)=\tgt(x),
					\qandq
					\vass{\tgt(x)-\tgt(y)}\leq\vass{x-y},
				\end{equation}
				for every $j\in\{1,2,\ldots,N\}$ let $f_j\colon\R\to\R$ satisfy for all $x\in[0,2\pi)$, $y\in(-\infty,-2^n\pi)\cup[2^n\pi,\infty)$,  $k\in\Z\cap[-2^{n-1},2^{n-1})$ that $f_j(y)=0$ and
				\begin{equation}
					\label{setup:hat:functions}
					f_j(x+2k\pi)=f_j(x)=
					\begin{cases}
						(\tgt(c_j)-\tgt(0))\edgy\prb{\frac{(x-c_{j-1})(N+1)}{2\pi}} &\colon x\in [c_{j-1},c_{j+1}]\\
						0 &\colon x\notin [c_{j-1},c_{j+1}],
					\end{cases}
				\end{equation}
				and let $F\colon\R\to\R$ satisfy for all $x\in\R$ that
				\begin{equation}
					\label{def:sum:of:base:hats}
					F(x)=\tgt(0)+\sum_{j=1}^N f_j(x).
				\end{equation}
				\cfload[.]Then
				\begin{enumerate}[(i)]
					\item{
						it holds that
						\label{approx:prec}
						$
						\sup_{x\in[-2^n\pi,2^n\pi]}\vass*{\tgt(x)-F(x)}\leq \eps
						$ and
					}
					\item{
						\label{lipschitz}
						it holds for all $x,y\in\R$ that
						$
						\vass{F(x)-F(y)}\leq \vass{x-y}
						$.
					}
				\end{enumerate}
			\end{athm}
			
			\begin{proof}[Proof of \cref{sin:approximator}]
				\Nobs that \eqref{setup:hat:functions} and \eqref{def:sum:of:base:hats} \prove\ that for all $x\in(2^n \pi,\infty)$ it holds that
				\begin{equation}
					\label{eq:outside:vanishing:tails}
					F(-x)=\tgt(0)=F(x).
				\end{equation}
				\Moreover \eqref{setup:hat:functions} and the fact that $\edgy(0)=0=\edgy(2)$ \prove\ that for all $j\in\{1,2,\ldots,N\}$ it holds that
				\begin{equation}
					\label{eq:welldefined:boundaries}
					f_j(c_{j-1})=(\tgt(c_j)-\tgt(0))\edgy(0)=0=(\tgt(c_j)-\tgt(0))\edgy(2)=f_j(c_{j+1}).
				\end{equation}
				This, \eqref{setup:hat:functions}, \eqref{def:sum:of:base:hats}, and the fact that for all $x\in[0,1]$ it holds that $\edgy(x)=x$ \prove\ that for all $x\in\PRb{c_0,c_1}$ it holds that
				\begin{equation}
					\begin{split}
						\label{linear:in:the:beginning}
						F(x)=\tgt(0)+\sum_{k=1}^N f_k(x)
						=\tgt(0)+f_{1}(x)
						&=\tgt(0)+(\tgt(c_1)-\tgt(0))\edgy\pr*{\frac{(x-c_0)(N+1)}{2\pi}}
						\\&=\tgt(0)+(\tgt(c_1)-\tgt(0))\pr*{\frac{(N+1)x}{2\pi}}.
					\end{split}
				\end{equation}
				\Moreover \eqref{setup:hat:functions}, \eqref{def:sum:of:base:hats}, \eqref{eq:welldefined:boundaries}, and the fact that for all $x\in[1,2]$ it holds that $\edgy(x)=2-x$ \prove\ that for all $x\in\PRb{c_N,c_{N+1}}$ it holds that
				\begin{equation}
					\begin{split}
						\label{linear:in:the:end}
						F(x)=\tgt(0)+\sum_{k=1}^N f_k(x)=&\tgt(0)+f_{N}(x)
						\\=&\tgt(0)+(\tgt(c_N)-\tgt(0))\edgy\pr*{\frac{(x-c_{N-1})(N+1)}{2\pi}}
						\\=&\tgt(0)+(\tgt(c_N)-\tgt(0))\pr*{2-\frac{(x-c_{N-1})(N+1)}{2\pi}}.
					\end{split}
				\end{equation}
				\Moreover \eqref{setup:hat:functions}, \eqref{def:sum:of:base:hats}, \eqref{eq:welldefined:boundaries}, and the fact that for all $x\in[0,1]$ it holds that $\edgy(x)=x$ and $\edgy(x+1)=2-(x+1)$ \prove\ that for all $j\in\N\cap(1,N]$, $x\in\PRb{c_{j-1},c_j}$ it holds that
				\begin{equation}
					\begin{split}
						\label{linear:in:the:middle}
						F(x)
						&=\tgt(0)+\sum_{k=1}^N f_k(x)
						\\&=\tgt(0)+f_{j-1}(x)+f_j(x)
						\\&=\tgt(0)+(\tgt(c_{j-1})-\tgt(0))\edgy\pr*{\frac{(x-c_{j-2})(N+1)}{2\pi}}
						\\&\quad+(\tgt(c_{j})-\tgt(0))\edgy\pr*{\frac{(x-c_{j-1})(N+1)}{2\pi}}
						\\&=\tgt(0)+(\tgt(c_{j-1})-\tgt(0))\pr*{2-\frac{(x-c_{j-2})(N+1)}{2\pi}}
						\\&\quad+(\tgt(c_{j})-\tgt(0))\pr*{\frac{(x-c_{j-1})(N+1)}{2\pi}}.
					\end{split}
				\end{equation}
				Combining this and \eqref{linear:in:the:beginning} with \eqref{linear:in:the:end} \proves\ that for all $j\in\{1,2,\ldots,N+1\}$ it holds that
				\begin{equation}
					\begin{split}
						\label{linear:every:inerval}
						F|_{\PR{c_{j-1},c_j}}\in\Lin\pr*{\PRb{c_{j-1},c_j}}
					\end{split}
				\end{equation}
				\cfload. \Moreover 
				\eqref{setup:hat:functions}, 
				\eqref{def:sum:of:base:hats},
				\eqref{eq:welldefined:boundaries},
				and the fact that $\edgy(1)=1$
				\prove\ that for all $j\in\{1,2,\ldots,N\}$ it holds that
				\begin{equation}
					\label{eq:equal:on:points}
					F(c_{j})=\tgt(0)+\sum_{k=1}^N f_k(c_{j})= \tgt(0)+f_j(c_{j})=\tgt(0)+(\tgt(c_{j})-\tgt(0))\edgy(1)=\tgt(c_{j}).
				\end{equation}
				This, 
				\eqref{linear:in:the:beginning},
				\eqref{linear:in:the:end},
				\eqref{linear:in:the:middle},
				\eqref{linear:every:inerval}, and the fact that $F(0)=\tgt(0)$ and $F(2\pi)=\tgt(0)=\tgt(2\pi)$ \prove\ that for all $j\in\{1,2,\ldots,N+1\}$, $x\in\PRb{c_{j-1},c_j}$ it holds that
				\begin{equation}
					\begin{split}
						F(x)&=\pr*{\frac{(N+1)c_j-(N+1)x}{2\pi}}F(c_{j-1})+\pr*{\frac{(N+1)x-(N+1)c_{j-1}}{2\pi}}F(c_{j})
						\\&=\pr*{\frac{2j\pi-(N+1)x}{2\pi}}F(c_{j-1})+\pr*{\frac{(N+1)x-2(j-1)\pi}{2\pi}}F(c_{j}).
					\end{split}
				\end{equation}
				Combining this and \eqref{eq:equal:on:points} with the assumption that for all $x,y\in\R$ it holds that $\vass{\tgt(x)-\tgt(y)}\leq \vass{x-y}$ and $\eps (N+1) \geq 2\pi$ \proves\ that for all $j\in\{1,2,\ldots,N+1\}$, $x\in\PRb{c_{j-1},c_j}$ it holds that
				\begin{equation}
					\label{eq:approx:prec:on:start}
					\begin{split}
						\vass{\tgt(x)-F(x)}
						&=\vass*{\tgt(x)-\pr*{\tfrac{2j\pi-(N+1)x}{2\pi}}F(c_{j-1})-\pr*{\tfrac{(N+1)x-2(j-1)\pi}{2\pi}}F(c_{j})}
						\\&=\vass*{\pr*{\tfrac{2j\pi-(N+1)x}{2\pi}}\pr*{\tgt(x)-F(c_{j-1})}
							+\pr*{\tfrac{(N+1)x-2(j-1)\pi}{2\pi}}\pr*{\tgt(x)-F(c_{j})}}
						\\&=\vass*{\pr*{\tfrac{2j\pi-(N+1)x}{2\pi}}\pr*{\tgt(x)-\tgt(c_{j-1})}
							+\pr*{\tfrac{(N+1)x-2(j-1)\pi}{2\pi}}\pr*{\tgt(x)-\tgt(c_{j})}}
						\\&\leq\vass*{\pr*{\tfrac{2j\pi-(N+1)x}{2\pi}}\pr*{x-c_{j-1}}}
						+\vass*{\pr*{\tfrac{(N+1)x-2(j-1)\pi}{2\pi}}\pr*{x-c_{j}}}
						\\&\leq \vass*{\pr*{\tfrac{2j\pi-(N+1)x}{2\pi}}\tfrac{2\pi}{N+1}}
						+\vass*{\pr*{\tfrac{(N+1)x-2(j-1)\pi}{2\pi}}\tfrac{2\pi}{N+1}}
						=\tfrac{2\pi}{N+1}\leq \eps.
					\end{split}
				\end{equation}
				\Moreover \eqref{linear:every:inerval}, \eqref{eq:equal:on:points}, and the assumption that for all $j\in\{1,2,\ldots,N\}$, $x,y\in\R$ it holds that $\vass{\tgt(x)-\tgt(y)}\leq \vass{x-y}$ \prove\ that for all $x,y\in[0,2\pi]=[c_0,c_{N+1}]$ it holds that
				\begin{equation}
					\vass{F(x)-F(y)}\leq \vass{x-y} \max_{j\in\{1,2,\ldots,N+1\}}\frac{\vass{\tgt(c_{j})-\tgt(c_{j-1})}}{c_j-c_{j-1}}\leq \vass{x-y}.
				\end{equation}
				Combining this and
				\eqref{eq:outside:vanishing:tails}
				with \eqref{eq:approx:prec:on:start} and the fact that for all $x\in[0,2\pi]$, $k\in\Z\cap[-2^{n-1},2^{n-1})$ it holds that $\tgt(x+2k\pi)=\tgt(x)$ and $F(x+2k\pi)=F(x)$  \proves[ei]\ \cref{approx:prec,lipschitz}.
				The proof of \cref{sin:approximator} is thus complete.
				\cfload
			\end{proof}

			The following result, \cref{propo:sin:approximation} below,
			builds on the linear interpolation construction in
			Beneventano et al.~\cite[Subsection~4.1]{beneventano21}.
			We adapt this construction to the periodic target functions
			considered and establish the stated bounds on the
			absolute values of the parameters and the network architecture.

			\cfclear
			\begin{athm}{lemma}{propo:sin:approximation}
				Let $\varepsilon \in (0,1)$, $N \in \N$ satisfy $\eps (N+1) \geq 2\pi$ and let $\tgt\colon\R\to\R$ satisfy for all $x,y\in\R$, $k\in\Z$ that $\vass{\tgt(0)}\leq 2$, $\tgt(x+2k\pi)=\tgt(x)$, and $\vass{\tgt(x)-\tgt(y)}\leq\vass{x-y}$.
				Then there exists $\mathscr{f} \in \ANNSET$ such that
				\begin{enumerate}[(i)]
					\item {\label{propo:sin_layers:prop} it holds 
						that $\lengthANN(\f)=6$,
					}
					\item {\label{propo:sin_dims:prop}
						it holds that $\dims(\f)=(1,2,2N+3,2,2,2,1)$,}
					\item{\label{propo:sin_lipschitz:prop} it holds for all
						$x,y \in \R$ that 
						$\vass[\big]{\pr{\realisation(\f)} (x) - \pr{\realisation(\f)} (y)} \leq \vass{x-y}$,}
					\item{\label{propo:sin_approx:prop} it holds that
						$\sup_{x\in[-2\pi,2\pi]}\vass[\big]{\tgt(x) - \pr{\realisation(\f)} (x)} \leq \varepsilon$,}
					\item{\label{propo:sin_start:prop} it holds for all $x\in[2\pi,\infty)$ that 
						$(\realisation(\f)) (-x) =(\realisation(\f)) (x) =\tgt(0)$,}
					\item{\label{propo:sin_size:prop} it holds 
						that
						$\size(\mathscr{f}) \leq 2
						$, and}
					\item{\label{propo:sin_cost:prop} it holds 
						that
						$\param(\mathscr{f}) \leq  4N^2$}
				\end{enumerate}
				\cfout. 
		\end{athm}

		\begin{proof}[Proof of \cref{propo:sin:approximation}]
			Throughout this proof 
			let $M\in\N$, $d\in(0,\eps]$ satisfy $M=2N+2$ and $d=\tfrac{4\pi}{M}$,
			let $\xi_0,\xi_1,\ldots,\xi_{M}\in[-2\pi,2\pi]$, $\alpha_0,\alpha_1,\ldots,\alpha_M\in[-2,2]$ satisfy for all $k\in\N_0\cap[0,M]$ that
			\begin{equation}
				\label{eq:prop:flat:approx:param:def}
				\xi_k=kd-2\pi
				\qandq
				\alpha_k=\frac{\tgt(\xi_{\min\{k+1,M\}})-\tgt(\xi_{k})}{d}-\frac{\tgt(\xi_{k})-\tgt(\xi_{\max\{k-1,0\}})}{d},
			\end{equation}
			and let $\g_1\in\prb{(\R^{(M+1) \times 1} \times \R^{M+1}) \times (\R^{1 \times (M+1)} \times \R^1) }\subseteq\ANNSET$ satisfy
			\begin{equation}
				\label{eq:prop:flat:approx:fdef}
				\g_1 =
				\pr*{\!
					\pr*{\!
						\begin{pmatrix} 4^{-1}\\4^{-1}\\\vdots\\4^{-1} \end{pmatrix},\!
						\,
						\begin{pmatrix} -4^{-1}\xi_0 \\-4^{-1}\xi_1 \\\vdots\\-4^{-1}\xi_M \end{pmatrix}
						\!}
					,
					\pr*{
						\begin{pmatrix} \alpha_0 & \alpha_1 & \cdots & \alpha_M \end{pmatrix},\!
						\,
						4^{-1}\tgt(0) 
					}
					\!}
			\end{equation}
			\cfload.
			\Nobs that \cref{gen:scaling:networks} (applied with
			$\beta \curvearrowleft  4$,
			$L \curvearrowleft  2$
			in the notation of \cref{gen:scaling:networks}) \proves\ that there exists $\g_2\in\ANNSET$ which satisfies that
			\begin{enumerate}[(I)]
				\item{
					\label{gen:scaling:realisation:flat:approx}
					it holds for all $x\in\R$ that $(\realisation(\g_2))(x)=4 x$,
				}
				\item{
					\label{gen:scaling:dims:flat:approx}
					it holds that $\dims(\g_2)=(1,2,2,1)\in\N^{4}$, and
				}
				\item{
					\label{gen:scaling:size:flat:approx}
					it holds that $\size(\g_2)\leq 2$
				}
			\end{enumerate}
			\cfload.
			Next let $\f\in\ANNSET$ satisfy
			\begin{equation}
				\label{eq:flat:approx:netw:def}
				\f=\compANN{\g_2}{\ReLUidANN{1}}\bullet\g_1\bullet\ReLUidANN{1}
			\end{equation}
			\cfload[.]\Nobs that
			\cref{Lemma:PropertiesOfCompositions_n2},
			\cref{Prop:identity_representation}, and
			\eqref{eq:flat:approx:netw:def} \prove\ that
			\begin{equation}
				\label{eq:flat:approx:length}
				\lengthANN(\f)=\lengthANN(\g_1)+\lengthANN(\g_2)+2\lengthANN(\ReLUidANN{1})-3=2+3+4-3=6.
			\end{equation}
			\Nobs that 
			\eqref{eq:prop:flat:approx:fdef} \proves\ that for all $x\in\R$ it holds that
			\begin{equation}
				\label{eq:prop:flat:approx:realisation}
				(\realisation(\g_1))(x)
				=\tfrac{\tgt(0)}{4}+\sum_{k=0}^M \alpha_k\RELU\pr*{\tfrac{x}{4}-\tfrac{\xi_k}{4}}
				=4^{-1}\PR*{\tgt(0)+\sum_{k=0}^M \alpha_k\RELU\pr*{x-\xi_k}}
			\end{equation}
			\cfload.
			This and \eqref{eq:prop:flat:approx:param:def} \prove\ that for all $x\in(-\infty,\xi_0]$ it holds that
			\begin{equation}
				\label{eq:prop:flat:approx:lower:part}
				(\realisation(\g_1))(x)=\frac{\tgt(0)}{4}=\frac{\tgt(\xi_0)}{4}
				\qandq
				(\realisation(\g_1))(\xi_1)=\frac{\tgt(0)+ \alpha_0\RELU\pr*{\xi_1-\xi_0}}{4}=\frac{\tgt(\xi_1)}{4}.
			\end{equation}
			\Moreover
			\eqref{eq:prop:flat:approx:param:def} and
			\eqref{eq:prop:flat:approx:realisation} \prove[de]\ that for all $k\in\N_0\cap[0,M)$ with $(\realisation(\g_1))(\xi_k)=\frac{\tgt(\xi_k)}{4}$ it holds that
			\begin{equation}
				\label{eq:prop:flat:approx:ind:step}
				\begin{split}
					(\realisation(\g_1))(\xi_{k+1})
					&=4^{-1}\PR*{\tgt(0)+\sum_{j=0}^{k} \alpha_j\RELU\pr*{\xi_{k+1}-\xi_j}}
					\\&=4^{-1}\PR*{\tgt(0)+\PR*{\sum_{j=0}^{k} \alpha_j\RELU\pr*{\xi_{k}-\xi_j}}+d\PR*{\sum_{j=0}^{k} \alpha_j}}
					\\&=(\realisation(\g_1))(\xi_{k})+4^{-1}\PR*{\tgt(\xi_{k+1})-\tgt(\xi_{k})}
					\\&=\frac{\tgt(\xi_{k+1})}{4}.
				\end{split}
			\end{equation}
			\Moreover 
			\eqref{eq:prop:flat:approx:param:def} and
			\eqref{eq:prop:flat:approx:realisation} \prove\ that for all $x\in[\xi_M,\infty)$ it holds that
			\begin{equation}
				\label{eq:prop:flat:approx:upper:part}
				\begin{split}
					(\realisation(\g_1))(x)
					&=4^{-1}\PR*{\tgt(0)+\sum_{j=0}^{M} \alpha_j\RELU\pr*{x-\xi_j}}
					\\&=4^{-1}\PR*{\tgt(0)+\PR*{\sum_{j=0}^{M} \alpha_j\RELU\pr*{\xi_{M}-\xi_j}}+\pr*{x-\xi_{M}}\PR*{\sum_{j=0}^{M} \alpha_j}}
					\\&=(\realisation(\g_1))(\xi_{M}).
				\end{split}
			\end{equation}
			Combining this,
			\eqref{eq:prop:flat:approx:lower:part}, and
			\eqref{eq:prop:flat:approx:ind:step}
			with induction \proves\ that for all $k\in\N_0\cap[0,M]$, $x\in(-\infty,\xi_0]$, $y\in[\xi_M,\infty)$ it holds that
			\begin{equation}
				\label{eq:flat:approx:unscaled}
				(\realisation(\g_1))(\xi_k)=\frac{\tgt(\xi_k)}{4},
				\qquad
				(\realisation(\g_1))(x)=\frac{\tgt(\xi_0)}{4},
				\qandq
				(\realisation(\g_1))(y)=\frac{\tgt(\xi_M)}{4}.
			\end{equation}
			\Nobs that 
			\cref{Lemma:PropertiesOfCompositions_n2},
			\cref{Prop:identity_representation},
			\eqref{eq:prop:flat:approx:fdef}, and
			\eqref{eq:flat:approx:netw:def}
			\prove\ that for all $x\in\R$ it holds that
			\begin{equation}
				\label{eq:prop:flat:approx:realisation2}
				(\realisation(\f))(x)
				=(\realisation(\g_2))\prb{\realisation(\g_1)(x)}
				=4(\realisation(\g_1))(x)
			\end{equation}
			This and \eqref{eq:prop:flat:approx:realisation} \prove\ that for all $k\in\{1,2,\ldots,M\}$, $x\in[\xi_{k-1},\xi_k]$ it holds that
			\begin{equation}
				\begin{split}
					(\realisation(\f))(x)
					&=\tgt(0)+\sum_{j=0}^{M} \alpha_j\RELU\pr*{x-\xi_j}
					=\tgt(0)+\sum_{j=0}^{k-1} \alpha_j\pr*{x-\xi_j}.
				\end{split}
			\end{equation}
			This,
			\eqref{eq:flat:approx:unscaled},
			\eqref{eq:prop:flat:approx:realisation}, and 
			the fact that $\tgt(\xi_0)=\tgt(0)=\tgt(\xi_M)$ \prove\ that for all $x\in[2\pi,\infty)$, $k\in\{1,2,\ldots,M\}$ it holds that
			\begin{equation}
				\label{eq:flat:approx:linearity:and:outside}
				\realisation(\f)|_{[\xi_{k-1},\xi_k]}\in\Lin([\xi_{k-1},\xi_k])
				\qandq
				(\realisation(\f))(x)=\tgt(0)=(\realisation(\f))(-x)
			\end{equation}
			\cfload. Combining this and
			\eqref{eq:flat:approx:unscaled} with 
			\eqref{eq:prop:flat:approx:realisation} \proves\ that for all $k\in\{1,2,\ldots,M\}$, $x\in[\xi_{k-1},\xi_k]$ it holds that
			\begin{equation}
				\begin{split}
					\label{eq:flat:approx:linear:parts}
					(\realisation(\f))(x)
					&=\pr*{\tfrac{\xi_k-x}{\xi_{k}-\xi_{k-1}}}(\realisation(\f))(\xi_{k-1})+\pr*{\tfrac{x-\xi_{k-1}}{\xi_{k}-\xi_{k-1}}}(\realisation(\f))(\xi_{k})
					\\&=\pr*{\tfrac{\xi_k-x}{\xi_{k}-\xi_{k-1}}}\tgt(\xi_{k-1})+\pr*{\tfrac{x-\xi_{k-1}}{\xi_{k}-\xi_{k-1}}}\tgt(\xi_{k}).
				\end{split}
			\end{equation}
			This \proves\ that for all $k\in\{1,2,\ldots,M\}$, $x\in[\xi_{k-1},\xi_k]$ it holds that
			\begin{equation}
				\label{eq:flat:approx}
				\begin{split}
					\vass{(\realisation(\f))(x)-\tgt(x)}
					&=\vass*{\pr*{\tfrac{\xi_k-x}{\xi_{k}-\xi_{k-1}}}\tgt(\xi_{k-1})+\pr*{\tfrac{x-\xi_{k-1}}{\xi_{k}-\xi_{k-1}}}\tgt(\xi_{k})-\tgt(x)}
					\\&=\vass*{\pr*{\tfrac{\xi_k-x}{\xi_{k}-\xi_{k-1}}}(\tgt(\xi_{k-1})-\tgt(x))+\pr*{\tfrac{x-\xi_{k-1}}{\xi_{k}-\xi_{k-1}}}(\tgt(\xi_{k})-\tgt(x))}
					\\&\leq\pr*{\tfrac{\xi_k-x}{\xi_{k}-\xi_{k-1}}}\vass*{\tgt(\xi_{k-1})-\tgt(x)}+\pr*{\tfrac{x-\xi_{k-1}}{\xi_{k}-\xi_{k-1}}}\vass{\tgt(\xi_{k})-\tgt(x)}
					\\&\leq\pr*{\tfrac{\xi_k-x}{\xi_{k}-\xi_{k-1}}}\vass*{\xi_{k-1}-x}+\pr*{\tfrac{x-\xi_{k-1}}{\xi_{k}-\xi_{k-1}}}\vass{x-\xi_{k}}
					\\&\leq \vass{\xi_{k}-\xi_{k-1}}=d\leq \eps
				\end{split}
			\end{equation}
			\Moreover 
			\eqref{eq:flat:approx:linearity:and:outside} and
			\eqref{eq:flat:approx:linear:parts} \prove\ that for all $x,y\in\R$ it holds that
			\begin{equation}
				\label{eq:flat:approx:lipschitz}
				\vass{(\realisation(\f))(x)-(\realisation(\f))(y)}\leq \vass{x-y}\max_{k\in\{1,2,\ldots,M\}}\frac{\vass{\tgt(\xi_k)-\tgt(\xi_{k-1})}}{\xi_k-\xi_{k-1}}\leq \vass{x-y}.
			\end{equation}
			\Moreover
			\dimsbene,
			\cref{Prop:identity_representation},
			\cref{gen:scaling:dims:flat:approx},
			\eqref{eq:prop:flat:approx:fdef}, and
			\eqref{eq:flat:approx:netw:def} \prove\ that
			\begin{equation}
				\label{eq:flat:approx:dims}
				\dims(\f)=(1,2,M+1,2,2,2,1)\in\N^7.
			\end{equation}
			This and the fact that $M=2N+2$ and $N\geq 2\pi-1>5$ \prove\ that
			\begin{equation}
				\label{eq:flat:approx:param}
				\begin{split}
					\param(\f)
					&=\sum_{k=1}^6 \singledims_k(\f)(\singledims_{k-1}(\f)+1)
					\\&=2(1+1)+(M+1)(2+1)+2(M+1+1)+2(2(2+1))+1(2+1)
					\\&=26+5M
					=36+10N\leq 4N^2.
				\end{split}
			\end{equation}
			\Moreover 
			\cref{Prop:identity_representation:prop},
			\cref{gen:scaling:size:flat:approx},
			\eqref{eq:prop:flat:approx:param:def},
			\eqref{eq:prop:flat:approx:fdef}, and
			\eqref{eq:flat:approx:netw:def} \prove\ that
			\begin{equation}
				\size(\f)
				=\max\{\size(\g_2),\size(\g_1\bullet\ReLUidANN{1})\}
				=\max\{\size(\g_2),\size(\g_1),1\}
				\leq\max\pRb{2,\tfrac{2\pi}{4},2,1}=2.
			\end{equation}
			Combining this,
			\eqref{eq:flat:approx:length},
			\eqref{eq:flat:approx:linearity:and:outside},
			\eqref{eq:flat:approx},
			\eqref{eq:flat:approx:lipschitz}, and
			\eqref{eq:flat:approx:dims} with
			\eqref{eq:flat:approx:param}
			 \proves[ei]\
			\cref{propo:sin_cost:prop,propo:sin_size:prop,propo:sin_start:prop,propo:sin_lipschitz:prop,propo:sin_layers:prop,propo:sin_dims:prop,propo:sin_approx:prop}.
			The proof of \cref{propo:sin:approximation} is thus complete.
		\end{proof}
		
		\cfclear
		\begin{athm}{lemma}{sin:approximation}
			Let $\varepsilon \in (0,1)$, $n,N \in \N$ satisfy $\eps (N+1) \geq 2\pi$ and let $\tgt\colon\R\to\R$ satisfy for all $x,y\in\R$, $k\in\Z$ that $\vass{\tgt(0)}\leq 2$, $\tgt(x+2k\pi)=\tgt(x)$, and $\vass{\tgt(x)-\tgt(y)}\leq\vass{x-y}$.
			Then there exists $\mathscr{f} \in \ANNSET$ such that
			\begin{enumerate}[(i)]
				\item {\label{sin_layers:prop} it holds 
					that $\lengthANN(\f)\leq n+11$,
				}
				\item {\label{sin_dims:prop}
					it holds that $\singledims_0(\f)=\singledims_{\lengthANN(\f)}(\f)=1$, 
					$\singledims_1(\f)\leq 2N$, and $\singledims_{\hidlengthANN(\f)}(\f)=2$,}
				\item{\label{sin_lipschitz:prop} it holds for all
					$x,y \in \R$ that 
					$\vass[\big]{\pr{\realisation(\f)} (x) - \pr{\realisation(\f)} (y)} \leq \vass{x-y}$,}
				\item{\label{sin_approx:prop} it holds that
					$\sup_{x \in [-2^n\pi,2^n\pi]}\vass[\big]{\tgt(x) - \pr{\realisation(\f)} (x)} \leq \varepsilon$,}
				\item{\label{sin_start:prop} it holds for all $x\in[2^n\pi,\infty)$ that 
					$(\realisation(\f)) (-x) =(\realisation(\f)) (x) =\tgt(0)$,}
				\item{\label{sin_size:prop} it holds 
					that
					$\size(\mathscr{f}) \leq 2
					$, and}
				\item{\label{sin_cost:prop} it holds 
					that
					$\param(\mathscr{f}) \leq  \pr*{24+18n}N^2$}
			\end{enumerate}
			\cfout. 
	\end{athm}
	
	\begin{proof}[Proof of \cref{sin:approximation}]
		Throughout this proof 
		assume without loss of generality that $n\geq 2$ \cfadd{propo:sin:approximation}\cfload, let $\edgy\colon\R\to\R$ satisfy for all $k\in\Z$, $x\in[2k,2k+2)$ that $\edgy(x)=1-\vass{x-2k-1}$,
		for every $j\in\{0,1,\ldots,N+1\}$ let 
		$
		c_j=\frac{2j\pi}{N+1}
		$,
		and for every $j\in\{1,2,\ldots,N\}$ 
		let $f_j\colon\R\to\R$ satisfy for all $x\in[0,2\pi)$, $y\in(-\infty,-2^n\pi)\cup[2^n\pi,\infty)$,  $k\in\Z\cap[-2^{n-1},2^{n-1})$ that $f_j(y)=0$ and
		\begin{equation}
			\label{eq:setup:hatfunctions2}
			f_j(x+2k\pi)=f_j(x)=
			\begin{cases}
				(\tgt(c_j)-\tgt(0))\edgy\pr*{\frac{(x-c_{j-1})(N+1)}{2\pi}} &\colon x\in [c_{j-1},c_{j+1}]\\
				0 &\colon x\notin [c_{j-1},c_{j+1}].
			\end{cases}
		\end{equation}
		\Nobs that the assumption that for all $x,y\in\R$, $k\in\Z$ it holds that $\vass{\tgt(0)}\leq 2$, $\tgt(x+2k\pi)=\tgt(x)$, and $\vass{\tgt(x)-\tgt(y)}\leq\vass{x-y}$ \proves\ that
		for all $j\in\{1,2,\dots,N\}$ it holds that
		\begin{equation}
			\textstyle
			\vass{g(c_j)-g(0)}
			=\inf_{k\in\Z}	\vass{g(c_j+2k\pi)-g(0)}
			\leq \inf_{k\in\Z}	\vass{c_j+2k\pi}\leq \pi<4.
		\end{equation}
		This and \cref{sum:hat:networks:eps} (applied with
		$n \curvearrowleft  n$,
		$N \curvearrowleft  N$,
		$C \curvearrowleft  4$,
		$\shft \curvearrowleft  \tgt(0)$,
		$a_j \curvearrowleft  c_{j-1}$,
		$b_j \curvearrowleft  c_{j+1}$,
		$c_j \curvearrowleft  \tgt(c_j)-\tgt(0)$,
		$f_j \curvearrowleft  f_j$ for $j\in\{1,2,\ldots,N\}$
		in the notation of \cref{sum:hat:networks:eps}) \prove\ that there exists $\f\in\ANNSET$ which satisfies that 
		\begin{enumerate}[(I)]
			\item{
				\label{it:sum:hat:networks1}
				it holds for all $x\in\R$ that
				$
				(\realisation(\f))(x)=\tgt(0)+\sum_{j=1}^N f_j(x)
				$,
			}
			\item{
				\label{it:sum:hat:networks3}
				it holds that $\lengthANN(\f)\leq n+\log_2(4)+9=n+11$,
			}
			\item{
				\label{it:sum:hat:networks4}
				it holds that $\singledims_0(\f)=\singledims_{\lengthANN(\f)}(\f)=1$, $\singledims_1(\f)=2N$, and $\singledims_{\hidlengthANN(\f)}(\f)=2$,
			}
			\item{
				\label{it:sum:hat:networks5}
				it holds that $\param(\f)\leq \pr*{14+18n+5\log_2(4)}N^2
				=\pr*{24+18n}N^2$, and
			}
			\item{
				\label{it:sum:hat:networks6}
				it holds that $\size(\f)\leq 2$
			}
		\end{enumerate}
		\cfload.
		\Nobs that \cref{it:sum:hat:networks1} and \cref{sin:approximator} (applied with
		$\eps \curvearrowleft  \eps$,
		$n \curvearrowleft  n$,
		$N \curvearrowleft  N$,
		$c_j \curvearrowleft  c_j$,
		$\tgt \curvearrowleft \tgt$,
		$f_j \curvearrowleft  f_j$,
		$F \curvearrowleft  \realisation(\f)$
		for $j\in\{1,2,\ldots,N\}$
		in the notation of \cref{sin:approximator}) \prove\ that for all $x\in[-2^n\pi,2^n\pi]$, $y,z\in\R$ it holds that
		\begin{equation}
			\label{eq:lipschitz}
			\vass*{\tgt(x)-(\realisation(\f))(x)}\leq \eps
			\qandq
			\vass*{(\realisation(\f))(y)-(\realisation(\f))(z)}\leq \vass{y-z}.
		\end{equation}
		Combining this and
		the fact that $\realisation(\f)\in C(\R,\R)$ with
		\cref{it:sum:hat:networks1,it:sum:hat:networks3,it:sum:hat:networks4,it:sum:hat:networks5,it:sum:hat:networks6}
		 \proves[ei]\ \cref{sin_cost:prop,sin_dims:prop,sin_layers:prop,sin_size:prop,sin_start:prop,sin_approx:prop,sin_lipschitz:prop}.
		The proof of \cref{sin:approximation} is thus complete.
		\cfload
	\end{proof}

	\begin{athm}{cor}{1:dim:scaled:test}
		Let $R\in(0,\infty)$, $\gamma\in(0,1]$, $\scl\in[1,\infty)$,  $\varepsilon \in (0,1)$ and let $\tgt\colon\R\to\R$ satisfy for all $x,y\in\R$, $k\in\Z$ that $\vass{\tgt(0)}\leq 2$, $\tgt(x+2k\pi)=\tgt(x)$, and $\vass{\tgt(x)-\tgt(y)}\leq\vass{x-y}$.
		Then there exists $\mathscr{f} \in \ANNSET$ such that
		
		\begin{enumerate}[(i)]
			\item \label{1:dim:scaled:test:continuous} it holds 
			that $\realisation(\mathscr{f}) \in C(\R,\R)$,

			\item\label{1:dim:scaled:test:approx} it holds that 
			$\sup_{x  \in [-R,R]}\vass[\big]{\tgt(\gamma\scl x)- (\realisation\pr{\mathscr{f}}) (x)} \leq \varepsilon$,
			
			\item\label{1:dim:scaled:test:hid} it holds 
			that
			$\singledims_{1}(\f)=\singledims_{\hidlengthANN(\mathscr{f})}(\f)=2$,
			
			\item\label{1:dim:scaled:test:length} it holds 
			that
			$\lengthANN(\mathscr{f}) \leq 10\max\{1,\ceil{\log_2(\scl)},\ceil{\log_2(R)}\}
			$,
			
			\item\label{1:dim:scaled:test:param} it holds 
			that
			$\paramANN(\mathscr{f}) \leq  1935 \max\{1,\ceil{\log_2(\scl)},\ceil{\log_2(R)}\}\eps^{-2}
			$, and

			\item\label{1:dim:scaled:test:size} it holds 
			that
			$\size(\mathscr{f}) \leq 2
			$
		\end{enumerate}
		\cfout. 
\end{athm}

\cfclear
\begin{proof}[Proof of \cref{1:dim:scaled:test}]
	Throughout this proof assume without loss of generality that $R\geq 2$ and let $M,N\in\N$ satisfy $M=\max\{1,\ceil{\log_2(\scl)},\ceil{\log_2(R)}\}$ and $N=\ceil*{\tfrac{2\pi}{\eps}}-1$.
	\Nobs that \cref{gen:scaling:networks} (applied with 
	$\beta \curvearrowleft \gamma\scl$,
	$L \curvearrowleft \max\{0,\ceil{\log_2(\gamma\scl)}\}$
	in the notation of \cref{gen:scaling:networks}) \proves\ that there exists $\g_1\in\ANNSET$ which satisfies that
	\begin{enumerate}[(A)]
		\item{
			\label{1:dim:gen:scaling:realisation}
			it holds for all $x\in\R$ that $(\realisation(\g_1))(x)=\gamma\scl x$,
		}
		\item{
			\label{1:dim:gen:scaling:dims}
			it holds that
			 $\dims(\g_1)=(1,2,2,\ldots,2,1)\in\N^{\max\{0,\ceil{\log_2(\gamma\scl)}\}+2}$, and
		}
		\item{
			\label{1:dim:gen:scaling:size}
			it holds that $\size(\g_1)\leq 2$}
	\end{enumerate}
	\cfload.
	In our proof we distinguish between the case $M=1$ and the case $M\geq 2$. 
	We first \prove[pes]\
	\cref{1:dim:scaled:test:continuous,1:dim:scaled:test:approx,1:dim:scaled:test:length,1:dim:scaled:test:param,1:dim:scaled:test:size,1:dim:scaled:test:hid}
	 in the case of $M=1$.
	 \Nobs that \cref{propo:sin:approximation} (applied with 
	 $\eps \curvearrowleft \eps$,
	 $N \curvearrowleft N$,
	 $g \curvearrowleft g$
	 in the notation of \cref{propo:sin:approximation}) \proves\ that there exists $\g_2\in\ANNSET$ which satisfies that
	 \begin{enumerate}[(I)]
	 	\item {\llabel{1:dim:sin_layers:prop} it holds 
	 		that $\lengthANN(\g_2)= 6$,
	 	}
	 	\item {\llabel{1:dim:sin_dims:prop}
	 		it holds that $\dims(\g_2)=(1,2,2N+3,2,2,2,1)$,}
	 	\item{\llabel{1:dim:sin_approx:prop} it holds that
	 		$\sup_{x \in [-2\pi ,2\pi]}\vass[\big]{\tgt(x) - \pr{\realisation(\g_2)} (x)} \leq \varepsilon$,}
	 	%
	 	%
	 	\item{\llabel{1:dim:sin_size:prop} it holds 
	 		that
	 		$\size(\g_2) \leq 2
	 		$, and}
	 	\item{\llabel{1:dim:sin_cost:prop} it holds 
	 		that
	 		$\param(\g_2) \leq  4N^2$\ifnocf.
	 	}
	 \end{enumerate}
	 \Nobs that 
	 \cref{Lemma:PropertiesOfCompositions_n2},
	 \cref{Prop:identity_representation:prop},
	 \lref{1:dim:sin_dims:prop}, and
	 \cref{1:dim:gen:scaling:dims}
	 \prove\ that 
	 \begin{equation}
	 	\llabel{eq:1:dim:function}
	 	\realisation(\compANN{\g_2\bullet\ReLUidANN{1}}{\g_1})=\PR{\realisation(\g_2)}\circ\PR{\realisation(\g_1)}\in C(\R,\R)
	 \end{equation}
	 \cfload.
	 This,
	 \lref{1:dim:sin_approx:prop}, 
	 \cref{1:dim:gen:scaling:realisation}, and the fact that $\fa{x}[-R,R]\colon\gamma\scl x\in[-\scl R,\scl R]\subseteq[-2\pi,2\pi]$ \prove\ that for all $x\in[-R,R]$ it holds that
	 \begin{equation}
	 	\llabel{eq:1:dim:approx}
	 	\vass{\tgt(\gamma\scl x)-(\realisation(\compANN{\g_2\bullet\ReLUidANN{1}}{\g_1}))(x)}
	 	=\vass{\tgt(\gamma\scl x)-(\realisation(\g_2))(\gamma\scl x)}\leq \eps.
	 \end{equation}
	 \Nobs that 
	 \cref{Lemma:PropertiesOfCompositions_n2},
	 \cref{Prop:identity_representation:prop},
	 \lref{1:dim:sin_layers:prop}, 
	 \cref{1:dim:gen:scaling:dims}, and the assumption that $\gamma\leq 1$ 
	 \prove\ that 
	 \begin{equation}
	 	\llabel{eq:1:dim:length}
	 	\begin{split}
	 		\lengthANN(\compANN{\g_2\bullet\ReLUidANN{1}}{\g_1})
	 		=\lengthANN(\g_2)+\lengthANN(\g_1)
	 		&\leq (\max\{0,\ceil{\log_2(\gamma\scl)}\}+1)+6
	 		\\&\leq (\ceil{\log_2(\scl)}+1)+6
	 		\leq 10M.
	 	\end{split}
	 \end{equation}
	 This,
	 \dimsbene, and
	 \lref{1:dim:sin_dims:prop}
	 \prove\ that for all $k\in\N_0\cap[0,\lengthANN(\g_2)+\lengthANN(\g_1)]$ it holds that
	 \begin{equation}
	 	\llabel{eq:1:dim:dims}
	 	\begin{split}
	 		\singledims_k(\compANN{\g_2\bullet\ReLUidANN{1}}{\g_1})
	 		=
	 		\begin{cases}
	 			1 &\colon k=0\\
	 			2 &\colon k\in\N\cap(0,\lengthANN(\g_1)]\\
	 			\singledims_{k-\lengthANN(\g_1)}(\g_2) &\colon k\in\N\cap(\lengthANN(\g_1),\lengthANN(\g_2)+\lengthANN(\g_1)].
	 		\end{cases}
	 	\end{split}
	 \end{equation}
	 Combining this and \lref{1:dim:sin_dims:prop} with \lref{eq:1:dim:length} \proves\ that
	 \begin{equation}
	 	\llabel{eq:1:dim:param}
	 	\begin{split}
	 		&\param(\compANN{\g_2\bullet\ReLUidANN{1}}{\g_1})
	 		\\&=\ssum_{k=1}^{\lengthANN(\g_2)+\lengthANN(\g_1)}\singledims_k(\compANN{\g_2\bullet\ReLUidANN{1}}{\g_1})(\singledims_{k-1}(\compANN{\g_2\bullet\ReLUidANN{1}}{\g_1})+1)
	 		\\&=2(1+1)+\PR*{\ssum_{k=2}^{\lengthANN(\g_1)}2(2+1)}+\singledims_{1}(\g_2)(2+1)
	 		+\PR*{\ssum_{k=2}^{\lengthANN(\g_2)}\singledims_{k}(\g_2)(\singledims_{k-1}(\g_2)+1)}
	 		\\&=4+6(\lengthANN(\g_1)-1)+\singledims_{1}(\g_2)+\PR*{\ssum_{k=1}^{\lengthANN(\g_2)}\singledims_{k}(\g_2)(\singledims_{k-1}(\g_2)+1)}
	 		\\&\leq 4+6\max\{0,\ceil{\log_2(\gamma\scl)}\}+2+\param(\g_2)
	 		\\&\leq 6+6\max\{0,\ceil{\log_2(\gamma\scl)}\}+4N^2
	 		\\&\leq 6+6\max\{0,\ceil{\log_2(\scl)}\}+16\pi^2\eps^{-2}
	 		\leq 172M\eps^{-2}\leq 1935M\eps^{-2}
	 		.
	 	\end{split}
	 \end{equation}
	 \Moreover
	 \cref{Prop:identity_representation:prop},
	 \lref{1:dim:sin_size:prop}, and
	 \cref{1:dim:gen:scaling:size} \prove\ that
	 \begin{equation}
	 	\size(\compANN{\g_2\bullet\ReLUidANN{1}}{\g_1})=\max\{\size(\g_2),\size(\g_1)\}\leq\max\{2,2\}=2.
	 \end{equation}
	 This,
	 \lref{eq:1:dim:function},
	 \lref{eq:1:dim:approx},
	 \lref{eq:1:dim:length},
	 \lref{eq:1:dim:dims}, and
	 \lref{eq:1:dim:param}
	 \prove[ei]\
	 \cref{1:dim:scaled:test:continuous,1:dim:scaled:test:approx,1:dim:scaled:test:length,1:dim:scaled:test:param,1:dim:scaled:test:size,1:dim:scaled:test:hid}.
	 Next we \prove[pes]\
	 \cref{1:dim:scaled:test:continuous,1:dim:scaled:test:approx,1:dim:scaled:test:length,1:dim:scaled:test:param,1:dim:scaled:test:size,1:dim:scaled:test:hid}
	 in the case of $M\geq 2$.
	\Nobs that \cref{sin:approximation} (applied with 
	$\eps \curvearrowleft \eps$,
	$n \curvearrowleft \ceil{\log_2(\scl R)}$,
	$N \curvearrowleft N$,
	$g \curvearrowleft g$
	in the notation of \cref{sin:approximation}) \proves\ that there exists $\g_2\in\ANNSET$ which satisfies that
	\begin{enumerate}[(a)]
		\item {\label{1:dim:sin_layers:prop} it holds 
			that $\lengthANN(\g_2)\leq \ceil{\log_2(\scl R)}+11$,
		}
		\item {\label{1:dim:sin_dims:prop}
			it holds that $\singledims_0(\g_2)=\singledims_{\lengthANN(\g_2)}(\g_2)=1$, 
			$\singledims_1(\g_2)\leq 14\eps^{-1}$, and $\singledims_{\hidlengthANN(\g_2)}(\g_2)=2$,}
		%
		%
		\item{\label{1:dim:sin_approx:prop} it holds that
			$\sup_{x \in [-\scl R ,\scl R]}\vass[\big]{\tgt(x) - \pr{\realisation(\g_2)} (x)} \leq \varepsilon$,}
		%
		%
		\item{\label{1:dim:sin_size:prop} it holds 
			that
			$\size(\g_2) \leq 2
			$, and}
		\item{\label{1:dim:sin_cost:prop} it holds 
			that
			$\param(\g_2) \leq  \pr*{24+18\ceil{\log_2(\scl R)}}\prb{2\pi\eps^{-1}}^2$\ifnocf.
		}
	\end{enumerate}
	\Nobs that 
	\cref{Lemma:PropertiesOfCompositions_n2},
	\cref{Prop:identity_representation:prop},
	\cref{1:dim:sin_dims:prop}, and
	\cref{1:dim:gen:scaling:dims}
	\prove\ that 
	\begin{equation}
		\label{eq:1:dim:function}
		\realisation(\compANN{\g_2\bullet\ReLUidANN{1}}{\g_1})=\PR{\realisation(\g_2)}\circ\PR{\realisation(\g_1)}\in C(\R,\R).
	\end{equation}
	This,
	\cref{1:dim:sin_approx:prop}, 
	\cref{1:dim:gen:scaling:realisation}, and the fact that $\fa{x}[-R,R]\colon\gamma\scl x\in[-\scl R,\scl R]$ \prove\ that for all $x\in[-R,R]$ it holds that
	\begin{equation}
		\label{eq:1:dim:approx}
		\vass{\tgt(\gamma\scl x)-(\realisation(\compANN{\g_2\bullet\ReLUidANN{1}}{\g_1}))(x)}
		=\vass{\tgt(\gamma\scl x)-(\realisation(\g_2))(\gamma\scl x)}\leq \eps.
	\end{equation}
	\Nobs that 
	\cref{Lemma:PropertiesOfCompositions_n2},
	\cref{Prop:identity_representation:prop},
	\cref{1:dim:sin_layers:prop}, 
	\cref{1:dim:gen:scaling:dims}, and the assumption that $\gamma\leq 1$ 
	\prove\ that 
	\begin{equation}
		\label{eq:1:dim:length}
		\begin{split}
			\lengthANN(\compANN{\g_2\bullet\ReLUidANN{1}}{\g_1})
			=\lengthANN(\g_2)+\lengthANN(\g_1)
			&\leq (\max\{0,\ceil{\log_2(\gamma\scl)}\}+1)+(\ceil{\log_2(\scl R)}+11)
			\\&\leq (\ceil{\log_2(\scl)}+1)+(\ceil{\log_2(\scl )}+\ceil{\log_2(R)}+11)
			\\&\leq 10M.
		\end{split}
	\end{equation}
	This,
	\dimsbene, and
	\cref{1:dim:sin_dims:prop}
	\prove\ that for all $k\in\N_0\cap[0,\lengthANN(\g_2)+\lengthANN(\g_1)]$ it holds that
	\begin{equation}
		\label{eq:1:dim:dims}
		\begin{split}
			\singledims_k(\compANN{\g_2\bullet\ReLUidANN{1}}{\g_1})
			=
			\begin{cases}
				1 &\colon k=0\\
				2 &\colon k\in\N\cap(0,\lengthANN(\g_1)]\\
				\singledims_{k-\lengthANN(\g_1)}(\g_2) &\colon k\in\N\cap(\lengthANN(\g_1),\lengthANN(\g_2)+\lengthANN(\g_1)].
			\end{cases}
		\end{split}
	\end{equation}
	Combining this and \cref{1:dim:sin_dims:prop} with \eqref{eq:1:dim:length} \proves\ that
	\begin{equation}
		\label{eq:1:dim:param}
		\begin{split}
			&\param(\compANN{\g_2\bullet\ReLUidANN{1}}{\g_1})
			\\&=\ssum_{k=1}^{\lengthANN(\g_2)+\lengthANN(\g_1)}\singledims_k(\compANN{\g_2\bullet\ReLUidANN{1}}{\g_1})(\singledims_{k-1}(\compANN{\g_2\bullet\ReLUidANN{1}}{\g_1})+1)
			\\&=2(1+1)+\PR*{\ssum_{k=2}^{\lengthANN(\g_1)}2(2+1)}+\singledims_{1}(\g_2)(2+1)
			+\PR*{\ssum_{k=2}^{\lengthANN(\g_2)}\singledims_{k}(\g_2)(\singledims_{k-1}(\g_2)+1)}
			\\&=4+6(\lengthANN(\g_1)-1)+\singledims_{1}(\g_2)+\PR*{\ssum_{k=1}^{\lengthANN(\g_2)}\singledims_{k}(\g_2)(\singledims_{k-1}(\g_2)+1)}
			\\&\leq 4+6\max\{0,\ceil{\log_2(\gamma\scl)}\}+14\eps^{-1}+\param(\g_2)
			\\&\leq 4+6\ceil{\log_2(\scl)}+14\eps^{-1}+\pr*{24+18\ceil{\log_2(\scl R)}}\prb{2\pi\eps^{-1}}^2
			\\&\leq 4+6\ceil{\log_2(\scl)}+14\eps^{-1}+4\pi^2\pr*{24+18\ceil{\log_2(\scl )}+18\ceil{\log_2(R)}}\eps^{-2}
			\\&\leq 2M+6M+7M\eps^{-1}+40\pr*{12M+18M+18M}\eps^{-2}
			\\&\leq 15M\eps^{-1}+1920M\eps^{-2}\leq1935M\eps^{-2}
			.
		\end{split}
	\end{equation}
	\Moreover
	\cref{Prop:identity_representation:prop},
	\cref{1:dim:sin_size:prop}, and
	\cref{1:dim:gen:scaling:size} \prove\ that
	\begin{equation}
		\size(\compANN{\g_2\bullet\ReLUidANN{1}}{\g_1})=\max\{\size(\g_2),\size(\g_1)\}\leq\max\{2,2\}=2.
	\end{equation}
	This,
	\eqref{eq:1:dim:function},
	\eqref{eq:1:dim:approx},
	\eqref{eq:1:dim:length},
	\eqref{eq:1:dim:dims}, and
	\eqref{eq:1:dim:param}
	 \prove[ei]\
	\cref{1:dim:scaled:test:continuous,1:dim:scaled:test:approx,1:dim:scaled:test:length,1:dim:scaled:test:param,1:dim:scaled:test:size,1:dim:scaled:test:hid}.
	The proof of \cref{1:dim:scaled:test} is thus complete.
\end{proof}

\subsection{Upper bounds for approximations of compositions of periodic and product functions}
\label{subsection:upper_sine_prod}

\cfclear
\begin{athm}{lemma}{final:approximation}
    Let $d \in \N$, $\varepsilon \in (0,1)$, $R\in(0,\infty)$, $\gamma\in(0,1]$, $\scl\in[1,\infty)$,
		$\tgt\in C(\R,\R)$ satisfy for all $x,y\in\R$, $k\in\Z$ that $\vass{\tgt(0)}\leq 2$, $\tgt(x+2k\pi)=\tgt(x)$, and $\vass{\tgt(x)-\tgt(y)}\leq\vass{x-y}$.
    Then there exists $\mathscr{f} \in \ANNSET$ such that

    \begin{enumerate}[(i)]
    \item \label{prop:new_continuous_d} it holds 
    that $\realisation(\mathscr{f}) \in C(\R^d,\R)$,

    \item\label{prop:new_approx_d} it holds that
    $
    \sup_{x = (x_1,\ldots, x_{d}) \in [-\radius, \radius]^{d}}\vass[\big]{g\prb{\gamma\scl^d\sprod_{i = 1}^d x_i} - \realisation\pr{\mathscr{f}} (x)} \leq \varepsilon$,

    \item\label{prop:new_hid_d} it holds 
	that
	$\singledims_{\hidlengthANN(\mathscr{f})}(\f)=2$,

    \item\label{prop:new_length_d} it holds 
    that
    $\lengthANN(\mathscr{f}) \leq  16 \max\{1,\ceil{\log_2( R)},\ceil{\log_2(\scl)}\}d^2\max\{1,\log_2(\eps^{-1})\}
    $,

    \item\label{prop:new_cost_d} it holds 
    that
    $\paramANN(\mathscr{f}) \leq   9105\max\{1,\ceil{\log_2( R)},\log_2(\scl)\}d^3\eps^{-2}
    $, and

		\item\label{prop:new_size_d} it holds 
    that
    $\size(\mathscr{f}) \leq 4
    $
    \end{enumerate}
    \cfout. 
\end{athm}

\begin{proof}[Proof of \cref{final:approximation}]
Throughout this proof assume without loss of generality that $R\geq 2$ and $d>1$ (cf.~\cref{1:dim:scaled:test}),
let $n,N\in\N$ satisfy $n =\ceil*{d\log_2(\scl R)}-1$ and $N =\ceil{\frac{4\pi}{\eps}}-1$.
\Nobs that \cref{lemma:d_product} (applied with 
$d \curvearrowleft d$,
$\eps \curvearrowleft \frac{\varepsilon}{2}$,
$R \curvearrowleft R$,
$\gamma \curvearrowleft \gamma$, 
$\scl \curvearrowleft \scl$
in the notation of \cref{lemma:d_product}) \proves\ that there exists $\g_1\in\ANNSET$ which satisfies that
\begin{enumerate}[(I)]
    \item \label{it:d_product_continuous} it holds 
    that $\realisation(\g_1) \in C(\R^{d},\R)$,

    \item\label{it:d_product_approx} it holds that
    $
    \sup_{x = (x_1,\ldots, x_{d}) \in [-\radius, \radius]^{d}}\vass[\big]{\gamma\scl^d\sprod_{i=1}^{d}x_i - (\realisation(\g_1)) (x)} 
                \leq 
                \frac{\varepsilon}{2}$,
		

    \item\label{it:d_product_length} it holds 
    that
    $\lengthANN(\g_1) \leq 8d^2+2d^2\ceil{\log_2(R)}+d(\log_2(\eps^{-1})+1)+d^2\ceil{\log_2(\scl)}+2
    $,

		\item\label{it:d_product_dims} it holds that
				$\singledims_1(\g_1)\leq 4d
		$ and
				$\singledims_{\hidlengthANN(\g_1)}(\g_1)=2
		$,

    \item\label{it:d_product_cost} it holds 
    that
    $\paramANN(\g_1) \leq 4107 d^3+2048 d^3\ceil{\log_2(R)}+512 d^2 (\log_2(\eps^{-1})+1)+514d^3\log_2(\scl)
    $, and
		
		\item\label{it:d_product_size} it holds 
    that
    $\size(\g_1) \leq 4
    $
\end{enumerate}
		\cfload. \Moreover 
		the fact that $\sup_{x\in\R}\vass{\tgt(x)}\leq 2+\pi \leq 6$ and
		\cref{sin:approximation} (applied with 
$\eps \curvearrowleft \frac{\varepsilon}{2}$,
$n \curvearrowleft n$,
$N \curvearrowleft N$,
$g \curvearrowleft g$
in the notation of \cref{sin:approximation}) \prove\ that there exists $\g_2\in\ANNSET$ which satisfies that
		     \begin{enumerate}[(A)]
        \item {\label{it:sin_layers:prop} it holds that $\lengthANN(\g_2)\leq n+11$,}
		\item {\label{it:sin_dims:prop}
    it holds that $\singledims_0(\g_2)=\singledims_{\lengthANN(\g_2)}(\g_2)=1$, 
		$\singledims_1(\g_2)\leq 2N$, and $\singledims_{\hidlengthANN(\g_2)}(\g_2)=2$,}
    \item{\label{it:sin_lipschitz:prop} it holds for all
    $x,y \in \R$ that 
    $\vass[\big]{\pr{\realisation(\g_2)} (x) - \pr{\realisation(\g_2)} (y)} \leq \vass{x-y}$,}
    \item{\label{it:sin_approx:prop} it holds that
    $\sup_{x \in [-2^n\pi,2^n\pi]}\vass[\big]{\tgt(x) - \pr{\realisation(\g_2)} (x)} \leq \tfrac{\varepsilon}{2}$,}
    \item{\label{it:sin_start:prop} it holds for all $x\in[2^n\pi,\infty)$ that 
    $\pr{\realisation(\g_2)} (-x)=\pr{\realisation(\g_2)} (x) =\tgt(0)$,}
		\item{\label{it:sin_size:prop} it holds 
    that
    $\size(\g_2) \leq 2
    $, and}
    \item{\label{it:sin_cost:prop} it holds 
    that
    $\param(\g_2) \leq (24+18n)N^2
    $.}
    \end{enumerate}
		Next let $\f\in\ANNSET$ satisfy
		\begin{equation}
		\label{definition:of:final:imp:approximator}
		\f=\compANN{\g_2}{\ReLUidANN{1}}\bullet\g_1
		\end{equation}
\cfload[.]\Nobs that 
\eqref{definition:of:final:imp:approximator}, 
\cref{it:d_product_length}, 
\cref{it:sin_layers:prop}, 
\cref{Lemma:PropertiesOfCompositions_n2}, 
\cref{Prop:identity_representation}, and
the fact that $n\leq d\ceil{\log_2(\scl R)}-1\leq 2d\max\{\ceil{\log_2(\scl)},\ceil{\log_2( R)}\}-1$ \prove\ that
\begin{equation}
	\llabel{eq:loc:length:ref}
\begin{split}
&\lengthANN(\f)
\\&=\lengthANN(\g_2)+\lengthANN(\ReLUidANN{1})+\lengthANN(\g_1)-2
\\&=\lengthANN(\g_2)+\lengthANN(\g_1)
\\&\leq(n+11)+\prb{8d^2+2d^2\ceil{\log_2(R)}+d(\log_2(\eps^{-1})+1)+d^2\ceil{\log_2(\scl)}+2}
\\&\leq (d\ceil{\log_2(\scl R)}+10)+\prb{8d^2+2d^2\ceil{\log_2(R)}+d\log_2(\eps^{-1})+d+d^2\ceil{\log_2(\scl)}+2}
\\&\leq (2d+10+8d^2+2d^2+d+d+d^2+2) \max\{1,\ceil{\log_2( R)},\log_2(\eps^{-1}),\ceil{\log_2(\scl)}\}
\\&= (11d^2+4d+12) \max\{\ceil{\log_2( R)},\log_2(\eps^{-1}),\ceil{\log_2(\scl)}\}
\\&\leq 16d^2 \max\{\ceil{\log_2( R)},\log_2(\eps^{-1}),\ceil{\log_2(\scl)}\}
\\&\leq 16 \max\{\ceil{\log_2( R)},\ceil{\log_2(\scl)}\}d^2\max\{1,\log_2(\eps^{-1})\}
.
 \end{split}
\end{equation}
This, 
\eqref{definition:of:final:imp:approximator}, 
\cref{it:d_product_continuous}, 
\cref{it:sin_dims:prop}, and
\dimsbene\ \prove\ that $\singledims_{\hidlengthANN(\mathscr{f})}(\f)=2$ and
\begin{multline}
\label{eq:comp:dims2}
\dims(\f)\\
=(\singledims_0(\g_1),\singledims_1(\g_1),\dots,\singledims_{\hidlengthANN(\g_1)}(\g_1),
			2\singledims_{\lengthANN(\g_1)}(\g_1),\singledims_1(\g_2),\singledims_2(\g_2),\dots,\singledims_{\lengthANN(\g_2)}(\g_2))
			.
\end{multline}
\Nobs that 
\eqref{definition:of:final:imp:approximator}, 
\cref{Lemma:PropertiesOfCompositions_n2}, 
and 
\cref{Prop:identity_representation} \prove\ that for all $x\in\R^d$ it holds that
\begin{equation}
\begin{split}
\label{eq:specific:realisation}
(\realisation(\f))(x)
&=\pr*{\PR*{\realisation(\g_2)}\circ\PR*{\realisation(\g_1)}}(x)
=\PR*{\realisation(\g_2)}\pr{(\realisation(\g_1))(x)}.
\end{split}
\end{equation}
This, 
\cref{it:d_product_approx}, 
\cref{it:sin_lipschitz:prop}, 
\cref{it:sin_approx:prop}, 
 and the fact that for all $x_1,x_2,\ldots,x_d\in [-R,R]$ it holds that $\vass{\gamma\scl^d\sprod_{i = 1}^d x_i}\leq (\scl R)^d \leq 2^{d\log_2(\scl R)-1}\pi\leq 2^n\pi$ \prove\ that for all $x=(x_1,\ldots,x_d)\in [-R,R]^d$ it holds that
\begin{align}
\label{eq:approx:g:prod}
\nonumber&\vass[\bigg]{\tgt\prbb{\gamma\scl^d\sprod_{i = 1}^d x_i} - \realisation\pr{\mathscr{f}} (x)}
\\&=\vass[\bigg]{\tgt\prbb{\gamma\scl^d\sprod_{i = 1}^d x_i} - (\realisation(\g_2))\pr{(\realisation(\g_1))(x)}}
\\\nonumber&=\vass[\bigg]{\tgt\prbb{\gamma\scl^d\sprod_{i = 1}^d x_i} - \pr*{\realisation\pr{\g_2}}\prbb{\gamma\scl^d\sprod_{i = 1}^d x_i}+\pr*{\realisation\pr{\g_2}}\prbb{\gamma\scl^d\sprod_{i = 1}^d x_i}-\pr*{\realisation\pr{\g_2}}\pr*{(\realisation\pr{\g_1}) (x)}}
\\\nonumber&\leq \vass[\bigg]{\tgt\prbb{\gamma\scl^d\sprod_{i = 1}^d x_i} - \pr*{\realisation\pr{\g_2}}\prbb{\gamma\scl^d\sprod_{i = 1}^d x_i}}
+
\vass[\bigg]{\pr*{\realisation\pr{\g_2}}\prbb{\gamma\scl^d\sprod_{i = 1}^d x_i}-\pr*{\realisation\pr{\g_2}}\pr*{(\realisation\pr{\g_1}) (x)}}
\\\nonumber&\leq \frac{\varepsilon}{2}+\vass[\big]{\gamma\scl^d\sprod_{i = 1}^d x_i -  (\realisation\pr{\g_1}) (x)}\leq \frac{\varepsilon}{2}+\frac{\varepsilon}{2}=\eps.
\end{align}
\Nobs that 
\eqref{eq:comp:dims2}, 
\cref{it:sin_layers:prop}, 
\cref{it:sin_dims:prop},  
\cref{it:d_product_dims}, and the fact that $\singledims_{\lengthANN(\g_1)}(\g_1)=\singledims_{0}(\g_2)$ \prove\ that
\begin{equation}
\begin{split}
\param\pr{\f}
&=\sum_{k=1}^{\lengthANN(\f)} \singledims_k(\f)(\singledims_{k-1}(\f)+1)
\\&=\PR*{\sum_{k=1}^{\lengthANN(\g_1)-1} \singledims_k(\g_1)(\singledims_{k-1}(\g_1)+1)}
+2\singledims_{\lengthANN(\g_1)}(\g_1)(\singledims_{\lengthANN(\g_1)-1}(\g_1)+1)
\\&\quad+\singledims_{1}(\g_2)(2\singledims_{0}(\g_2)+1)
+\PR*{\sum_{k=2}^{\lengthANN(\g_2)} \singledims_k(\g_2)(\singledims_{k-1}(\g_2)+1)}
\\&=\PR*{\sum_{k=1}^{\lengthANN(\g_1)} \singledims_k(\g_1)(\singledims_{k-1}(\g_1)+1)}
+\singledims_{\lengthANN(\g_1)}(\g_1)(\singledims_{\lengthANN(\g_1)-1}(\g_1)+1)
\\&\quad+\singledims_{1}(\g_2)\singledims_{0}(\g_2)
+\PR*{\sum_{k=1}^{\lengthANN(\g_2)} \singledims_k(\g_2)(\singledims_{k-1}(\g_2)+1)}
\\&\leq\param(\g_1)
+1(2+1)
+2N
+\param(\g_2)
.
\end{split}
\end{equation}
Combining this, \cref{it:d_product_cost} and \cref{it:sin_cost:prop} with the fact that for all $x\in(0,1)$ it holds that $-\log_2\pr*{x}\leq x^{-2}$ and $\pi^2\leq 10$ with the assumption that $\varepsilon \in (0,1)$, $R\geq 2$, $d\geq 2$, $n\leq \max\{d\log_2(\scl R),1\}$, and $N \leq\frac{4\pi}{\eps}$ \proves\ that
\begin{equation}
\label{eq:d-dim:prod:param}
\begin{split}
\param\pr{\f}
&\leq
\param(\g_1)+\param(\g_2)+3+2N
\\&\leq \pr*{4107 d^3+2048 d^3\ceil{\log_2(R)}+512 d^2 (\log_2(\eps^{-1})+1)+514d^3\log_2(\scl)}
\\&\quad+\pr*{(18n+24)N^2}+3+2N
\\&\leq 4107 d^3+2048 d^3\ceil{\log_2(R)}+512 d^2 \log_2(\eps^{-1})+512d^2+514d^3\log_2(\scl)
\\&\quad+\pr*{18\max\{d\log_2(\scl R),1\}+24}16\pi^2\eps^{-2}+3+8\pi\eps^{-1}
\\&\leq 4107 d^3+2048 d^3\ceil{\log_2(R)}+512 d^2 \log_2(\eps^{-1})+512d^2+514d^3\log_2(\scl)
\\&\quad+\pr*{36\max\{\log_2( R),\log_2(\scl),1\}d+24}160\eps^{-2}+3+26\eps^{-1}
\\&= 4107 d^3+2048 d^3\ceil{\log_2(R)}+512 d^2 \log_2(\eps^{-1})+512d^2+514d^3\log_2(\scl)
\\&\quad+5760\max\{\log_2( R),\log_2(\scl)\}d\eps^{-2}+3840\eps^{-2}+3+26\eps^{-1}
\\&\leq \prb{4107 d^3+2048 d^3+512 d^2+512d^2+514d^3
\\&\quad
+5760d+3840+3+26}\max\{\ceil{\log_2( R)},\log_2(\scl)\}\eps^{-2}
\\&= \prb{6669 d^3+1024 d^2+5760d+3869}\max\{\ceil{\log_2( R)},\log_2(\scl)\}\eps^{-2}
\\&\leq \prb{6669+512+1440+484}d^3\max\{\ceil{\log_2( R)},\log_2(\scl)\}\eps^{-2}
\\&\leq 9105d^3\max\{\ceil{\log_2( R)},\log_2(\scl)\}\eps^{-2}
.
\end{split}
\end{equation}
\Moreover 
\cref{it:d_product_size}, 
\cref{it:sin_size:prop},
and 
\cref{Prop:identity_representation:prop} \prove\ that
\begin{equation}
\begin{split}
\size(\f)=\max\pR*{\size(\g_2),\size(\g_1)}\leq\max\pR{2,4}=4.
\end{split}
\end{equation}
Combining this, 
\lref{eq:loc:length:ref},
\eqref{eq:comp:dims2}, and 
\eqref{eq:d-dim:prod:param} with 
\eqref{eq:approx:g:prod}
 \proves[ei]\ \cref{prop:new_continuous_d,prop:new_approx_d,prop:new_length_d,prop:new_cost_d,prop:new_size_d,prop:new_hid_d}.
The proof of \cref{final:approximation} is thus complete.
\cfload
\end{proof}

\cfclear
\begin{athm}{lemma}{final:approximation:streched}
    Let $d \in \N$, $\kappa,R\in(0,\infty)$, $\varepsilon \in (0,\kappa)$, $\gamma\in(0,1]$, $\scl\in[1,\infty)$,
		$\tgt\in C(\R,\R)$ satisfy for all $x,y\in\R$, $k\in\Z$ that $\vass{\tgt(0)}\leq 2\kappa$, $\tgt(x+2k\pi)=\tgt(x)$, and $\vass{\tgt(x)-\tgt(y)}\leq\kappa\vass{x-y}$.
    Then there exists $\mathscr{f} \in \ANNSET$ such that

    \begin{enumerate}[(i)]
    \item \label{prop:new_continuous_d:streched} it holds 
    that $\realisation(\mathscr{f}) \in C(\R^d,\R)$,

    \item\label{prop:new_approx_d:streched} it holds that
    $
    \sup_{x = (x_1,\ldots, x_{d}) \in [-\radius, \radius]^{d}}\vass[\big]{ g\prb{\gamma\scl^d\sprod_{i = 1}^d x_i} - \realisation\pr{\mathscr{f}} (x)} \leq \varepsilon$,

    \item\label{prop:new_hid_d:streched} it holds 
that
$\singledims_{\hidlengthANN(\f)}(\f)=2$,

    \item\label{prop:new_cost_d:streched} it holds 
    that
    $\paramANN(\mathscr{f}) \leq   9123\max\{1,\ceil{\log_2( R)},\log_2(\scl)\}\max\{1,\kappa^3\}d^3\eps^{-2}
    $,
		
   \item\label{prop:new_length_d:streched} it holds 
	that
	$\lengthANN(\mathscr{f}) \leq 19 \max\{1,\ceil{\log_2( R)},\ceil{\log_2(\scl)}\}\max\{1,\kappa^2\}d^2\eps^{-1}
	$, and

		\item\label{prop:new_size_d:streched} it holds 
    that
    $\size(\mathscr{f}) \leq 4
    $
    \end{enumerate}
    \cfout. 
\end{athm}

\begin{proof}[Proof of \cref{final:approximation:streched}]
Throughout this proof 
let $L\in\N$ satisfy $L=\max\{1,\ceil{\log_2(\kappa)}\}$, and
let $f\in C(\R,\R)$ satisfy for all $x\in\R$ that $f(x)=\kappa^{-1}\tgt(x)$.
\Nobs that \cref{final:approximation} (applied with 
$d \curvearrowleft d$,
$\eps \curvearrowleft \frac{\varepsilon}{\kappa}$,
$R \curvearrowleft R$,
$\gamma \curvearrowleft \gamma$, 
$\scl \curvearrowleft \scl$, 
$g \curvearrowleft f$
in the notation of \cref{final:approximation}) \proves\ that there exists $\g_1\in\ANNSET$ which satisfies that
    \begin{enumerate}[(I)]
    \item \label{it:to:stretch:prop:new_continuous_d} it holds 
    that $\realisation(\g_1) \in C(\R^d,\R)$,

    \item\label{it:to:stretch:prop:new_approx_d} it holds that
    $
    \sup_{x = (x_1,\ldots, x_{d}) \in [-\radius, \radius]^{d}}\vass[\big]{f\prb{\gamma\scl^d\sprod_{i = 1}^d x_i} - \realisation\pr{\g_1} (x)} \leq \frac{\varepsilon}{\kappa}$,

    \item\label{it:to:stretch:prop:new_hid_d} it holds 
    that
		$\singledims_{\hidlengthANN(\g_1)}(\g_1)=2$,

    \item\label{it:to:stretch:prop:new_length_d} it holds 
that
$\lengthANN(\g_1) \leq  16 \max\{1,\ceil{\log_2( R)},\ceil{\log_2(\scl)}\}d^2\max\{1,\log_2(\kappa\eps^{-1})\}
$,

    \item\label{it:to:stretch:prop:new_cost_d} it holds 
    that
    $\paramANN(\g_1) \leq    9105\max\{1,\ceil{\log_2( R)},\log_2(\scl)\}d^3\kappa^2\eps^{-2}
    $, and

		\item\label{it:to:stretch:prop:new_size_d} it holds 
    that
    $\size(\g_1) \leq 4
    $
    \end{enumerate}
\cfload. \Nobs that \cref{gen:scaling:networks} (applied with 
$\beta \curvearrowleft \kappa$,
$L \curvearrowleft L$
in the notation of \cref{gen:scaling:networks}) \proves\ that there exists $\g_2\in\ANNSET$ which satisfies that
\begin{enumerate}[(A)]
\item{
\label{it:to:stretch:gen:scaling:realisation}
it holds for all $x\in\R$ that $(\realisation(\g_2))(x)=\kappa x$,
}
\item{
\label{it:to:stretch:gen:scaling:dims}
it holds that $\dims(\g_2)=(1,2,2,\ldots,2,1)\in\N^{L+2}$, and
}
\item{
\label{it:to:stretch:gen:scaling:size}
it holds that $\size(\g_2)\leq 2$.
}
\end{enumerate}
\Nobs that
\cref{it:to:stretch:prop:new_approx_d},
\cref{it:to:stretch:gen:scaling:realisation},
\cref{Lemma:PropertiesOfCompositions_n2}, and
\cref{Prop:identity_representation} \prove\ that for all $x=(x_1,\ldots,x_d)\in[-R,R]^d$ it holds that $\realisation(\g_2\bullet\compANN{\ReLUidANN{1}}{\g_1})\in C(\R^d,\R)$ and
\begin{equation}
\label{eq:stratched:prod:approx}
\begin{split}
\textstyle
\vass{\pr{\realisation(\g_2\bullet\compANN{\ReLUidANN{1}}{\g_1})}(x)-g\prb{\gamma\scl^d\prod_{i = 1}^d x_i}}
&=\vass{\pr{\PR{\realisation(\g_2)}\circ\PR{\realisation(\g_1)}}(x)-g\prb{\gamma\scl^d\textstyle\prod_{i = 1}^d x_i}}
\\&=\kappa\vass{\pr{\realisation(\g_1)}(x)-f\prb{\gamma\scl^d\textstyle\prod_{i = 1}^d x_i}}
\leq\eps.
\end{split}
\end{equation}
\Nobs that 
\cref{it:to:stretch:prop:new_length_d}, 
\cref{it:to:stretch:gen:scaling:dims}, 
\cref{Lemma:PropertiesOfCompositions_n2}, 
\cref{Prop:identity_representation}, and the fact that for all $x\in(1,\infty)$, $y\in(0,\infty)$ it holds that $\max\{1,\log_2(x)\}\leq x$ and $\ceil{\log_2(y)}\leq y+1$ \prove\
\begin{equation}
\label{eq:stratched:prod:length}
\begin{split}
&\lengthANN(\g_2\bullet\compANN{\ReLUidANN{1}}{\g_1})
\\&=\lengthANN(\g_2)+2+\lengthANN(\g_1)-2
\\&=\lengthANN(\g_1)+L+1
\\&\leq 16 \max\{1,\ceil{\log_2( R)},\ceil{\log_2(\scl)}\}d^2\max\{1,\log_2(\kappa\eps^{-1})\} + \max\{1,\ceil{\log_2(\kappa)}\}+1
\\&\leq 16 \max\{1,\ceil{\log_2( R)},\ceil{\log_2(\scl)}\}d^2\kappa\eps^{-1} +\kappa+2
\\&\leq 16 \max\{1,\ceil{\log_2( R)},\ceil{\log_2(\scl)}\}d^2\kappa\eps^{-1} +\kappa^2\eps^{-1}+2\kappa\eps^{-1}
\\&\leq 19 \max\{1,\ceil{\log_2( R)},\ceil{\log_2(\scl)}\}\max\{1,\kappa^2\}d^2\eps^{-1}
\end{split}
\end{equation}
\Nobs that 
\cref{it:to:stretch:prop:new_continuous_d},
\cref{it:to:stretch:prop:new_length_d}, 
\cref{it:to:stretch:gen:scaling:dims}, and
\cref{Lemma:PropertiesOfCompositions_n2} \hence \prove\ that for all $k\in\N_0\cap[0,\lengthANN(\g_2\bullet\compANN{\ReLUidANN{1}}{\g_1})]$ it holds that
\begin{equation}
\label{eq:stratched:prod:dims}
\singledims_k(\g_2\bullet\compANN{\ReLUidANN{1}}{\g_1})=
\begin{cases}
\singledims_k(\g_1) &\colon k\in\N_0\cap[0,\lengthANN(\g_1))\\
2 &\colon k\in\N\cap[\lengthANN(\g_1),\lengthANN(\g_1)+L+1)\\
1 &\colon k=\lengthANN(\g_1)+L+1.
\end{cases}
\end{equation}
This,
 \cref{it:to:stretch:prop:new_hid_d}, 
 and the fact that for all $x\in(0,\infty)$ it holds that $\ceil{\log_2(x)}\leq x+1$ \prove\ that
\begin{equation}
\label{eq:stratched:prod:param}
\begin{split}
\param(\g_2\bullet\compANN{\ReLUidANN{1}}{\g_1})
&=\sum_{k=1}^{\lengthANN(\g_2\bullet\compANN{\ReLUidANN{1}}{\g_1})} \singledims_k(\g_2\bullet\compANN{\ReLUidANN{1}}{\g_1})(\singledims_{k-1}(\g_2\bullet\compANN{\ReLUidANN{1}}{\g_1})+1)
\\&=\PR*{\sum_{k=1}^{\lengthANN(\g_1)-1} \singledims_k(\g_1)(\singledims_{k-1}(\g_1)+1)}
+\PR*{\sum_{k=\lengthANN(\g_1)}^{\lengthANN(\g_1)+L} 2(2+1)}
+1(2+1)
\\&=\param(\g_1)-\singledims_{\lengthANN(\g_1)}(\g_1)(\singledims_{\hidlengthANN(\g_1)}(\g_1)+1)+6(L+1)+3
\\&=\param(\g_1)+6(L+1)
\\&\leq 9105\max\{1,\ceil{\log_2( R)},\log_2(\scl)\}d^3\kappa^2\eps^{-2}+6(\max\{1,\ceil{\log_2(\kappa)}\}+1)
\\&\leq 9105\max\{1,\ceil{\log_2( R)},\log_2(\scl)\}d^3\kappa^2\eps^{-2}+6\kappa+12
\\&\leq 9117\max\{1,\ceil{\log_2( R)},\log_2(\scl)\}d^3\kappa^2\eps^{-2}+6\kappa^3\eps^{-2}
\\&\leq 9123\max\{1,\ceil{\log_2( R)},\log_2(\scl)\}\max\{\kappa^3,1\}d^3\eps^{-2}
.
\end{split}
\end{equation}
\Moreover 
\cref{it:to:stretch:prop:new_size_d},
\cref{it:to:stretch:gen:scaling:size}, and
\cref{Prop:identity_representation:prop} \prove\ that
\begin{equation}
\size(\g_2\bullet\compANN{\ReLUidANN{1}}{\g_1})
=\max\{\size(\g_2),\size(\g_1)\}\leq\max\{2,4\}=4
.
\end{equation}
Combining this, 
\eqref{eq:stratched:prod:approx},
\eqref{eq:stratched:prod:length}, 
\eqref{eq:stratched:prod:dims}, and
\eqref{eq:stratched:prod:param}
 \proves[ei]\
\cref{prop:new_continuous_d:streched,prop:new_length_d:streched,prop:new_size_d:streched,prop:new_cost_d:streched,prop:new_hid_d:streched,prop:new_approx_d:streched}. 
The proof of \cref{final:approximation:streched} is thus complete.
\end{proof}

\cfclear
\begin{athm}{cor}{final:approximation:imp}
Let $a\in\R$, $b\in[a,\infty)$, $d\in\N$, $\kappa,\mathfrak{c}\in(0,\infty)$, $\eps\in(0,\kappa)$, $\gamma\in(0,1]$, $\scl\in[1,\infty)$ satisfy 
$\mathfrak{c}\geq
9885 \ceil{\log_2\pr*{\max\{2,\abs{a},\abs{b},\scl\}}}\max\{1,\kappa^3\}$, 
let $g\colon\R\to\R$ satisfy for all $x,y\in\R$, $k\in\Z$ that 
$\vass{\tgt(0)}\leq 2\kappa$, 
$\tgt(x+2k\pi)=\tgt(x)$, 
and $\vass{\tgt(x)-\tgt(y)}\leq\kappa\vass{x-y}$,
and let $f\colon\R^d\to\R$ satisfy for all $x=(x_1,\ldots,x_d)\in[a,b]^d$ that $f(x)=g\prb{\gamma\scl^d\sprod_{i = 1}^d x_i}$
\cfload. Then 
\begin{equation}
\label{eq:final:prod:cost:representation}
\begin{split}
	&\min \pr*{ \pR*{ p \in \N \colon \PR*{\!\!
	\begin{array}{c}
	    \exists\, \mathscr{f} \in \ANNSET \colon
			(\paramANN(\mathscr{f})\leq p)\land
			(\lengthANN(\mathscr{f})\leq \mathfrak{c}d^2\eps^{-1})
			\land{}\\
			(\size(\mathscr{f})\leq 1)\land{}(\realisation(\mathscr{f}) \in C(\R^d,\R))\land{}\\
			(\sup_{x\in[a,b]^d}\vass{(\realisation(\mathscr{f}))(x)-f(x)} \leq \varepsilon)\\
    \end{array} \!\! }
}
\cup\{\infty\}
 }\leq  \mathfrak{c}d^3\eps^{-2}
\end{split}
\end{equation}
    \cfout. 
\end{athm}

\begin{proof}[Proof of \cref{final:approximation:imp}]
Throughout this proof 
assume without loss of generality that $0<\max\{\abs{a},\abs{b}\}$ and let $c\in\R$ satisfy $c=\ceil{\log_2\pr*{\max\{2,\abs{a},\abs{b},\scl\}}}$.
\Nobs that \cref{final:approximation:streched} (applied with 
$d \curvearrowleft d$,
$R \curvearrowleft \max\{\abs{a},\abs{b}\}$,
$\kappa \curvearrowleft \kappa$,
$\eps \curvearrowleft \eps$,
$\gamma \curvearrowleft \gamma$
 in the notation of \cref{final:approximation:streched}) \proves\ that there exists $\g\in\ANNSET$ which satisfies that
\begin{enumerate}[(i)]
    \item  it holds 
    that $\realisation(\g) \in C(\R^d,\R)$,

    \item it holds that 
    $\sup_{x  \in [a,b]^d}\vass[\big]{(\realisation\pr{\g}) (x)-f(x) } \leq \varepsilon$,

    \item it holds 
that 
$\singledims_{\hidlengthANN(\g)}(\g)=2$,

    \item it holds 
    that 
		$
		\lengthANN(\g) \leq 19 c\max\{1,\kappa^2\} d^2\eps^{-1}
		$
    ,

    \item it holds 
    that
    $\paramANN(\g) \leq    9123 c\max\{1,\kappa^3\}d^3\eps^{-2}$, and
		
		\item it holds 
    that
    $\size(\g) \leq 4
    $
    \end{enumerate}
		\cfload. \Nobs that \cref{lemma:network:halfed:paramsize} (applied with 
$d \curvearrowleft d$,
$\f \curvearrowleft \g$,
$\dnsc \curvearrowleft 4^{-1}$
 in the notation of \cref{lemma:network:halfed:paramsize}) \hence \proves\ that there exists $\h\in\ANNSET$ which satisfies that 
\begin{enumerate}[(I)]
    \item  \label{it:final:approximation:imp:1} it holds 
    that $\realisation(\h)=\realisation(\g) \in C(\R^d,\R)$,

    \item \label{it:final:approximation:imp:2} it holds that 
    $\sup_{x  \in [a,b]^d}\vass[\big]{(\realisation\pr{\h}) (x)-f(x)} \leq \varepsilon$,

    \item \label{it:final:approximation:imp:3} it holds 
    that 
		$
		\lengthANN(\h) \leq 3(19 c\max\{1,\kappa^2\} d^2\eps^{-1})+1\leq (57c\max\{1,\kappa^2\}+\kappa)d^2\eps^{-1}
		$,

    \item \label{it:final:approximation:imp:4} it holds 
    that
    $\paramANN(\h) \leq 9123 c\max\{1,\kappa^3\}d^3\eps^{-2}+2+40(19 c \max\{1,\kappa^2\} d^2\eps^{-1})$, and
		
	\item \label{it:final:approximation:imp:5} it holds 
    that
    $\size(\h) \leq 1
    $.
    \end{enumerate}
    \Nobs that \cref{it:final:approximation:imp:4} and the assumption that $\eps<\kappa$ and that $\mathfrak{c}\geq 9885 c\max\{1,\kappa^3\}$ \prove\ that
    \begin{equation}
    	\begin{split}
    	\paramANN(\h) &\leq 9123 c\max\{1,\kappa^3\}d^3\eps^{-2}+2+40(19 c \max\{1,\kappa^2\} d^2\eps^{-1})
    	\\&\leq
    	\pr{9123 c\max\{1,\kappa^3\}+2\kappa^2+760c\max\{1,\kappa^3\}}d^3\eps^{-2}
    	\\&\leq
    	9885 c\max\{1,\kappa^3\}d^3\eps^{-2}
    	\leq \mathfrak{c}d^3\eps^{-2}.
    	\end{split}
    \end{equation}
Combining this with the fact that $57c\max\{1,\kappa^2\}+\kappa\leq\mathfrak{c}$ and \cref{it:final:approximation:imp:1,it:final:approximation:imp:2,it:final:approximation:imp:3,it:final:approximation:imp:5}
\proves\ \eqref{eq:final:prod:cost:representation}. 
The proof of \cref{final:approximation:imp} is thus complete.
\end{proof}

\subsection{Upper bounds for approximations of certain smooth and bounded functions}
\label{subsection:upper_high_osz}

\cfclear
\begin{athm}{lemma}{final:approximation2}
    Let $a\in\R$, $b\in[a,\infty)$, $d \in \N$, $\gamma\in(0,\infty)$, $\varepsilon \in (0,1)$,  $\tgt\in C(\R,\R)$ satisfy for all $x,y\in\R$, $k\in\Z$ that $\vass{\tgt(0)}\leq 2$, $\tgt(x+2k\pi)=\tgt(x)$, and $\vass{\tgt(x)-\tgt(y)}\leq\vass{x-y}$.
    Then there exists $\mathscr{f} \in \ANNSET$ such that

    \begin{enumerate}[(i)]
    \item \label{prop:sum_d} it holds 
    that $\realisation(\mathscr{f}) \in C(\R^d,\R)$,

    \item\label{prop:sum_approx_d} it holds that
    $\sup_{x = (x_1, \ldots, x_d) \in [a,b]^d}\vass[\big]{\tgt\prb{\gamma 2^d\ssum_{i = 1}^d x_i} - (\realisation\pr{\mathscr{f}}) (x)} \leq \varepsilon$,

				\item\label{prop:sum_hid_d} it holds that
		$\singledims_1(\f)=\singledims_{\hidlengthANN(\mathscr{f})}(\f)=2$,
		
		\item\label{prop:sum_length_d} it holds that
    $\lengthANN(\mathscr{f}) \leq 10d\ceil{\log_2\pr*{\max\{1,\gamma\}\max\{2,\vass{a},\vass{b}\}}}
    $,

    \item\label{prop:sum_cost_d} it holds 
    that
    $\paramANN(\mathscr{f}) \leq 2239 \ceil{\log_2\pr*{\max\{1,\gamma\}\max\{2,\vass{a},\vass{b}\}}}d\eps^{-2}
    $, and
		
		\item\label{prop:sum_size_d} it holds 
    that
    $\size(\mathscr{f}) \leq 2
    $
    \end{enumerate}
    \cfout. 
\end{athm}

\begin{proof}[Proof of \cref{final:approximation2}]
Throughout this proof 
assume without loss of generality that $d> 1$ (cf.~\cref{1:dim:scaled:test} applied with 
$R\curvearrowleft\max\{2,\vass{a},\vass{b}\}$,
$g\curvearrowleft g$,
$\gamma\curvearrowleft\gamma\pr{\max\{1,\gamma\}}^{-1}$,
$\scl\curvearrowleft 2\max\{1,\gamma\}$,
$\eps\curvearrowleft\eps$
in the notation of \cref{1:dim:scaled:test}),
let $n,N\in\N$ satisfy $N=\ceil{\tfrac{2\pi}{\eps}}-1$ and $n =2d\ceil{\log_2\pr*{\max\{1,\gamma\}\max\{2,\vass{a},\vass{b}\}}}$,
 and
let $\g_1\in\ANNSET$ satisfy
\begin{equation}
\label{def:summing:netw:post2:split}
\begin{split}
	&\g_1 = 
		\pr*{
			\begin{pmatrix}
			1 & 1 & \cdots & 1
			\end{pmatrix},\!
			\,0
		}
		\in (\R^{1 \times d} \times \R^1)
\end{split}
\end{equation}
\cfload.
\Nobs that \cref{gen:scaling:networks} (applied with
$\beta \curvearrowleft \gamma 2^d$,
$L \curvearrowleft d+\ceil{\log_2(\max\{1,\gamma\})}$
in the notation of \cref{gen:scaling:networks}) \proves\ that there exists $\g_2\in\ANNSET$ which satisfies that
\begin{enumerate}[(I)]
\item{
\label{it:scaling:realisation:sum}
it holds for all $x\in\R$ that $(\realisation(\g_2))(x)=\gamma 2^d x$,
}
\item{
\label{it:scaling:dims:sum}
it holds that $\dims(\g_2)=(1,2,2,\ldots,2,1)\in\N^{d+\ceil{\log_2(\max\{1,\gamma\})}+2}$,
}
\item{
\label{it:scaling:size:sum}
it holds that  $\insize(\g_2)\leq1$, $\outsize(\g_2)\leq2$, and $\size(\g_2)\leq2$, and
}
\item{
\label{it:scaling:param:sum}
it holds that $\param(\g_2)\leq6(d+\ceil{\log_2(\max\{1,\gamma\})})+1$.
}
\end{enumerate}
		\cfload. \Moreover 
		\cref{sin:approximation} (applied with 
$\eps \curvearrowleft \eps$,
$n \curvearrowleft n$,
$N \curvearrowleft N$,
$g \curvearrowleft g$
in the notation of \cref{sin:approximation}) \proves\ that there exists $\g_3\in\ANNSET$ which satisfies that
		     \begin{enumerate}[(A)]
        \item {\label{it:sin_layers:prop:sum} it holds that $\lengthANN(\g_3)\leq n+11$,
}
    \item {\label{it:sin_dims:prop:sum}
    it holds that $\singledims_0(\g_3)=\singledims_{\lengthANN(\g_3)}(\g_3)=1$, 
		$\singledims_1(\g_3)\leq 2N$, and $\singledims_{\hidlengthANN(\g_3)}(\g_3)=2$,}
    \item{\label{it:sin_lipschitz:prop:sum} it holds for all
    $x,y \in \R$ that 
    $\vass[\big]{\pr{\realisation(\g_3)} (x) - \pr{\realisation(\g_3)} (y)} \leq \vass{x-y}$,}
    \item{\label{it:sin_approx:prop:sum} it holds that
    $\sup_{x \in [-2^n\pi,2^n\pi]}\vass[\big]{\tgt(x) - \pr{\realisation(\g_3)} (x)} \leq \eps$,}
    \item{\label{it:sin_start:prop:sum} it holds for all $x\in[2^n\pi,\infty)$ that 
    $\pr{\realisation(\g_3)} (-x)=\pr{\realisation(\g_3)} (x) =\tgt(0)$,}
		\item{\label{it:sin_size:prop:sum} it holds 
    that
    $\size(\g_3) \leq 2
    $, and}
    \item{\label{it:sin_cost:prop:sum} it holds 
    that
    $\param(\g_3) \leq  (18n+24)N^2
    $.}
    \end{enumerate}
		Next let $\f\in\ANNSET$ satisfy
		\begin{equation}
		\label{definition:of:final:imp:approximator:sum}
		\f=\compANN{\g_3}{\ReLUidANN{1}}\bullet\g_2\bullet\g_1
		\end{equation}
\cfload[.]\Nobs that 
\eqref{def:summing:netw:post2:split},
\eqref{definition:of:final:imp:approximator:sum}, 
\cref{it:scaling:dims:sum}, 
\cref{it:sin_layers:prop:sum}, 
\cref{Lemma:PropertiesOfCompositions_n2}, and
\cref{Prop:identity_representation} \prove\ that
\begin{equation}
\label{eq:length:g:sum}
\begin{split}
\lengthANN(\f)
&=\lengthANN(\g_3)+\lengthANN(\ReLUidANN{1})+\lengthANN(\g_2)+\lengthANN(\g_1)-3
\\&=\lengthANN(\g_3)+\lengthANN(\g_2)
\\&\leq(n+11)+(d+\ceil{\log_2(\max\{1,\gamma\})}+1)
\\&= d+2d\ceil{\log_2\pr*{\max\{1,\gamma\}\max\{2,\vass{a},\vass{b}\}}}+\ceil{\log_2(\max\{1,\gamma\})}+12
\\&\leq (3d+1+12)\ceil{\log_2\pr*{\max\{1,\gamma\}\max\{2,\vass{a},\vass{b}\}}}
\\&\leq 10d\ceil{\log_2\pr*{\max\{1,\gamma\}\max\{2,\vass{a},\vass{b}\}}}
.
\end{split}
\end{equation}
\Nobs that 
\eqref{definition:of:final:imp:approximator:sum}, 
\eqref{def:summing:netw:post2:split},
\cref{it:scaling:realisation:sum},
\cref{Lemma:PropertiesOfCompositions_n2}, 
and 
\cref{Prop:identity_representation} \prove\ that for all $x=(x_1,\ldots,x_d)\in\R^d$ it holds that
\begin{equation}
\begin{split}
\label{eq:specific:realisation2}
(\realisation(\f))(x)
=\pr*{\PR*{\realisation(\g_3)}\circ\PR*{\realisation(\g_2)}\circ\PR*{\realisation(\g_1)}}(x)
&=\pr*{\PR*{\realisation(\g_3)}\circ\PR*{\realisation(\g_2)}}\prbb{\ssum_{i=1}^d x_i}
\\&=(\realisation(\g_3))\prbb{\gamma 2^d\ssum_{i=1}^d x_i}.
\end{split}
\end{equation}
This, 
\cref{it:sin_approx:prop:sum}, 
and the fact that for all $x_1,x_2,\ldots,x_d\in [a,b]$ it holds that $\vass{\gamma 2^d\ssum_{i = 1}^d x_i}\leq \gamma d2^d \max\{2,\vass{a},\vass{b}\}\leq 2^{2d\ceil{\log_2\pr*{\max\{1,\gamma\}\max\{2,\vass{a},\vass{b}\}}}}\pi= 2^n\pi$ \prove\ that for all $x=(x_1,\ldots\allowbreak,x_d)\in [a,b]^d$ it holds that
\begin{equation}
\begin{split}
\label{eq:approx:g:sum}
&\vass[\bigg]{\tgt\prbb{\gamma 2^d\ssum_{i = 1}^d x_i} - (\realisation\pr{\mathscr{f}}) (x)}
=\vass[\bigg]{\tgt\prbb{ \gamma2^d\ssum_{i = 1}^d x_i} - (\realisation(\g_3))\prbb{\gamma 2^d\ssum_{i=1}^d x_i}}
\leq \eps.
\end{split}
\end{equation}
\Moreover 
\eqref{def:summing:netw:post2:split},
\eqref{definition:of:final:imp:approximator:sum}, 
\eqref{eq:length:g:sum}, 
\cref{it:scaling:dims:sum}, 
\cref{it:sin_dims:prop:sum}, and
\dimsbene\ \prove\ that for all $k\in\N_0\cap[0,\lengthANN(\f)]$ it holds that
\begin{equation}
\label{eq:comp:dims2:sum}
\begin{split}
\singledims_k(\f)=
\begin{cases}
d &\colon k=0\\
\singledims_k(\g_2) &\colon k\in\N\cap(0,\lengthANN(\g_2))\\
2 &\colon k=\lengthANN(\g_2)\\
\singledims_{k-\lengthANN(\g_2)}(\g_3) &\colon k\in\N\cap(\lengthANN(\g_2), \lengthANN(\f)].\\
\end{cases}
\end{split}
\end{equation}
This,
\cref{it:scaling:dims:sum}, and 
\cref{it:sin_dims:prop:sum} \prove\ that
\begin{equation}
\begin{split}
&\param\pr{\f}
\\&=\sum_{k=1}^{\lengthANN(\f)} \singledims_k(\f)(\singledims_{k-1}(\f)+1)
\\&=
\singledims_1(\g_2)(d+1)
+\PR*{\sum_{k=2}^{\lengthANN(\g_2)-1} \singledims_k(\g_2)(\singledims_{k-1}(\g_2)+1)}
+2(\singledims_{\lengthANN(\g_2)-1}(\g_2)+1)
\\&\quad
+\singledims_{1}(\g_3)(2+1)
+\PR*{\sum_{k=2}^{\lengthANN(\g_3)} \singledims_k(\g_3)(\singledims_{k-1}(\g_3)+1)}
\\&=
\singledims_1(\g_2)(d-1)
+\PR*{\sum_{k=1}^{\lengthANN(\g_2)} \singledims_k(\g_2)(\singledims_{k-1}(\g_2)+1)}
+\singledims_{\lengthANN(\g_2)-1}(\g_2)+1
\\&\quad
+\singledims_{1}(\g_3)
+\PR*{\sum_{k=1}^{\lengthANN(\g_3)} \singledims_k(\g_3)(\singledims_{k-1}(\g_3)+1)}
\\&\leq 2d-2+\param(\g_2)+3+2N+\param(\g_3)=\param(\g_2)+\param(\g_3)+2N+2d+1
.
\end{split}
\end{equation}
Combining this, 
\cref{it:scaling:param:sum}, 
and \cref{it:sin_cost:prop:sum} with the fact that for all $x\in(0,1)$ it holds that $-\log_2\pr*{x}\leq x^{-2}$ and $\pi^2\leq 10$ with the assumption that  
$n= 2d\ceil{\log_2\pr*{\max\{1,\gamma\}\max\{2,\vass{a},\vass{b}\}}}$, 
$\varepsilon \in (0,1)$, 
$d\geq 2$, and $N \leq\frac{2\pi}{\eps}$ \proves\ that
\begin{equation}
\label{eq:d-dim:sum:param}
\begin{split}
&\param\pr{\f}
\\&\leq
\param(\g_2)+\param(\g_3)+2N+2d+1
\\&\leq
(6(d+\ceil{\log_2(\max\{1,\gamma\})})+1)+(18n+24)N^2+2N+2d+1
\\&\leq
(18n+34)4\pi^2\eps^{-2}+6\ceil{\log_2(\max\{1,\gamma\})}+4\pi\eps^{-1}+8d+2
\\&\leq
((18n+34)40+6\ceil{\log_2(\max\{1,\gamma\})}+13+8d+2)\eps^{-2}
\\&=
(1440d\ceil{\log_2\pr*{\max\{1,\gamma\}\max\{2,\vass{a},\vass{b}\}}}+6\ceil{\log_2(\max\{1,\gamma\})}+8d+1375)\eps^{-2}
\\&\leq
(1440d\ceil{\log_2\pr*{\max\{1,\gamma\}\max\{2,\vass{a},\vass{b}\}}}+3d\ceil{\log_2(\max\{1,\gamma\})}+8d+688d)\eps^{-2}
\\&\leq
(1440+3+8+688)\ceil{\log_2\pr*{\max\{1,\gamma\}\max\{2,\vass{a},\vass{b}\}}}d\eps^{-2}
\\&=2139 \ceil{\log_2\pr*{\max\{1,\gamma\}\max\{2,\vass{a},\vass{b}\}}}d\eps^{-2}
.
\end{split}
\end{equation}
\Moreover 
\eqref{def:summing:netw:post2:split},
\eqref{definition:of:final:imp:approximator:sum}, 
\cref{it:scaling:size:sum},
\cref{it:sin_size:prop:sum}, 
\cref{lem:sizecomp},
and 
\cref{Prop:identity_representation:prop} \prove\ that
\begin{equation}
\begin{split}
\size(\f)
=\max\pR*{\size(\g_3),\size(\g_2\bullet\g_1)}
&\leq \max\pR*{\size(\g_3),\size(\g_2),\size(\g_1),\insize(\g_2)(\outsize(\g_1)+1)}
\\&\leq\max\pR{2,2,1,1(1+1)}=2.
\end{split}
\end{equation}
Combining this with 
\eqref{eq:length:g:sum}, 
\eqref{eq:comp:dims2:sum},
\eqref{eq:d-dim:sum:param}, and 
\eqref{eq:approx:g:sum}  \proves[ei]\ 
\cref{prop:sum_d,prop:sum_approx_d,prop:sum_size_d,prop:sum_length_d,prop:sum_cost_d,prop:sum_hid_d}.
The proof of \cref{final:approximation2} is thus complete.
\end{proof}

\cfclear
\begin{athm}{cor}{cor:rescaled:final:approximation2}
    Let $a\in\R$, $b\in[a,\infty)$, $d \in \N$, $\gamma,\kappa\in(0,\infty)$, $\varepsilon \in (0,\kappa)$, $\tgt\in C(\R,\R)$ satisfy for all $x,y\in\R$, $k\in\Z$ that $\vass{\tgt(0)}\leq 2\kappa$, $\tgt(x+2k\pi)=\tgt(x)$, and $\vass{\tgt(x)-\tgt(y)}\leq\kappa\vass{x-y}$.
    Then there exists $\mathscr{f} \in \ANNSET$ such that

    \begin{enumerate}[(i)]
    \item \label{cor:rescaled:prop:sum_d} it holds 
    that $\realisation(\mathscr{f}) \in C(\R^d,\R)$,

    \item\label{cor:rescaled:prop:sum_approx_d} it holds that
    $\sup_{x = (x_1,\ldots, x_d) \in [a,b]^d}\vass[\big]{g\prb{\gamma 2^d\ssum_{i = 1}^d x_i}- (\realisation\pr{\mathscr{f}}) (x)} \leq \varepsilon$,

				\item\label{cor:rescaled:prop:sum_hid_d} it holds that
		$
		\singledims_{\hidlengthANN(\f)}(\f)=2
		$,
		
		\item\label{cor:rescaled:prop:sum_length_d} it holds that
    $\lengthANN(\mathscr{f}) 
		\leq 17\ceil{\log_2\pr*{\max\{1,\gamma\}\max\{2,\vass{a},\vass{b}\}}}\max\{\kappa,1\}d
    $,

    \item\label{cor:rescaled:prop:sum_cost_d} it holds 
    that
    $\paramANN(\mathscr{f}) \leq 2251\ceil{\log_2\pr*{\max\{1,\gamma\}\max\{2,\vass{a},\vass{b}\}}}\max\{\kappa^3,1\} d\eps^{-2}
    $, and   
		
		\item\label{cor:rescaled:prop:sum_size_d} it holds 
    that
    $\size(\mathscr{f}) \leq 2
    $
    \end{enumerate}
    \cfout. 
\end{athm}
\begin{proof}[Proof of \cref{cor:rescaled:final:approximation2}]
Throughout this proof let $n,R\in\N$, $f\in C(\R,\R)$ satisfy for all $x\in\R$ that $n=\ceil*{\log_2(\max\{\kappa,2\})}$, $R= \ceil{\log_2\pr*{\max\{1,\gamma\}\max\{2,\vass{a},\vass{b}\}}}$, and $f(x)=\kappa^{-1}\tgt(x)$ \cfload.
\Nobs that \cref{final:approximation2} (applied with
$a \curvearrowleft a$,
$b \curvearrowleft b$,
$d \curvearrowleft d$,
$\gamma \curvearrowleft \gamma$,
$\eps \curvearrowleft \frac{\eps}{\kappa}$, 
$g \curvearrowleft f$
in the notation of \cref{final:approximation2}) \proves\ that there exists $\g_1\in\ANNSET$ which satisfies that
\begin{enumerate}[(I)]
    \item \label{it:rescaled:prop:sum_d} it holds 
    that $\realisation(\g_1) \in C(\R^d,\R)$,

    \item\label{it:rescaled:prop:sum_approx_d} it holds that
    $\sup_{x = (x_1, \ldots, x_d) \in [a,b]^d}\vass[\big]{f\prb{\gamma 2^d\ssum_{i = 1}^d x_i}- (\realisation\pr{\g_1}) (x)} \leq \frac{\varepsilon}{\kappa}$,
		
		\item\label{it:rescaled:prop:sum_hid_d} it holds that 
    	$\singledims_1(\g_1)=\singledims_{\hidlengthANN(\g_1)}(\g_1)=2$,
    
    		\item\label{it:rescaled:prop:sum_length_d} it holds that 
    $\lengthANN(\g_1) \leq 10Rd
    $,

    \item\label{it:rescaled:prop:sum_cost_d} it holds 
    that
    $\paramANN(\g_1) \leq 2239 R\kappa^2d\eps^{-2}
    $, and
		
		\item\label{it:rescaled:prop:sum_size_d} it holds 
    that
    $\size(\g_1) \leq 2
    $
    \end{enumerate}
		\cfload. \Nobs that \cref{gen:scaling:networks} (applied with
$\beta \curvearrowleft \kappa$,
$L \curvearrowleft n$
in the notation of \cref{gen:scaling:networks}) and the fact that $ \kappa\leq 2^n$ \prove\ that there exists $\g_2\in\ANNSET$ which satisfies that
\begin{enumerate}[(A)]
\item{
\label{it:scaling:realisation:rescaling:kappa}
it holds for all $x\in\R$ that $(\realisation(\g_2))(x)=\kappa x$,
}
\item{
\label{it:scaling:dims:rescaling:kappa}
it holds that $\dims(\g_2)=(1,2,2,\ldots,2,1)\in\N^{n+2}$,
}
\item{
\label{it:scaling:size:rescaling:kappa}
it holds that  $\insize(\g_2)\leq1$ and $\size(\g_2)\leq 2$, and
}
\item{
\label{it:scaling:param:rescaling:kappa}
it holds that $\param(\g_2)\leq6n+1$.
}
\end{enumerate}
\Nobs that \cref{it:rescaled:prop:sum_d}, \cref{it:scaling:dims:rescaling:kappa}, \cref{Lemma:PropertiesOfCompositions_n2}, and \cref{Prop:identity_representation} \prove\ that 
\begin{equation}
\label{cor:rescaled:final:approximation2:eq:1}
\realisation(\compANN{\g_2}{\ReLUidANN{1}}\bullet\g_1)\in C(\R^d,\R)
\qandq
\lengthANN(\compANN{\g_2}{\ReLUidANN{1}}\bullet\g_1)=\lengthANN(\g_2)+\lengthANN(\g_1)
\end{equation}
\cfload. Combining this with \cref{it:rescaled:prop:sum_length_d}, \cref{it:scaling:dims:rescaling:kappa}, and the fact that for all $x\in[2,\infty)$ it holds that $\ceil*{\log_2(x)}\leq \log_2(x)+1\leq x$ \proves\ that
\begin{equation}
\label{cor:rescaled:final:approximation2:eq:2}
\lengthANN(\compANN{\g_2}{\ReLUidANN{1}}\bullet\g_1)= \lengthANN(\g_2)+\lengthANN(\g_1)\leq 10Rd+6\ceil*{\log_2(\max\{\kappa,2\})}+1\leq 17R\max\{\kappa,1\}d.
\end{equation}
\Moreover \cref{it:rescaled:prop:sum_approx_d}, \cref{it:scaling:realisation:rescaling:kappa}, \cref{Lemma:PropertiesOfCompositions_n2}, and \cref{Prop:identity_representation} \prove\ that for all $x = (x_1, \ldots, x_d) \in [a,b]^d$ it holds that
\begin{equation}
\label{cor:rescaled:final:approximation2:eq:3}
\begin{split}
\vass[\bigg]{g\prbb{\gamma 2^d\ssum_{i = 1}^d x_i}- (\realisation\pr{\compANN{\g_2}{\ReLUidANN{1}}\bullet\g_1}) (x)}
&=
\kappa\vass[\bigg]{f\prbb{\gamma 2^d\ssum_{i = 1}^d x_i}- (\realisation\pr{\g_1}) (x)}\leq \eps.
\end{split}
\end{equation}
\Moreover \cref{it:rescaled:prop:sum_size_d}, \cref{it:scaling:size:rescaling:kappa}, \cref{Prop:identity_representation:prop}, and the fact that $\kappa^{\nicefrac{1}{n}}\leq 2$ \prove\ that
\begin{equation}
\label{cor:rescaled:final:approximation2:eq:4}
\size(\compANN{\g_2}{\ReLUidANN{1}}\bullet\g_1)\leq\max\{\size(\g_2),\size(\g_1)\}\leq \max\pRb{2,2}=2.
\end{equation}
\Moreover 
\eqref{cor:rescaled:final:approximation2:eq:2},
\cref{it:scaling:dims:rescaling:kappa},
\dimsbene, and
\cref{Prop:identity_representation}
 \prove\ that for all $k\in\N_0\cap[0,\lengthANN(\compANN{\g_2}{\ReLUidANN{1}}\bullet\g_1)]$ it holds that
\begin{equation}
\label{cor:rescaled:final:approximation2:eq:5}
\singledims_k(\compANN{\g_2}{\ReLUidANN{1}}\bullet\g_1)=
\begin{cases}
\singledims_k(\g_1) &\colon k\in\N_0\cap[0,\lengthANN(\g_1))\\
2 &\colon k\in\N\cap[\lengthANN(\g_1), \lengthANN(\g_1)+\lengthANN(\g_2))\\
1 &\colon k=\lengthANN(\g_1)+\lengthANN(\g_2).
\end{cases}
\end{equation}
Combining this, 
\cref{it:rescaled:prop:sum_hid_d}, and
\cref{it:rescaled:prop:sum_cost_d}
with \cite[Proposition~2.19]{beneventano21}
and the fact that for all $x\in[2,\infty)$ it holds that $\ceil*{\log_2(x)}\leq \log_2(x)+1\leq x$ \proves\ that
\begin{align}
\nonumber&\param(\compANN{\g_2}{\ReLUidANN{1}}\bullet\g_1)
\\\nonumber&=\sum_{k=1}^{\lengthANN(\compANN{\g_2}{\ReLUidANN{1}}\bullet\g_1)} \singledims_k(\compANN{\g_2}{\ReLUidANN{1}}\bullet\g_1)(\singledims_{k-1}(\compANN{\g_2}{\ReLUidANN{1}}\bullet\g_1)+1)
\\\nonumber&=\PR*{\sum_{k=1}^{\lengthANN(\g_1)-1} \singledims_k(\g_1)(\singledims_{k-1}(\g_1)+1)}
+2(\singledims_{\lengthANN(\g_1)-1}(\g_1)+1)
+\PR*{\sum_{k=1}^{\lengthANN(\g_2)-1} 2(2+1)}
+1(2+1)
\\\nonumber&=\PR*{\sum_{k=1}^{\lengthANN(\g_1)} \singledims_k(\g_1)(\singledims_{k-1}(\g_1)+1)}
+\singledims_{\lengthANN(\g_1)-1}(\g_1)+1
+6(\lengthANN(\g_2)-1)
+3
\\&=\param(\g_1)
+6(\lengthANN(\g_2)-1)
+6
\\\nonumber&\leq 2239 R\kappa^2d\eps^{-2}
+6n
+6
\\\nonumber&= 2239 R\kappa^2d\eps^{-2}+6\max\{\ceil*{\log_2(\kappa)},1\}+6
\\\nonumber&\leq 2245 R\kappa^2d\eps^{-2}+6\max\{\kappa,1\}
\\\nonumber&\leq 2251 R\max\{\kappa^3,1\}d\eps^{-2}
.
\end{align}
This, 
\eqref{cor:rescaled:final:approximation2:eq:1}, 
\eqref{cor:rescaled:final:approximation2:eq:2}, 
\eqref{cor:rescaled:final:approximation2:eq:3}, 
\eqref{cor:rescaled:final:approximation2:eq:4}, and
\eqref{cor:rescaled:final:approximation2:eq:5}
 \prove[ei]\ \cref{cor:rescaled:prop:sum_d,cor:rescaled:prop:sum_approx_d,cor:rescaled:prop:sum_length_d,cor:rescaled:prop:sum_cost_d,cor:rescaled:prop:sum_size_d,cor:rescaled:prop:sum_hid_d}.
The proof of \cref{cor:rescaled:final:approximation2} is thus complete.
\end{proof}

\cfclear
\begin{athm}{cor}{cor:rescaled:final:approximation3}
    Let $a\in\R$, $b\in[a,\infty)$, $d\in\N$, $\mathfrak{c},\gamma,\kappa\in(0,\infty)$, $\eps\in(0,\kappa)$ satisfy $\mathfrak{c}\geq 4504\max\{\kappa^3,1\}\ceil{\log_2\pr*{\max\{1,\gamma\}\max\{2,\vass{a},\vass{b}\}}}$,
		and let $\tgt\colon\R\to\R$ satisfy for all $x,y\in\R$, $k\in\Z$ that $\vass{\tgt(0)}\leq 2\kappa$, $\tgt(x+2k\pi)=\tgt(x)$, and $\vass{\tgt(x)-\tgt(y)}\leq\kappa\vass{x-y}$. 
    Then there exists $\f\in\ANNSET$ such that

    \begin{enumerate}[(i)]
    \item \label{cor:rescaled3:prop:sum_d} it holds
    that $\realisation(\f) \in C(\R^d,\R)$,

    \item\label{cor:rescaled3:prop:sum_approx_d} it holds that
    $\sup_{x = (x_1, \ldots, x_d) \in [a,b]^d}
		\textstyle
    \vass[\big]{g\prb{\gamma 2^d\ssum_{i = 1}^d x_i}- (\realisation\pr{\f}) (x)} \leq \varepsilon
		$,
		
		\item\label{cor:rescaled3:prop:sum_length_d} it holds that 
    $\lengthANN(\f) 
		\leq \mathfrak{c}d
    $,

    \item\label{cor:rescaled3:prop:sum_cost_d} it holds
    that
    $\paramANN(\f) \leq \mathfrak{c} d\eps^{-2}
    $, and
		
		\item\label{cor:rescaled3:prop:sum_size_d} it holds
    that
    $\size(\f) \leq 1
    $
    \end{enumerate}
    \cfout. 
\end{athm}

\begin{proof}[Proof of \cref{cor:rescaled:final:approximation3}]
Throughout this proof let $c=2251\ceil{\log_2\pr*{\max\{1,\gamma\}\max\{2,\vass{a},\vass{b}\}}}\in\N$.
\Nobs that \cref{cor:rescaled:final:approximation2} (applied with
$a \curvearrowleft a$,
$b \curvearrowleft b$,
$d \curvearrowleft d$,
$\kappa \curvearrowleft \kappa$,
$\eps \curvearrowleft \eps$
in the notation of \cref{cor:rescaled:final:approximation2}) \proves\ that there exists $\g\in\ANNSET$ which satisfies that
\begin{enumerate}[(I)]
    \item \label{it:cor:rescaled:prop:sum_d} it holds
    that $\realisation(\g) \in C(\R^d,\R)$,

    \item\label{it:cor:rescaled:prop:sum_approx_d} it holds that
    $\sup_{x = (x_1, \ldots, x_d) \in [a,b]^d}
		\textstyle
    \vass[\big]{g\prb{\gamma 2^d\ssum_{i = 1}^d x_i}- (\realisation\pr{\g}) (x)} \leq \varepsilon
		$,
		
		\item\label{it:cor:rescaled:prop:sum_hid_d} it holds that 
		$
		\singledims_{\hidlengthANN(\g)}(\g)=2
		$,
		
		\item\label{it:cor:rescaled:prop:sum_length_d} it holds that
    $\lengthANN(\g) 
		\leq 20^{-1}\max\{\kappa,1\}cd
    $,

    \item\label{it:cor:rescaled:prop:sum_cost_d} it holds
    that
    $\paramANN(\g) \leq \max\{\kappa^3,1\}cd\eps^{-2}
    $, and
		
		\item\label{it:cor:rescaled:prop:sum_size_d} it holds
    that
    $\size(\g) \leq 2
    $
    \end{enumerate}
\cfload. This,
the fact that $\max\{\kappa,1\}\leq\max\{\kappa^3,1\}\min\{\eps^{-2},1\}$, 
 and \cref{lemma:network:halfed:paramsize} (applied with
$\f \curvearrowleft \g$,
$d \curvearrowleft d$,
$\dnsc \curvearrowleft 2^{-1}$
 in the notation of \cref{lemma:network:halfed:paramsize}) \prove\ that there exists $\f\in\ANNSET$ which satisfies that
\begin{enumerate}[(A)]
    \item \label{it2:cor:rescaled:prop:sum_d} it holds
    that $\realisation(\f)=\realisation(\g) \in C(\R^d,\R)$,

    \item\label{it2:cor:rescaled:prop:sum_approx_d} it holds that
    $\sup_{x = (x_1, \ldots, x_d) \in [a,b]^d}
		\textstyle
    \vass[\big]{g\prb{\gamma 2^d\ssum_{i = 1}^d x_i}- (\realisation\pr{\f}) (x)} \leq \varepsilon
		$,
		
		\item\label{it2:cor:rescaled:prop:sum_length_d} it holds that 
    $\lengthANN(\f)\leq 2 (20^{-1}\max\{\kappa,1\}cd)+1 \leq  \mathfrak{c}d
    $,

    \item\label{it2:cor:rescaled:prop:sum_cost_d} it holds
    that
    \begin{equation}
		\paramANN(\f) \leq \max\{\kappa^3,1\}cd\eps^{-2}+2+\max\{\kappa,1\}cd\leq \max\{\kappa^3,1\}(2c+2)d\eps^{-2}\leq\mathfrak{c}d\eps^{-2},
    \end{equation} and
		
		\item\label{it2:cor:rescaled:prop:sum_size_d} it holds
    that
    $\size(\f) \leq 1
    $.
    \end{enumerate}
\Nobs that 
\cref{it2:cor:rescaled:prop:sum_d,it2:cor:rescaled:prop:sum_length_d,it2:cor:rescaled:prop:sum_approx_d,it2:cor:rescaled:prop:sum_cost_d,it2:cor:rescaled:prop:sum_size_d}
 \prove[ei]\
\cref{cor:rescaled3:prop:sum_d,cor:rescaled3:prop:sum_length_d,cor:rescaled3:prop:sum_approx_d,cor:rescaled3:prop:sum_cost_d,cor:rescaled3:prop:sum_size_d}.
The proof of \cref{cor:rescaled:final:approximation3} is thus complete.
\end{proof}

\newpage
\section{Lower and upper bounds for the minimal number of ANN parameters in the approximation of certain high-dimensional functions}
\label{Upper and Lower Bounds}

 In this section we establish in \cref{Thm}, \cref{Thm1}, \cref{Thm6}, and \cref{Thm6.1} below that certain families of functions can be approximated without the curse of dimensionality by deep \ANNs\ but neither by shallow nor insufficiently deep \ANNs. 
 %

Specifically, \cref{Thm} proves that the plain vanilla \textit{product functions}
can neither be approximated without the curse of dimensionality by means of shallow \ANNs\ nor insufficiently deep \ANNs\ even if the absolute values of the \ANN\ parameters may be arbitrarily large but can be approximated without the curse of dimensionality by sufficiently deep \ANNs\ even if the absolute values of the \ANN\ parameters are assumed to be uniformly bounded by $1$.
Our proof of \cref{Thm} employs 
\begin{itemize}
	\item{
\textit{the lower bound results} for the minimal number of non-zero parameters of \ANNs\ to approximate the product functions in  \cref{prop:product:unrestricted} and}
\item{\textit{the upper bound result} for the minimal number of parameters of \ANNs\ to approximate the product functions in \cref{cor:downsized:products}.}
\end{itemize}
\Nobs that \cref{Thm:prod} in the introduction is a direct consequence of \cref{Cor:oneplus:for:prod} which is implied by \cref{Thm}.

\cref{Thm1} proves that \textit{compositions of certain periodic functions and certain scaled product functions} can neither be approximated without the curse of dimensionality by means of shallow \ANNs\ nor insufficiently deep \ANNs\ even if the \ANN\ parameters may be arbitrarily large but can be approximated without the curse of dimensionality by sufficiently deep \ANNs\ even if the absolute values of the \ANN\ parameters are assumed to be uniformly bounded by $1$.
Our proof of \cref{Thm1} employs 
\begin{itemize}
	\item{\textit{the lower bound result} for the minimal number of non-zero parameters of \ANNs\ to approximate the considered compositions in \cref{LowerBound:prod} and }
	\item{\textit{the upper bound result} for the minimal number of parameters of \ANNs\ to approximate the considered compositions in \cref{final:approximation:imp}.}
	\end{itemize}
\Nobs that \cref{Thm:sin} in the introduction follows immediately from \cref{cor:sin:prod:oneplus} and \cref{cor:sin:prod:oneplus}, in turn, follows from \cref{Thm1}.

\cref{Thm6} proves that \textit{compositions of certain periodic functions and certain scaled sum functions} can neither be approximated without the curse of dimensionality by means of shallow \ANNs\ nor insufficiently deep \ANNs\ even if the \ANN\ parameters may be arbitrarily large but can be approximated without the curse of dimensionality by sufficiently deep \ANNs\ even if the absolute values of the \ANN\ parameters are assumed to be uniformly bounded by $1$.
Our proof of \cref{Thm6} employs 
\begin{itemize}
	\item{
\textit{the lower bound result} for the minimal number of non-zero parameters of \ANNs\ to approximate the considered compositions in \cref{LowerBound:sum} and}
\item{\textit{the upper bound result} for the minimal number of parameters of \ANNs\ to approximate the considered compositions in \cref{cor:rescaled:final:approximation3}.}
\end{itemize}
\cref{Thm6} and the elementary result regarding multidimensional localizing functions in \cref{cor:localized:sum} \prove\ \cref{Thm6.1}.
\cref{Thm6.1} establishes the existence of \emph{smooth
	and uniformly globally bounded functions with compact support} which 
can neither be approximated without the curse of dimensionality by means of shallow \ANNs\ nor insufficiently deep \ANNs\ even if the \ANN\ parameters may be arbitrarily large but which can be approximated without the curse of dimensionality by sufficiently deep \ANNs\ even if the absolute values of the \ANN\ parameters are assumed to be uniformly bounded by $1$.
\Nobs that \cref{Thm:gen} in the introduction is a direct consequence of \cref{Cor:Thm6.1} which is implied by \cref{Thm6.1}.

\subsection{ANN approximations regarding high-dimensional product functions}
\label{subsection:prod_results}
\cfclear

\begin{definition}
	[Cost of ANN approximations]
	\label{def:cost}
	\cfconsiderloaded{def:cost}
	Let $a\in\R$, $b\in[a,\infty)$.
	Then we denote by  $\COST_{a,b}\colon\pr{\cup_{d\in\N}C(\R^d,\R)}\times[0,\infty]^3\to[0,\infty]$ and $\cost_{a,b}\colon\pr{\cup_{d\in\N}C(\R^d,\R)}\times[0,\infty]^3\to[0,\infty]$
	the functions which satisfy for all $d\in\N$, $f\in C(\R^d,\R)$,  $L,S,\eps\in[0,\infty]$ that
	\begin{equation}
		\label{Cost:def}
		\begin{split}
			\COST_{a,b}(f,L,S,\eps)=
			&\min \pr*{ \pR*{ p \in \N \colon \PR*{\!\!
						\begin{array}{c}
							\exists\, \mathscr{f} \in \ANNSET \colon
							(\paramANN(\mathscr{f})\leq p)\land
							(\lengthANN(\mathscr{f})\leq L)
							\\\land{}
							(\size\pr{\mathscr{f}}\leq S)	\land{}
							(\realisation(\mathscr{f}) \in C(\R^d,\R))
							\\\land{}
							(\sup_{x\in[a,b]^d}\vass{(\realisation(\mathscr{f}))(x)-f(x)} \leq \varepsilon)\\
						\end{array} \!\! }
				}
				\cup\{\infty\}
			}
		\end{split}
	\end{equation}
and
	\begin{equation}
	\label{cost:def}
	\begin{split}
		\cost_{a,b}(f,L,S,\eps)=
		&\min \pr*{ \pR*{ p \in \N \colon \PR*{\!\!
					\begin{array}{c}
						\exists\, \mathscr{f} \in \ANNSET \colon
						(\paramnotzero(\mathscr{f})\leq p)\land
						(\lengthANN(\mathscr{f})\leq L)
						\\\land{}
						(\size\pr{\mathscr{f}}\leq S)	\land{}
						(\realisation(\mathscr{f}) \in C(\R^d,\R))
						\\\land{}
						(\sup_{x\in[a,b]^d}\vass{(\realisation(\mathscr{f}))(x)-f(x)} \leq \varepsilon)\\
					\end{array} \!\! }
			}
			\cup\{\infty\}
		}
	\end{split}
\end{equation}
\cfload. 
\end{definition}

\cfclear
\begin{athm}{cor}{cost:inequalities}
	Let 
	$a\in\R$, 
	$b\in[a,\infty)$, 
	$d\in\N$,
	$f\in C(\R^d,\R)$, 
	$L,S,\eps\in[0,\infty]$ \cfload. Then 
	\begin{equation}
\label{eq:cost:inequalities}
		1
		\leq
		\cost_{a,b}(f,L,S,\eps)
		\leq
		\COST_{a,b}(f,L,S,\eps)
	\end{equation}
\cfout.
\end{athm}

\begin{proof}[Proof of \cref{cost:inequalities}]
	\Nobs that \eqref{Cost:def}, \eqref{cost:def}, and the fact that for all $\f\in\ANNSET$ it holds that $\paramnotzero(\f)\leq\paramANN(\f)$ \prove[pe]\ \eqref{eq:cost:inequalities}.
	\finishproofthus
\end{proof}

\cfclear
\begin{athm}{theorem}{Thm}
Let $a\in\R$, $b\in(a,\infty)$, $\dom\in(1,\infty)$ satisfy $\dom=\max\{\vass{a},\vass{b}\}$ and
for every $d\in\N$ let $f_d\colon\R^d\to\R$ satisfy for all $x=(x_1,\ldots,x_d)\in\R^d$ that $f_d(x)=\prod_{i=1}^d x_i$. Then 
\begin{enumerate}[(i)]
\item{
\label{it:Thm:eq1}
it holds for all $c\in[3,\infty)$, $d,L\in\N$, $\eps\in(0,\dom)$ that 
\begin{equation}
\begin{split}
	&\cost_{a,b}(f_d,L,cd^c,\eps)\geq
(\dom^{\nicefrac{d}{2}}-\eps^{\nicefrac{1}{2}})^{\nicefrac{1}{L}}\tfrac{\ln(\dom)^c}{\dom c^{c+1} L^c},
\end{split}
\end{equation}
}
\item{
	\label{it:Thm:eq1.5}
	it holds for all $d,L\in\N\cap[2,\infty)$, $\eps\in(0,\infty)$ that
	\begin{equation}
		\begin{split}
			&\cost_{a,b}(f_d,L,\infty,\eps)\geq
			(L-1)\prb{
					\prb{(b-a)\dom^{(d-2)/2}\pr{8\eps}^{-\nicefrac{1}{2}}}^{1/(L-1)}-1}
		\end{split}
	\end{equation}
}
\item{
\label{it:Thm:eq2}
it holds for all $c\in\big[2^{15}\ln(\max\{2,\dom\}),\infty\big)$, $d\in\N$, $\eps\in(0,\nicefrac{1}{2})$ that
\begin{equation}
\begin{split}
	&\COST_{a,b}(f_d,cd^2\vass{\ln(\eps)},1,\eps)\leq
cd^3\vass{\ln(\eps)}
\end{split}
\end{equation}
}
\end{enumerate}
\cfout.
\end{athm}

\begin{proof}[Proof of \cref{Thm}]
\Nobs that \cref{cor:lower:bound:prod:sized} (applied with
$a \curvearrowleft a$,
$b \curvearrowleft b$,
$c \curvearrowleft c$,
$\eps \curvearrowleft \eps$,
$d \curvearrowleft d$,
$L \curvearrowleft L$,
$f \curvearrowleft f_d$
for $c\in[3,\infty)$, $d,L\in\N$, $\eps\in(0,\dom)$ in the notation of \cref{cor:lower:bound:prod:sized}) \proves\ that for all $c\in[3,\infty)$, $d,L\in\N$, $\eps\in(0,\dom)$ it holds that
\begin{equation}
\begin{split}
	&\cost_{a,b}(f_d,L,cd^c,\eps)\geq
(\dom^{\nicefrac{d}{2}}-\eps^{\nicefrac{1}{2}})^{\nicefrac{1}{L}}\tfrac{\ln(\dom)^c}{\dom c^{c+1} L^c}
\end{split}
\end{equation}
\cfload. This \proves[ei]\ \cref{it:Thm:eq1}.
\Nobs that \cref{prop:product:unrestricted} (applied with
$a \curvearrowleft a$,
$b \curvearrowleft b$,
$\dom \curvearrowleft \dom$,
$d \curvearrowleft  d$,
$\eps \curvearrowleft  \eps$,
$f \curvearrowleft f_d|_{[a,b]^d}$
for $d\in\N\cap[2,\infty)$, $\eps\in\pr{0,\infty}$ in the notation of \cref{prop:product:unrestricted}) \proves\ that for all $d,L\in\N\cap[2,\infty)$, $\eps\in(0,\infty)$ it holds that
\begin{equation}
	\cost_{a,b}(f_d,L,\infty,\eps)\geq
	(L-1)\prb{
			\prb{(b-a)\dom^{(d-2)/2}\pr{8\eps}^{-\nicefrac{1}{2}}}^{1/(L-1)}-1}
	.
\end{equation}
\Nobs that \cref{cor:downsized:products} (applied with
		$d \curvearrowleft  d$,
		$\eps \curvearrowleft  \eps$,
		$a \curvearrowleft  a$,
		$b \curvearrowleft  b$,
		$\gamma \curvearrowleft  1$,
		$\scl \curvearrowleft 1$
for $d\in\N$, $\eps\in\pr{0,\nicefrac{1}{2}}$ in the notation of \cref{cor:downsized:products}) \proves\ that for all $c\in\big[7743\max\{1,\ceil{\log_2(\dom)}\}\ln(2)^{-1},\infty\big)\supseteq \big[2^{15}\ln(\max\{2,\dom\}),\infty\big)$, $d\in\N$, $\eps\in(0,\nicefrac{1}{2})$ there exists $\f\in\ANNSET$ such that
\begin{enumerate}[(I)]
    \item \label{it2:d_downsized:product_continuous} it holds 
    that $\realisation(\mathscr{f}) \in C(\R^{d},\R)$,

    \item\label{it2:d_downsized:product_approx} it holds that
    $\sup_{x = (x_1, \ldots, x_d) \in [a,b]^d}\vass[\big]
    {\sprod_{i=1}^{d}x_i - (\realisation(\f)) (x)} 
                \leq 
                \eps$,

    \item\label{it2:d_downsized:product_length} it holds 
    that
    $\lengthANN(\mathscr{f}) \leq cd^2\vass{\ln(\eps)}
    $,

    \item\label{it2:d_downsized:product_cost} it holds 
    that
    $\paramANN(\mathscr{f}) \leq cd^3\vass{\ln(\eps)}
    $, and
		
		\item\label{it2:d_downsized:product_size} it holds 
    that
    $\size(\mathscr{f}) \leq 1
    $.
\end{enumerate}
\Nobs that
\cref{it2:d_downsized:product_continuous,it2:d_downsized:product_approx,it2:d_downsized:product_length,it2:d_downsized:product_size,it2:d_downsized:product_cost}  \prove[ei]\ \cref{it:Thm:eq2}.
\finishproofthus
\end{proof}

\cfclear
\begin{athm}{cor}{Cor:oneplus:for:prod}
	Let $a\in\R$, $b\in(a,\infty)$ satisfy $\max\{\vass{a},\vass{b}\}>1$ and
	for every $d\in\N$ let $f_d\colon\R^d\to\R$ satisfy for all $x=(x_1,\ldots,x_d)\in\R^d$ that $f_d(x)=\prod_{i=1}^d x_i$. Then there exists $\mathfrak{c}\in(0,\infty)$ such that for all $c\in[\mathfrak{c},\infty)$, $d,L\in\N$, $\eps\in(0,1)$ it holds that
			\begin{equation}
				\label{it:Cor:oneplus:eq1,it:Cor:oneplus:eq2.1}
					\cost_{a,b}(f_d,L,cd^c,\eps)\geq
					\pr{cL}^{-3c}(1+c^{-1})^{\nicefrac{d}{L}},
					\qquad
					\COST_{a,b}(f_d,L,\infty,\varepsilon)
					\ge (1+c^{-1})^{d/L},
			\end{equation}
			\begin{equation}
				\label{it:Cor:oneplus:eq1,it:Cor:oneplus:eq2}
				\qandq
				\COST_{a,b}(f_d,cd^2\eps^{-1},1,\eps)\leq
				cd^3\eps^{-1}
			\end{equation}
	\cfout.
\end{athm}

\begin{proof}[Proof of \cref{Cor:oneplus:for:prod}]
	Throughout this proof let $c,\dom,B\in\R$ satisfy $\dom=\max\{\vass{a},\vass{b}\}$,
	$B=\sqrt{8}\dom(b-a)^{-1}$, 
	and
	\begin{equation}
		\label{eq:oneplus:onlyprod:setup}
		c\geq\max
		\pRb{
			4\ln(2B)\pr{\ln(2)\ln(\dom)}^{-1},
			\dom\ln(\dom)^{-1},
			\pr{
		\dom^{\nicefrac{1}{4}}-1}^{-1},
		2^{16}\ln(\max\{2,\dom\})
	    }
.
	\end{equation} 
\Nobs that \cref{cost:inequalities} and \eqref{eq:oneplus:onlyprod:setup} \prove\ that for all $\eps\in(0,1)$, $d,L\in\N$ with $d< 4\ln(2)\ln(\dom)^{-1}$ it holds that
\begin{equation}
	\label{eq:oneplus:onlyprod:lowdim}
\cost_{a,b}(f_d,L,cd^c,\eps)
\geq 1
=(cL)^{-3c} (cL)^{3c} 
\geq (cL)^{-3c} \pr{1+c^{-1}}^{c}
\geq (cL)^{-3c} \pr{1+c^{-1}}^{\nicefrac{d}{L}}
\end{equation}
\cfload.
\Moreover \eqref{eq:oneplus:onlyprod:setup} \proves\ that for all $\eps\in(0,1)$, $d,L\in\N$ with $d\geq 4\ln(2)\ln(\dom)^{-1}$ it holds that
\begin{equation}
	\prb{\dom^{\nicefrac{d}{2}}-\eps^{\nicefrac{1}{2}}}^{\nicefrac{1}{L}}
	\geq
	\prb{\dom^{\nicefrac{d}{4}}\dom^{\nicefrac{d}{4}}-1}^{\nicefrac{1}{L}}
	\geq
	\prb{2\dom^{\nicefrac{d}{4}}-1}^{\nicefrac{1}{L}}
	\geq
	\prb{\dom^{\nicefrac{d}{4}}+2-1}^{\nicefrac{1}{L}}
	\geq
	\dom^{\nicefrac{d}{4L}}.
\end{equation}
This and \eqref{eq:oneplus:onlyprod:setup} \prove\ that for all $\eps\in(0,1)$, $d,L\in\N$ with $d\geq 4\ln(2)\ln(\dom)^{-1}$ it holds that
\begin{equation}
	\begin{split}
	\prb{\dom^{\nicefrac{d}{2}}-\eps^{\nicefrac{1}{2}}}^{\nicefrac{1}{L}}\tfrac{\ln(\dom)^c}{\dom c^{c+1} L^c}
	\geq
    \dom^{\nicefrac{d}{4L}}\prb{\tfrac{\ln(\dom)}{\dom }}^{c}\pr{\tfrac{1}{cL}}^{c+1}
    \geq
    \pr{cL}^{-3c}(1+c^{-1})^{\nicefrac{d}{L}}.
	\end{split}
\end{equation}
This and \cref{Thm} (applied with
$a \curvearrowleft  a$,
$b \curvearrowleft  b$,
$\dom \curvearrowleft  \dom$,
$f_d \curvearrowleft  f_d$
for $d\in\N$ in the notation of \cref{Thm}) \prove\ that
\begin{enumerate}[(I)]
	\item{
		\llabel{it:first:LB}
		it holds for all $d,L\in\N$, $\eps\in(0,1)$ with $d\geq 4\ln(2)\ln(\dom)^{-1}$ that 
			\begin{equation}
			\cost_{a,b}(f_d,L,cd^c,\eps)\geq  
				(\dom^{\nicefrac{d}{2}}-\eps^{\nicefrac{1}{2}})^{\nicefrac{1}{L}}\tfrac{\ln(\dom)^c}{\dom c^{c+1} L^c}
				\geq 
				    \pr{cL}^{-3c}(1+c^{-1})^{\nicefrac{d}{L}},
			\end{equation}
	}
	\item{
		\llabel{it:Thm:eq1.5}
		it holds for all $d,L\in\N\cap[2,\infty)$, $\eps\in(0,\infty)$ that
		\begin{equation}
			\begin{split}
				&\cost_{a,b}(f_d,L,\infty,\eps)\geq
				(L-1)\prb{
					\prb{\dom^{d/2}B^{-1}\eps^{-\nicefrac{1}{2}}}^{1/(L-1)}-1},
			\end{split}
		\end{equation}
		and
	}
	\item{
		it holds for all $d\in\N$, $\eps\in(0,\nicefrac{1}{2})$ that
		\begin{equation}
			\llabel{it:loc:UB}
			\begin{split}
			\COST_{a,b}(f_d,cd^2\vass{\ln(\eps)},1,\eps)
			\leq  cd^3\vass{\ln(\eps)}.
			\end{split}
		\end{equation}
	}
\end{enumerate}
\Nobs that \eqref{eq:def:ANN:operators}, \eqref{Cost:def}, and the fact that
$\ln(1+c^{-1})\leq c^{-1}$ \prove\ that for all $d,L\in\N$,
$\eps\in(0,1)$ with $d<cL\ln(2)$ it holds that
\begin{equation}
	\llabel{eq:unrestricted:lowdim}
	\COST_{a,b}(f_d,L,\infty,\eps)
	\geq 2
	\geq \exp\prb{d\pr{cL}^{-1}}
	\geq (1+c^{-1})^{\nicefrac{d}{L}}.
\end{equation}
\Moreover \eqref{eq:oneplus:onlyprod:setup} \proves\ that for all
$d,L\in\N$ with $d\geq cL\ln(2)$ it holds that $d\geq2$ and
\begin{equation}
	\dom^{\nicefrac{d}{4L}}
	\geq \dom^{c\ln(2)/4}
	\geq 2B.
\end{equation}
Combining this and \lref{it:Thm:eq1.5} 
with \cref{cost:inequalities}, \eqref{cost:def}, and the fact that
$B\geq1$ \proves\ that for all $d,L\in\N$, $\eps\in(0,1)$
with $d\geq cL\ln(2)$ it holds that
\begin{equation}
	\begin{split}
		\COST_{a,b}(f_d,L,\infty,\eps)
		\geq \cost_{a,b}(f_d,L+1,\infty,\eps)
		&\geq L\prb{\pr*{\dom^{\nicefrac{d}{2}}B^{-1}\eps^{-\nicefrac{1}{2}}}^{\nicefrac{1}{L}}-1}
		\\&\geq L\pr*{B^{-1}\dom^{\nicefrac{d}{2L}}-1}\\
		&\geq L\pr*{2\dom^{\nicefrac{d}{4L}}-1}
		\geq \dom^{\nicefrac{d}{4L}}
		\geq (1+c^{-1})^{\nicefrac{d}{L}}.
	\end{split}
\end{equation}
This, \cref{eq:oneplus:onlyprod:lowdim}, \lref{eq:unrestricted:lowdim}, and \lref{it:first:LB} \prove\ 
\eqref{it:Cor:oneplus:eq1,it:Cor:oneplus:eq2.1}.
\Nobs that \lref{it:loc:UB} \proves\ that for all $d\in\N$, $\eps\in[\nicefrac{1}{2},1)$ it holds that
\begin{equation}
	\COST_{a,b}(f_d,cd^2\eps^{-1},1,\eps)
	\leq
	\COST_{a,b}(f_d,cd^2\vass{\ln(e)},1,e^{-1})
	\leq  cd^3\vass{\ln(e)}
	=
	cd^3
	\leq
	cd^3\eps^{-1}.
\end{equation}
Combining this and \lref{it:loc:UB} with the fact that for all $\eps\in(0,\nicefrac{1}{2})$ it holds that $\vass{\ln(\eps)}\leq\eps^{-1}$ \proves[ei]\ \eqref{it:Cor:oneplus:eq1,it:Cor:oneplus:eq2}.
\finishproofthus
\end{proof}

\cfclear
\begin{athm}{cor}{prod:old}
	Let $a\in\R$, $b\in[a+4,\infty)$ and for every $d\in\N$ let $f_d\colon\R^d\to\R$ satisfy for all $x=(x_1,\ldots,x_d)\in\R^d$ that $f_d(x)=\sprod_{i=1}^d x_i$.
	Then 
	\begin{enumerate}[(i)]
	\item{
		\label{it:prod:old:eq1}
		it holds for all $c\in[3,\infty)$, $d,L\in\N$, $\eps\in(0,1)$ that 
		\begin{equation}
			\begin{split}
				\cost_{a,b}(f_d,L,cd^c,\eps)\geq  	\PRb{\pr{cL}^{-3c}}2^{\nicefrac{d}{3L}}
			\end{split}
		\end{equation}
		and
	}
	\item{
		\label{it:prod:old:eq2}
		it holds for all $c\in[2^{15}\ln(\max\{\vass{a},\vass{b}\}),\infty)$, $d\in\N$, $\eps\in(0,\nicefrac{1}{2})$ that 
		\begin{equation}
			\COST_{a,b}(f_d,cd^2\vass{\ln(\eps)},1,\eps)
			\leq  cd^3\vass{\ln(\eps)}
		\end{equation}
	}
\end{enumerate}
\cfout.
\end{athm}

\begin{proof}[Proof of \cref{prod:old}]
Throughout this proof 
let $\dom\in\R$ satisfy $\dom=\max\{\vass{a},\vass{b}\}$.
\Nobs that \cref{cost:inequalities} \proves\ that for all $c\in[2,\infty)$, $L\in\N$, $\eps\in(0,1)$, $d\in\N\cap[1,12L]$ it holds that 
		\begin{equation}
		\label{eq:lowerdimensions:parameter:bound}
	\begin{split}
  \cost_{a,b}(f_d,L,cd^c,\eps)
	&\geq  1
	\geq  16 c^{-3c}
	\geq  \PRb{(cL)^{-3c}}2^{\nicefrac{12L}{3L}}
	\geq  \PRb{(cL)^{-3c}}2^{\nicefrac{d}{3L}}
	\end{split}
\end{equation}
\cfload.
\Moreover the fact that $\dom\geq 2$ \proves\ that for all $c\in[2,\infty)$, $L\in\N$, $\eps\in(0,1)$, $d\in\N\cap(12L,\infty)$ it holds that
\begin{equation}
	\begin{split}
(\dom^{\nicefrac{d}{2}}-\eps^{\nicefrac{1}{2}})^{\nicefrac{1}{L}}\tfrac{\ln(\dom)^c}{\dom c^{c+1} L^c}
\geq
(\dom^{\nicefrac{d}{12}}\dom^{\nicefrac{d}{12}}\dom^{\nicefrac{d}{3}}-1)^{\nicefrac{1}{L}}\tfrac{\ln(\dom)^c}{\dom c^{c+1} L^c}
&\geq
\dom^{\nicefrac{d}{12L}}\dom^{\nicefrac{d}{3L}}\tfrac{\ln(\dom)^c}{\dom c^{c+1} L^c}
\\&\geq
\dom^{\nicefrac{d}{3L}}\tfrac{\ln(2)^c}{ c^{c+1} L^c}
\\&\geq
\dom^{\nicefrac{d}{3L}}(2cL)^{-c-1}
\\&\geq
\PRb{(cL)^{-3c}}2^{\nicefrac{d}{3L}}
	    .
	\end{split}
\end{equation}
Combining this and \eqref{eq:lowerdimensions:parameter:bound} with \cref{Thm} (applied with
$a \curvearrowleft  a$,
$b \curvearrowleft  b$,
$\dom \curvearrowleft  \dom$,
$f_d \curvearrowleft  f_d$
for $d\in\N$ in the notation of \cref{Thm}) \proves\ that
\begin{enumerate}[(I)]
	\item{
		it holds for all $c\in[3,\infty)$, $d,L\in\N$, $\eps\in(0,1)$ that 
		\begin{equation}
			\begin{split}
				&\cost_{a,b}(f_d,L,cd^c,\eps)\geq  
				\PRb{\pr{cL}^{-3c}}2^{\nicefrac{d}{3L}}
			\end{split}
		\end{equation}
		and
	}
	\item{
		it holds for all $c\in\big[2^{15}\ln(\max\{\vass{a},\vass{b}\}),\infty\big)$, $d\in\N$, $\eps\in(0,\nicefrac{1}{2})$ that
		\begin{equation}
			\begin{split}
				&\COST_{a,b}(f_d,cd^2\vass{\ln(\eps)},1,\eps)\leq  cd^3\vass{\ln(\eps)}.
			\end{split}
		\end{equation}
	}
\end{enumerate}
	\finishproofthus
\end{proof}

\newcommand{\oneplus}{\mathscr{c}}

\cfclear
\begin{athm}{theorem}{Thm1}

Let $\varphi,\dom\in\R$, $\gamma\in(0,1]$, $\beta\in[1,\infty)$, $a\in\R$, $b\in[a,\infty)$, $\mathfrak{c},\kappa\in(0,\infty)$ satisfy
$\mathfrak{c}\geq 19770 \ceil{\ln\pr*{\max\{2,\abs{a},\abs{b},\scl\}}} \max\{1,\kappa^3\}$
and
$\dom=\max\{\vass{a},\vass{b}\}>\scl^{-1}$,
 and
for every $d\in\N$ let
$f_d\colon\R^d\to\R$ satisfy for all $x=(x_1,\ldots,x_d)\in\R^d$ that $f_d(x)=\kappa\sin\prb{\gamma\beta^d\prb{\sprod_{i = 1}^d x_i}+\varphi}$. Then
\begin{enumerate}[(i)]
\item{
\label{it:Thm1:eq1}
it holds for all $d,L\in\N$, $\eps\in(0,\kappa)$ that 
\begin{equation}
\cost_{a,b}(f_d,L,\infty,\eps)
\geq 	\vass[\Big]{\frac{\floor{(b-a)\gamma(\dom\scl)^{d-1}\pi^{-1}-1}}{2}}^{\nicefrac{1}{\max\{1,L-1\}}}
\end{equation}
and
}
\item{
\label{it:Thm1:eq2}
it holds for all $d\in\N$, $\eps\in(0,\kappa)$ that
\begin{equation}
\begin{split}
\COST_{a,b}(f_d,\mathfrak{c}d^2\eps^{-1},1,\eps)
\leq  \mathfrak{c}d^3\eps^{-2}
\end{split}
\end{equation}
}
\end{enumerate}
\cfout.
\end{athm}

\begin{proof}[Proof of \cref{Thm1}]
	Throughout this proof let $g\colon\R\to\R$ satisfy for all $x\in\R$ that $\tgt(x)=\kappa\sin(x+\varphi)$.
	\Nobs that the fact that $\fa{x}[2,\infty)\colon  \ceil{\log_2\pr{x}}\leq 2\ceil{\ln(x)}$
	and \cref{final:approximation:imp} (applied with
	$a \curvearrowleft  a$,
	$b \curvearrowleft  b$,
	$d \curvearrowleft  d$,
	$\kappa \curvearrowleft  \kappa$,
	$\mathfrak{c} \curvearrowleft \mathfrak{c}$,
	$\eps \curvearrowleft  \eps$,
	$\gamma \curvearrowleft  \gamma$,
	$g \curvearrowleft  g$,
	$f \curvearrowleft  f_d$
	for $d\in\N$, $\eps\in(0,\kappa)$ in the notation of \cref{final:approximation:imp}) \prove\ that for all $d\in\N$, $\eps\in(0,\kappa)$ it holds that
	\begin{equation}
		\label{eq:Thm1:prove}
		\begin{split}
		\COST_{a,b}(f_d,\mathfrak{c}d^2\eps^{-1},1,\eps)
		\leq 9885\ceil{\log_2\pr*{\max\{2,\abs{a},\abs{b},\scl\}}}\max\{1,\kappa^3\}d^3\eps^{-2}
		\leq  \mathfrak{c}d^3\eps^{-2}
		\end{split}
	\end{equation}
\cfload.
\Moreover \cref{LowerBound:prod} (applied with
$\varphi \curvearrowleft  \varphi$,
$\gamma \curvearrowleft  \gamma$,
$\scl \curvearrowleft  \scl$,
$a \curvearrowleft  a$,
$b \curvearrowleft  b$,
$\kappa \curvearrowleft  \kappa$,
$f_d \curvearrowleft f_d$
for $d\in\N$ in the notation of \cref{LowerBound:prod}) \proves\ that for all $d,L\in\N$, $\eps\in(0,\kappa)$ it holds that 
\begin{equation}
	\cost_{a,b}(f_d,L,\infty,\eps)
	\geq 
	\vass[\Big]{\frac{\floor{(b-a)\gamma(\dom\scl)^{d-1}\pi^{-1}-1}}{2}}^{\nicefrac{1}{\max\{1,L-1\}}}.
\end{equation}
This and \eqref{eq:Thm1:prove}  \prove[ei]\ \cref{it:Thm1:eq1,it:Thm1:eq2}.
The proof of \cref{Thm1} is thus complete.
\end{proof}

\cfclear
\begin{athm}{cor}{cor:sin:prod:oneplus}
	Let $\varphi\in\R$, $\gamma\in(0,1]$, $\beta\in[1,\infty)$, $a\in\R$, $b\in(a,\infty)$, $\kappa\in(0,\infty)$ satisfy
	$\scl\max\{\vass{a},\vass{b}\}>1$,
	and
	for every $d\in\N$ let
	$f_d\colon\R^d\to\R$ satisfy for all $x=(x_1,\ldots,x_d)\in\R^d$ that $f_d(x)=\kappa\sin\prb{\gamma\scl^d\prb{\sprod_{i = 1}^d x_i}+\varphi}$. Then there exists $\mathfrak{c}\in(0,\infty)$ such that for all $c\in[\mathfrak{c},\infty)$, $d,L\in\N$, $\eps\in(0,\kappa)$ it holds that
	\begin{equation}
		\label{it:cor:oneplus:eq2}
		\COST_{a,b}(f_d,L,\infty,\eps)\geq  (1+c^{-1})^{\nicefrac{d}{\max\{1,L-1\}}}
		\qandq
		\COST_{a,b}(f_d,cd^2\eps^{-1},1,\eps)\leq cd^3\eps^{-2}
	\end{equation}
	\cfout.
\end{athm}

\begin{proof}[Proof of \cref{cor:sin:prod:oneplus}]
	Throughout this proof let $D\in\N$, $\dom\in\R$ satisfy
	\begin{equation}
		\label{eq:cor:oneplus:setup1}
		\dom=\max\{\vass{a},\vass{b}\}
		\qandq
		D\geq  \ln\prb{\tfrac{6\pi}{(b-a)\gamma}}\PR{\ln(\dom\scl)}^{-1}
		\end{equation}
	and let $c\in\R$ satisfy
	\begin{equation}
		\label{eq:cor:oneplus:setup2}
		c\geq
		\max\pRb{
			\pr{\dom\scl-1}^{-1},
			\prb{2^{\nicefrac{1}{D+1}}-1}^{-1},
			19770 \ceil{\ln\pr*{\max\{2,\dom,\scl\}}} \max\{1,\kappa^3\}
		}.
	\end{equation}
\Nobs that 
\eqref{eq:cor:oneplus:setup2} \proves\ that for all $d\in\N\cap[1,D]$, 
$L\in\N$, 
$\eps\in(0,\kappa)$, 
$\f\in\ANNSET$ it holds that
\begin{equation}
	\label{eq:cor:oneplus:lowdim}
	\paramANN(\f)
	\geq 2
	\geq (1+c^{-1})^D
	\geq(1+c^{-1})^{\nicefrac{d}{\max\{1,L-1\}}}
\end{equation}
\cfload.
\Moreover 
\eqref{eq:cor:oneplus:setup1} and 
\eqref{eq:cor:oneplus:setup2} \prove\ that 
\begin{equation}
	\label{eq:cor:oneplus:ind:base}
	\frac{(b-a)\gamma(\dom\scl)^{D}}{2\pi}-1
	\geq
	\frac{(b-a)\gamma(6\pi(b-a)^{-1}\gamma^{-1})}{2\pi}-1
	=2\geq
	(1+c^{-1})^{D+1}.
\end{equation}
\Moreover \eqref{eq:cor:oneplus:setup2} \proves\ that for all $d\in\N$ with $ \frac{(b-a)\gamma(\dom\scl)^{d-1}}{2\pi}-1\geq (1+c^{-1})^d$ it holds that
\begin{equation}
	\begin{split}
		\frac{(b-a)\gamma(\dom\scl)^{d+1-1}}{2\pi}-1
		&=
		\dom\scl\prbb{\frac{(b-a)\gamma(\dom\scl)^{d-1}}{2\pi}-\frac{1}{\dom\scl}+1-1}
		\\&\geq
		\dom\scl\prb{(1+c^{-1})^d+1-(\dom\scl)^{-1}}
		\\&>
		\dom\scl(1+c^{-1})^d
		\geq (1+c^{-1})^{d+1}.
	\end{split}
\end{equation}
Combining this and \eqref{eq:cor:oneplus:ind:base} with induction \proves\ that for all $d\in\N\cap(D,\infty)$ it holds that
\begin{equation}
	\begin{split}
	\vass[\Big]{\frac{\floor{(b-a)\gamma(\dom\scl)^{d-1}\pi^{-1}-1}}{2}}
	=
	\frac{\floor{(b-a)\gamma(\dom\scl)^{d-1}\pi^{-1}-1}}{2}
	&\geq
	\frac{(b-a)\gamma(\dom\scl)^{d-1}\pi^{-1}-2}{2}
	\\&=
	\frac{(b-a)\gamma(\dom\scl)^{d-1}}{2\pi}-1
	\\&\geq (1+c^{-1})^{d}.
	\end{split}
\end{equation}
\Nobs that \eqref{eq:cor:oneplus:setup2} and
 \cref{Thm1} (applied with
$\varphi \curvearrowleft  \varphi$,
$\gamma \curvearrowleft  \gamma$,
$\scl \curvearrowleft  \scl$,
$a \curvearrowleft  a$,
$b \curvearrowleft  b$,
$\mathfrak{c} \curvearrowleft c$,
$\kappa \curvearrowleft  \kappa$,
$f_d \curvearrowleft  f_d$
for $d\in\N$ in the notation of \cref{Thm1}) \hence \prove\ that
\begin{enumerate}[(I)]
	\item{
		\label{it:cor:fromThm:forCor:eq1}
		it holds for all $d\in\N\cap(D,\infty)$, $L\in\N$, $\eps\in(0,\kappa)$ that 
		\begin{equation}
			\cost_{a,b}(f_d,L,\infty,\eps)\geq 		\vass[\Big]{\frac{\floor{(b-a)\gamma(\dom\scl)^{d-1}\pi^{-1}-1}}{2}}^{\nicefrac{1}{\max\{1,L-1\}}}\geq (1+c^{-1})^{\nicefrac{d}{\max\{1,L-1\}}}
		\end{equation}
		and
	}
	\item{
		\label{it:cor:fromThm:forCor:eq2}
		it holds for all $d\in\N$, $\eps\in(0,\kappa)$ that
		\begin{equation}
			\begin{split}
			\COST_{a,b}(f_d,cd^2\eps^{-1},1,\eps)\leq  cd^3\eps^{-2}.
			\end{split}
		\end{equation}
	}
\end{enumerate}
Combining this and \eqref{eq:cor:oneplus:lowdim} with \cref{cost:inequalities}  \proves[ei]\ \eqref{it:cor:oneplus:eq2}.
	\finishproofthus
\end{proof}

\cfclear
\begin{athm}{cor}{cor:sin:prod:oneplus:nonzc}
	Let $\varphi\in\R$, $\gamma\in(0,1]$, $\beta\in[1,\infty)$, $a\in\R$, $b\in(a,\infty)$, $\kappa\in(0,\infty)$ satisfy
	$\scl\max\{\vass{a},\vass{b}\}>1$,
	and
	for every $d\in\N$ let
	$f_d\colon\R^d\to\R$ satisfy for all $x=(x_1,\ldots,x_d)\in\R^d$ that $f_d(x)=\kappa\sin\prb{\gamma\scl^d\prb{\sprod_{i = 1}^d x_i}+\varphi}$. Then there exists $\mathfrak{c}\in(0,\infty)$ such that for all $c\in[\mathfrak{c},\infty)$, $d,L\in\N$, $\eps\in(0,\kappa)$ it holds that
	\begin{equation}
		\label{it:cor:oneplus:eq2:nonzc}
		\cost_{a,b}(f_d,L,\infty,\eps)\geq c^{-1} (1+c^{-1})^{\nicefrac{d}{\max\{1,L-1\}}}
		\quad\text{and}\quad
		\COST_{a,b}(f_d,cd^2\eps^{-1},1,\eps)\leq cd^3\eps^{-2}
	\end{equation}
	\cfout.
\end{athm}

\begin{proof}[Proof of \cref{cor:sin:prod:oneplus:nonzc}]
	Throughout this proof let $D\in\N$, $\dom\in\R$ satisfy
	\begin{equation}
		\label{eq:cor:oneplus:setup1:nonzc}
		\dom=\max\{\vass{a},\vass{b}\}
		\qandq
		D\geq  \ln\prb{\tfrac{5\pi}{(b-a)\gamma}}\PR{\ln(\dom\scl)}^{-1}
	\end{equation}
	and let $c\in\R$ satisfy
	\begin{equation}
		\label{eq:cor:oneplus:setup2:nonzc}
		c\geq
		\max\pRb{
			\prb{\dom\scl-1}^{-1},
			2^D,
			19770 \ceil{\ln\pr*{\max\{2,\dom,\scl\}}} \max\{1,\kappa^3\}}
		.
	\end{equation}
	\Nobs that 
	\eqref{eq:cor:oneplus:setup2:nonzc} \proves\ that for all $d\in\N\cap[1,D]$, 
	$L\in\N$, 
	$\eps\in(0,\kappa)$, 
	$\f\in\ANNSET$ it holds that
	\begin{equation}
		\label{eq:cor:oneplus:lowdim:nonzc}
		\cost_{a,b}(f_d,L,\infty,\eps)
		\geq 1
		\geq c^{-1}2^D
		\geq c^{-1}(1+c^{-1})^D
		\geq c^{-1}(1+c^{-1})^{\nicefrac{d}{\max\{1,L-1\}}}
	\end{equation}
	\cfload.
	\Moreover 
	\eqref{eq:cor:oneplus:setup1:nonzc} and \eqref{eq:cor:oneplus:setup2:nonzc} \prove\ that for all $d\in\N\cap[D,\infty)$ it holds that
			\begin{equation}
		\begin{split}
			\frac{\floor{(b-a)\gamma(\dom\scl)^{(d+1)-1}\pi^{-1}-1}}{2}
			&=
			\frac{\floor{(b-a)\gamma(\dom\scl)^{D}\pi^{-1}(\dom\scl)^{d-D}-1}}{2}
			\\&\geq
			\frac{\floor{5(\dom\scl)^{d-D}-1}}{2}
						\\&\geq
			\frac{5(\dom\scl)^{d-D}-2}{2}
						\\&\geq
			\prb{\tfrac{3}{2}}(\dom\scl)^{d-D}
			\\&\geq
			2^{-D}\prb{\tfrac{3}{2}}^{D+1}(\dom\scl)^{d-D}
			\\&\geq
			c^{-1}\pr{1+c^{-1}}^{D+1}\pr{1+c^{-1}}^{d-D}=c^{-1}\pr{1+c^{-1}}^{d+1}
			.
		\end{split}
	\end{equation}
	\Nobs that \eqref{eq:cor:oneplus:setup2:nonzc} and
	\cref{Thm1} (applied with
	$\varphi \curvearrowleft  \varphi$,
	$\gamma \curvearrowleft  \gamma$,
	$\scl \curvearrowleft  \scl$,
	$a \curvearrowleft  a$,
	$b \curvearrowleft  b$,
	$\mathfrak{c} \curvearrowleft c$,
	$\kappa \curvearrowleft  \kappa$,
	$f_d \curvearrowleft  f_d$
	for $d\in\N$ in the notation of \cref{Thm1}) \hence \prove\ that
	\begin{enumerate}[(I)]
		\item{
			\label{it:cor:fromThm:forCor:eq1:nonzc}
			it holds for all $d\in\N\cap(D,\infty)$, $L\in\N$, $\eps\in(0,\kappa)$ that 
			\begin{equation}
				\begin{split}
				\cost_{a,b}(f_d,L,\infty,\eps)
				&\geq 	\vass[\Big]{\frac{\floor{(b-a)\gamma(\dom\scl)^{d-1}\pi^{-1}-1}}{2}}^{\nicefrac{1}{\max\{1,L-1\}}}
				\\&\geq c^{-1}(1+c^{-1})^{\nicefrac{d}{\max\{1,L-1\}}}
				\end{split}
			\end{equation}
			and
		}
		\item{
			\label{it:cor:fromThm:forCor:eq2:nonzc}
			it holds for all $d\in\N$, $\eps\in(0,\kappa)$ that
			\begin{equation}
				\begin{split}
					\COST_{a,b}(f_d,cd^2\eps^{-1},1,\eps)\leq  cd^3\eps^{-2}.
				\end{split}
			\end{equation}
		}
	\end{enumerate}
	Combining this and \eqref{eq:cor:oneplus:lowdim:nonzc}  \proves[ei]\ \eqref{it:cor:oneplus:eq2:nonzc}.
	\finishproofthus
\end{proof}

\cfclear
\begin{athm}{cor}{old:sin:prod}
	Let $a\in\R$, $b\in[a+7,\infty)$ and for every $d\in\N$ let $f_d\colon\R^d\to\R$ satisfy for all $x=(x_1,\ldots,x_d)\in\R^d$ that $f_d(x)=\sin\prb{\sprod_{i=1}^d x_i}$.
	Then there exists $\mathfrak{c}\in(0,\infty)$ such that for all $d,L\in\N$, $\eps\in(0,1)$ it holds that
	\begin{equation}
		\label{it:old:sin:prod:eq1}
		\COST_{a,b}(f_d,L,\infty,\eps)\geq  2^{\nicefrac{d}{\max\{1,L-1\}}}
		\qandq
		\COST_{a,b}(f_d,\mathfrak{c}d^2\eps^{-1},1,\eps)
		\leq  \mathfrak{c}d^3\eps^{-2}
	\end{equation}
	\cfout.
\end{athm}

\begin{proof}[Proof of \cref{old:sin:prod}]
	\Nobs that the fact that $\fa{d}\N\cap[3,\infty)\colon (\tfrac{7}{2})^{d-1}\geq(\tfrac{3}{2})2^d$ \proves\ that for all $d\in\N\cap[3,\infty)$ it holds that
	\begin{equation}
		\label{eq:param:trans:from:thm:d:geq:3}
		\begin{split}
	\frac{\floor{(b-a)(\max\{\vass{a},\vass{b}\})^{d-1}\pi^{-1}-1}}{2}
	\geq
	\frac{\floor{7(\tfrac{7}{2})^{d-1}\pi^{-1}-1}}{2}
	&\geq
	\frac{\floor{(\tfrac{7}{2})^{d-1}2-1}}{2}
	\\&\geq
	\frac{\floor{2^{d} 3-1}}{2}=(\tfrac{3}{2})2^d-\tfrac{1}{2}\geq 2^d.
	\end{split}
	\end{equation}
	\Moreover \cref{prop:low:dim:gen:low:len:cases} \proves\ that for all $d\in\{1,2\}$, $L\in\N$, $\eps\in(0,1)$ it holds that
	\begin{equation}
		\begin{split}
			\COST_{a,b}(f_d,L,\infty,\eps) \geq 
			2^{\nicefrac{d}{\max\{1,L-1\}}}
		\end{split}
	\end{equation}
	\cfload.
	Combining this, \cref{cost:inequalities}, and \eqref{eq:param:trans:from:thm:d:geq:3} with \cref{Thm1} (applied with
	$\varphi \curvearrowleft  0$,
	$\gamma \curvearrowleft  1$,
	$\scl \curvearrowleft  1$,
	$a \curvearrowleft  a$,
	$b \curvearrowleft  b$,
	$\kappa \curvearrowleft  1$,
	$f_d \curvearrowleft  f_d$
	for $d\in\N$ in the notation of \cref{Thm1}) \proves\ that there exists $\mathfrak{c}\in(0,\infty)$ such that for all $d,L\in\N$, $\eps\in(0,1)$ it holds that
	\begin{equation}
		\COST_{a,b}(f_d,L,\infty,\eps)\geq  2^{\nicefrac{d}{\max\{1,L-1\}}}
		\qandq
		\COST_{a,b}(f_d,\mathfrak{c}d^2\eps^{-1},1,\eps)
		\leq  \mathfrak{c}d^3\eps^{-2}.
	\end{equation}
	\finishproofthus
\end{proof}

\cfclear
\begin{athm}{cor}{old:sin:prod:nozp}
	Let $a\in\R$, $b\in[a+7,\infty)$ and for every $d\in\N$ let $f_d\colon\R^d\to\R$ satisfy for all $x=(x_1,\ldots,x_d)\in\R^d$ that $f_d(x)=\sin\prb{\sprod_{i=1}^d x_i}$.
	Then there exists $\mathfrak{c}\in(0,\infty)$ such that for all $d\in\N\cap(2,\infty)$, $L\in\N$, $\eps\in(0,1)$ it holds that
	\begin{equation}
		\label{it:old:sin:prod:eq1:nozp}
		\cost_{a,b}(f_d,L,\infty,\eps)\geq  2^{\nicefrac{d}{\max\{1,L-1\}}}
		\qandq
		\COST_{a,b}(f_d,\mathfrak{c}d^2\eps^{-1},1,\eps)
		\leq  \mathfrak{c}d^3\eps^{-2}
	\end{equation}
\cfout.
\end{athm}

\begin{proof}[Proof of \cref{old:sin:prod:nozp}]
\Nobs that the fact that $\fa{d}\N\cap[3,\infty)\colon (\tfrac{7}{2})^{d-1}\geq(\tfrac{3}{2})2^d$ \proves\ that for all $d\in\N\cap[3,\infty)$ it holds that
\begin{equation}
	\label{eq:param:trans:from:thm:d:geq:3:nozp}
	\begin{split}
\frac{\floor{(b-a)(\max\{\vass{a},\vass{b}\})^{d-1}\pi^{-1}-1}}{2}
	\geq
\frac{7(\tfrac{7}{2})^{d-1}}{2\pi}-1
	\geq
(\tfrac{7}{2})^{d-1}-1
&\geq
(\tfrac{3}{2})2^d-1
\geq
2^d.
\end{split}
\end{equation}
Combining this and \cref{cost:inequalities} with \cref{Thm1} (applied with
$\varphi \curvearrowleft  0$,
$\gamma \curvearrowleft  1$,
$\scl \curvearrowleft  1$,
$a \curvearrowleft  a$,
$b \curvearrowleft  b$,
$\kappa \curvearrowleft  1$,
$f_d \curvearrowleft  f_d$
for $d\in\N$ in the notation of \cref{Thm1}) \proves[ps]\ that there exists $\mathfrak{c}\in(0,\infty)$ such that for all $d\in\N\cap(2,\infty)$, $L\in\N$, $\eps\in(0,1)$ it holds that
\begin{equation}
	\cost_{a,b}(f_d,L,\infty,\eps)\geq  2^{\nicefrac{d}{\max\{1,L-1\}}}
	\qandq
	\COST_{a,b}(f_d,\mathfrak{c}d^2\eps^{-1},1,\eps)
	\leq  \mathfrak{c}d^3\eps^{-2}.
\end{equation}
	\finishproofthus
\end{proof}

\subsection{Localizing functions}
\label{subsection:localizing_function}
\renewcommand{\shft}{\delta}

\begin{athm}{lemma}{density:one_side:base}
	Let $f\colon\R\to\R$ satisfy for all $x\in\R$ that
	\begin{equation}
		\label{eq:base:smooth}
		f(x)=
		\begin{cases}
			0 & \colon x\leq 0\\
			e^{-\frac{1}{x}} & \colon x>0.
		\end{cases}
	\end{equation}
	Then  
	\begin{enumerate}[(i)]
		\item{
			\label{it:localizer:base2}
			it holds that $f\in C^{\infty}(\R,\R)$ and}
		\item{
			\label{it:localizer:base1}
			it holds for all $x\in\R$ that $\vass{f'(x)}\leq 1$.}
	\end{enumerate}
\end{athm}

\begin{proof}[Proof of \cref{density:one_side:base}]
	\Nobs that \eqref{eq:base:smooth} \proves\ that
	\begin{equation}
		\lim_{h\searrow 0}\frac{f(h)-f(0)}{h}=\lim_{h\searrow 0}\frac{1}{h e^{\nicefrac{1}{h}}}=0=\lim_{h\nearrow 0}\frac{f(h)-f(0)}{h}.
	\end{equation}
	Combining this with \eqref{eq:base:smooth} \proves\ that for all $x\in\R$ it holds that
	\begin{equation}
		\label{eq:base:smooth1}
		\pr{\R\ni x\mapsto f(x)\in[0,1]}\in C^1(\R,[0,1])
		\qandq
		f'(x)=
		\begin{cases}
			0 & \colon x\leq 0\\
			\frac{1}{x^2}e^{-\frac{1}{x}} & \colon x>0.
		\end{cases}
	\end{equation}
	This \proves\ that
	\begin{equation}
		\label{eq:first:der:at:0}
		\lim_{h\searrow 0}\frac{f'(h)-f'(0)}{h}=\lim_{h\searrow 0}\frac{1}{h^3 e^{\nicefrac{1}{h}}}=0=\lim_{h\nearrow 0}\frac{f'(h)-f'(0)}{h}.
	\end{equation}
	\Moreover the chain and the product rule \prove\ that for all $g\in C^1(\R,\R)$, $x\in(0,\infty)$ with $\fa{y}(0,\infty)\colon \tgt(y)=\tfrac{1}{y^2}e^{-\frac{1}{y}}$ it holds that
	\begin{equation}
		\begin{split}
			g'(x)=-\tfrac{2}{x^3}e^{-\frac{1}{x}}+\tfrac{1}{x^4}e^{-\frac{1}{x}}
			=\pr*{\tfrac{1}{x}-2}\tfrac{1}{x^3}e^{-\frac{1}{x}}.
		\end{split}
	\end{equation}
	Combining this, \eqref{eq:base:smooth1}, and \eqref{eq:first:der:at:0} \proves\ that
	\begin{equation}
		\begin{split}
			\textstyle
		\sup_{x\in\R}\vass{f'(x)}
		&\textstyle=\max\pRb{\PRb{\sup_{x\in(-\infty,0]}\vass{f'(x)}},\PRb{\sup_{x\in(0,\infty)}\vass{f'(x)}}}	\\&\textstyle=\max\{0,\vass{f'\pr*{\tfrac{1}{2}}}\}
		=4e^{-2}\leq 1.
		\end{split}
	\end{equation}
	This  \proves[ei]\ \cref{it:localizer:base1}. \Moreover for all $n\in\N$, $ p\in\Z[X]$ with 
	\begin{equation}
		\forall\,x\in\R\colon f^{(n)}(x)=
		\begin{cases}
			0 & \colon x\leq 0\\
			p\pr*{\frac{1}{x}}e^{-\frac{1}{x}} & \colon x>0
		\end{cases}
	\end{equation}
	it holds that
	\begin{equation}
		\label{eq:differentialquotient:n-th:der}
		\lim_{h\searrow 0}\frac{f^{(n)}(h)-f^{(n)}(0)}{h}=\lim_{h\searrow 0}\frac{p\pr*{\frac{1}{h}}}{h e^{\nicefrac{1}{h}}}=0=\lim_{h\nearrow 0}\frac{f^{(n)}(h)-f^{(n)}(0)}{h}.
	\end{equation}
	\Moreover for all $p\in\Z[X]$, $g\in C^1(\R,\R)$, $x\in(0,\infty)$ with $\fa{y}(0,\infty)\colon \tgt(y)=p\pr{\tfrac{1}{y}}e^{-\frac{1}{y}}$ it holds that
	\begin{equation}
		\begin{split}
			g'(x)
			=-\tfrac{1}{x^2}p'\pr*{\tfrac{1}{x}}e^{-\frac{1}{x}}+p\pr*{\tfrac{1}{x}}\tfrac{1}{x^2}e^{-\frac{1}{x}}
			&=\pr*{p\pr*{\tfrac{1}{x}}\tfrac{1}{x^2}-\tfrac{1}{x^2}p'\pr*{\tfrac{1}{x}}}e^{-\frac{1}{x}}.
		\end{split}
	\end{equation}
	Combining this, 
	\eqref{eq:base:smooth},
	\eqref{eq:base:smooth1}, and
	\eqref{eq:differentialquotient:n-th:der} with 
	the fact that for all $p\in\Z[X]$ it holds that $p'\in\Z[X]$
	and induction \proves\ that for all $n\in\N_0$ there exists $p\in\Z[X]$ such that for all $x\in\R$ it holds that
	\begin{equation}
		f^{(n)}\in C^{1}(\R,\R)
		\qandq
		f^{(n)}(x)=
		\begin{cases}
			0 & \colon x\leq 0\\
			p\pr*{\frac{1}{x}}e^{-\frac{1}{x}} & \colon x>0.
		\end{cases}
	\end{equation}
	This  \proves[ei]\ \cref{it:localizer:base2}.
	The proof of \cref{density:one_side:base} is thus complete.
\end{proof}

\begin{athm}{lemma}{density:one_side}
	Let $\shft\in(0,\infty)$.
	Then there exists $\varphi\in C^{\infty}(\R,[0,1])$ such that
	\begin{enumerate}[(i)]
		\item{
			\label{density:lower}
			it holds for all $x\in(-\infty,0]$ that $\varphi(x)=0$,
		}
		\item{
			\label{density:middle}
			it holds for all $x\in(0,\shft)$ that $\varphi(x)\in(0,1)$,
		}
		\item{
			\label{density:upper}
			it holds for all $x\in[\shft,\infty)$ that $\varphi(x)=1$, and
		}
		\item{
			\label{density:der}
			it holds for all $x\in\R$ that $\vass{\varphi '(x)}\leq \frac{48}{\shft}$.
		}
	\end{enumerate}
\end{athm}

\begin{proof}[Proof of \cref{density:one_side}]
	Throughout this proof let $f\colon\R\to\R$ and $\varphi\colon\R\to\R$ satisfy for all $x\in\R$ that
	\begin{equation}
		\label{eq:base:smooth:3}
		f(x)=
		\begin{cases}
			0 & \colon x\leq 0\\
			e^{-\frac{\shft}{x}} & \colon x>0
		\end{cases}
		\qandq
		\varphi(x)=\frac{f(x)}{f(x)+f(\shft-x)}.
	\end{equation}
	\Nobs that \eqref{eq:base:smooth:3}, \cref{density:one_side:base}, and the fact that for all $x\in\R$ it holds that $f(x)+f(\shft-x)\geq f\pr*{\frac{\shft}{2}}=e^{-2}\geq \pr*{\tfrac{4}{11}}^2=\tfrac{16}{121}\geq \tfrac18$ \prove\ that for all $x\in(0,\shft)$ it holds that
	\begin{equation}
		\label{eq:for:item:well:base:localizing:function}
		2 f(\shft)\geq f(x)+f(\shft-x)\geq \tfrac18
		\qandq
		\varphi\in C^{\infty}(\R,\R).
	\end{equation}
	\Moreover \cref{density:one_side:base}, \eqref{eq:base:smooth:3}, and the chain rule \prove\ that for all $x\in\R$ it holds that
	\begin{equation}
		\label{eq:bound:der:in:n}
		\vass{f'(x)}
		\leq \tfrac{1}{\shft}.
	\end{equation}
	\Moreover \eqref{eq:base:smooth:3} \proves\ that for all $x\in(-\infty,0]$ it holds that
	\begin{equation}
		\label{eq:for:item:lower:base:localizing:function}
		\varphi(x)=\frac{f(x)}{f(x)+f(\shft-x)}=\frac{0}{0+e^{-\frac{\shft}{\shft-x}}}=0.
	\end{equation}
	\Moreover \eqref{eq:base:smooth:3} \proves\ that for all $x\in(0,\shft)$ it holds that
	\begin{equation}
		\label{eq:for:item:middle:base:localizing:function}
		0=\frac{0}{e^{-\frac{\shft}{x}}+e^{-\frac{\shft}{\shft-x}}}< \frac{e^{-\frac{\shft}{x}}}{e^{-\frac{\shft}{x}}+e^{-\frac{\shft}{\shft-x}}}=\varphi(x)= \frac{e^{-\frac{\shft}{x}}}{e^{-\frac{\shft}{x}}+e^{-\frac{\shft}{\shft-x}}}< \frac{e^{-\frac{\shft}{x}}}{e^{-\frac{\shft}{x}}}=1.
	\end{equation}
	\Moreover \eqref{eq:base:smooth:3} \proves\ that for all $x\in[\shft,\infty)$ it holds that
	\begin{equation}
		\label{eq:for:item:upper:base:localizing:function}
		\varphi(x)=\frac{f(x)}{f(x)+f(\shft-x)}=\frac{e^{-\frac{\shft}{x}}}{e^{-\frac{\shft}{x}}+0}=1.
	\end{equation}
	\Moreover \eqref{eq:base:smooth:3}, 
	\eqref{eq:for:item:well:base:localizing:function}, 
	\eqref{eq:bound:der:in:n},
	\cref{density:one_side:base},
	the quotient rule, and the fact that $e^{-1}\leq\frac{10}{27}$ \prove\ that for all $x\in(0,\shft)$ it holds that
	\begin{equation}
		\label{eq:for:item:der:bound:base:localizing:function}
		\begin{split}
			\vass{\varphi'(x)}
			&=
			\vass[\bigg]{\frac{f'(x)\pr*{f(x)+f(\shft-x)}-\pr*{f'(x)-f'(\shft-x)}f(x)}{\pr*{f(x)+f(\shft-x)}^2}}
			\\&=
			\vass[\bigg]{\frac{f'(x)f(\shft-x)+f'(\shft-x)f(x)}{\pr*{f(x)+f(\shft-x)}^2}}
			\\&\leq
			\frac{\vass{f'(x)f(\shft-x)}+\vass{f'(\shft-x)f(x)}}{\pr*{f(x)+f(\shft-x)}^2}
			\\&\leq
			\tfrac{8^2}{\shft}\pr*{\vass{f(\shft-x)}+\vass{f(x)}}\leq \tfrac{64\pr*{2f(\shft)}}{\shft}=\tfrac{128e^{-1}}{\shft}\leq \tfrac{1280}{27\shft}\leq \tfrac{48}{\shft}
			.
		\end{split}
	\end{equation}
	This, \eqref{eq:for:item:well:base:localizing:function}, \eqref{eq:for:item:lower:base:localizing:function}, \eqref{eq:for:item:middle:base:localizing:function}, and \eqref{eq:for:item:upper:base:localizing:function}  \prove[ei]\ \cref{density:lower,density:middle,density:upper,density:der}.
	The proof of \cref{density:one_side} is thus complete.
\end{proof}

\begin{athm}{lemma}{density:two_side}
	Let $a\in\R$, $b\in(a,\infty)$, $\shft\in(0,\infty)$. Then there exists $\varphi\in C^{\infty}(\R,[0,1])$ such that
	\begin{enumerate}[(i)]
		\item{
			\label{density:lower2}
			it holds for all $x\in(-\infty,a-\shft]\cup[b+\shft,\infty)$ that $\varphi(x)=0$,
		}
		\item{
			\label{density:middle2}
			it holds for all $x\in(a-\shft,a)\cup(b,b+\shft)$ that $\varphi(x)\in(0,1)$,
		}
		\item{
			\label{density:upper2}
			it holds for all $x\in[a,b]$ that $\varphi(x)=1$, and
		}
		\item{
			\label{density:der2}
			it holds for all $x\in\R$ that $\vass{\varphi'(x)}\leq \frac{48}{\shft}$.
		}
	\end{enumerate}
\end{athm}

\begin{proof}[Proof of \cref{density:two_side}]
	\Nobs that \cref{density:one_side} \proves\ that there exists ${}f{}\in C^{\infty}(\R,[0,1])$ which satisfies that
	\begin{enumerate}[(I)]
		\item{
			\label{it:density:lower}
			it holds for all $x\in(-\infty,0]$ that ${}f{}(x)=0$,
		}
		\item{
			\label{it:density:middle}
			it holds for all $x\in(0,\shft)$ that ${}f{}(x)\in(0,1)$,
		}
		\item{
			\label{it:density:upper}
			it holds for all $x\in[\shft,\infty)$ that ${}f{}(x)=1$, and
		}
		\item{
			\label{it:density:der}
			it holds for all $x\in\R$ that $\vass{{}f{} '(x)}\leq \tfrac{48}{\shft}$.
		}
	\end{enumerate}
	Next let $\varphi\colon\R\to[0,1]$ satisfy for all $x\in\R$ that
	\begin{equation}
		\label{fdef:density:two:sides}
		\varphi(x)=
		\begin{cases}
			f(x-a+\shft)&\colon x< a\\
			f(b-x+\shft)&\colon x\geq a
		\end{cases}.
	\end{equation}
	\Nobs that \eqref{fdef:density:two:sides} and \cref{it:density:lower} \prove\ that for all $x\in(-\infty,a-\shft]$, $y\in[b+\shft,\infty)$ it holds that
	\begin{equation}
		\label{eq:two:sides:outer}
		\varphi(x)=f(x-a+\shft)=0
		\qandq
		\varphi(y)=f(b-y+\shft)=0.
	\end{equation}
	\Moreover \eqref{fdef:density:two:sides} and \cref{it:density:middle} \prove\ that for all $x\in(a-\shft,a)$, $y\in(b,b+\shft)$ it holds that
	\begin{equation}
		\label{eq:two:sides:inbtw}
		\varphi(x)=f(x-a+\shft)\in(0,1)
		\qandq
		\varphi(y)=f(b-y+\shft)\in(0,1).
	\end{equation}
	\Moreover \eqref{fdef:density:two:sides} and \cref{it:density:upper} \prove\ that for all $x\in[a,b]$ it holds that
	\begin{equation}
		\label{eq:two:sides:inner}
		\varphi(x)=f(b-x+\shft)=1
		.
	\end{equation}
	\Moreover 
	\eqref{fdef:density:two:sides},
	\cref{it:density:upper},
	and the fact that for all $k\in\N_0$ it holds that $f\in C^{\infty}(\R,[0,1])$ and $(-1)^k f^{(k)}(b-a+\shft)=f^{(k)}(\shft)$ \prove\ that for all $k\in\N_0$ with 
	$\varphi\in C^k(\R,[0,1])$, 
	$\fa{x}(-\infty,a)\colon\varphi^{(k)}(x)=f^{(k)}(x-a+\shft)$, and 
	$\fa{x}[a,\infty)\colon\varphi^{(k)}(x)=(-1)^k f^{(k)}(b-x+\shft)$ it holds that
	\begin{equation}
		\label{eq:two:sides:diff}
		\begin{split}
			\lim_{h\nearrow 0}\frac{\varphi^{(k)}(a+h)-\varphi^{(k)}(a)}{h}
			&=\lim_{h\nearrow 0}\frac{f^{(k)}(h+\shft)-(-1)^k f^{(k)}(b-a+\shft)}{h}
			\\&=\lim_{h\nearrow 0}\frac{f^{(k)}(\shft+h)-f^{(k)}(\shft)}{h}
			\\&=f^{(k+1)}(\shft)
			\\&=(-1)^{k+1}f^{(k+1)}(b-a+\shft)
			\\&=\lim_{h\searrow 0}\frac{\varphi^{(k)}(a+h)-\varphi^{(k)}(a)}{h}
			.
		\end{split}
	\end{equation}
	Combining this,
	\eqref{fdef:density:two:sides}, 
	and the fact that $\lim_{x\nearrow a}\varphi(x)=1=\lim_{x\searrow a}\varphi(x)$
	with induction \proves\ that for all $x\in(-\infty,a)$, $y\in[a,\infty)$ it holds that
	\begin{equation}
		\label{eq:two:sides:smooth}
		\varphi\in C^{\infty}(\R,[0,1]),
		\qquad
		\varphi'(x)=f'(x-a+\shft)
		\qandq
		\varphi'(y)=-f'(b-y+\shft)
		.
	\end{equation}
	This and \cref{it:density:der} \prove\ that for all $x\in\R$ it holds that
	\begin{equation}
		\vass{\varphi'(x)}\leq\max\{\vass{f'(x-a+\shft)},\vass{f'(b-x+\shft)}\}\leq \tfrac{48}{\shft}.
	\end{equation}
	This, 
	\eqref{eq:two:sides:outer}, 
	\eqref{eq:two:sides:inbtw}, 
	\eqref{eq:two:sides:inner}, and
	\eqref{eq:two:sides:smooth}  \prove[ei]\ \cref{density:lower2,density:middle2,density:upper2,density:der2}.
	The proof of \cref{density:two_side} is thus complete.
\end{proof}

\cfclear
\begin{athm}{lemma}{density:two_side_d}
	Let $a\in\R$, $b\in(a,\infty)$, $\shft\in(0,\infty)$, $d\in\N$. Then there exists $\varphi\in C^{\infty}(\R^d,[0,1])$ such that
	\begin{enumerate}[(i)]
		\item{
			\label{density:lower2_d}
			it holds for all $x\in\R^d\backslash(a-\shft,b+\shft)^d$ that $\varphi(x)=0$,
		}
		\item{
			\label{density:upper2_d}
			it holds for all $x\in[a,b]^d$ that $\varphi(x)=1$, and
		}
		\item{
			\label{density:der2_d}
			it holds for all $x,y\in\R^d$ that $\vass{\varphi(x)-\varphi(y)}\leq \frac{48d}{\shft}\norm{x-y}$
		}
	\end{enumerate}
	\cfout.
\end{athm}

\begin{proof}[Proof of \cref{density:two_side_d}]
	\Nobs that \cref{density:two_side} (applied with 
	$a \curvearrowleft a$,
	$b \curvearrowleft b$,
	$\shft \curvearrowleft \shft$
	in the notation of \cref{density:two_side}) \proves\ that there exists $f\in C^{\infty}(\R,[0,1])$ which satisfies that
	\begin{enumerate}[(I)]
		\item{
			\label{it:density:lower2}
			it holds for all $x\in(-\infty,a-\shft]\cup[b+\shft,\infty)$ that ${}f{}(x)=0$,
		}
		\item{
			\label{it:density:middle2}
			it holds for all $x\in(a-\shft,a)\cup(b,b+\shft)$ that ${}f{}(x)\in(0,1)$,
		}
		\item{
			\label{it:density:upper2}
			it holds for all $x\in[a,b]$ that ${}f{}(x)=1$, and
		}
		\item{
			\label{it:density:der2}
			it holds for all $x\in\R$ that $\vass{{}f{}'(x)}\leq \frac{48}{\shft}$.
		}
	\end{enumerate}
	Next let $\varphi\in C^{\infty}(\R^d,[0,1])$ satisfy for all $x=(x_1,\ldots,x_d)\in\R^d$ that
	\begin{equation}
		\label{eq:fdef:product:density}
		\varphi(x)=\textstyle\prod_{i=1}^d f(x_i).
	\end{equation}
	\Nobs that \eqref{eq:fdef:product:density} and \cref{it:density:lower2} \prove\ that for all $x=(x_1,\ldots,x_d)\in\R^d\setminus(a-\shft,b+\shft)^d$ it holds that
	\begin{equation}
		\label{eq:prduct:density:outer}
		\varphi(x)=\textstyle\prod_{i=1}^d f(x_i)=0.
	\end{equation}
	\Moreover \eqref{eq:fdef:product:density} and \cref{it:density:upper2} \prove\ that for all $x=(x_1,\ldots,x_d)\in[a,b]^d$ it holds that
	\begin{equation}
		\label{eq:prduct:density:inner}
		\varphi(x)=\textstyle\prod_{i=1}^d f(x_i)=\textstyle\prod_{i=1}^d 1=1.
	\end{equation}
	\Moreover \eqref{eq:fdef:product:density}, \cref{it:density:der2}, and the fact that $f\in C^{\infty}(\R,[0,1])$ \prove\ that for all $x=(x_1,\ldots,x_d)$, $y=(y_1,\ldots,y_d)\in\R^d$ it holds that
	\begin{equation}
		\label{eq:prduct:density:lipschitz}
		\begin{split}
			\vass{\varphi(x)-\varphi(y)}
			&=\vass*{\PR*{\textstyle\prod_{i=1}^d f(x_i)}-\PR*{\textstyle\prod_{i=1}^d f(y_i)}}
			\\&=\vass*{\textstyle\sum_{j=1}^d\pr*{\PR*{\textstyle\prod_{i=1}^{j-1} f(y_i)}\PR*{\textstyle\prod_{i=j}^d f(x_i)}-\PR*{\textstyle\prod_{i=1}^{j} f(y_i)}\PR*{\textstyle\prod_{i=j+1}^d f(x_i)}}}
			\\&=\vass*{\textstyle\sum_{j=1}^d(f(x_j)-f(y_j))\PR*{\textstyle\prod_{i=1}^{j-1} f(y_i)}\PR*{\textstyle\prod_{i=j+1}^d f(x_i)}}
			\\&\leq\textstyle\sum_{j=1}^d\vass*{f(x_j)-f(y_j)}\PR*{\textstyle\prod_{i=1}^{j-1} \vass*{f(y_i)}}\PR*{\textstyle\prod_{i=j+1}^d \vass*{f(x_i)}}
			\\&\leq\textstyle\sum_{j=1}^d\vass*{f(x_j)-f(y_j)}
			\leq\textstyle\sum_{j=1}^d \tfrac{48}{\shft}\vass*{x_j-y_j}
			\leq \tfrac{48d}{\shft}\norm*{x-y}
			.
		\end{split}
	\end{equation}
	Combining this, \eqref{eq:prduct:density:outer}, and \eqref{eq:prduct:density:inner}  \proves[ei]\ \cref{density:lower2_d,density:upper2_d,density:der2_d}.
	The proof of \cref{density:two_side_d} is thus complete.
\end{proof}

\newcommand{\lip}{L}
\begin{athm}{cor}{cor:localized:sum}
	Let $a\in\R$, $b\in(a,\infty)$, $d\in\N$, $\kappa,\shft,\lip\in(0,\infty)$ and $g\in C^{\infty}(\R^d,\R)$ satisfy for all $x,y\in[a-\shft,b+\shft]^d$ that $\vass{\tgt(x)-\tgt(y)}\leq \lip \norm{x-y}$ and $\vass{\tgt(x)}\leq\kappa$.
	Then there exists $f\in C^{\infty}(\R^d,\R)$ such that
	\begin{enumerate}[(i)]
		\item{
			\label{cor:localized:sum:prop:funct}
			it holds for all $x\in[a,b]^d$ that $f(x)=\tgt(x)$,
		}
		
		\item{
			\label{cor:localized:sum:prop:funct:outside}
			it holds for all $x\in\R^d$ that $\vass{f(x)}\leq \kappa\prb{\Indfct{[a-\shft,b+\shft]^d}(x)}$, and
		}
		
		\item{
			\label{cor:localized:sum:prop:funct:lipschitz}
			it holds for all $x,y\in\R^d$ that $\vass{f(x)-f(y)}\leq (\tfrac{48 \kappa d}{\shft}+\lip)\norm{x-y}$
		}
	\end{enumerate}
	\cfout.
\end{athm}

\begin{proof}[Proof of \cref{cor:localized:sum}]
	\Nobs that \cref{density:two_side_d} (applied with
	$a \curvearrowleft a$,
	$b \curvearrowleft b$,
	$d \curvearrowleft d$,
	$\shft \curvearrowleft \shft$
	in the notation of \cref{density:two_side_d}) \proves\ that there exists $\varphi\in C^{\infty}(\R^d,[0,1])$ which satisfies that
	\begin{enumerate}[(I)]
		\item{
			\label{it:density:lower2_d}
			it holds for all $x\in\R^d\backslash(a-\shft,b+\shft)^d$ that $\varphi(x)=0$,
		}
		\item{
			\label{it:density:upper2_d}
			it holds for all $x\in[a,b]^d$ that $\varphi(x)=1$, and
		}
		\item{
			\label{it:density:der2_d}
			it holds for all $x,y\in\R^d$ that $\vass{\varphi(x)-\varphi(y)}\leq \frac{48 d}{\shft}\norm{x-y}$.
		}
	\end{enumerate}
	Next let $f\in C^{\infty}(\R^d,\R)$ satisfy for all $x\in\R^d$ that 
	\begin{equation}
		\label{eq:fdef:localized:testfunction:Thm6}
		f(x)=\varphi(x) \tgt(x).
	\end{equation}
	\Nobs that \eqref{eq:fdef:localized:testfunction:Thm6} and \cref{it:density:lower2_d} \prove\ that for all $x\in\R^d\setminus(a-\shft,b+\shft)^d$ it holds that
	\begin{equation}
		\label{eq:outer:localized:testfunction:Thm6}
		f(x)=\varphi(x) \tgt(x)=0.
	\end{equation}
	\Moreover \eqref{eq:fdef:localized:testfunction:Thm6} and \cref{it:density:upper2_d} \prove\ that for all $x\in[a,b]^d$ it holds that
	\begin{equation}
		\label{eq:inner:localized:testfunction:Thm6}
		f(x)=\varphi(x) \tgt(x)= \tgt(x).
	\end{equation}
	\Moreover 
	\eqref{eq:fdef:localized:testfunction:Thm6},
	\eqref{eq:outer:localized:testfunction:Thm6}, 
	\cref{it:density:der2_d}, and the fact that $\varphi\in C^{\infty}(\R^d,[0,1])$ \prove\ that for all $x$, $y\in[a-\shft,b+\shft]^d$ it holds that $\vass{f(x)}\leq \kappa$ and
	\begin{equation}
		\llabel{eq:lipsch:inner}
		\begin{split}
			\vass{f(x)-f(y)}
			&=\vass{\varphi(x)\tgt(x)-\varphi(y)\tgt(y)}
			\\&\leq\vass{
				\varphi(x)\tgt(x)
				-\varphi(y)\tgt(x)}
			+\vass{
				\varphi(y)\tgt(x)
				-\varphi(y)\tgt(y)}
			\\&\leq
			\vass*{\varphi(x)-\varphi(y)}\vass{\tgt(x)}
			+\vass{\varphi(y)}\vass{\tgt(x)-\tgt(y)}
			\\&\leq
			\tfrac{48 \kappa d}{\shft}\norm{x-y}+\lip\norm{x-y}
			\leq
			(\tfrac{48 \kappa d}{\shft}+\lip)\norm{x-y}
			.
		\end{split}
	\end{equation}
	This, 
	\eqref{eq:outer:localized:testfunction:Thm6}, and the fact that $f\in C(\R^d,\R)$ \prove\ that for all $x,y\in\R^d$ it holds that
	\begin{equation}
		\vass{f(x)-f(y)}\leq	(\tfrac{48 \kappa d}{\shft}+\lip)\norm{x-y}.
	\end{equation}
	Combining this,
	\eqref{eq:fdef:localized:testfunction:Thm6},
	\eqref{eq:outer:localized:testfunction:Thm6}, and
	\eqref{eq:inner:localized:testfunction:Thm6}
	 \proves[ei]\ \cref{cor:localized:sum:prop:funct,cor:localized:sum:prop:funct:outside,cor:localized:sum:prop:funct:lipschitz}.
	The proof of \cref{cor:localized:sum} is thus complete.
\end{proof}

\subsection{ANN approximations for classes of smooth and bounded functions}
\label{subsection:sum_results}

\cfclear
\begin{athm}{theorem}{Thm6}
Let $\varphi\in\R$, $\gamma,\kappa\in(0,\infty)$, $a\in\R$, $b\in[a+\pi\gamma^{-1},\infty)$, $\mathfrak{c}\in(0,\infty)$ satisfy $\mathfrak{c}\geq4504\ceil{\log_2\pr*{\max\{1,\gamma\}\max\{\abs{a},\abs{b},2\}}}\max\{\kappa^3,1\}$ and
for every $d\in\N$ let
$f_d\colon\R^d\to\R$ satisfy for all $x=(x_1,\ldots,x_d)\in\R^d$ that $f_d(x)=\kappa\sin\prb{\gamma 2^d\prb{\ssum_{i = 1}^d x_i}+\varphi}$.
Then it holds for all $d,L\in\N$, $\eps\in(0,\kappa)$ that 
\begin{equation}
\label{it:Thm6:eq1}
\begin{split}
\cost_{a,b}(f_d,L,\infty,\eps)
\geq  \pr{d2^{d-1}}^{\nicefrac{1}{\max\{1,L-1\}}}
\qandq
\COST_{a,b}(f_d,\mathfrak{c}d,1,\eps) 
\leq  \mathfrak{c}d\eps^{-2}
\end{split}
\end{equation}
\cfout.
\end{athm}

\begin{proof}[Proof of \cref{Thm6}]
	Throughout this proof let $g\colon\R\to\R$ satisfy for all $x\in\R$ that $\tgt(x)=\kappa\sin(x+\varphi)$.
	\Nobs that \cref{cor:rescaled:final:approximation3} (applied with
	$a \curvearrowleft  a$,
	$b \curvearrowleft  b$,
	$d \curvearrowleft  d$,
	$\kappa \curvearrowleft  \kappa$,
	$\eps \curvearrowleft  \eps$,
	$\gamma \curvearrowleft  \gamma$,
	$\mathfrak{c} \curvearrowleft \mathfrak{c}$,
	$g \curvearrowleft  g$
	for $d\in\N$, $\eps\in(0,\kappa)$ in the notation of \cref{cor:rescaled:final:approximation3}) \proves\ that for all $d\in\N$, $\eps\in(0,\kappa)$ it holds that
	\begin{equation}
		\label{eq:Thm6:prove}
		\begin{split}
\COST_{a,b}(f_d,\mathfrak{c}d,1,\eps) \leq  \mathfrak{c}d\eps^{-2}.
		\end{split}
	\end{equation}
\Moreover \cref{LowerBound:sum} (applied with
$\varphi \curvearrowleft  \varphi$,
$\kappa \curvearrowleft  \kappa$,
$\gamma \curvearrowleft  \gamma$,
$a \curvearrowleft  a$,
$b \curvearrowleft  b$,
$f_d \curvearrowleft f_d$
for $d\in\N$ in the notation of \cref{LowerBound:sum}) \proves\ that for all $d,L\in\N$, $\eps\in(0,\kappa)$ it holds that
\begin{equation}
	\begin{split}
\cost_{a,b}(f_d,L,\infty,\eps)\geq  \pr{d2^{d-1}}^{\nicefrac{1}{\max\{1,L-1\}}}.
	\end{split}
\end{equation} 
This and \eqref{eq:Thm6:prove}  \prove[ei]\ \eqref{it:Thm6:eq1}.
The proof of \cref{Thm6} is thus complete.
\end{proof}

\cfclear
\begin{athm}{cor}{Thm6.1}
Let $a\in\R$, $b\in(a,\infty)$, $\kappa,\delta\in(0,\infty)$. 
Then there exist $\mathfrak{c}\in(0,\infty)$ and $f_d\in C^{\infty}(\R^d,\R)$, $d\in\N$, with $\fa{d}\N, x\in\R^d\colon\vass{f_d(x)}\leq\kappa\Indfct{[a-\delta,b+\delta]^d}(x)$ such that for all $d,L\in\N$, $\eps\in(0,\kappa)$ it holds that
\begin{equation}
	\label{it:Thm6.1:eq1}
	\begin{split}
		\cost_{a,b}(f_d,L,\infty,\eps)
		\geq  \prb{d2^{d-1}}^{\nicefrac{1}{\max\{1,L-1\}}}
		\qandq
		\COST_{a,b}(f_d,\mathfrak{c}d,1,\eps) 
		\leq  \mathfrak{c}d\eps^{-2}
	\end{split}
\end{equation}
\cfout.
\end{athm}

\begin{proof}[Proof of \cref{Thm6.1}]
Throughout this proof let $g_d\in C^{\infty}(\R^d,\R)$, $d\in\N$, satisfy for all $d\in\N$, $x=(x_1,\ldots,x_d)\in\R^d$ that 
\begin{equation}
\label{eq:Thm6.1:fdef}
g_d(x)=\kappa\sin\prb{ \tfrac{2^d\pi}{b-a}\prb{\textstyle\ssum_{i = 1}^d x_i}}.
\end{equation}
\Nobs that \cref{cor:localized:sum} (applied with
$a \curvearrowleft a$,
$b \curvearrowleft b$,
$d \curvearrowleft d$,
$\kappa \curvearrowleft \kappa$,
$\delta \curvearrowleft \delta$,
$L\curvearrowleft \kappa 2^d d\pi(b-a)^{-1}$,
$g \curvearrowleft g_d$
for $d\in\N$ in the notation of \cref{cor:localized:sum}) \proves\ that there exist $f_d\in C^{\infty}(\R^d,\R)$, $d\in\N$, which satisfy that
\begin{enumerate}[(I)]
\item{
\label{cor:localized:sum:prop:funct:Thm6.1}
it holds for all $d\in\N$, $x=(x_1,\ldots,x_d)\in[a,b]^d$ that $f_d(x)=\kappa\sin\prb{ \tfrac{2^d\pi}{b-a}\prb{\textstyle\ssum_{i = 1}^d x_i}}$, and
}

\item{
\label{cor:localized:sum:prop:funct:outside:Thm6.1}
it holds for all $d\in\N$, $x\in\R^d$ that $\vass{f_d(x)}\leq \kappa\prb{\Indfct{[a-\delta,b+\delta]^d}(x)}$.
}

\end{enumerate}
\Nobs that 
\cref{cor:localized:sum:prop:funct:Thm6.1},
the fact that for all $d,L\in\N$ it holds that $\max\{L-1,1\}\leq L$, $f_d|_{[a,b]^d}=g_d|_{[a,b]^d}$, and
\cref{Thm6} (applied with
$\varphi \curvearrowleft 0$,
$\gamma \curvearrowleft \tfrac{\pi}{b-a}$,
$a \curvearrowleft a$,
$b \curvearrowleft b$,
$\kappa \curvearrowleft \kappa$,
$f_d \curvearrowleft g_d$
for $d\in\N$ in the notation of \cref{Thm6}) \prove\ that there exists $\mathfrak{c}\in(0,\infty)$ such that for all $d,L\in\N$, $\eps\in(0,\kappa)$ it holds that 
\begin{equation}
	\label{prop:it:Thm6:eq1}
	\begin{split}
		\cost_{a,b}(f_d,L,\infty,\eps)
		\geq \prb{d2^{d-1}}^{\nicefrac{1}{\max\{1,L-1\}}}
		\qandq
		\COST_{a,b}(f_d,\mathfrak{c}d,1,\eps) 
		\leq  \mathfrak{c}d\eps^{-2}
	\end{split}
\end{equation}
\cfload.
Combining 
\cref{cor:localized:sum:prop:funct:Thm6.1},
\cref{cor:localized:sum:prop:funct:outside:Thm6.1}, and
\eqref{prop:it:Thm6:eq1}
  \proves[ei]\ \eqref{it:Thm6.1:eq1}.
The proof of \cref{Thm6.1} is thus complete.
\end{proof}

\cfclear
\begin{athm}{cor}{Cor:Thm6.1}
	Let $a\in\R$, $b\in(a,\infty)$, $\kappa,\delta\in(0,\infty)$. 
	Then there exist $\mathfrak{c}\in(0,\infty)$ and $f_d\in C^{\infty}(\R^d,\R)$, $d\in\N$, with $\fa{d}\N, x\in\R^d\colon\vass{f_d(x)}\leq\kappa\Indfct{[a-\delta,b+\delta]^d}(x)$ such that for all $d,L\in\N$, $\eps\in(0,\kappa)$ it holds that
	\begin{equation}
		\label{it:Cor:Thm6.1:eq1}
		\begin{split}
			\COST_{a,b}(f_d,L,\infty,\eps)
			\geq  2^{\nicefrac{d}{\max\{1,L-1\}}}
			\qandq
			\COST_{a,b}(f_d,\mathfrak{c}d,1,\eps) 
			\leq  \mathfrak{c}d\eps^{-2}
		\end{split}
	\end{equation}
	\cfout.
\end{athm}

\begin{proof}[Proof of \cref{Cor:Thm6.1}]
\Nobs that \cref{Thm6.1} and the fact that for all $d\in\N\setminus\{1\}$, $\f\in\ANNSET$ it holds that $\paramANN\pr{\f}\geq 2$ and $d2^{d-1}\geq 2^d$ \prove[ep]\ \eqref{it:Cor:Thm6.1:eq1}.
\finishproofthus
\end{proof}

\newcommand{\pot}{5^7}
\newcommand{\B}{B}
\newcommand{\cnt}{\mathfrak{B}}
\cfclear

\section{Lower and upper bounds for the amount of memory needed to store ANN approximations of certain high-dimensional functions}

\label{subsection:entropy_results}
 In this section we study the amount of memory needed on a computer to store \ANNs\ that approximate suitable globally bounded and infinitely often differentiable functions.
 More specifically, in \cref{cor:bit:both:bounds} below we show that the number of bits needed on a computer to store the approximating \ANN\ must grow at least exponentially in the input dimension of the considered target function for shallow or insufficiently deep \ANNs\ but grows at most polynomially in the input dimension of the considered target function and the reciprocal of the approximation accuracy for sufficiently deep \ANNs\ (see \cref{rmk:1} below). 
 Our proof of \cref{cor:bit:both:bounds} employs 
 \begin{itemize}
 	\item{
 		\textit{the lower bound result} for the number of bits needed to store \ANNs\ approximating the considered globally bounded and infinitely often differentiable target functions in \cref{cor:bit:lower:bound:v2} and}
 	\item{\textit{the upper bound result} for the number of bits needed to store \ANNs\ approximating the considered globally bounded and infinitely often differentiable target functions in \cref{cor:bit:upper:bound}.}
 \end{itemize}
In \cref{cor:bit:lower:bound:v2} we bound the number of bits needed to store the approximating \ANN\ from below by the number of non-zero \ANN\ parameters 
and
in \cref{cor:bit:upper:bound} we bound the number of bits to store the approximating \ANN\ from above by the product of the number of \ANN\ parameters and the maximal number of bits needed to store each parameter of the \ANN.

\Cref{cor:bit:lower:bound:v2} follows from \cref{cor:bit:lower:bound:v1} which provides lower bounds for the minimal number of non-zero \ANN\ parameters to approximate compositions of certain periodic functions and the product functions. \Cref{cor:bit:lower:bound:v1}, in turn, is a direct consequence of the lower bound result in \cref{it:Cor:Thm1:eq1} in \cref{Cor:Thm1}. The upper and lower bound results in \cref{Cor:Thm1} follow directly from  \cref{Thm1}.

Our proof of \cref{cor:bit:upper:bound} is based on an application of \cref{lem:bit:upper:bound} and an application of the upper bound result for the minimal number of \ANN\ parameters to approximate compositions of certain periodic functions and the product functions in \cref{final:approximation:imp}. \Cref{lem:bit:upper:bound}, in turn, follows from \cref{lem:bit:approx:error} which provides, among other things, upper bounds for the maximal number of bits needed to store the parameters of a suitably discretized \ANN\ approximating a non-discretized \ANN.

\subsection{Upper bounds for the number of bits to store ANN approximations}

\cfclear
\begin{athm}{lemma}{lem:bit:approx:error}
	Let $\f\in\ANNSET$, $\eps\in(0,1]$, $R\in[1,\infty)$, $\B\in\N\setminus\{1\}$ satisfy $\size(\f)\leq 1=\outDimANN(\f)$.
	Then there exist $\cnt\in\N$, $\g\in\ANNSET$
	such that 
	\begin{enumerate}[(i)]
		\item \label{it:label:dims}it holds that $\dims(\f)=\dims(\g)$,
		\item \label{it:label:approx}it holds that $\sup_{x\in[-R,R]^{\inDimANN(\f)}}\vass{(\realisation(\f))(x)-(\realisation(\g))(x)}\leq\eps$,
		\item \label{it:label:bits} it holds that $\cnt\leq \lengthANN(\f)\log_{\B}(4\paramANN(\f))+\log_{\B}(R\eps^{-1})+2		$, and
		\item \label{it:label:discretized}it holds for all $p\in\{1,2,\ldots,\param(\g)\}$ that there exist $\beta_1,\beta_2,\ldots,\beta_{\cnt}\in\N_0\cap[0,\B)$ such that
		\begin{equation}
			\form{p}{\g}
			=	(-1)^{\beta_\cnt}\textstyle\sum_{i=1}^{\cnt-1} \beta_i \B^{-i},
		\end{equation}
	\end{enumerate}
	\cfout.
\end{athm}

\begin{proof}[Proof of \cref{lem:bit:approx:error}]
	Throughout this proof 
	let $\delta\in\R$, $\cnt\in\N$ satisfy $\delta=B^{1-\cnt}$ and $\cnt=\ceil{\lengthANN(\f)\log_{\B}(4\paramANN(\f))+\log_{\B}(R\eps^{-1})}+1$,
	let $\g\in\ANNSET$ satisfy for all $p\in\{1,2,\ldots,\param(\f)\}$ that there exist $\beta_1,\beta_2,\ldots,\beta_{\cnt}\in\N_0\cap[0,\B)$ such that
	\begin{equation}
		\label{eq:setup:disc:netw}
		\dims(\g)=\dims(\f),
		\qquad
		\vass{\form{p}{\g}-\form{p}{\f}}\leq\delta,
		\qandq
		\form{p}{\g}
		=	(-1)^{\beta_{\cnt}}\textstyle\sum_{i=1}^{\cnt-1} \beta_i \B^{-i},
	\end{equation}
	let $L,l_0,l_1,\ldots,l_L\in\N$ satisfy $\dims(\f)=(l_0,l_1,\ldots,l_L)$,
	let $x_j\in\R^{l_j}$, $j\in\N_0\cap[0,L]$, 
	and
	let $y_j\in\R^{l_j}$, $j\in\N_0\cap[0,L]$,
	satisfy for all $j\in\N\cap(0,L]$ that
	\begin{equation}
		\label{eq:setup:perturbated}
		x_0=y_0\in[-R,R]^{l_0},
		\qquad
		x_j=\RELU(\weight{j}{\f}x_{j-1}+\bias{j}{\f}),
		\qquad\text{and}\qquad
		y_j=\RELU(\weight{j}{\g}y_{j-1}+\bias{j}{\g})
	\end{equation}
	\cfload. \Nobs that \eqref{eq:setup:perturbated} and induction \prove\ that for all $j\in\N_0\cap[0,L]$ it holds that
	\begin{equation}
		\llabel{eq:sump:bd}
		\supn{x_j}\leq \max\{\size(\f),1\}^j\prb{\textstyle\prod_{k=0}^{j-1}(l_k+1)}R.
	\end{equation}
	This, \eqref{eq:setup:perturbated}, and the fact that for all $d\in\N$, $a,b\in\R^d$ it holds that $\supn{\RELU(a)-\RELU(b)}\leq \supn{a-b}$ \prove\ that for all $j\in\N_0\cap[0,L)$ with $\supn{x_j-y_j}\leq \max\{\size(\f)+\delta,1\}^{j}
	\pr{\prod_{k=0}^{j-1}(2l_k+1)}R\delta$ it holds that
	\begin{equation}
		\label{eq:distance:j+1st:iteration:perturbated:netw}
		\begin{split}
			&\supn{x_{j+1}-y_{j+1}}
			\\&\leq\supn{(\weight{j+1}{\f}x_{j}+\bias{j+1}{\f})-(\weight{j+1}{\g}y_{j}+\bias{j+1}{\g})}
			\\&=\supn{(\weight{j+1}{\f}-\weight{j+1}{\g})x_{j}+\weight{j+1}{\g}(x_j-y_j)+\bias{j+1}{\f}-\bias{j+1}{\g}}
			\\&\leq \supn{\weight{j+1}{\f}-\weight{j+1}{\g}}l_j\supn{x_{j}}+\supn{\weight{j+1}{\g}}l_j\supn{x_{j}-y_j}+\supn{\bias{j+1}{\f}-\bias{j+1}{\g}}
			\\&\leq \delta l_j\max\{\size(\f),1\}^j\prb{\textstyle\prod_{k=0}^{j-1}(l_k+1)}R
			\\&\quad+(\size(\f)+\delta)l_j\max\{\size(\f)+\delta,1\}^{j}
			\prb{\textstyle\prod_{k=0}^{j-1}(2l_k+1)}R\delta+\delta
			\\&\leq \prb{ l_j\max\{\size(\f),1\}^j+l_j \max\{\size(\f)+\delta,1\}^{j+1}+1}
			\prb{\textstyle\prod_{k=0}^{j-1}(2l_k+1)}R\delta
			\\&\leq \prb{l_j + l_j+1}\max\{\size(\f)+\delta,1\}^{j+1}
			\prb{\textstyle\prod_{k=0}^{j-1}(2l_k+1)}R\delta
			\\&\leq \max\{\size(\f)+\delta,1\}^{j+1}
			\prb{\textstyle\prod_{k=0}^{j}(2l_k+1)}R\delta.
		\end{split}
	\end{equation}
	Combining this, \eqref{eq:setup:perturbated}, and \lref{eq:sump:bd} with induction \proves\ that 
	\begin{equation}
		\llabel{eq:1step:bd}
		\begin{split}
			&\vass{(\realisation(\f))(x_0)-(\realisation(\g))(x_0)}
			\\&=
			\vass{\weight{L}{\f}x_{L-1}+\bias{L}{\f}-\pr{\weight{L}{\g}y_{L-1}+\bias{L}{\g}}}
			\\&\leq
			\vass{\pr{\weight{L}{\f}-\weight{L}{\g}}x_{L-1}}
			+\vass{\weight{L}{\g}(y_{L-1}-x_{L-1})}
			+
			\vass{\bias{L}{\f}-\bias{L}{\g}}
			\\&\leq
			\supn{\weight{L}{\f}-\weight{L}{\g}}l_{L-1}\supn{x_{L-1}}
			+\supn{\weight{L}{\g}}l_{L-1}\supn{x_{L-1}-y_{L-1}}
			+\supn{\bias{L}{\f}-\bias{L}{\g}}
			\\&\leq
			\delta l_{L-1}\max\{\size(\f),1\}^{L-1}\prb{\textstyle\prod_{k=0}^{L-2}(l_k+1)}R
			\\&\quad+ \pr{\size(\f)+\delta}l_{L-1} \max\{\size(\f)+\delta,1\}^{L-1}\prb{\textstyle\prod_{k=0}^{L-2}(2l_k+1)}R\delta
			+\delta
			\\&\leq \max\{\size(\f)+\delta,1\}^{L}\prb{\textstyle\prod_{k=0}^{L-1}(2l_k+1)}R\delta
			\\&\leq \max\{\size(\f)+\delta,1\}^{\lengthANN(\f)}
			\prb{\textstyle\prod_{k=0}^{\hidlengthANN(\f)}(2l_k+1)}R\delta
			.
		\end{split}
	\end{equation}
		\Moreover the inequality of arithmetic and geometric means 
	\proves\ that  
	\begin{equation}
		\sprod_{k=0}^{\hidlengthANN(\f)}(2l_k+1)
		\leq \PR*{\frac{\ssum_{k=0}^{\hidlengthANN(\f)}(2l_k+1)}{\lengthANN(\f)}}^{\lengthANN(\f)}
		\leq
		\PRbbb{\ssum_{k=0}^{\hidlengthANN(\f)}(2l_k+1)}^{\lengthANN(\f)}
		\leq
		(2\param(\f))^{\lengthANN(\f)}.
	\end{equation}
	Combining this and \lref{eq:1step:bd} with the fact that $\delta=B^{1-\cnt}$ and $\cnt=\ceil{\lengthANN(\f)\log_{\B}(4\paramANN(\f))+\log_{\B}(R\eps^{-1})}+1$ \proves\ that
		\begin{equation}
		\begin{split}
			\vass{(\realisation(\f))(x_0)-(\realisation(\g))(x_0)}
			&\leq \max\{\size(\f)+\delta,1\}^{\lengthANN(\f)}
			\prb{\textstyle\prod_{k=0}^{\hidlengthANN(\f)}(2l_k+1)}R\delta
			\\&\leq 2^{\lengthANN(\f)}
			\prb{\textstyle\prod_{k=0}^{\hidlengthANN(\f)}(2l_k+1)}R\delta
			\\&= 2^{\lengthANN(\f)}\B^{1-\cnt}
			\prb{\textstyle\prod_{k=0}^{\hidlengthANN(\f)}(2l_k+1)}R
			\\&= 2^{\lengthANN(\f)}\B^{-\ceil{\lengthANN(\f)\log_{\B}(4\paramANN(\f))+\log_{\B}(R\eps^{-1})}}
			\prb{\textstyle\prod_{k=0}^{\hidlengthANN(\f)}(2l_k+1)}R
			\\&\leq 2^{\lengthANN(\f)}\B^{-\lengthANN(\f)\log_{\B}(4\paramANN(\f))-\log_{\B}(R\eps^{-1})}
			\prb{\textstyle\prod_{k=0}^{\hidlengthANN(\f)}(2l_k+1)}R
			\\&= 2^{\lengthANN(\f)}(4\paramANN(\f))^{-\lengthANN(\f)}R^{-1}\eps
			\prb{\textstyle\prod_{k=0}^{\hidlengthANN(\f)}(2l_k+1)}R
			\\&= (2\paramANN(\f))^{-\lengthANN(\f)}
			\prb{\textstyle\prod_{k=0}^{\hidlengthANN(\f)}(2l_k+1)}\eps
			\\&\leq (2\paramANN(\f))^{-\lengthANN(\f)}
			(2\paramANN(\f))^{\lengthANN(\f)}\eps=\eps
			.
		\end{split}
	\end{equation}
	\Nobs that \eqref{eq:setup:disc:netw} hence \proves[ep]\ \cref{it:label:dims,it:label:approx,it:label:discretized,it:label:bits}.
	\finishproofthus
\end{proof}

\cfclear
\begin{athm}{prop}{lem:bit:upper:bound}[Approximations through ANNs with discretized parameters]
	Let $a,H\in\R$, $b\in[a,\infty)$, $k,l,m,n\in\N_0$, $d\in\N$, $c\in[2,\infty)$, $\B\in\N\setminus\{1\}$, $\eps\in(0,1]$ satisfy $H=\max\{1,\log_B(\max\{1,\vass{a},\vass{b}\})\}$ and
	let $f\in C(\R^d,\R)$ satisfy for all $\epsilon\in(0,1)$ that
	\begin{equation}
		\label{eq:assumption:cost:upper:bound}
		\COST_{a,b}(f,cd^k\epsilon^{-l},1,\epsilon) \leq cd^{m}\epsilon^{-n}
	\end{equation}\cfload.
	Then there exist $\cnt\in\N$, $\f\in\ANNSET$
	such that
	\begin{enumerate}[(i)]
		\item \label{it:perturbated:statement:eq.2}it holds that $\realisation(\f)\in C(\R^d,\R)$,
		
		\item \label{it:perturbated:statement:eq.4}it holds that $\sup_{x\in[a,b]^d}\vass{(\realisation(\f))(x)-f(x)}\leq\eps$,
				
		\item \label{it:perturbated:statement:eq.3aa}it holds that 
		$\paramANN(\f)\leq 2^n c d^{m}\eps^{-n}$,
		
		\item \label{it:perturbated:statement:eq.3a}it holds that 
		$\cnt\leq (c+m+n+H)c^{l+2} d^{k+1}\eps^{-(l+1)}$,

		\item \label{it:perturbated:statement:eq.3b}it holds that 
		$\paramANN(\f)\cnt\leq (c+m+n+H) 2^n c^{l+3}d^{m+k+1}\eps^{-(n+l+1)}$, and
		
		\item \label{it:perturbated:statement:eq.5}it holds for all $p\in\{1,2,\ldots,\param(\f)\}$ that
		there exist $\beta_1,\beta_2,\ldots,\beta_{\cnt}\in\N_0\cap[0,\B)$ such that
		\begin{equation}
			\form{p}{\f}
			=	(-1)^{\beta_{\cnt}}\textstyle\sum_{i=1}^{\cnt-1} \beta_i \B^{-i}
		\end{equation}
	\end{enumerate}
	\cfout.
\end{athm}

\begin{proof}[Proof of \cref{lem:bit:upper:bound}]
	Throughout this proof 
	let $R\in\R$ satisfy $R=\max\{1,\vass{a},\vass{b}\}$.
	\Nobs that the fact that $\COST_{a,b}(f,cd^k\pr{\tfrac{\eps}{2}}^{-l},1,\tfrac{\eps}{2}) \leq cd^{m}\pr{\tfrac{\eps}{2}}^{-n}$ \proves[pis]\ that there exists $\g\in\ANNSET$ which satisfies that
	\begin{enumerate}[(I)]
		\item \label{it:cont:net:dims}it holds that $\realisation(\g)\in C(\R^d,\R)$,
		\item \label{it:cont:net:approx}it holds that $\sup_{x\in[a,b]^d}\vass{(\realisation(\g))(x)-f(x)}\leq\tfrac{\eps}{2}$,
		\item \label{it:cont:net:length}it holds that $\lengthANN(\g)\leq 2^l c d^k \eps^{-l}$,
		\item \label{it:cont:net:param}it holds that $\paramANN(\g)\leq 2^n c d^m \eps^{-n}$, and
		\item \label{it:cont:net:size}it holds that $\size(\g)\leq 1$
	\end{enumerate}
	\cfload. \Moreover \cref{lem:bit:approx:error} (applied with
	$\f \curvearrowleft \g$,
	$\eps \curvearrowleft \tfrac{\eps}{2}$,
	$R \curvearrowleft R$,
	$\B \curvearrowleft \B$
	in the notation of \cref{lem:bit:approx:error}) \proves[pis]\ that there exist $\cnt\in\N$, $\f\in\ANNSET$ which satisfy that
	\begin{enumerate}[(A)]
		\item \label{it:disc:net:param:dims}it holds that $\dims(\f)=\dims(\g)$,
		\item \label{it:disc:net:param:approx}it holds that $\sup_{x\in[-R,R]^d}\vass{(\realisation(\f))(x)-(\realisation(\g))(x)}\leq\tfrac{\eps}{2}$,
		\item \label{it:disc:net:nbbits}it holds that $\cnt\leq \lengthANN(\g)\log_{\B}(4\paramANN(\g))+\log_{\B}(2R\eps^{-1})+2$, and
		\item \label{it:disc:net:exbits}it holds for all $p\in\{1,2,\ldots,\param(\f)\}$ that there exist $\beta_1,\beta_2,\ldots,\beta_{\cnt}\in\N_0\cap[0,\B)$ such that
		\begin{equation}
			\form{p}{\f}
			=	(-1)^{\beta_{\cnt}}\textstyle\sum_{i=1}^{\cnt-1} \beta_i \B^{-i}
		\end{equation}
	\end{enumerate}
	\cfload.
	Combining \cref{it:cont:net:param,it:cont:net:length,it:disc:net:nbbits} 
	and the fact that for all $x\in[2,\infty)$ it holds that $\log_{\B}( x)+1\leq x$ and $c\geq2$
	\proves\ that
	\begin{equation}
		\label{eq:bits:per:param:est}
		\begin{split}
			&\cnt
			\\&\leq \lengthANN(\g)\log_{\B}(4\param(\g))+\log_{\B}(2R\eps^{-1})+2
			\\&\leq 2^l cd^k\eps^{-l}\log_{\B}(2^{n+2} c d^m \eps^{-n})+\log_{\B}(R)+\log_{\B}(2\eps^{-1})+2
			\\&\leq ((n+2)\log_{\B}(2)+\log_{\B}(c)+m\log_{\B}(d)-n\log_{\B}(\eps))2^l cd^k\eps^{-l}+H+\log_{\B}(2\eps^{-1})+2
			\\&\leq (n+1+c+md+n\eps^{-1})2^l cd^k\eps^{-l}+H+2\eps^{-1}+1
			\\&\leq (1+c+m+2n)2^l cd^{k+1}\eps^{-(l+1)}+H+2\eps^{-1}+1
			\\&\leq (2+c+m+2n+H)2^l cd^{k+1}\eps^{-(l+1)}
			\\&\leq (c+m+n+H)c^{l+2}d^{k+1}\eps^{-(l+1)}
			.
		\end{split}
	\end{equation}
	\Moreover \cref{it:cont:net:approx,it:disc:net:param:approx} and the fact that $[a,b]^d\subseteq [-\dom,\dom]^d$ \prove\ that for all $x\in[a,b]^d$ it holds that
	\begin{equation}
		\vass{(\realisation(\f))(x)-f(x)}
		\leq
		\vass{(\realisation(\f))(x)-(\realisation(\g))(x)}+\vass{(\realisation(\g))(x)-f(x)}
		\leq \tfrac{\eps}{2}+\tfrac{\eps}{2}=\eps.
	\end{equation}
	This,
	\eqref{eq:bits:per:param:est},
	\cref{it:cont:net:param},
	\cref{it:disc:net:param:dims},
	and \cref{it:disc:net:exbits} \prove[ipe]\ \cref{it:perturbated:statement:eq.2,it:perturbated:statement:eq.3aa,it:perturbated:statement:eq.3a,it:perturbated:statement:eq.3b,it:perturbated:statement:eq.4,it:perturbated:statement:eq.5}.
	\finishproofthus
\end{proof}

\cfclear
\begin{athm}{cor}{lem:bit:upper:2:bound}[Bitwise approximation through ANNs]
	Let $a\in\R$, $b\in[a,\infty)$, $k,l,m,n\in\N_0$, $c\in[2,\infty)$, $d\in\N$, $f\in C(\R^d,\R)$ satisfy for all $\eps\in(0,1)$ that
	\begin{equation}
		\label{eq:assumption:cost:upper:bound:2}
		\COST_{a,b}(f,cd^k\eps^{-l},1,\eps) \leq cd^{m}\eps^{-n}
	\end{equation} and let $\eps\in(0,1]$ \cfload.
	Then there exist $\cnt\in\N$, $\f\in\ANNSET$
	such that
	\begin{enumerate}[(i)]
		\item \label{perturbated:statement:2:eq.2}it holds that $\realisation(\f)\in C(\R^d,\R)$,
		
		\item \label{perturbated:statement:2:eq.4}it holds that $\sup_{x\in[a,b]^d}\vass{(\realisation(\f))(x)-f(x)}\leq\eps$,
		
		\item \label{perturbated:statement:2:eq.3}it holds that 
		$\paramANN(\f)\cnt\leq (c+m+n+\max\{\vass{a},\vass{b},1\}) c^{n+l+3}d^{m+k+1}\eps^{-(n+l+1)}$, and
		
		\item \label{perturbated:statement:2:eq.5}it holds for all $p\in\{1,2,\ldots,\param(\f)\}$ that
		there exist $\beta_1,\beta_2,\ldots,\beta_{\cnt}\in\{0,1\}$ such that
		\begin{equation}
			\form{p}{\f}
			=	(-1)^{\beta_{\cnt}}\textstyle\sum_{i=1}^{\cnt-1} \beta_i 2^{-i}
		\end{equation}
	\end{enumerate}
	\cfout.
\end{athm}

\begin{proof}[Proof of \cref{lem:bit:upper:2:bound}]
	\Nobs that \cref{lem:bit:upper:bound}  (applied with
			$B \curvearrowleft 2$, $H\curvearrowleft\max\{1,\log_2(\allowbreak\max\{1,\vass{a},\vass{b}\})\}$ in the notation of \cref{lem:bit:upper:bound}) \proves[ipe]\ \cref{perturbated:statement:2:eq.2,perturbated:statement:2:eq.3,perturbated:statement:2:eq.4,perturbated:statement:2:eq.5}.
	\finishproofthus
\end{proof}

\subsection{Lower and upper bounds for the number of bits to store ANN approximations}

\begin{definition}
	[Entropy for ANN approximations]
	\label{def:ent}
	\cfconsiderloaded{def:ent}
	Let $a\in\R$, $b\in[a,\infty)$, $\B\in(1,\infty)$.
	Then we denote by  $\ENT_{\B}^{a,b}\colon\pr{\cup_{d\in\N}C(\R^d,\R)}\times[0,\infty]^2\to[0,\infty]$ 
	and
	$\ent_{\B}^{a,b}\colon\pr{\cup_{d\in\N}C(\R^d,\R)}\times[0,\infty]^2\to[0,\infty]$
	the functions which satisfy for all $d\in\N$, $f\in C(\R^d,\R)$,  $L,\eps\in[0,\infty]$ that
	\begin{equation}
		\label{Ent:def}
		\begin{split}
			\ENT_{\B}^{a,b}(f,L,\eps)=
			&\inf \pr*{ \pR*{ p\in\N \colon \PR*{\!\!
			\begin{array}{c}
				\exists\, \mathscr{f} \in \ANNSET \colon
				(\lengthANN(\mathscr{f})\leq L)
				\land{}
				(\realisation(\mathscr{f}) \in C(\R^d,\R))
				\land\\
				((\log_{\B}(\max\{1,\size(\f)\})+1)\param(\f)\leq p)
				\land{}\\
				(\sup_{x\in[a,b]^d}\vass{(\realisation(\mathscr{f}))(x)-f(x)} \leq \varepsilon)\\
			\end{array} \!\! }
	}
	\cup\{\infty\}
}
		\end{split}
	\end{equation}
	and
	\begin{equation}
		\label{ent:def}
		\begin{split}
			\ent_{\B}^{a,b}(f,L,\eps)=
			&\inf \pr*{ \pR*{ p\in\N \colon \PR*{\!\!
						\begin{array}{c}
							\exists\, \mathscr{f} \in \ANNSET \colon
							(\lengthANN(\mathscr{f})\leq L)
							\land{}
							(\realisation(\mathscr{f}) \in C(\R^d,\R))
							\land{}\\
							(\textstyle\sum_{i=1}^{\param(\f)}\PR{\log_{\B}\pr{\max\{1,\vass{\form{i}{\f}}\}}+\Indfct{\R\setminus\{0\}}(\form{i}{\f})}\leq p)\land\\
							(\sup_{x\in[a,b]^d}\vass{(\realisation(\mathscr{f}))(x)-f(x)} \leq \varepsilon)\\
						\end{array} \!\! }
				}
				\cup\{\infty\}
			}.
		\end{split}
	\end{equation}
	\cfload.
\end{definition}

\cfclear
\begin{athm}{cor}{cost:inequalities:2}
	Let $a\in\R$, $b\in[a,\infty)$, $\B\in(1,\infty)$, $d\in\N$, $f\in C(\R^d,\R)$, $L,\eps\in[0,\infty]$ \cfload. Then 
	\begin{equation}
		\label{eq:cost:inequalities:2}
		\cost_{a,b}(f,L,\infty,\eps)
		\leq
		\ent_{\B}^{a,b}(f,L,\eps)
		\leq
		\ENT_{\B}^{a,b}(f,L,\eps)
		\leq
		\COST_{a,b}(f,L,1,\eps)
	\end{equation}
	\cfout.
\end{athm}

\begin{proof}[Proof of \cref{cost:inequalities:2}]
	\Nobs that \eqref{Cost:def}, \eqref{cost:def}, \eqref{Ent:def}, and \eqref{ent:def} \prove\ that for all $\f\in\ANNSET$ it holds that
	\begin{equation}
		\begin{split}
		\paramnotzero(\f)
		=\sum_{i=1}^{\param(\f)}\Indfct{\R\setminus\{0\}}(\form{i}{\f})
		&\leq\sum_{i=1}^{\param(\f)}\PR{\log_{\B}\pr{\max\{1,\vass{\form{i}{\f}}\}}+\Indfct{\R\setminus\{0\}}(\form{i}{\f})}
		\\&\leq\sum_{i=1}^{\param(\f)}\PR{\Indfct{\R\setminus\{0\}}(\size(\f))}\PR{\log_{\B}(\max\{1,\size(\f)\})+1}
		\\&\leq \param(\f)\PR{\log_{\B}(\max\{1,\size(\f)\})+1}
	\end{split}
	\end{equation}
\cfload. Combining this with the fact that for all $\f\in\ANNSET$ with $\size(\f)\leq 1$ it holds that $\param(\f)\PR{1+\log_{\B}\pr{\max\{1,\size(\f)\}}}=\param(\f)$ \proves\ \eqref{eq:cost:inequalities:2}.
	\finishproofthus
\end{proof}

\cfclear
\begin{athm}{cor}{Cor:Thm1}
	Let $\varphi,\dom\in\R$, $\gamma\in(0,1]$, $\beta\in[1,\infty)$, $a\in\R$, $b\in(a,\infty)$, $\mathfrak{c},\kappa\in(0,\infty)$, $\B\in(1,\infty)$ satisfy
	$\mathfrak{c}\geq \pot \ln\pr*{\max\{2,\abs{a},\abs{b},\scl\}} \max\{1,\kappa^3\}$ 
	and
	$\dom=\max\{\vass{a},\vass{b}\}>\scl^{-1}$,
	and
	for every $d\in\N$ let
	$f_d\colon\R^d\to\R$ satisfy for all $x=(x_1,\ldots,x_d)\in\R^d$ that $f_d(x)=\kappa\sin\prb{\gamma\beta^d\prb{\sprod_{i = 1}^d x_i}+\varphi}$. Then
	\begin{enumerate}[(i)]
		\item{
			\label{it:Cor:Thm1:eq1}
			it holds for all $d,L\in\N$, $\eps\in(0,\kappa)$ that 
			\begin{equation}
				\ent_{\B}^{a,b}(f_d,L,\eps)
				\geq 
				\cost_{a,b}(f_d,L,\infty,\eps)
				\geq	\vass[\Big]{\frac{\floor{(b-a)\gamma(\dom\scl)^{d-1}\pi^{-1}-1}}{2}}^{\nicefrac{1}{\max\{1,L-1\}}}
			\end{equation}
			and
		}
		\item{
			\label{it:Cor:Thm1:eq2}
			it holds for all $d\in\N$, $\eps\in(0,\kappa)$ that
			\begin{equation}
				\begin{split}
					\ENT_{\B}^{a,b}(f_d,\mathfrak{c}d^2\eps^{-1},\eps)
					\leq  
					\COST_{a,b}(f_d,\mathfrak{c}d^2\eps^{-1},1,\eps)
					\leq
					\mathfrak{c}d^3\eps^{-2}
				\end{split}
			\end{equation}
		}
	\end{enumerate}
	\cfout.
\end{athm}

\begin{proof}[Proof of \cref{Cor:Thm1}]
	\Nobs that \cref{Thm1}, \cref{cost:inequalities:2}, and the fact that for all $x\in[2,\infty)$ it holds that $19770\ceil{\ln(x)}\leq \pot \ln(x)$ \prove[pe]\ \cref{it:Cor:Thm1:eq1,it:Cor:Thm1:eq2}.
	\finishproofthus
\end{proof}

\cfclear

\begin{athm}{cor}{Cor:Thm6}
	Let $\varphi\in\R$, $\gamma,\kappa\in(0,\infty)$, $a\in\R$,  $b\in[a+\pi\gamma^{-1},\infty)$, $\mathfrak{c}\in(0,\infty)$, $\B\in(1,\infty)$ satisfy $\mathfrak{c}\geq4504\ceil{\log_2\pr*{\max\{1,\gamma\}\max\{\abs{a},\abs{b},2\}}}\max\{\kappa^3,1\}$ and
	for every $d\in\N$ let
	$f_d\colon\R^d\to\R$ satisfy for all $x=(x_1,\ldots,x_d)\in\R^d$ that $f_d(x)=\kappa\sin\prb{\gamma 2^d\pr{\ssum_{i = 1}^d x_i}+\varphi}$.
	Then 
	\begin{enumerate}[(i)]
	\item{
		\label{it:Cor:Thm6:eq1}
		it holds for all $d,L\in\N$, $\eps\in(0,\kappa)$ that 
		\begin{equation}
			\ent_{\B}^{a,b}(f_d,L,\eps)
			\geq 
			\cost_{a,b}(f_d,L,\infty,\eps)
			\geq	
			\pr{d2^{d-1}}^{\nicefrac{1}{\max\{1,L-1\}}}
		\end{equation}
		and
	}

	\item{
		\label{it:Cor:Thm6:eq2}
		it holds for all $d\in\N$, $\eps\in(0,\kappa)$ that
		\begin{equation}
			\begin{split}
				\ENT_{\B}^{a,b}(f_d,\mathfrak{c}d,\eps)
				\leq  
				\COST_{a,b}(f_d,\mathfrak{c}d,1,\eps)
				\leq
				\mathfrak{c}d\eps^{-2}
			\end{split}
		\end{equation}
	}
\end{enumerate}
	\cfout.
\end{athm}

\begin{proof}[Proof of \cref{Cor:Thm6}]
	\Nobs that \cref{Thm6} and \cref{cost:inequalities:2} \prove[pe]\ \cref{it:Cor:Thm6:eq1,it:Cor:Thm6:eq2}.
	\finishproofthus
\end{proof}

\cfclear
\begin{athm}{cor}{cor:bit:lower:bound:v1}
	Let $\varphi,a\in\R$, $b\in[a+4,\infty)$, $d,L\in\N$, $\kappa\in(0,\infty)$, $\eps\in(0,\kappa)$
	and let
	$f\colon\R^d\to\R$ satisfy for all $x=(x_1,\ldots,x_d)\in\R^d$ that $f(x)=\kappa\sin\prb{\prb{\sprod_{i = 1}^d x_i}+\varphi}$. 
	Then 
\begin{equation}
 \label{it:cost:discretized:lb}  				
 		\ent_{2}^{a,b}(f,L,\eps)
		\geq 
		\cost_{a,b}(f,L,\infty,\eps)
		\geq
		2^{\frac{d-2}{L}}
	\end{equation}
	\cfout.
\end{athm}

\begin{proof}[Proof of \cref{cor:bit:lower:bound:v1}]
	\Nobs that \cref{it:Cor:Thm1:eq1} in \cref{Cor:Thm1} (applied with
	$\varphi \curvearrowleft \varphi$,
	$\dom \curvearrowleft \max\{\vass{a},\vass{b}\}$,
	$\gamma \curvearrowleft 1$,
	$\beta \curvearrowleft 1$,
	$a \curvearrowleft a$,
	$b \curvearrowleft b$,
	$\kappa \curvearrowleft \kappa$,
	$\B \curvearrowleft 2$
	in the notation of \cref{Cor:Thm1}) \proves\ that
	\begin{equation}
		\label{eq:cor:bit:lower:bound:v1:left}
		\begin{split}
			\ent_{2}^{a,b}(f,L,\eps)
			\geq 
			\cost_{a,b}(f,L,\infty,\eps)
			&\geq
			\vass[\Big]{\frac{\floor{(b-a)(\max\{\vass{a},\vass{b}\})^{d-1}\pi^{-1}-1}}{2}}^{\nicefrac{1}{\max\{1,L-1\}}}
			\\&\geq
			\vass[\Big]{\frac{\floor{4(2)^{d-1}\pi^{-1}-1}}{2}}^{\nicefrac{1}{\max\{1,L-1\}}}
			\\&\geq
			\vass[\Big]{\frac{\floor{\tfrac{5}{4}(2)^{d-1}-1}}{2}}^{\nicefrac{1}{\max\{1,L-1\}}}
		\end{split}
	\end{equation}
	\cfload.
\Moreover for all $D\in\N\cap[3,\infty)$ it holds that
\begin{equation}
	\frac{\floor{\tfrac{5}{4}(2)^{D-1}-1}}{2}
	=
	\frac{\tfrac{5}{4}(2)^{D-1}-1}{2}
	=
	\tfrac{5}{4}(2)^{D-2}-\tfrac{1}{2}
	\geq
	2^{D-2}.
\end{equation}
This,
\eqref{eq:cor:bit:lower:bound:v1:left}, and the fact that for all $D\in\{1,2\}$ it holds that $\cost_{a,b}(f,L,\infty,\eps)\geq 1 \geq 2^{\frac{D-2}{L}}$ \prove\ that
		\begin{equation}
		\begin{split}
			\ent_{2}^{a,b}(f,L,\eps)
			\geq 
			\cost_{a,b}(f,L,\infty,\eps)
			\geq
			\prb{\max\pRb{1,2^{d-2}}}^{\nicefrac{1}{\max\{1,L-1\}}}
			&\geq	
			2^{\frac{d-2}{L}}.		
		\end{split}
		\end{equation}
	\finishproofthus
\end{proof}

\cfclear
\begin{athm}{cor}{cor:bit:lower:bound:v2}
	Let $\varphi,a\in\R$, $b\in[a+4,\infty)$, $d\in\N$, $\kappa\in(0,\infty)$, $\eps\in(0,\kappa)$, $\f\in\ANNSET$
	satisfy  
	\begin{equation}
		\label{eq:cor:bit:lower:bound:v2}
		\realisation(\f)\in C(\R^d,\R)
		\qandq
		\textstyle\sup_{x=(x_1,\ldots,x_d)\in[a,b]^d}
		\vass[\big]{(\realisation(\f))(x)-\kappa\sin\prb{\prb{\textstyle\prod_{i = 1}^d x_i}+\varphi}}\leq\eps
	\end{equation}
\cfload.
	Then 
	\begin{equation}
		\textstyle\sum_{i=1}^{\param(\f)}\PR{\log_{2}\pr{\max\{1,\vass{\form{i}{\f}}\}}+\Indfct{\R\setminus\{0\}}(\form{i}{\f})}
		\geq
		\paramnotzero(\f)
		\geq
		2^{\frac{d-2}{\lengthANN(\f)}}
	\end{equation}
	\cfout.
\end{athm}

\begin{proof}[Proof of \cref{cor:bit:lower:bound:v2}]
	Throughout this proof let $f\colon\R^d\to\R$ satisfy for all $x=(x_1,\dots,x_d)\in\R^d$ that $f(x)=\kappa\sin\prb{\prb{\textstyle\prod_{i = 1}^d x_i}+\varphi}$.
	\Nobs that \eqref{eq:cor:bit:lower:bound:v2} and the fact that $\sup_{x\in[a,b]^d}\vass{f(x)}=\kappa>\eps$ \prove\ that $\paramnotzero(\f)\geq 1$.
	\Nobs that \eqref{eq:cor:bit:lower:bound:v2} and \cref{cor:bit:lower:bound:v1} (applied with
	$\varphi \curvearrowleft \varphi$,
	$a \curvearrowleft a$,
	$b \curvearrowleft b$,
	$L \curvearrowleft \lengthANN(\f)$,
	$\kappa \curvearrowleft \kappa$,
	$\eps \curvearrowleft \eps$,
	$f\curvearrowleft f$
	in the notation of \cref{cor:bit:lower:bound:v1}) hence
	\prove\ that 
	\begin{equation}
		\begin{split}
		\textstyle\sum_{i=1}^{\param(\f)}\PR{\log_{2}\pr{\max\{1,\vass{\form{i}{\f}}\}}+\Indfct{\R\setminus\{0\}}(\form{i}{\f})}
		&\geq
		\textstyle\sum_{i=1}^{\param(\f)}\Indfct{\R\setminus\{0\}}(\form{i}{\f})
		\\&=
		\max\{1,\paramnotzero(\f)\}
		\\&\geq
		\cost_{a,b}(f,\lengthANN(\f),\infty,\eps)
		\geq
		2^{\frac{d-2}{\lengthANN(\f)}}	
		\end{split}
	\end{equation}	
	\cfload.
	\finishproofthus
\end{proof}

\cfclear
\begin{athm}{cor}{cor:bit:upper:bound}[Bitwise approximation through ANNs]
	Let $\varphi,a\in\R$, $b\in[a,\infty)$, $\kappa\in(0,\infty)$, $\eps\in(0,\kappa)$, $c\in [\pot\ln(\max\{2,\vass{a},\vass{b}\}) \max\{1,\kappa^3\},\infty)$, $d\in\N$ and let
$f\colon\R^d\to\R$ satisfy for all $x=(x_1,\ldots,x_d)\in\R^d$ that $f(x)=\kappa\sin\prb{\prb{\sprod_{i = 1}^d x_i}+\varphi}$. 
Then there exist $\cnt\in\N$, $\f\in\ANNSET$ such that 
	\begin{enumerate}[(i)]
		
		\item \label{it:perturbated:statement:final:eq.2}it holds that $\realisation(\f)\in C(\R^d,\R)$,
		
		\item \label{it:perturbated:statement:final:eq.4}it holds that $\sup_{x\in[a,b]^d}\vass{(\realisation(\f))(x)-f(x)}\leq\eps$,
		
		\item \label{it:perturbated:statement:final:eq.3}it holds that 
		$\paramANN(\f)\cnt\leq c^8d^6\eps^{-4}$, and
		
		\item \label{it:perturbated:statement:final:eq.5}it holds for all  $p\in\{1,2,\ldots,\param(\f)\}$ that
		there exist $\beta_1,\beta_2,\ldots,\beta_{\cnt}\in\{0,1\}$ such that
		\begin{equation}
			\form{p}{\f}
			=	(-1)^{\beta_{\cnt}}\textstyle\sum_{i=1}^{\cnt-1} \beta_i 2^{-i}
		\end{equation}
	\end{enumerate}
	\cfout.
\end{athm}

\begin{proof}[Proof of \cref{cor:bit:upper:bound}]
Throughout this proof 
let $\dom,H\in\R$ satisfy $\dom=\max\{1,\abs{a},\abs{b}\}$ and $H=\max\{1,\log_2(\dom)\}$ and
let $g\colon\R\to\R$ satisfy for all $x\in\R$ that $\tgt(x)=\kappa\sin(x+\varphi)$. 
\Nobs that
\begin{equation}
	\llabel{eq:setup:ineq}
	\max\pR{5,H,4\max\{1,\kappa^2\}}\leq c
	\qandq
	\max\{1,\eps^{-2}\}\leq \eps^{-2}\max\{1,\kappa^2\}\leq \frac{c}{4\eps^2}.
\end{equation}
\Moreover the fact that $c\geq9885\ceil{\log_2\pr*{\max\{2,\abs{a},\abs{b},1\}}}\max\{1,\kappa^3\}$ and \cref{final:approximation:imp} (applied for every $\epsilon\in(0,\kappa)$ with
$a \curvearrowleft a$,
$b \curvearrowleft b$,
$d \curvearrowleft d$,
$\kappa \curvearrowleft \kappa$,
$\mathfrak{c} \curvearrowleft c$,
$\eps \curvearrowleft \epsilon$,
$\gamma \curvearrowleft 1$,
$\beta \curvearrowleft 1$,
$g \curvearrowleft g$,
$f \curvearrowleft f$
in the notation of \cref{final:approximation:imp}) \prove[sip]\ that for all $\epsilon\in(0,\kappa)$ it holds that
\begin{equation}
	\llabel{eq:constr:UB}
	\COST_{a,b}(f,cd^2\epsilon^{-1},1,\epsilon)
	\leq
	cd^3\epsilon^{-2}
\end{equation}
\cfload.
\Moreover the fact that there exists $\f\in\ANNSET$ with $\dims(\f)=(d,1)$, $\param(\f)=d+1$, and $\paramnotzero(\f)=0$ and the fact that $\sup_{x\in[a,b]^d}\vass{f(x)}\leq\kappa$ \prove\ that for all $\epsilon\in[\kappa,\infty)\cap(0,1]$ it holds that
\begin{equation}
	\llabel{eq:constr:UB:triv}
	\COST_{a,b}(f,cd^2\epsilon^{-1},1,\epsilon)
	\leq
	d+1
	\leq
	cd^3\epsilon^{-2}
	.
\end{equation}
This,
\lref{eq:setup:ineq},
\lref{eq:constr:UB}, and \cref{lem:bit:upper:bound} (applied with
$k \curvearrowleft 2$,
$l \curvearrowleft 1$,
$m \curvearrowleft 3$,
$n \curvearrowleft 2$,
$\eps \curvearrowleft \min\{1,\eps\}$,
$c\curvearrowleft c$,
$B\curvearrowleft 2$,
$H\curvearrowleft H$
in the notation of \cref{lem:bit:upper:bound}) \prove\ that
there exist $\cnt\in\N$, $\f\in\ANNSET$
which satisfy that
\begin{enumerate}[(I)]
	\item \label{it:ref:perturbated:statement:eq.2}it holds that $\realisation(\f)\in C(\R^d,\R)$,
	
	\item \label{it:ref:perturbated:statement:eq.4}it holds that $\sup_{x\in[a,b]^d}\vass{(\realisation(\f))(x)-f(x)}\leq\eps$,
	
	\item \label{it:ref:perturbated:statement:eq.3aa}it holds that 
	$\paramANN(\f)\leq 4cd^{3}\max\{1,\eps^{-2}\}\leq c^2d^{3}\eps^{-2}$,
	
	\item \label{it:ref:perturbated:statement:eq.3a}it holds that 
	$\cnt\leq (c+3+2+H)c^{1+2} d^{2+1}\max\{1,\eps^{-(1+1)}\}\leq c^6d^3\eps^{-2}$,

	\item \label{it:ref:perturbated:statement:eq.3b}it holds that 
	$\paramANN(\f)\cnt\leq (c+3+2+H)4c^{1+3}d^{2+3+1}\max\{1,\eps^{-(1+2+1)}\}
	\leq  c^8d^6\eps^{-4}$, and
	
	\item \label{it:ref:perturbated:statement:eq.5}it holds for all $p\in\{1,2,\ldots,\param(\f)\}$ that
	there exist $\beta_1,\beta_2,\ldots,\beta_{\cnt}\in\{0,1\}$ such that
	\begin{equation}
		\form{p}{\f}
		=	(-1)^{\beta_{\cnt}}\textstyle\sum_{i=1}^{\cnt-1} \beta_i 2^{-i}
	\end{equation}
\end{enumerate}
\cfload. \Nobs that \cref{it:ref:perturbated:statement:eq.2,it:ref:perturbated:statement:eq.3b,it:ref:perturbated:statement:eq.4,it:ref:perturbated:statement:eq.5} \prove\ \cref{it:perturbated:statement:final:eq.2,it:perturbated:statement:final:eq.3,it:perturbated:statement:final:eq.4,it:perturbated:statement:final:eq.5}.
\finishproofthus
\end{proof}

\cfclear
\begin{athm}{cor}{cor:bit:both:bounds}[Bitwise approximation through ANNs]
	Let $\varphi,a\in\R$, $b\in[a+4,\infty)$, $\eps\in(0,1)$, $c\in [\pot\ln(\max\{\vass{a},\vass{b}\}),\infty)$, $d\in\N$
	and let
	$f\colon\R^d\to\R$ satisfy for all $x=(x_1,\ldots,x_d)\in\R^d$ that $f(x)=\sin\prb{\prb{\sprod_{i = 1}^d x_i}+\varphi}$. 
	Then 
	\begin{enumerate}[(i)]
		\item \label{cor:bit:both:bounds:lower}it holds for all $L\in\N$, $\f\in\ANNSET$ with $\realisation(\f)\in C(\R^d,\R)$,
$
\sup_{x\in[a,b]^d}\vass{(\realisation(\f))(x)-f(x)}\leq\eps
$, $
\lengthANN(\f)\leq L$
that
\begin{equation}
	\textstyle\sum_{i=1}^{\param(\f)}\PR{\log_{2}\pr{\max\{1,\vass{\form{i}{\f}}\}}+\Indfct{\R\setminus\{0\}}(\form{i}{\f})}
	\geq
	\paramnotzero(\f)
	\geq
	2^{\frac{d-2}{L}}
\end{equation}
	and
			\item \label{cor:bit:both:bounds:upper}there exist $\cnt\in\N$, $\f\in\ANNSET$ with $\realisation(\f)\in C(\R^d,\R)$, $\sup_{x\in[a,b]^d}\vass{(\realisation(\f))(x)-f(x)}\leq\eps$,
	$			\pr{\cup_{p=1}^{\paramANN(\f)}\pR{\form{p}{\f}}}
	\subseteq
	\pr{\cup_{ \beta_1,\beta_2,\ldots,\beta_{\cnt}\in\{0,1\} }
		\pR{(-1)^{\beta_{\cnt}}\sum_{i=1}^{\cnt-1} \beta_i 2^{-i}
	}}$
	such that 
	\begin{equation}
		\paramANN(\f)\cnt\leq c^8d^6\eps^{-4}
	\end{equation}
	\end{enumerate}
	\cfout.
\end{athm}

\begin{proof}[Proof of \cref{cor:bit:both:bounds}]
 \Nobs that \cref{cor:bit:lower:bound:v2} \proves\ that for all $L\in\N$, $\f\in\ANNSET$ with $\realisation(\f)\in C(\R^d,\R)$,
 $
 \sup_{x\in[a,b]^d}\vass{(\realisation(\f))(x)-f(x)}\leq\eps
 $, $
 \lengthANN(\f)\leq L$
 it holds that	
 \begin{equation}
 	\textstyle\sum_{i=1}^{\param(\f)}\PR{\log_{2}\pr{\max\{1,\vass{\form{i}{\f}}\}}+\Indfct{\R\setminus\{0\}}(\form{i}{\f})}
 	\geq
 	\paramnotzero(\f)
 	\geq
 	2^{\frac{\max\{0,d-2\}}{\lengthANN(\f)}}
 	\geq
 	2^{\frac{\max\{0,d-2\}}{L}}.
 \end{equation}
 \cfload. This and the fact that for all $\f\in\ANNSET$ with $\realisation(\f)\in C(\R,\R)$ and $\paramnotzero(\f)=0$ it holds that $\sup_{x\in[a,b]}\vass{(\realisation(\f))(x)-\sin(x+\varphi)}=1>\eps$ \prove\ \cref{cor:bit:both:bounds:lower}.
 	\Nobs that \cref{cor:bit:upper:bound} \proves\ \cref{cor:bit:both:bounds:upper}.
	\finishproofthus
\end{proof}

\begin{athm}{remark}{rmk:1}
	Let us add two comments 
	concerning the conclusion of \cref{cor:bit:both:bounds} in concrete implementations on a computer: 
	\begin{enumerate}[(I)]
		\item 
	\Nobs that 
	\cref{cor:bit:both:bounds:lower} in \cref{cor:bit:both:bounds} shows that for any arbitrarily large fixed \ANN\ depth $L \in \N$ we have that the number of bits needed to approximate the considered globally bounded and infinitely often differentiable functions via \ANNs\ with the depth being bounded by $L$ \emph{must grow at least exponentially} in the dimension $d \in \N$ of the function. In this sense the amount of memory needed to store such shallow or insufficiently deep \ANN\ approximations \emph{must grow at least exponentially} in the problem dimension.
	\item
	\Nobs that \cref{cor:bit:both:bounds:upper} in \cref{cor:bit:both:bounds} reveals that it is possible to approximate the considered globally bounded and infinitely often differentiable functions via sufficiently deep \ANNs\ with the number of bits needed to store the approximating \ANNs\ \emph{growing at most polynomially} in both the dimension $d\in\N$ of the function and the reciprocal $\eps^{-1}$ of the prescribed approximation accuracy $\eps\in(0,1)$.
	In this sense the amount of memory needed to store such deep \ANN\ approximations \emph{grows at most polynomially} in the problem dimension and the reciprocal of the prescribed approximation accuracy.
	\end{enumerate}
\end{athm}

\subsection*{Acknowledgements}
This work has been partially supported by the National Natural Science Foundation of China (NSFC) under grant number W2531010. 
This work has been partially funded by the European Union (ERC, MONTECARLO, 101045811). The views and the opinions expressed in this work are however those of the authors only and do not necessarily reflect those of the European Union or the European Research Council. Neither the European Union nor the granting authority can be held responsible for them.
In addition, the third author gratefully acknowledges the Cluster of Excellence EXC 2044-390685587, Mathematics Münster: Dynamics-Geometry-Structure funded by the Deutsche Forschungsgemeinschaft (DFG, German Research Foundation).
Most of the specific formulations in the proofs of this work have been created using \cite{KuckuckTools}.

\begin{figure}[h]
	\centering
	\includegraphics[width=70mm]{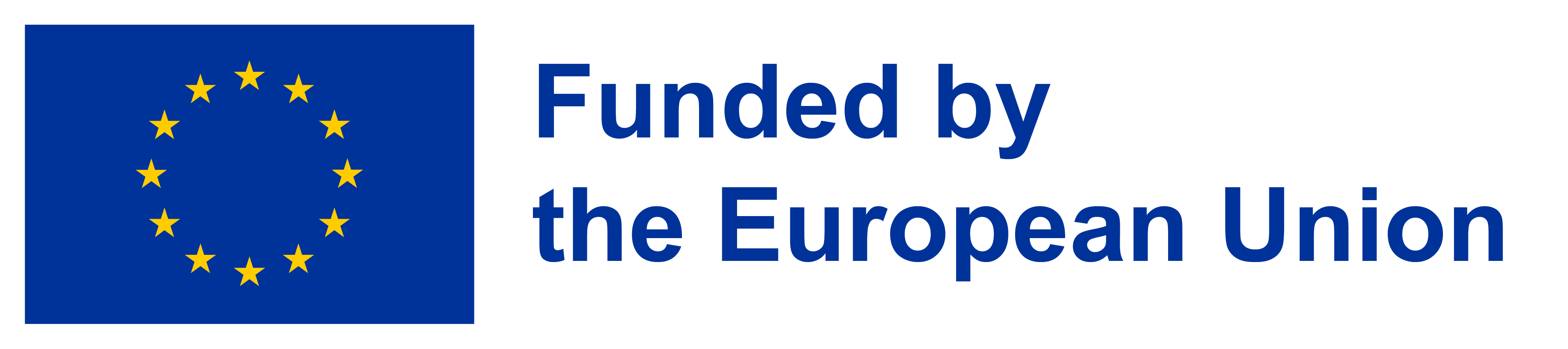}
	\label{fig:meine-grafik}
\end{figure}

\subsection*{Use of large language models}
\label{intro:llms}
The large language models GPT 6.2 and Claude 5.5 have helped to identify several errors and inaccuracies within this article.
Moreover, GPT 6.2 has helped organize and extend the literature review in \cref{intro:literature}.
The ideas underlying the mathematical results and their proofs originate from the authors of this article, who take full responsibility for the content and the statements in this article.

\newpage
\bibliographystyle{acm}
\bibliography{bibfile}

@article {Shapira2023,
	AUTHOR = {Shapira, Itai},
	TITLE = {Expressivity of shallow and deep neural networks
	for polynomial approximation},
	JOURNAL = {\href{https://arxiv.org/abs/2303.03544}{arXiv:2303.03544}},
	YEAR = {2023},
	DOI = {10.48550/arXiv.2303.03544},
	URL = {https://arxiv.org/abs/2303.03544},
}

@inproceedings {BlanchardBennouna2022,
	AUTHOR = {Blanchard, Mo{\"i}se and Bennouna, Amine},
	TITLE = {Shallow and deep networks are near-optimal
	approximators of {Korobov} functions},
	BOOKTITLE = {The Tenth International Conference on Learning
	Representations},
	YEAR = {2022},
	URL = {https://openreview.net/forum?id=AV8FPoMTTa},
}

@article {Safran2026EveryLayerCounts,
	AUTHOR = {Safran, Itay},
	TITLE = {Every layer counts: an exponential {$L_2$}
	depth hierarchy for {ReLU} networks},
	JOURNAL = {\href{https://arxiv.org/abs/2608.23877}{arXiv:2608.23877}},
	YEAR = {2026},
	DOI = {10.48550/arXiv.2608.23877},
	URL = {https://arxiv.org/abs/2608.23877},
}

@article {Siegel2023Optimal,
	AUTHOR = {Siegel, Jonathan W.},
	TITLE = {Optimal approximation rates for deep {ReLU} neural
	networks on {Sobolev} and {Besov} spaces},
	JOURNAL = {J. Mach. Learn. Res.},
	FJOURNAL = {Journal of Machine Learning Research},
	VOLUME = {24},
	YEAR = {2023},
	PAGES = {Paper No. 357, 52},
	ISSN = {1532-4435},
	URL = {https://www.jmlr.org/papers/v24/23-0025.html},
}

@article {GuehringRaslan2021,
	AUTHOR = {G{\"u}hring, Ingo and Raslan, Mones},
	TITLE = {Approximation rates for neural networks with encodable
	weights in smoothness spaces},
	JOURNAL = {Neural Netw.},
	FJOURNAL = {Neural Networks},
	VOLUME = {134},
	YEAR = {2021},
	PAGES = {107--130},
	ISSN = {0893-6080},
	DOI = {10.1016/j.neunet.2020.11.010},
	URL = {https://doi.org/10.1016/j.neunet.2020.11.010},
}

@article {OuSchenkelBolcskei2024,
	AUTHOR = {Ou, Weigutian and Schenkel, Philipp and
	B{\"o}lcskei, Helmut},
	TITLE = {Three quantization regimes for {ReLU} networks},
	JOURNAL = {\href{https://arxiv.org/abs/2405.01952}{arXiv:2405.01952}},
	YEAR = {2024},
	DOI = {10.48550/arXiv.2405.01952},
	URL = {https://arxiv.org/abs/2405.01952},
}

@inproceedings {SafranReichmanValiant2025,
	AUTHOR = {Safran, Itay and Reichman, Daniel and Valiant, Paul},
	TITLE = {Depth separations in neural networks: separating the
	dimension from the accuracy},
	BOOKTITLE = {Proceedings of Thirty Eighth Conference on Learning Theory},
	EDITOR = {Haghtalab, Nika and Moitra, Ankur},
	SERIES = {Proceedings of Machine Learning Research},
	VOLUME = {291},
	YEAR = {2025},
	PAGES = {5108--5142},
	PUBLISHER = {PMLR},
	URL = {https://proceedings.mlr.press/v291/safran25a.html},
}

@article {BortolottoGanacimOliveira2026,
	AUTHOR = {Bortolotto, Cynthia and Ganacim, Francisco and
	Oliveira, Roberto I.},
	TITLE = {Depth separation for neural networks under
	{Gaussian} measure},
	JOURNAL = {\href{https://arxiv.org/abs/2610.28800}{arXiv:2610.28800}},
	YEAR = {2026},
	DOI = {10.48550/arXiv.2610.28800},
	URL = {https://arxiv.org/abs/2610.28800},
}

@inproceedings {Yarotsky2018Optimal,
	AUTHOR = {Yarotsky, Dmitry},
	TITLE = {Optimal approximation of continuous functions by
	very deep {ReLU} networks},
	BOOKTITLE = {Proceedings of the 31st Conference on Learning Theory},
	EDITOR = {Bubeck, S{\'e}bastien and Perchet, Vianney and
	Rigollet, Philippe},
	SERIES = {Proceedings of Machine Learning Research},
	VOLUME = {75},
	YEAR = {2018},
	PAGES = {639--649},
	PUBLISHER = {PMLR},
	URL = {https://proceedings.mlr.press/v75/yarotsky18a.html},
}

@inproceedings {YarotskyZhevnerchuk2020,
	AUTHOR = {Yarotsky, Dmitry and Zhevnerchuk, Anton},
	TITLE = {The phase diagram of approximation rates for deep
	neural networks},
	BOOKTITLE = {Advances in Neural Information Processing Systems},
	EDITOR = {Larochelle, Hugo and Ranzato, Marc'Aurelio and
	Hadsell, Raia and Balcan, Maria-Florina and
	Lin, Hsuan-Tien},
	VOLUME = {33},
	YEAR = {2020},
	PAGES = {13005--13015},
	PUBLISHER = {Curran Associates, Inc.},
	ADDRESS = {Red Hook, NY},
	URL = {https://proceedings.neurips.cc/paper/2020/hash/979a3f14bae523dc5101c52120c535e9-Abstract.html},
}

@article {CoDHinrichs,
	AUTHOR = {Hinrichs, Aicke and Novak, Erich and Wo\'zniakowski, Henryk},
	TITLE = {The curse of dimensionality for the class of monotone
	functions and for the class of convex functions},
	JOURNAL = {J. Approx. Theory},
	FJOURNAL = {Journal of Approximation Theory},
	VOLUME = {163},
	YEAR = {2011},
	NUMBER = {8},
	PAGES = {955--965},
	ISSN = {0021-9045,1096-0430},
	MRCLASS = {41A63 (41A55 65D30 65Y20)},
	MRNUMBER = {2832759},
	DOI = {10.1016/j.jat.2011.02.009},
	URL = {https://doi.org/10.1016/j.jat.2011.02.009},
}

@InProceedings{He_2016_CVPR,
author = {He, Kaiming and Zhang, Xiangyu and Ren, Shaoqing and Sun, Jian},
title = {Deep Residual Learning for Image Recognition},
booktitle = {Proceedings of the IEEE Conference on Computer Vision and Pattern Recognition (CVPR)},
month = {June},
year = {2016},
 pages={770--778}
}

@inproceedings {BERT,
	AUTHOR = {Devlin, Jacob and Chang, Ming-Wei and Lee, Kenton
	and Toutanova, Kristina},
	TITLE = {{BERT}: Pre-training of deep bidirectional transformers
	for language understanding},
	BOOKTITLE = {Proceedings of the 2019 Conference of the North
	American Chapter of the Association for Computational
	Linguistics: Human Language Technologies,
	Volume 1 (Long and Short Papers)},
	EDITOR = {Burstein, Jill and Doran, Christy and Solorio, Thamar},
	YEAR = {2019},
	PAGES = {4171--4186},
	PUBLISHER = {Association for Computational Linguistics},
	ADDRESS = {Minneapolis, Minnesota},
	DOI = {10.18653/v1/N19-1423},
	URL = {https://aclanthology.org/N19-1423/},
}

@InProceedings{Huang_2017_CVPR,
author = {Huang, Gao and Liu, Zhuang and van der Maaten, Laurens and Weinberger, Kilian Q.},
title = {Densely Connected Convolutional Networks},
booktitle = {Proceedings of the IEEE Conference on Computer Vision and Pattern Recognition (CVPR)},
month = {July},
pages={4700--4708},
year = {2017}
}

@misc{Petersen2020NeuralNT,
	AUTHOR       = {Petersen, Philipp Christian},
	TITLE        = {Neural network theory},
	YEAR         = {2020},
	HOWPUBLISHED = {Lecture notes},
	URL          = {http://pc-petersen.eu/Neural_Network_Theory.pdf},
	NOTE         = {\url{http://pc-petersen.eu/Neural_Network_Theory.pdf},
	accessed 6 November 2026},
}

@article {Venturi21,
	AUTHOR = {Venturi, Luca and Jelassi, Samy and Ozuch, Tristan
	and Bruna, Joan},
	TITLE = {Depth separation beyond radial functions},
	JOURNAL = {J. Mach. Learn. Res.},
	FJOURNAL = {Journal of Machine Learning Research},
	VOLUME = {23},
	YEAR = {2022},
	PAGES = {Paper No. 122, 56},
	URL = {https://jmlr.org/papers/v23/21-1109.html},
}

@article{Telgarsky15,
  author    = {Matus Telgarsky},
  title     = {Representation Benefits of Deep Feedforward Networks},
JOURNAL = {\href{https://arxiv.org/abs/1509.08101}{arXiv:1509.08101}},
eprint = {1509.08101},
  year      = {2015},
  url       = {http://arxiv.org/abs/1509.08101},
  bibsource = {dblp computer science bibliography, https://dblp.org}
}

@inproceedings{Pascanu13,
	AUTHOR    = {Pascanu, Razvan and Mont{\'u}far, Guido
	and Bengio, Yoshua},
	TITLE     = {On the number of inference regions of deep feed
	forward networks with piece-wise linear activations},
	BOOKTITLE = {Proceedings of the 2nd International Conference
	on Learning Representations},
	YEAR      = {2014},
	ADDRESS   = {Banff, Canada},
	URL       = {https://arxiv.org/abs/1312.6098},
}

@article {Goujon22,
	AUTHOR = {Goujon, Alexis and Etemadi, Arian and Unser, Michael},
	TITLE = {On the number of regions of piecewise linear neural
	networks},
	JOURNAL = {J. Comput. Appl. Math.},
	FJOURNAL = {Journal of Computational and Applied Mathematics},
	VOLUME = {441},
	YEAR = {2024},
	PAGES = {Paper No. 115667, 22},
	ISSN = {0377-0427},
	DOI = {10.1016/j.cam.2023.115667},
	URL = {https://doi.org/10.1016/j.cam.2023.115667},
}

@inproceedings {TakaiYuukiSannai20,
	AUTHOR = {Takai, Yuuki and Sannai, Akiyoshi and Cordonnier,
	Matthieu},
	TITLE = {On the number of linear functions composing deep neural
	network: towards a refined definition of neural networks
	complexity},
	BOOKTITLE = {Proceedings of the 24th International Conference on
	Artificial Intelligence and Statistics},
	EDITOR = {Banerjee, Arindam and Fukumizu, Kenji},
	SERIES = {Proceedings of Machine Learning Research},
	VOLUME = {130},
	YEAR = {2021},
	PAGES = {3799--3807},
	PUBLISHER = {PMLR},
	URL = {https://proceedings.mlr.press/v130/takai21a.html},
}

@article{Opschoor2021ExponentialRD,
    AUTHOR = {Opschoor, J. A. A. and Schwab, Ch. and Zech, J.},
     TITLE = {Exponential {R}e{LU} {DNN} expression of holomorphic maps in
              high dimension},
   JOURNAL = {Constr. Approx.},
  FJOURNAL = {Constructive Approximation. An International Journal for
              Approximations and Expansions},
    VOLUME = {55},
      YEAR = {2022},
    NUMBER = {1},
     PAGES = {537--582},
      ISSN = {0176-4276},
   MRCLASS = {41A46 (41A63 68T07)},
  MRNUMBER = {4376568},
MRREVIEWER = {Agnieszka Lisowska},
       DOI = {10.1007/s00365-021-09542-5},
       URL = {https://doi.org/10.1007/s00365-021-09542-5},
}

@article{YuAnnan21,
  doi = {10.48550/ARXIV.2109.11354},
  
  url = {https://arxiv.org/abs/2109.11354},

JOURNAL = {\href{https://arxiv.org/abs/2109.11354}{arXiv:2109.11354}},
eprint = {2109.11354},
  
  author = {Yu, Annan and Becquey, Chloé and Halikias, Diana and Mallory, Matthew Esmaili and Townsend, Alex},
  
  title = {Arbitrary-Depth Universal Approximation Theorems for Operator Neural Networks},
  
  publisher = {arXiv},
  
  year = {2021},
  
  copyright = {Creative Commons Attribution 4.0 International}
}

@inproceedings {VardiShamir20,
	AUTHOR = {Vardi, Gal and Shamir, Ohad},
	TITLE = {Neural networks with small weights and depth-separation
	barriers},
	BOOKTITLE = {Advances in Neural Information Processing Systems},
	EDITOR = {Larochelle, Hugo and Ranzato, Marc'Aurelio and
	Hadsell, Raia and Balcan, Maria-Florina and
	Lin, Hsuan-Tien},
	VOLUME = {33},
	YEAR = {2020},
	PAGES = {19433--19442},
	PUBLISHER = {Curran Associates, Inc.},
	ADDRESS = {Red Hook, NY},
	URL = {https://proceedings.neurips.cc/paper/2020/hash/e1fe6165cad3f7f3f57d409f78e4415f-Abstract.html},
}

@article {SafranEldanShamir19,
	AUTHOR = {Safran, Itay and Eldan, Ronen and Shamir, Ohad},
	TITLE = {Depth separations in neural networks: what is actually
	being separated?},
	JOURNAL = {Constr. Approx.},
	FJOURNAL = {Constructive Approximation. An International Journal for
	Approximations and Expansions},
	VOLUME = {55},
	YEAR = {2022},
	NUMBER = {1},
	PAGES = {225--257},
	ISSN = {0176-4276},
	DOI = {10.1007/s00365-021-09532-7},
	URL = {https://doi.org/10.1007/s00365-021-09532-7},
}

@inproceedings {Telgarsky16,
	AUTHOR = {Telgarsky, Matus},
	TITLE = {Benefits of depth in neural networks},
	BOOKTITLE = {29th Annual Conference on Learning Theory},
	EDITOR = {Feldman, Vitaly and Rakhlin, Alexander and Shamir, Ohad},
	SERIES = {Proceedings of Machine Learning Research},
	VOLUME = {49},
	YEAR = {2016},
	PAGES = {1517--1539},
	PUBLISHER = {PMLR},
	URL = {https://proceedings.mlr.press/v49/telgarsky16.html},
}

@inproceedings {RolnickT17,
	AUTHOR = {Rolnick, David and Tegmark, Max},
	TITLE = {The power of deeper networks for expressing natural
	functions},
	BOOKTITLE = {The Sixth International Conference on Learning
	Representations},
	YEAR = {2018},
	URL = {https://openreview.net/forum?id=SyProzZAW},
}

@article{Lin_2017,
	doi = {10.1007/s10955-017-1836-5},
  
	url = {https://doi.org/10.1007%2Fs10955-017-1836-5},
  
	year = 2017,
	month = {jul},
  
	publisher = {Springer Science and Business Media {LLC}
},
  
	volume = {168},
  
	number = {6},
  
	pages = {1223--1247},
  
	author = {Henry W. Lin and Max Tegmark and David Rolnick},
  
	title = {Why Does Deep and Cheap Learning Work So Well?},
  
	journal = {Journal of Statistical Physics}
}

@article{beneventano21,
    doi = {10.48550/ARXIV.2112.14523},
  
  url = {https://arxiv.org/abs/2112.14523},
JOURNAL = {\href{https://arxiv.org/abs/2112.14523}{arXiv:2112.14523}},
eprint = {2112.14523},
  
  author = {Beneventano, Pierfrancesco and Cheridito, Patrick and Graeber, Robin and Jentzen, Arnulf and Kuckuck, Benno},
  
  title = {Deep neural network approximation theory for high-dimensional functions},
  
  publisher = {arXiv},
  
  year = {2021},
  
  copyright = {arXiv.org perpetual, non-exclusive license}
}

@InProceedings{safran17a,
  title={Depth-width tradeoffs in approximating natural functions with neural networks},
  author =       {Itay Safran and Ohad Shamir},
  booktitle = 	 {Proceedings of the 34th International Conference on Machine Learning},
  pages = 	 {2979--2987},
  year = 	 {2017},
  month = 	 {06--11 Aug},
  editor = 	 {Precup, Doina and Teh, Yee Whye},
  volume = 	 {70},
  series = 	 {Proceedings of Machine Learning Research},
  publisher =    {PMLR},
  url = 	 {https://proceedings.mlr.press/v70/safran17a.html}
}

@article {MaioPinkus1999,
AUTHOR = {Maiorov, Vitaly and Pinkus, Allan},
TITLE = {Lower bounds for approximation by {MLP} neural networks},
JOURNAL = {Neurocomputing},
VOLUME = {25},
YEAR = {1999},
NUMBER = {1},
PAGES = {81--91},
ISSN = {0925-2312},
DOI = {10.1016/S0925-2312(98)00111-8},
URL = {http://www.sciencedirect.com/science/article/pii/S0925231298001118},
PUBLISHER = {Elsevier},
}

@article {GULIYEV2018296,
AUTHOR = {Namig J. Guliyev and Vugar E. Ismailov},
TITLE = {{On the approximation by single hidden layer feedforward neural networks with fixed weights}},
JOURNAL = {Neural Netw.},
FJOURNAL = {Neural Networks},
VOLUME = {98},
YEAR = {2018},
PAGES = {296--304},
ISSN = {0893-6080},
DOI = {https://doi.org/10.1016/j.neunet.2017.12.007},
URL = {https://www.sciencedirect.com/science/article/pii/S0893608017302927},
}

@article {GuliIsm2018a,
AUTHOR = {Guliyev, Namig J. and Ismailov, Vugar E.},
TITLE = {Approximation capability of two hidden layer feedforward neural networks with fixed weights},
JOURNAL = {Neurocomputing},
VOLUME = {316},
YEAR = {2018},
PAGES = {262--269},
ISSN = {0925-2312},
DOI = {10.1016/j.neucom.2018.07.075},
URL = {http://www.sciencedirect.com/science/article/pii/S0925231218309111},
PUBLISHER = {Elsevier},
}

@article {LuShenYangZhang2020,
	AUTHOR = {Lu, Jianfeng and Shen, Zuowei and Yang, Haizhao and
	Zhang, Shijun},
	TITLE = {Deep network approximation for smooth functions},
	JOURNAL = {SIAM J. Math. Anal.},
	FJOURNAL = {SIAM Journal on Mathematical Analysis},
	VOLUME = {53},
	YEAR = {2021},
	NUMBER = {5},
	PAGES = {5465--5506},
	ISSN = {0036-1410},
	DOI = {10.1137/20M134695X},
	URL = {https://doi.org/10.1137/20M134695X},
}

@article{mhaskar1993,
    AUTHOR = {Mhaskar, H. N.},
     TITLE = {Approximation properties of a  multilayered feedforward artificial neural network},
   JOURNAL = {Adv. Comput. Math.},
   FJOURNAL = {Advances in Computational Mathematics},
    VOLUME = {1},
      YEAR = {1993},
     PAGES = {61--80},
	URL = {https://doi.org/10.1007/BF02070821}
}

@article{bach2017breaking,
    AUTHOR = {Bach, Francis},
     TITLE = {Breaking the curse of dimensionality with convex neural networks},
   JOURNAL = {J. Mach. Learn. Res.},
  FJOURNAL = {Journal of Machine Learning Research (JMLR)},
    VOLUME = {18},
      YEAR = {2017},
      NUMBER = {19},
     pages={1--53} ,
      ISSN = {1532-4435},
   MRCLASS = {92B20 (62G05 90C25)},
  MRNUMBER = {3634886},
}

@article {barron1993universal,
    AUTHOR = {Barron, Andrew R.},
     TITLE = {Universal approximation bounds for superpositions of a sigmoidal function},
   JOURNAL = {IEEE Trans. Inform. Theory},
  FJOURNAL = {Institute of Electrical and Electronics Engineers. Transactions on Information Theory},
    VOLUME = {39},
      YEAR = {1993},
    NUMBER = {3},
     PAGES = {930--945},
      ISSN = {0018-9448},
   MRCLASS = {92B20 (68T05)},
  MRNUMBER = {1237720},
       DOI = {10.1109/18.256500},
       URL = {https://doi.org/10.1109/18.256500},
}

@article {barron1994approximation,
    AUTHOR = {Barron, Andrew R.},
     TITLE = {Approximation and estimation bounds for artificial neural networks},
   JOURNAL = {Mach. Learn.},
  FJOURNAL = {Machine Learning},
    VOLUME = {14},
      YEAR = {1994},
    NUMBER = {1},
     PAGES = {115--133},
      ISSN = {0885-6125},
       DOI = {10.1007/BF00993164},
       URL = {https://doi.org/10.1007/BF00993164},
 PUBLISHER = {Springer},
}

@article {BeckJentzenKuckuck2019,
  AUTHOR = {Beck, Christian and Jentzen, Arnulf and Kuckuck, Benno},
     TITLE = {Full error analysis for the training of deep neural networks},
   JOURNAL = {Infin. Dimens. Anal. Quantum Probab. Relat. Top.},
  FJOURNAL = {Infinite Dimensional Analysis, Quantum Probability and Related
              Topics},
    VOLUME = {25},
      YEAR = {2022},
    NUMBER = {2},
     PAGES = {Paper No. 2150020, 76},
      ISSN = {0219-0257},
   MRCLASS = {65G99 (68T07)},
  MRNUMBER = {4440215},
       DOI = {10.1142/S021902572150020X},
       URL = {https://doi.org/10.1142/S021902572150020X},
}

@book {Bellman1957,
    AUTHOR = {Bellman, Richard},
     TITLE = {Dynamic programming},
    SERIES = {Princeton Landmarks in Mathematics},
      NOTE = {Reprint of the 1957 edition,
              With a new introduction by Stuart Dreyfus},
 PUBLISHER = {Princeton University Press, Princeton, NJ},
      YEAR = {2010},
     PAGES = {xxx+340},
      ISBN = {978-0-691-14668-3},
   MRCLASS = {90-01 (49L20 90C39)},
  MRNUMBER = {2641641},
}

@article {BolcskeiGrohsKutyniokPetersen2019OptimalApproximation,
    AUTHOR = {B\"{o}lcskei, Helmut and Grohs, Philipp and Kutyniok, Gitta and Petersen, Philipp},
     TITLE = {Optimal approximation with sparsely connected deep neural networks},
   JOURNAL = {SIAM J. Math. Data Sci.},
  FJOURNAL = {SIAM Journal on Mathematics of Data Science},
    VOLUME = {1},
      YEAR = {2019},
    NUMBER = {1},
     PAGES = {8--45},
      ISSN = {2577-0187},
   MRCLASS = {41A25 (41A46 42C15 42C40 68T05 82C32 94A12 94A34)},
  MRNUMBER = {3949699},
MRREVIEWER = {Asuman Aksoy},
       DOI = {10.1137/18M118709X},
       URL = {https://doi.org/10.1137/18M118709X},
}

@article {ChenWu2019,
   AUTHOR = {Chen, Liang and Wu, Congwei},
     TITLE = {A note on the expressive power of deep rectified linear unit
              networks in high-dimensional spaces},
   JOURNAL = {Math. Methods Appl. Sci.},
  FJOURNAL = {Mathematical Methods in the Applied Sciences},
    VOLUME = {42},
      YEAR = {2019},
    NUMBER = {9},
     PAGES = {3400--3404},
      ISSN = {0170-4214},
   MRCLASS = {41A30 (41A63)},
  MRNUMBER = {3949573},
MRREVIEWER = {Aurelian Bejancu, Jr.},
       DOI = {10.1002/mma.5575},
       URL = {https://doi.org/10.1002/mma.5575},
}

@article {Cybenko1989,
AUTHOR = {Cybenko, G.},
     TITLE = {Approximation by superpositions of a sigmoidal function},
   JOURNAL = {Math. Control Signals Systems},
  FJOURNAL = {Mathematics of Control, Signals, and Systems},
    VOLUME = {2},
      YEAR = {1989},
    NUMBER = {4},
     PAGES = {303--314},
      ISSN = {0932-4194},
   MRCLASS = {41A30 (92A09 92A25)},
  MRNUMBER = {1015670},
MRREVIEWER = {A. Haimovici},
       DOI = {10.1007/BF02551274},
       URL = {https://doi.org/10.1007/BF02551274},
}

@article {wang2018exponential,
    AUTHOR = {E, Weinan and Wang, Qingcan},
     TITLE = {Exponential convergence of the deep neural network approximation for analytic functions},
   JOURNAL = {Science China Mathematics},
    VOLUME = {61},
      YEAR = {2018},
    NUMBER = {10},
     PAGES = {1733--1740},
      ISSN = {1869-1862},
       DOI = {10.1007/s11425-018-9387-x},
       URL = {https://doi.org/10.1007/s11425-018-9387-x},
}

@article {ElbraechterSchwab2018,
AUTHOR = {Elbr\"{a}chter, Dennis and Grohs, Philipp and Jentzen, Arnulf and
              Schwab, Christoph},
     TITLE = {D{NN} expression rate analysis of high-dimensional {PDE}s:
              application to option pricing},
   JOURNAL = {Constr. Approx.},
  FJOURNAL = {Constructive Approximation. An International Journal for
              Approximations and Expansions},
    VOLUME = {55},
      YEAR = {2022},
    NUMBER = {1},
     PAGES = {3--71},
      ISSN = {0176-4276},
   MRCLASS = {65D30 (35K25 68T07 91G20 91G60)},
  MRNUMBER = {4376559},
       DOI = {10.1007/s00365-021-09541-6},
       URL = {https://doi.org/10.1007/s00365-021-09541-6},
}

@inproceedings {eldan2016power,
  title = 	 {The Power of Depth for Feedforward Neural Networks},
  author = 	 {Eldan, Ronen and Shamir, Ohad},
  booktitle = 	 {Proceedings of the 29th Annual Conference on Learning Theory},
  pages = 	 {907--940},
  year = 	 {2016},
  editor = 	 {Feldman, Vitaly and Rakhlin, Alexander and Shamir, Ohad},
  volume = 	 {49},
  series = 	 {Proceedings of Machine Learning Research},
  address = 	 {Columbia University, New York, New York, USA},
  month = 	 {23--26 Jun},
  publisher =    {PMLR},
  url = 	 {https://proceedings.mlr.press/v49/eldan16.html}
}

@article {funahashi1989approximate,
    AUTHOR = {Funahashi, Ken-Ichi},
     TITLE = {On the approximate realization of continuous mappings by neural networks},
   FJOURNAL = {Neural Networks},
   JOURNAL = {Neural Netw.},
    VOLUME = {2},
      YEAR = {1989},
    NUMBER = {3},
     PAGES = {183--192},
      ISSN = {0893-6080},
       DOI = {10.1016/0893-6080(89)90003-8},
       URL = {http://www.sciencedirect.com/science/article/pii/0893608089900038},
 PUBLISHER = {Elsevier},
}

@article {GrohsHornungJentzen2019,
    AUTHOR = {Grohs, Philipp and Hornung, Fabian and Jentzen, Arnulf and Zimmermann, Philipp},
     TITLE = {Space-time error estimates for deep neural network approximations for differential equations},
   JOURNAL = {Adv. Comput. Math.},
  FJOURNAL = {Advances in Computational Mathematics},
    VOLUME = {49},
      YEAR = {2023},
    PAGES = {Paper No. 4, 78},
       DOI = {10.1007/s10444-022-09970-2},
       URL = {https://doi.org/10.1007/s10444-022-09970-2},
}

@article {GrohsIbrgimovJentzen2021,
	AUTHOR = {Grohs, Philipp and Ibragimov, Shokhrukh and Jentzen,
	Arnulf and Koppensteiner, Sarah},
	TITLE = {Lower bounds for artificial neural network approximations:
	{A} proof that shallow neural networks fail to overcome
	the curse of dimensionality},
	JOURNAL = {J. Complexity},
	FJOURNAL = {Journal of Complexity},
	VOLUME = {77},
	YEAR = {2023},
	PAGES = {Paper No. 101746, 53},
	ISSN = {0885-064X},
	DOI = {10.1016/j.jco.2023.101746},
	URL = {https://doi.org/10.1016/j.jco.2023.101746},
}

@article {GrohsJentzenSalimova2019,
    AUTHOR = {Grohs, Philipp and Jentzen, Arnulf and Salimova, Diyora},
     TITLE = {Deep neural network approximations for solutions of {PDE}s
              based on {M}onte {C}arlo algorithms},
   JOURNAL = {Partial Differ. Equ. Appl.},
  FJOURNAL = {Partial Differential Equations and Applications},
    VOLUME = {3},
      YEAR = {2022},
    NUMBER = {4},
     PAGES = {Paper No. 45, 41},
      ISSN = {2662-2963},
   MRCLASS = {65C05 (35-04)},
  MRNUMBER = {4437476},
       DOI = {10.1007/s42985-021-00100-z},
       URL = {https://doi.org/10.1007/s42985-021-00100-z},
}

@article {guhring2019error,
    AUTHOR = {G\"{u}hring, Ingo and Kutyniok, Gitta and Petersen, Philipp},
     TITLE = {Error bounds for approximations with deep {R}e{LU} neural
              networks in {$W^{s,p}$} norms},
   JOURNAL = {Anal. Appl. (Singap.)},
  FJOURNAL = {Analysis and Applications},
    VOLUME = {18},
      YEAR = {2020},
    NUMBER = {5},
     PAGES = {803--859},
      ISSN = {0219-5305},
   MRCLASS = {41A30 (35A35 41A25 46E35 65N30 68T07)},
  MRNUMBER = {4131039},
MRREVIEWER = {Bao Huai Sheng},
       DOI = {10.1142/S0219530519410021},
       URL = {https://doi.org/10.1142/S0219530519410021},
}

@article {Hanin2017,
	AUTHOR = {Hanin, Boris},
	TITLE = {Universal function approximation by deep neural nets with
	bounded width and {ReLU} activations},
	JOURNAL = {Mathematics},
	FJOURNAL = {Mathematics},
	VOLUME = {7},
	YEAR = {2019},
	NUMBER = {10},
	PAGES = {Paper No. 992, 9},
	ISSN = {2227-7390},
	DOI = {10.3390/math7100992},
	URL = {https://doi.org/10.3390/math7100992},
}

@article {hornik1991approximation,
    AUTHOR = {Hornik, Kurt},
     TITLE = {Approximation capabilities of multilayer feedforward networks},
   JOURNAL = {Neural Netw.},
   FJOURNAL = {Neural Networks},
    VOLUME = {4},
      YEAR = {1991},
    NUMBER = {2},
     PAGES = {251--257},
      ISSN = {0893-6080},
       DOI = {10.1016/0893-6080(91)90009-T},
       URL = {http://www.sciencedirect.com/science/article/pii/089360809190009T},
 PUBLISHER = {Elsevier},
}

@article {hornik1989multilayer,
    AUTHOR = {Hornik, Kurt and Stinchcombe, Maxwell and White, Halbert},
     TITLE = {Multilayer feedforward networks are universal approximators},
   JOURNAL = {Neural Netw.},
   FJOURNAL = {Neural Networks},
    VOLUME = {2},
      YEAR = {1989},
    NUMBER = {5},
     PAGES = {359--366},
      ISSN = {0893-6080},
       DOI = {10.1016/0893-6080(89)90020-8},
       URL = {http://www.sciencedirect.com/science/article/pii/0893608089900208},
 PUBLISHER = {Elsevier},
}

@article {leshno1993multilayer,
    AUTHOR = {Leshno, Moshe and Lin, Vladimir Ya and Pinkus, Allan and Schocken, Shimon},
     TITLE = {Multilayer feedforward networks with a nonpolynomial activation function can approximate any function},
   JOURNAL = {Neural Netw.},
   FJOURNAL = {Neural Networks},
    VOLUME = {6},
      YEAR = {1993},
    NUMBER = {6},
     PAGES = {861--867},
      ISSN = {0893-6080},
       DOI = {10.1016/S0893-6080(05)80131-5},
       URL = {http://www.sciencedirect.com/science/article/pii/S0893608005801315},
 PUBLISHER = {Elsevier},
}

@book {NovakWozniakowski2008,    
AUTHOR = {Novak, Erich and Wo\'{z}niakowski, Henryk},
     TITLE = {Tractability of multivariate problems. {V}ol. 1: {L}inear
              information},
    SERIES = {EMS Tracts in Mathematics},
    VOLUME = {6},
 PUBLISHER = {European Mathematical Society (EMS), Z\"{u}rich},
      YEAR = {2008},
     PAGES = {xii+384},
      ISBN = {978-3-03719-026-5},
   MRCLASS = {46E22 (41A63 41A65 46N40 47N40)},
  MRNUMBER = {2455266},
MRREVIEWER = {Arthur G. Werschulz},
       DOI = {10.4171/026},
       URL = {https://doi.org/10.4171/026},
}

@book {NovakWozniakowski2010,
    AUTHOR = {Novak, Erich and Wo\'{z}niakowski, Henryk},
     TITLE = {Tractability of multivariate problems. {V}olume {II}:
              {S}tandard information for functionals},
    SERIES = {EMS Tracts in Mathematics},
    VOLUME = {12},
 PUBLISHER = {European Mathematical Society (EMS), Z\"{u}rich},
      YEAR = {2010},
     PAGES = {xviii+657},
      ISBN = {978-3-03719-084-5},
   MRCLASS = {46N40 (41A55 41A63 47A50 47N40 65Y20 68Q25)},
  MRNUMBER = {2676032},
MRREVIEWER = {Arthur G. Werschulz},
       DOI = {10.4171/084},
       URL = {https://doi.org/10.4171/084},
}

@article {perekrestenko2018universal,
    AUTHOR = {Perekrestenko, Dmytro and Grohs, Philipp and Elbr{\"a}chter, Dennis and B{\"o}lcskei, Helmut},
     TITLE = {The universal approximation power of finite-width deep {R}e{LU} networks},
   JOURNAL = {\href{https://arxiv.org/abs/1806.01528}{arXiv:1806.01528}},
   eprint = {1806.01528},
      YEAR = {2018},
     PAGES = {16 pages},
}

@article {petersen2018optimal,
    AUTHOR = {Petersen, Philipp and Voigtlaender, Felix},
     TITLE = {Optimal approximation of piecewise smooth functions using deep {ReLU} neural networks},
   JOURNAL = {Neural Netw.},
   FJOURNAL = {Neural Networks},
    VOLUME = {108},
      YEAR = {2018},
     PAGES = {296--330},
      ISSN = {0893-6080},
       DOI = {10.1016/j.neunet.2018.08.019},
       URL = {http://www.sciencedirect.com/science/article/pii/S0893608018302454},
 PUBLISHER = {Elsevier},
}

@article {ShahamCloningerCoifman2018,
    AUTHOR = {Shaham, Uri and Cloninger, Alexander and Coifman, Ronald R.},
     TITLE = {Provable approximation properties for deep neural networks},
   JOURNAL = {Appl. Comput. Harmon. Anal.},
  FJOURNAL = {Applied and Computational Harmonic Analysis. Time-Frequency and Time-Scale Analysis, Wavelets, Numerical Algorithms, and Applications},
    VOLUME = {44},
      YEAR = {2018},
    NUMBER = {3},
     PAGES = {537--557},
      ISSN = {1063-5203},
   MRCLASS = {94C15 (94A20)},
  MRNUMBER = {3768851},
MRREVIEWER = {Michael M. Dediu},
       DOI = {10.1016/j.acha.2016.04.003},
       URL = {https://doi.org/10.1016/j.acha.2016.04.003},
}

@inproceedings {voigtlaender2019approximation,
	AUTHOR = {Voigtlaender, Felix and Petersen, Philipp},
	TITLE = {Approximation in {$L^p(\mu)$} with deep {ReLU} neural
	networks},
	BOOKTITLE = {2019 13th International Conference on Sampling Theory
	and Applications ({SampTA})},
	YEAR = {2019},
	PAGES = {1--4},
	PUBLISHER = {IEEE},
	ADDRESS = {Piscataway, NJ},
	ISBN = {978-1-7281-3741-4},
	DOI = {10.1109/SampTA45681.2019.9030992},
	URL = {https://doi.org/10.1109/SampTA45681.2019.9030992},
}

@article {yarotsky2017error,
    AUTHOR = {Yarotsky, Dmitry},
     TITLE = {Error bounds for approximations with deep {ReLU} networks},
   JOURNAL = {Neural Netw.},
   FJOURNAL = {Neural Networks},
    VOLUME = {94},
      YEAR = {2017},
     PAGES = {103--114},
      ISSN = {0893-6080},
       DOI = {10.1016/j.neunet.2017.07.002},
       URL = {http://www.sciencedirect.com/science/article/pii/S0893608017301545},
 PUBLISHER = {Elsevier},
}

@ARTICLE{cheridito2021efficient,
    AUTHOR = {Cheridito, Patrick and Jentzen, Arnulf and Rossmannek,
              Florian},
     TITLE = {Efficient approximation of high-dimensional functions with
              neural networks},
   JOURNAL = {IEEE Trans. Neural Netw. Learn. Syst.},
  FJOURNAL = {IEEE Transactions on Neural Networks and Learning Systems},
    VOLUME = {33},
      YEAR = {2022},
    NUMBER = {7},
     PAGES = {3079--3093},
      ISSN = {2162-237X},
   MRCLASS = {65D15 (68T07)},
  MRNUMBER = {4449147},
}

@incollection {guehring2020expressivity,
	AUTHOR = {G{\"u}hring, Ingo and Raslan, Mones and Kutyniok, Gitta},
	TITLE = {Expressivity of deep neural networks},
	BOOKTITLE = {Mathematical aspects of deep learning},
	EDITOR = {Grohs, Philipp and Kutyniok, Gitta},
	PAGES = {149--199},
	PUBLISHER = {Cambridge University Press},
	ADDRESS = {Cambridge},
	YEAR = {2022},
	DOI = {10.1017/9781009025096.004},
	URL = {https://doi.org/10.1017/9781009025096.004},
}

@InProceedings{kidger2020universal,
  title = 	 {{Universal {A}pproximation with {D}eep {N}arrow {N}etworks}},
  author =       {Kidger, Patrick and Lyons, Terry},
  booktitle = 	 {Proceedings of Thirty Third Conference on Learning Theory},
  pages = 	 {2306--2327},
  year = 	 {2020},
  editor = 	 {Abernethy, Jacob and Agarwal, Shivani},
  volume = 	 {125},
  series = 	 {Proceedings of Machine Learning Research},
  month = 	 {09--12 Jul},
  publisher =    {PMLR},
  url = 	 {https://proceedings.mlr.press/v125/kidger20a.html}
}

@article {elbraechter2021deep,
    AUTHOR = {Elbr\"{a}chter, Dennis and Perekrestenko, Dmytro and Grohs,
              Philipp and B\"{o}lcskei, Helmut},
     TITLE = {Deep neural network approximation theory},
   JOURNAL = {IEEE Trans. Inform. Theory},
  FJOURNAL = {Institute of Electrical and Electronics Engineers.
              Transactions on Information Theory},
    VOLUME = {67},
      YEAR = {2021},
    NUMBER = {5},
     PAGES = {2581--2623},
      ISSN = {0018-9448},
   MRCLASS = {68T07},
  MRNUMBER = {4282380},
       DOI = {10.1109/TIT.2021.3062161},
       URL = {https://doi.org/10.1109/TIT.2021.3062161},
}

@article {chui2019deep,
    AUTHOR = {Chui, Charles K. and Lin, Shao-Bo and Zhou, Ding-Xuan},
     TITLE = {Deep neural networks for rotation-invariance approximation and
              learning},
   JOURNAL = {Anal. Appl. (Singap.)},
  FJOURNAL = {Analysis and Applications},
    VOLUME = {17},
      YEAR = {2019},
    NUMBER = {5},
     PAGES = {737--772},
      ISSN = {0219-5305},
   MRCLASS = {41A30 (41A25 68T05)},
  MRNUMBER = {3998499},
MRREVIEWER = {Michael M. Dediu},
       DOI = {10.1142/S0219530519400074},
       URL = {https://doi.org/10.1142/S0219530519400074},
}

@inproceedings{barron1992neural,
  author          = {Barron, A.},
  booktitle       = {Proceedings of the 7th Yale Workshop on Adaptive
and Learning Systems},
  title           = {Neural net approximation},
  year            = {1992},
  pages={69--72}
}

@article {jones1992simple,
    AUTHOR = {Jones, Lee K.},
     TITLE = {A simple lemma on greedy approximation in {H}ilbert space and
              convergence rates for projection pursuit regression and neural
              network training},
   JOURNAL = {Ann. Statist.},
  FJOURNAL = {The Annals of Statistics},
    VOLUME = {20},
      YEAR = {1992},
    NUMBER = {1},
     PAGES = {608--613},
      ISSN = {0090-5364},
   MRCLASS = {62H40 (62J02)},
  MRNUMBER = {1150368},
MRREVIEWER = {Lionel Weiss},
       DOI = {10.1214/aos/1176348546},
       URL = {https://doi.org/10.1214/aos/1176348546},
}

@article {klusowski2018approximation,
    AUTHOR = {Klusowski, Jason M. and Barron, Andrew R.},
     TITLE = {Approximation by combinations of {R}e{LU} and squared {R}e{LU}
              ridge functions with {$\ell^1$} and {$\ell^0$} controls},
   JOURNAL = {IEEE Trans. Inform. Theory},
  FJOURNAL = {Institute of Electrical and Electronics Engineers.
              Transactions on Information Theory},
    VOLUME = {64},
      YEAR = {2018},
    NUMBER = {12},
     PAGES = {7649--7656},
      ISSN = {0018-9448},
   MRCLASS = {41A30 (94A20)},
  MRNUMBER = {3883687},
MRREVIEWER = {Danilo Costarelli},
       DOI = {10.1109/tit.2018.2874447},
       URL = {https://doi.org/10.1109/tit.2018.2874447},
}

@InProceedings{daniely2017depth,
  title = 	 {Depth Separation for Neural Networks},
  author = 	 {Daniely, Amit},
  booktitle = 	 {Proceedings of the 2017 Conference on Learning Theory},
  pages = 	 {690--696},
  year = 	 {2017},
  editor = 	 {Kale, Satyen and Shamir, Ohad},
  volume = 	 {65},
  series = 	 {Proceedings of Machine Learning Research},
  month = 	 {07--10 Jul},
  publisher =    {PMLR},
  url = 	 {https://proceedings.mlr.press/v65/daniely17a.html}
}

@InProceedings{lee2017ability,
  title = 	 {On the {A}bility of {N}eural {N}ets to {E}xpress {D}istributions},
  author = 	 {Lee, Holden and Ge, Rong and Ma, Tengyu and Risteski, Andrej and Arora, Sanjeev},
  booktitle = 	 {Proceedings of the 2017 Conference on Learning Theory},
  pages = 	 {1271--1296},
  year = 	 {2017},
  editor = 	 {Kale, Satyen and Shamir, Ohad},
  volume = 	 {65},
  series = 	 {Proceedings of Machine Learning Research},
  month = 	 {07--10 Jul},
  publisher =    {PMLR},
  url = 	 {https://proceedings.mlr.press/v65/lee17a.html}
}

@article{pinkus1999approximation, 
title={Approximation theory of the {MLP} model in neural networks}, 
volume={8}, 
DOI={10.1017/S0962492900002919}, 
fjournal={Acta Numerica}, 
journal={Acta Numer.}, 
publisher={Cambridge University Press}, 
author={Pinkus, Allan}, 
year={1999}, 
pages={143--195}
}

@inproceedings {carroll1989construction,
	AUTHOR = {Carroll, S. M. and Dickinson, B. W.},
	TITLE = {Construction of neural nets using the {Radon} transform},
	BOOKTITLE = {Proceedings of the International Joint Conference
	on Neural Networks},
	VOLUME = {1},
	PAGES = {607--611},
	PUBLISHER = {IEEE},
	YEAR = {1989},
	DOI = {10.1109/IJCNN.1989.118639},
	URL = {https://doi.org/10.1109/IJCNN.1989.118639},
}

@inproceedings {gallant1988there,
	AUTHOR = {Gallant, A. Ronald and White, Halbert},
	TITLE = {There exists a neural network that does not make
	avoidable mistakes},
	BOOKTITLE = {Proceedings of the {IEEE} International Conference
	on Neural Networks},
	VOLUME = {1},
	PAGES = {657--664},
	PUBLISHER = {IEEE},
	YEAR = {1988},
	DOI = {10.1109/ICNN.1988.23903},
	URL = {https://doi.org/10.1109/ICNN.1988.23903},
}

@inproceedings {irie1988capabilities,
	AUTHOR = {Irie, Bunpei and Miyake, Sei},
	TITLE = {Capabilities of three-layered perceptrons},
	BOOKTITLE = {Proceedings of the {IEEE} International Conference
	on Neural Networks},
	VOLUME = {1},
	PAGES = {641--648},
	PUBLISHER = {IEEE},
	YEAR = {1988},
	DOI = {10.1109/ICNN.1988.23901},
	URL = {https://doi.org/10.1109/ICNN.1988.23901},
}

@article{makovoz1996random,
title = {Random {A}pproximants and {N}eural {N}etworks},
fjournal = {Journal of Approximation Theory},
journal = {J. Approx. Theory},
volume = {85},
number = {1},
pages = {98--109},
year = {1996},
issn = {0021-9045},
doi = {https://doi.org/10.1006/jath.1996.0031},
url = {https://www.sciencedirect.com/science/article/pii/S0021904596900313},
author = {Y. Makovoz},
}

@Article{shen2020deep,
author = {Shen, Zuowei and Yang, Haizhao and Zhang, Shijun},
title = {{D}eep {N}etwork {A}pproximation {C}haracterized by {N}umber of {N}eurons},
 JOURNAL = {Commun. Comput. Phys.},
  FJOURNAL = {Communications in Computational Physics},
  year = {2020},
volume = {28},
number = {5},
pages = {1768--1811},
issn = {1991-7120},
doi = {https://doi.org/10.4208/cicp.OA-2020-0149},
url = {http://global-sci.org/intro/article_detail/cicp/18396.html}}

@article {grohs2021proof,
	AUTHOR = {Grohs, Philipp and Voigtlaender, Felix},
	TITLE = {Proof of the theory-to-practice gap in deep learning via
	sampling complexity bounds for neural network
	approximation spaces},
	JOURNAL = {Found. Comput. Math.},
	FJOURNAL = {Foundations of Computational Mathematics},
	VOLUME = {24},
	YEAR = {2024},
	NUMBER = {4},
	PAGES = {1085--1143},
	ISSN = {1615-3375},
	MRNUMBER = {4783638},
	DOI = {10.1007/s10208-023-09607-w},
	URL = {https://doi.org/10.1007/s10208-023-09607-w},
}

@article{min2017deep,
  title={Deep learning in bioinformatics},
  author={Min, Seonwoo and Lee, Byunghan and Yoon, Sungroh},
  journal={Briefings in bioinformatics},
  volume={18},
  number={5},
  pages={851--869},
  year={2017},
  publisher={Oxford University Press}
}

@article {sidey2019machine,
	AUTHOR = {Sidey-Gibbons, Jenni A. M. and Sidey-Gibbons, Chris J.},
	TITLE = {Machine learning in medicine: a practical introduction},
	JOURNAL = {BMC Med. Res. Methodol.},
	FJOURNAL = {BMC Medical Research Methodology},
	VOLUME = {19},
	YEAR = {2019},
	NUMBER = {1},
	PAGES = {Paper No. 64, 18},
	ISSN = {1471-2288},
	DOI = {10.1186/s12874-019-0681-4},
	URL = {https://doi.org/10.1186/s12874-019-0681-4},
}

@article{PALTRINIERI2019475,
title = {Learning about risk: Machine learning for risk assessment},
journal = {Safety Science},
volume = {118},
pages = {475--486},
year = {2019},
issn = {0925-7535},
doi = {https://doi.org/10.1016/j.ssci.2019.06.001},
url = {https://www.sciencedirect.com/science/article/pii/S0925753518311184},
author = {Nicola Paltrinieri and Louise Comfort and Genserik Reniers}
}

@article{McCaneS17,
  author    = {Brendan McCane and
               Lech Szymanski},
  title     = {Deep Radial Kernel Networks: Approximating Radially Symmetric Functions
               with Deep Networks},
  journal   = {\href{https://arxiv.org/abs/1703.03470}{arXiv:1703.03470}},
  year      = {2017},
  url       = {http://arxiv.org/abs/1703.03470},
  eprinttype = {arXiv},
  eprint    = {1703.03470},
  bibsource = {dblp computer science bibliography, https://dblp.org}
}

@book {jovanovic2014analysis,
    AUTHOR = {Jovanovi\'{c}, Bo\v{s}ko S. and S\"{u}li, Endre},
     TITLE = {Analysis of finite difference schemes},
    SERIES = {Springer Series in Computational Mathematics},
    VOLUME = {46},
      NOTE = {For linear partial differential equations with generalized
              solutions},
 PUBLISHER = {Springer, London},
      YEAR = {2014},
     PAGES = {xiv+408},
      ISBN = {978-1-4471-5459-4; 978-1-4471-5460-0},
   MRCLASS = {65-02 (35B30 35J40 35K35 35L35 65M06 65N06)},
  MRNUMBER = {3136501},
MRREVIEWER = {Enrico Bertolazzi},
       DOI = {10.1007/978-1-4471-5460-0},
       URL = {https://doi.org/10.1007/978-1-4471-5460-0},
}

@article {TadmorEitan12,
    AUTHOR = {Tadmor, Eitan},
     TITLE = {A review of numerical methods for nonlinear partial
              differential equations},
   JOURNAL = {Bull. Amer. Math. Soc. (N.S.)},
  FJOURNAL = {American Mathematical Society. Bulletin. New Series},
    VOLUME = {49},
      YEAR = {2012},
    NUMBER = {4},
     PAGES = {507--554},
      ISSN = {0273-0979},
   MRCLASS = {65-02 (35J60 35K55 35L65 35L70 65Mxx 65Nxx)},
  MRNUMBER = {2958929},
MRREVIEWER = {Paul Andrew Martin},
       DOI = {10.1090/S0273-0979-2012-01379-4},
       URL = {https://doi.org/10.1090/S0273-0979-2012-01379-4},
}

@inproceedings {SaadSolla95,
	AUTHOR = {Saad, David and Solla, Sara A.},
	TITLE = {Dynamics of on-line gradient descent learning
	for multilayer neural networks},
	BOOKTITLE = {Advances in Neural Information Processing Systems},
	EDITOR = {Touretzky, David S. and Mozer, Michael C.
	and Hasselmo, Michael E.},
	VOLUME = {8},
	PAGES = {302--308},
	PUBLISHER = {MIT Press},
	YEAR = {1996},
	ISBN = {9780262201070},
	URL = {https://proceedings.neurips.cc/paper_files/paper/1995/hash/a1519de5b5d44b31a01de013b9b51a80-Abstract.html},
}

@inproceedings {wang2018glue,
	AUTHOR = {Wang, Alex and Singh, Amanpreet and Michael, Julian
	and Hill, Felix and Levy, Omer and Bowman, Samuel R.},
	TITLE = {{GLUE}: A multi-task benchmark and analysis platform
	for natural language understanding},
	BOOKTITLE = {International Conference on Learning Representations},
	YEAR = {2019},
	URL = {https://openreview.net/forum?id=rJ4km2R5t7},
}

@InProceedings{Williams18,
  author = {Williams, Adina
            and Nangia, Nikita
            and Bowman, Samuel},
  title = {{A Broad-Coverage Challenge Corpus for 
           Sentence Understanding through Inference}},
  booktitle = {Proceedings of the 2018 Conference of 
               the North American Chapter of the 
               Association for Computational Linguistics:
               Human Language Technologies, Volume 1 (Long
               Papers)},
  year = {2018},
  publisher = {Association for Computational Linguistics},
  pages = {1112--1122},
  location = {New Orleans, Louisiana},
  url = {http://aclweb.org/anthology/N18-1101}
}

@article{ILSVRC15,
 AUTHOR = {Olga Russakovsky and Jia Deng and Hao Su and Jonathan Krause and Sanjeev Satheesh and Sean Ma and Zhiheng Huang and Andrej Karpathy and Aditya Khosla and Michael Bernstein and Alexander C. Berg and Li Fei-Fei},
     TITLE = {Image{N}et large scale visual recognition challenge},
   JOURNAL = {Int. J. Comput. Vis.},
  FJOURNAL = {International Journal of Computer Vision},
    VOLUME = {115},
      YEAR = {2015},
    NUMBER = {3},
     PAGES = {211--252},
      ISSN = {0920-5691},
   MRCLASS = {94A08},
  MRNUMBER = {3422482},
       DOI = {10.1007/s11263-015-0816-y},
       URL = {https://doi.org/10.1007/s11263-015-0816-y},
}

@article{RieglerBiehl95,
doi = {10.1088/0305-4470/28/20/002},
url = {https://dx.doi.org/10.1088/0305-4470/28/20/002},
year = {1995},
month = {oct},
publisher = {},
volume = {28},
number = {20},
pages = {L507--L513},
author = {P Riegler and  M Biehl},
title = {On-line backpropagation in two-layered neural networks},
journal = {Journal of Physics A: Mathematical and General}
}

@article {ChenEtAl19,
    AUTHOR = {Chen, Tao and Lu, Shijian and Fan, Jiayuan},
     TITLE = {S{S}-{HCNN}: semi-supervised hierarchical convolutional neural
              network for image classification},
   JOURNAL = {IEEE Trans. Image Process.},
  FJOURNAL = {IEEE Transactions on Image Processing},
    VOLUME = {28},
      YEAR = {2019},
    NUMBER = {5},
     PAGES = {2389--2398},
      ISSN = {1057-7149},
   MRCLASS = {Expansion},
  MRNUMBER = {3915633},
       DOI = {10.1109/TIP.2018.2886758},
       URL = {https://doi.org/10.1109/TIP.2018.2886758},
}

@inproceedings {Tsung-Yiet2014,
	AUTHOR = {Lin, Tsung-Yi and Maire, Michael and Belongie, Serge
	and Hays, James and Perona, Pietro and Ramanan, Deva
	and Doll{\'a}r, Piotr and Zitnick, C. Lawrence},
	TITLE = {{Microsoft COCO}: Common objects in context},
	BOOKTITLE = {Computer vision---{ECCV} 2014. Part V},
	EDITOR = {Fleet, David and Pajdla, Tomas and Schiele, Bernt
	and Tuytelaars, Tinne},
	SERIES = {Lecture Notes in Comput. Sci.},
	VOLUME = {8693},
	YEAR = {2014},
	PAGES = {740--755},
	PUBLISHER = {Springer},
	ADDRESS = {Cham},
	DOI = {10.1007/978-3-319-10602-1_48},
	URL = {https://doi.org/10.1007/978-3-319-10602-1_48},
}

@inproceedings{raghu2017expressive,
	AUTHOR    = {Raghu, Maithra and Poole, Ben and Kleinberg, Jon
	and Ganguli, Surya and Sohl-Dickstein, Jascha},
	TITLE     = {On the expressive power of deep neural networks},
	BOOKTITLE = {Proceedings of the 34th International Conference
	on Machine Learning},
	EDITOR    = {Precup, Doina and Teh, Yee Whye},
	SERIES    = {Proceedings of Machine Learning Research},
	VOLUME    = {70},
	YEAR      = {2017},
	PAGES     = {2847--2854},
	PUBLISHER = {PMLR},
	URL       = {https://proceedings.mlr.press/v70/raghu17a.html},
}

@article {GononSchwab21,
    AUTHOR = {Gonon, Lukas and Schwab, Christoph},
     TITLE = {Deep {R}e{LU} network expression rates for option prices in
              high-dimensional, exponential {L}\'{e}vy models},
   JOURNAL = {Finance Stoch.},
  FJOURNAL = {Finance and Stochastics},
    VOLUME = {25},
      YEAR = {2021},
    NUMBER = {4},
     PAGES = {615--657},
      ISSN = {0949-2984,1432-1122},
   MRCLASS = {91G80 (60G51 60H35)},
  MRNUMBER = {4318896},
       DOI = {10.1007/s00780-021-00462-7},
       URL = {https://doi.org/10.1007/s00780-021-00462-7},
}

@article {GononGrohsJentzenKoflerSiska22,
    AUTHOR = {Gonon, Lukas and Grohs, Philipp and Jentzen, Arnulf and
              Kofler, David and \v{S}i\v{s}ka, David},
     TITLE = {Uniform error estimates for artificial neural network
              approximations for heat equations},
   JOURNAL = {IMA J. Numer. Anal.},
  FJOURNAL = {IMA Journal of Numerical Analysis},
    VOLUME = {42},
      YEAR = {2022},
    NUMBER = {3},
     PAGES = {1991--2054},
      ISSN = {0272-4979,1464-3642},
   MRCLASS = {65M99 (35K05)},
  MRNUMBER = {4454915},
       DOI = {10.1093/imanum/drab027},
       URL = {https://doi.org/10.1093/imanum/drab027},
}

@article{KuckuckTools,
	author  = {Kuckuck, Benno},
	title   = {Some useful {\LaTeX} commands},
	journal = {\href{https://latex.bennokuckuck.de}
	{https://latex.bennokuckuck.de}},
	year    = {Accessed on 22-August-2025}
}


\end{document}